\documentclass[12pt]{amsart}
\usepackage{amsmath, amsthm, amssymb,cite,enumerate}
\usepackage{fullpage}
\usepackage{color}
\usepackage{hyperref}
\usepackage[english]{babel}
\usepackage{framed}

\newtheorem{thm}{Theorem}
\newtheorem{lem}{Lemma}[section]
\newtheorem{prop}[lem]{Proposition}
\newtheorem{cor}[lem]{Corollary}
\newtheorem{rmkn}{Remark}[thm]
\newtheorem{rmk}{Remark}[lem]

\newtheorem*{defn}{Definition}
\numberwithin{equation}{section}

\newcommand{\R}{\mathbb{R}}

\newcommand{\F}{\mathcal{F}}
\newcommand{\eps}{\varepsilon}
\newcommand{\norm}[1]{\left\lVert#1\right\rVert}
\newcommand{\mcl}[1]{\mathcal{ #1 }}
\newcommand{\fm}[1]{\[\begin{aligned} #1 \end{aligned}\]}
\newcommand{\eq}[1]{\begin{equation}\begin{aligned} #1 \end{aligned}\end{equation}}
\newcommand{\lra}[1]{\langle #1 \rangle}

\newcommand{\wt}[1]{\widetilde{ #1 }}

\DeclareMathOperator{\tr}{tr}

\begin{document}

\title{Asymptotic completeness for a scalar quasilinear wave equation satisfying the weak null condition}

\author{Dongxiao Yu}
\address{Department of Mathematics, University of California at Berkeley}
\thanks{}
\email{yudx@math.berkeley.edu}

\providecommand{\keywords}[1]
{
  \small	
  \textbf{\textit{Keywords---}} #1
}

\begin{abstract}

In this paper, we prove the first asymptotic completeness result for a scalar quasilinear wave equation satisfying the weak null condition. The main tool we use in the study of this equation is the geometric reduced system introduced in \cite{yu2020}. Starting from a global solution $u$ to the quasilinear wave equation, we rigorously show that well chosen asymptotic variables solve the same reduced system with small error terms. This allows us to recover the scattering data for our system, as well as to construct a matching exact solution to the reduced system. 

\end{abstract}



\maketitle
\tableofcontents\addtocontents{toc}{\protect\setcounter{tocdepth}{1}}

\section{Introduction}
This paper is devoted to the study of the long time dynamics for a scalar quasilinear wave equation in $\R^{1+3}_{t,x}$, of the form
\eq{\label{qwe} g^{\alpha\beta}(u)\partial_\alpha \partial_\beta u=0}
with small initial data
\eq{\label{init} (u,u_t)|_{t=0}=(\eps u_0,\eps u_1)\in C_c^\infty(\R^3),\ 0<\eps\ll1.}
Here we use the Einstein summation convention, with the sum taken over $\alpha,\beta=0,1,2,3$ with  $\partial_0=\partial_t$, $\partial_i=\partial_{x_i}$, $i=1,2,3$. We assume that  $g^{\alpha\beta}(u)$ are smooth functions of $u$, such that $g^{\alpha\beta}=g^{\beta\alpha}$ and $g^{\alpha\beta}(0)=m^{\alpha\beta}$ where $(m^{\alpha\beta})=\text{diag}(-1,1,1,1)$ is the Minkowski metric. Here we can assume that $g^{00}\equiv -1$. In fact, since we expect $|u|\ll 1$, we have $g^{00}(u)<0$, so we can replace $(g^{\alpha\beta})$ with $(g^{\alpha\beta}/(-g^{00}))$ if necessary.

The study of global well-posededness theory of \eqref{qwe} started with Lindblad's paper~\cite{lind2}.  Lindblad conjectured that the equation \eqref{qwe} has a global solution if $\eps$ in \eqref{init} is sufficiently small. In the same paper, he proved the small data global existence for a special case 
\begin{equation}\label{qwespec}-\partial_t^2u+c(u)^2\Delta_xu=0,\hspace{1cm}\text{where }c(0)=1\end{equation}for radially symmetric data.
Later, Alinhac \cite{alin2} generalized  the result to general initial data for \eqref{qwespec}. The small data global existence result for the general case \eqref{qwe} was finally proved by Lindblad  \cite{lind}.

In the author's recent paper \cite{yu2020}, we have identified a new notion of asymptotic profile  and an associated notion of scattering data for the model equation, by deriving a new reduced system. With these new notions, we have proved the existence of the modified wave operators for \eqref{qwe}.

In this paper, we seek to continue the study of modified scattering by proving the asymptotic completeness for \eqref{qwe}. That is, given a global solution to the Cauchy problem \eqref{qwe} and \eqref{init}, we seek to find the corresponding asymptotic profile and scattering data associated to this global solution.

Given a global solution $u$, we start the proof with the construction of a global optical function $q=q(t,x)$. In other words, we solve the eikonal equation $g^{\alpha\beta}(u)q_\alpha q_\beta=0$ in a spacetime region $\Omega$ contained in $\{2r\geq t\geq \exp(\delta/\eps)\}$; this is where our evolution is expected to have a nonlinear behavior. Here $\delta>0$ is a small fixed parameter. We apply the method of characteristics and then follow the idea in Christodoulou-Klainerman \cite{chriklai}. By viewing $(g_{\alpha\beta})$, the inverse of the coefficient matrix $(g^{\alpha\beta}(u))$, as a Lorentzian metric in $[0,\infty)\times\R^3$, we construct a null frame $\{e_k\}_{k=1}^4$. Then, most importantly, we define the second fundamental form $\chi_{ab}$ for $a,b=1,2$ which are related to  the Levi-Civita connection and the null frame under the metric $(g_{\alpha\beta})$.  By  studying the Raychaudhuri equation and using a continuity argument, we can show that   $\tr\chi>0$  everywhere. This is the key step which guarantees that the solutions to the eikonal equation are global. In addition, we can prove that $q=q(t,x)$ is smooth in some weak sense (see Section \ref{sec2.5}). We refer our readers to Section \ref{s3} and Section \ref{s4} for more details in the proof.

Next, following \cite{yu2020}, we define our asymptotic variables  $(\mu,U)(t,x):=(q_t-q_r,\eps^{-1}ru)(t,x)$. The map \fm{\Omega\to[0,\infty)\times\R\times\mathbb{S}^2:\hspace{2em}(t,x)\mapsto (\eps\ln t-\delta,q(t,x),x/|x|):=(s,q,\omega)} is an injective smooth function with a smooth inverse, so a function $(\mu,U)(s,q,\omega)$ is obtained.  It can be proved that $(\mu,U)(s,q,\omega)$ is an approximate solution to the reduced system introduced in \cite{yu2020}, and that there is an exact solution $(\wt{\mu},\wt{U})(s,q,\omega)$ to the reduced system which matches $(\mu,U)(s,q,\omega)$ as $s\to\infty$. A key step  is to prove that  $A(q,\omega):=-\frac{1}{2}\lim_{s\to\infty}(\mu U_q)(s,q,\omega)$ is well-defined for each $(q,\omega)$. The function $A$ is called the \emph{scattering data} in this paper. We also show a gauge independence result, which  states that the scattering data for the solution $u$ is independent of  the choice of the optical function $q$ in a suitable sense. We refer our readers to Section \ref{s5} and Section \ref{ss5.7}.

Finally, starting from the scattering data $A$, we show that we can construct an approximate solution $\wt{u}$ to \eqref{qwe} in $\Omega$. The construction here is similar to that in Section 4 of \cite{yu2020}. That is, we construct a function $\wt{q}$ by solving \fm{\wt{q}_t-\wt{q}_r=\mu(\eps\ln t-\delta,\wt{q}(t,x),\omega)} by the method of characteristics, and then define \fm{\wt{u}(t,x):=\eps r^{-1}\wt{U}(\eps\ln t-\delta,\wt{q}(t,x),\omega).}
Then, in $\Omega$, $\wt{q}$ is an approximate optical function, and $\wt{u}$ is an approximate solution to \eqref{qwe}. In addition, near the light cone $t=r$, the difference $u-\wt{u}$, along with its derivatives, decays much faster than $\eps t^{-1+C\eps}$. Since $u$ and its derivatives is of size $O(\eps t^{-1+C\eps})$, we conclude that $\wt{u}$ offers a good approximation of $u$.

\subsection{Background}
Let us consider  a generalization of the scalar quasilinear wave equation \eqref{qwe} in~$\R^{1+3}_{t,x}$
\begin{equation}\label{nlw}
\square u=F(u,\partial u,\partial^2u).
\end{equation}
The nonlinear term is assumed to be smooth with the Taylor expansion
\begin{equation}\label{nonlinearity}F(u,\partial u,\partial^2u)=\sum a_{\alpha\beta}\partial^\alpha u\partial^\beta u+O(|u|^3+|\partial u|^3+|\partial^2 u|^3).\end{equation}
The sum  is taken over    all multiindices $\alpha,\beta$ with $|\alpha|\leq |\beta|\leq 2$, $|\beta|\geq1$ and $|\alpha|+|\beta|\leq 3$. Besides, the coefficients $a_{\alpha\beta}$'s are all universal constants. 

Since 1980's, several results  on the lifespan of the solutions to the Cauchy problem \eqref{nlw} with initial data \eqref{init} have been proved. For example, John \cite{john,john2} proved that \eqref{nlw} does not necessarily have a global solution; in fact, any nontrivial solution to  $\square u=u_t\Delta u$ or $\square u=u_t^2$ must blow up in finite time. In contrast, in $\R^{1+d}$ with $d\geq 4$, H\"{o}rmander \cite{horm3} proved  the small data global existence for \eqref{nlw}. For arbitrary nonlinearities in three space dimensions, the best result on the lifespan is the almost global existence: the solution exists for $t\leq \exp(c/\eps)$ where $\eps\ll 1$. The almost global existence for \eqref{nlw} was proved by Lindblad \cite{lind3}, and  we also refer to   \cite{johnklai,klai,horm2,horm} for some earlier work.

In contrast to the finite-time blowup in John's examples,  Klainerman \cite{klai3} and Christodoulou \cite{chri} proved that the null condition is sufficient for  small data global existence. The null condition, first introduced by Klainerman \cite{klai2}, states that for each $0\leq m\leq n\leq 2$ with $m+n\leq 3$, we have 
\begin{equation}\label{nullcond}A_{mn}(\omega):=\sum_{|\alpha|=m,|\beta|=n}a_{\alpha\beta}\widehat{\omega}^{\alpha}\widehat{\omega}^{\beta}=0,\hspace{1cm}\text{for all }\widehat{\omega}=(-1,\omega)\in\R\times\mathbb{S}^2.\end{equation}
Equivalently, we assume $A_{mn}\equiv 0$  on the null cone $\{m^{\alpha\beta}\xi_\alpha\xi_\beta=0\}$. The null condition leads to cancellations in the nonlinear terms \eqref{nonlinearity} so that the nonlinear effects of the equations are much weaker than the linear effects. 
However, note that the null condition is not necessary for  small data global existence. For example, the null condition fails for \eqref{qwe} in general, but \eqref{qwe} still has  small data global existence.  We also refer our readers to \cite{tata} for a general introduction on the null condition.

Later, Lindblad and Rodnianski \cite{lindrodn,lindrodn2} introduced the weak null condition. To state the weak null condition, we start with the asymptotic equations first introduced by H\"{o}rmander \cite{horm2,horm,horm3}. We make the ansatz\begin{equation}\label{ansatz}
u(t,x)\approx \eps r^{-1}U(s,q,\omega),\hspace{1cm}r=|x|,\ \omega_i=x_i/r,\ s=\eps\ln(t),\ q=r-t.
\end{equation}
Assuming that $t=r\to\infty$, we substitute this ansatz into \eqref{nlw} and compare the coefficients of  terms of order $\eps^2 t^{-2}$.  Nonrigorously, we can obtain the following asymptotic PDE for $U(s,q,\omega)$
\begin{equation}\label{asypde11}2\partial_s\partial_q U=\sum A_{mn}(\omega)\partial_q^mU\partial_q^nU.\end{equation}
Here $A_{mn}$ is defined in \eqref{nullcond} and the sum is taken over $0\leq m\leq n\leq 2$ with $m+n\leq 3$. We say that the \emph{weak null condition} is satisfied if \eqref{asypde11} has a global solution for all $s\geq 0$ and if the solution and all its derivatives grow at most exponentially in $s$, provided that the initial data decay sufficiently fast in $q$. In the same papers, Lindblad and Rodnianski  conjectured that the weak null condition is sufficient for small data global existence. To the best of the author's knowledge, this conjecture  remains open until today.

There are three remarks about the weak null condition and the corresponding conjecture. First, the  null condition implies the weak null condition. In fact, under the null condition, \eqref{asypde11} becomes $\partial_s\partial_q U=0$. Secondly, though the conjecture   remains open, there are many examples of \eqref{nlw} satisfying the weak null condition and admitting small data global existence at the same time. The equation \eqref{qwe} is one of several such examples: the small data global existence for  \eqref{qwe} has  been proved by Lindblad  \cite{lind}; meanwhile, the asymptotic equation \eqref{asypde11} now becomes
\begin{equation}\label{asypde1}2\partial_s\partial_q U=G(\omega)U\partial_q^2U,\end{equation}
where \eq{\label{gomedefn}G(\omega):=g_0^{\alpha\beta}\widehat{\omega}_\alpha\widehat{\omega}_\beta,\hspace{1cm}g_0^{\alpha\beta}=\frac{d}{du}g^{\alpha\beta}(u)|_{u=0},\ \widehat{\omega}=(-1,\omega)\in\R\times\mathbb{S}^2,}
whose solutions exist globally in $s$ and satisfy the decay requirements, so \eqref{qwe} satisfies the weak null condition. There are also many examples  violating the weak null condition and admitting finite-time blowup at the same time. Two  such examples are $\square u=u_t\Delta u$ and $\square u=u_t^2$: the corresponding asymptotic equations are  $(2\partial_s-U_q\partial_q)U_q=0$ (Burger's equation) and $\partial_s U_q=U_q^2$, respectively, whose solutions are known to blow up in finite time. Thirdly, in  recent years, Keir has made some further  progress. In \cite{keir}, he proved  the small data global existence  for a large class of quasilinear wave equations satisfying the weak null condition, significantly enlarging upon the class of equations for which global existence is known. His proof also applies to \eqref{qwe}. In \cite{keir2},  he proved that if the solutions to the asymptotic system are bounded (given small initial data) and  stable against rapidly decaying perturbations, then the corresponding system of nonlinear wave equations admits  small data global existence.

\subsection{The geometric reduced system}\label{sec1.2an} In \cite{yu2020}, we have constructed a new system of asymptotic equations. Our analysis starts as in H\"{o}rmander's derivation in \cite{horm2,horm,horm3}, but diverges at a key point: the choice of $q$ is different. One may contend from this work that this new system is more accurate than \eqref{asypde1}, in that it both describes the long time evolution and contains full information about it. In addition, if we choose the initial data appropriately, our reduced system will reduce to linear first order ODE's on $\mu$ and $U_q$, so it is easier to solve it than to solve \eqref{asypde1}. 

To derive the new equations, we still make the ansatz  \eqref{ansatz}, but now we replace $q=r-t$ with a solution $q(t,x)$ to the eikonal equation related to \eqref{qwe}\begin{equation}\label{eikeqn}g^{\alpha\beta}(u)\partial_\alpha q\partial_\beta q=0.\end{equation} In other words, $q(t,x)$ is an optical function. We remark that  the eikonal equations have been used in the previous works on the small data global existence for \eqref{qwe}. In  \cite{alin2}, Alinhac followed the method used in  Christodoulou and Klainerman \cite{chriklai}, and adapted the vector fields to the characteristic surfaces, i.e.\ the level surfaces of solutions to the eikonal equations. In \cite{lind}, Lindblad considered the radial eikonal equations when he derived the pointwise bounds of solutions to \eqref{qwe}. When they derived the energy estimates, both Alinhac and Lindblad considered a weight $w(q)$ where $q$ is an approximate solution to the eikonal equation. Their works suggest that the  eikonal equation  plays an important role when we study the long time behavior of solutions to \eqref{qwe}. Moreover, the eikonal equations have also been used in the study of the asymptotic behavior of solutions to the Einstein vacuum equations, an analogue of \eqref{qwe}; we refer our readers to  \cite{chriklai,lind4}. In addition, we also refer to \cite{smittata} where the eikonal equations are used to study the sharp local wellposedness for the nonlinear wave equations.

Since $u$ is  unknown, it is difficult to solve \eqref{eikeqn} directly. Instead,  we introduce a new auxiliary function $\mu=\mu(s,q,\omega)$ such that $q_t-q_r=\mu$. From \eqref{eikeqn}, we can express $q_t+q_r$ in terms of $\mu$ and $U$, and then solve for all partial derivatives of $q$, assuming that all the angular derivatives are negligible. Then from \eqref{qwe}, we can derive the following asymptotic equations for $\mu(s,q,\omega)$ and $U(s,q,\omega)$:
\begin{equation}\label{asypde2}\left\{
\begin{array}{l}\displaystyle
\partial_s\mu=\frac{1}{4}G(\omega)\mu^2 U_q,\\[.2cm]
\displaystyle\partial_sU_q=-\frac{1}{4}G(\omega) \mu U_q^2.
\end{array}\right.
\end{equation} Here $G(\omega)$ is defined by \eqref{gomedefn}. We call this new system of asymptotic equations the \emph{geometric reduced system} since it is related to the geometry of the null cone with respect to the Lorentzian metric $(g_{\alpha\beta})=(g^{\alpha\beta}(u))^{-1}$ instead of the Minkowski metric. For a derivation of \eqref{asypde2}, we refer our readers to Section 3 in \cite{yu2020}, or Chapter 2 in the author's PhD dissertation \cite{yuthesis}. Heuristically,   one expects the solution to the quasilinear wave equation \eqref{qwe} to correspond to an approximate solution to this geometric reduced system, and to be well approximated by an exact solution to the geometric reduced system.

Note that \eqref{asypde2} is a system of two ODE's for $(\mu,U_q)$. In addition, we have $\partial_s(\mu U_q)=0$ for each $(s,q,\omega)$. That is, if the initial data are given by \fm{(\mu,U_q)|_{s=0}(q,\omega)=(A_1,A_2)(q,\omega),} then we have $\mu U_q=A_1\cdot A_2$ at each $(s,q,\omega)$. In this paper, we define a function $A=A(q,\omega)$ for $(q,\omega)\in\R\times\mathbb{S}^2$ by
\fm{A(q,\omega):=-\frac{1}{2}A_1(q,\omega)\cdot A_2(q,\omega),}
and we call the function $A$ a \emph{scattering data} associated to a solution $u$ to the quasilinear wave equation \eqref{qwe}. Now \eqref{asypde2} reduces to a linear system of ODE's
\fm{\left\{\begin{array}{l}\displaystyle \partial_s\mu=-\frac{1}{2}G(\omega)A(q,\omega)\mu,\\[1em] \displaystyle \partial_s U_q=\frac{1}{2}G(\omega)A(q,\omega) U_q,\end{array}\right.} 
whose solutions are given by
\eq{\label{asypde2sol}\left\{\begin{array}{l}\displaystyle \mu(s,q,\omega)=A_1(q,\omega)\exp(-\frac{1}{2}G(\omega)A(q,\omega)s),\\[1em] \displaystyle U_q(s,q,\omega)=A_2(q,\omega)\exp(\frac{1}{2}G(\omega)A(q,\omega) s),\end{array}\right.}
To solve for $U(s,q,\omega)$ uniquely, we assume that \fm{\lim_{q\to-\infty}U(s,q,\omega)=0\qquad \text{or}\qquad\lim_{q\to\infty}U(s,q,\omega)=0,}
depending on which problem we are studying. For instance, in \cite{yu2020}, to guarantee that a solution to \eqref{qwe} is zero inside a certain light cone, we assume that $\lim_{q\to-\infty}U(s,q,\omega)=0$; in this paper, the global solution to \eqref{qwe} has localized initial data, so we assume that $\lim_{q\to\infty}U(s,q,\omega)=0$.

We end this subsection by proposing an alternative definition of the weak null condition. In the discussion above, we define $\mu=q_t-q_r$ and derive a geometric reduced system \eqref{asypde2} for $(\mu,U)(s,q,\omega)$. This method to derive a reduced system should not just work for \eqref{qwe}. A derivation of the geometric reduced systems for a system of general  quasilinear wave equations can be found in Chapter 2, \cite{yuthesis}. We can make the following definition.
\defn{\rm We say that a system of quasilinear wave equations satisfies the \emph{geometric weak null condition}, if for any initial data at $s=0$ decaying sufficiently fast in $q$, we have a global solution to the corresponding geometric reduced system for all $s\geq 0$, and if the solution and all the derivatives grow at most exponentially in
$s$.}\rm\bigskip

It is clear from \eqref{asypde2sol} that  \eqref{qwe} satisfies the geometric weak null condition. The author believes that it is interesting to study to what extent is the geometric weak null condition equivalent to the weak null condition, and whether this geometric weak null condition is sufficient for the global existence of general quasilinear wave equations with small and localized initial data.

\subsection{Modified scattering theory: an overview}
The objective of \cite{yu2020,yuthesis} and this paper is to study the long time dynamics, and more specifically, scattering theory for highly nonlinear dispersive equations. In other words, we would like to provide an accurate description of  asymptotic behavior of the global solutions. For many nonlinear dispersive PDE's, one can establish a linear scattering theory. That is, a global solution to a nonlinear PDE  scatters to a solution to the corresponding linear equation as time goes to infinity. Take the cubic defocusing NLS  \fm{iu_t+\Delta u=u|u|^2\hspace{2em} \text{in } \R^{1+3}_{t,x}} as an example. Its corresponding linear equation is the linear  Schr\"{o}dinger equation (LS) \fm{iw_t+\Delta w=0\hspace{2em} \text{in } \R^{1+3}_{t,x}.} One can prove that for each $u_0\in H^1$, there exists a unique $u_+\in H^1$ such that \fm{\norm{u(t)-w(t)}_{H^1}\to 0\hspace{2em}\text{as }t\to\infty} where $u$ (or $w$) is the global solution to NLS (or LS) with data $u_0$ (or $u_+$). This result is called the \emph{asymptotic completeness}. One can also prove that for each $u_+\in H^1$, there exists a unique $u_0\in H^1$ such that the same conclusion holds. This result is called the \emph{existence of wave operators}, where the wave operator is defined by $\Omega_+u_+=u_0$. We refer to Section 3.6 of  \cite{tao} for this result. Some other nonlinear PDE's have modified scattering instead of linear scattering. That is, each of their global solutions scatters to a suitable modification of a linear solution. Here the modification can be made in more than one way: we can add a phase correction term, an amplitude correction term, or a velocity correction term to the linear solution.  For example, in \cite{ifrim2015global}, when the authors study modified scattering for the cubic 1D NLS, they make use of the following asymptotic approximation:\fm{\widehat{u}(t,\xi)\approx e^{-it\xi^2}W(\xi)e^{i|W(\xi)|^2\ln t}.} That is,  a phase shift term is introduced. For nonlinear wave equations, the modification often corresponds to a change of the geometry of the light cone foliation of the space-time. This point is reflected in the ansatz used in Section \ref{sec1.2an}.

In general, the following steps are taken in order to study  modified scattering. Given a nonlinear dispersive PDE, we hope to identify a good notion of \emph{asymptotic profile} and an associated notion of \emph{scattering data} for the model equation. This can be achieved by introducing some type of asymptotic equations. Like linear scattering, the two main problems in modified scattering theory are as follows:
\begin{enumerate}[1.]

\item \emph{Asymptotic completeness}. Given an exact global solution to the model equation, can we  find the corresponding asymptotic profile and scattering data?

\item \emph{Existence of (modified) wave operators}. Given  an asymptotic profile constructed for a scattering data, can we construct a unique exact global solution to the model equation which matches  the asymptotic profile at infinite time?
\end{enumerate}

There have been only a few previous results on the (modified) scattering for general quasilinear wave equations and the Einstein's equations.  In \cite{dafermos2013scattering}, Dafermos, Holzegel and Rodnianski gave a scattering theory construction of nontrivial black hole solutions to the vacuum Einstein equations.
That is a backward scattering problem in General Relativity. In \cite{lindschl}, Lindblad and Schlue  proved the existence of the wave operators for the semilinear models of Einstein's equations. In \cite{dengpusa}, Deng and Pusateri used the original H\"{o}rmander's asymptotic system \eqref{asypde1} to prove a partial scattering result for \eqref{qwe}.  In their proof, they applied the spacetime resonance method; we refer to \cite{pusashat,pusa2} for some earlier applications of this method to the first order systems of wave equation. Recently, in \cite{yu2020}, by using a new reduced system,  the author proved  the existence of the modified wave operators for \eqref{qwe}.

\subsection{Construction of an optical function}
Let $u=u(t,x)$ be a global solution to \eqref{qwe} and \eqref{init} constructed in Lindblad \cite{lind}. Here we fix a constant $R>0$ such that $\text{supp }(u_0,u_1)\subset\{|x|\leq R\}$, so we have $u\equiv 0$ for $|x|\geq t+R$ by the finite speed of propagation. Our goal in this section is to construct an optical function, i.e.\ a solution to the eikonal equation \eqref{eikeqn}. Here we do not expect to solve \eqref{eikeqn} for all $(t,x)\in \R^{1+3}_{t,x}$. Instead, we solve it in a region $\Omega\subset  \R^{1+3}_{t,x}$ which is defined by 
\fm{\Omega:=\{(t,x):\ t>T_0,\ |x|>(t+T_0)/2+2R\}.}
Here $T_0=\exp(\delta/\eps)$ and $\delta>0$ is a fixed constant independent of $\eps$. We also assign the initial data by setting $q=r-t$ on $\partial\Omega$. It is then clear that $q=r-t$ in $\Omega\cap\{r-t>R\}$, so from now on we focus on the region $\Omega\cap\{r-t<2R\}$.

To construct an optical function, we apply the method of characteristics. In fact, the characteristics for \eqref{eikeqn} are the geodesics with respect to the Lorentzian metric $(g_{\alpha\beta})$ which is the inverse of the matrix $(g^{\alpha\beta}(u))$. Moreover, we only need to study those geodesics emanating from the cone
\fm{H:=\partial\Omega\cap\{t>T_0\}=\{(t,x):\ t>T_0,\ |x|=(t+T_0)/2+2R\}.}

Now we follow the idea in Christodoulou-Klainerman \cite{chriklai}. Fix $T>T_0$ and suppose that the optical function exists in $\Omega_T:=\Omega\cap\{t\leq T,r-t\leq 2R\}$. Then, every point in $\Omega_T$ can be reached by a unique characteristic emanating from $H$. We first define a null frame $\{e_k\}_{k=1}^4$ in $\Omega_T$, such that $e_4$ is tangent to the unique characteristic passing throught that point.   We then define the second fundamental form of the time slices of the null cones:
\fm{\chi_{ab}:=\lra{D_{e_a}e_4,e_b},\hspace{2em}a,b\in\{1,2\}.}
Here $D$ is the Levi-Civita connection associated to  the Lorentzian metric $(g_{\alpha\beta})$, and  $\lra{\cdot,\cdot}$ is the bilinear form asscociated to the  metric $(g_{\alpha\beta})$. We now use a continuity argument. Suppose that in $\Omega_T$ we have
\eq{\label{s1.4chiest} \max_{a,b=1,2}|\chi_{ab}-\delta_{ab}r^{-1}|\leq At^{-2+B\eps}.}
The positive constants $A$ and $B$ are both independent of $\eps$ and $T$. Our goal is to prove that \eqref{s1.4chiest} holds with $A$ replaced by $A/2$. It follows that  $\tr\chi:=\chi_{11}+\chi_{22}$, sometimes called the \emph{null mean curvature}\footnote{We will briefly explain the geometric meaning of $\tr\chi$ in Section \ref{s3}.} of the level sets of $q$, is positive everywhere, and that the characteristics emanating from $H$ will not intersect with each other. This allows us to extend the optical function to  $\Omega_{T+\epsilon}$ for a small $\epsilon>0$, such that \eqref{s1.4chiest} holds everywhere in $\Omega_{T+\epsilon}$. We conclude from this continuity argument that the optical function exists everywhere in $\Omega$.

In order to prove that \eqref{s1.4chiest} holds with $A$ replaced by $A/2$, we make use of the Raychaudhuri equation
\fm{ e_4(\chi_{ab})&=-\sum_{c=1,2}\chi_{ac}\chi_{cb}+\Gamma_{\alpha\beta}^0e_4^\alpha e_4^\beta\chi_{ab}+\lra{R(e_4,e_a)e_4,e_b},}
which describes the evolution of $\chi$ along the null geodesics foliating the light cones.
In this equation, $\Gamma_{**}^*$'s are the Christoffel symbols, and $\lra{R(X,Y)Z,W}$ is the curvature tensor, both with respect to the Lorentzian metric $(g_{\alpha\beta})$. Note that we have a decomposition
\fm{\lra{R(e_4,e_a)e_4,e_b}=e_4(f_1)+f_2} where $f_1=O(\eps t^{-2+C\eps})$ and $f_2=O(\eps t^{-3+C\eps})$; see Lemma \ref{l3curest} for a more accurate statement. We also refer our readers to Corollary 5.9 in \cite{smittata} for a similar decomposition of curvature tensors. Moreover, it follows from \eqref{qwe} that \fm{|e_4(e_3(u))+r^{-1}e_3(u)|\lesssim\eps At^{-3+B\eps},\hspace{2em}|e_4(e_3(u))|\lesssim \eps t^{-2}.} Combining all these estimates and the Gronwall's inequality, we are able to prove \eqref{s1.4chiest} with $A$ replaced by $A/2$.

So far, we have constructed a global optical function $q=q(t,x)$ in $\Omega$ which  is $C^2$ by the method of characteristics. In fact, the optical function $q=q(t,x)$ is smooth\footnote{See Section \ref{sec2.5}. In particular, a smooth function may not be $C^\infty$ in this paper. } in $\Omega$ in the followings sense: for each integer $N\geq 2$, there exists $\eps_N>0$ such that $q$ is a $C^N$ function in $\Omega$ for each $0<\eps<\eps_N$. Moreover, if $Z$ is one of the commuting vector fields: translations $\partial_\alpha$, scaling $t\partial_t+r\partial_r$, rotations $x_i\partial_j-x_j\partial_i$ and Lorentz boosts $x_i\partial_t+t\partial_i$, then in $\Omega$ we have $Z^Iq=O(\lra{q}t^{C\eps})$ and $Z^I\Omega_{ij}q=O(t^{C\eps})$ for each multiindex $I$ and $\eps\ll_I1$. To prove these estimates, we introduce the commutator coefficients $\{\xi_{k_1k_2}^l\}_{1\leq k_1,k_2,l\leq 4}$ for which we have $[e_{k_1},e_{k_2}]=\xi_{k_1k_2}^le_l$. We also introduce a weighted null frame \fm{(V_1,V_2,V_3,V_4):=(re_1,re_2,(3R-r+t)e_3,te_4)}
which combines the advantages of a usual null frame $\{e_k\}$ and the commuting vector fields $Z$'s. By computing $e_4(V^I\xi_{k_1k_2}^l)$ for each multiindex $I$ and applying the Gronwall's inequality, we are able to obtain several estimates for $V^I(\xi_{k_1k_2}^l)$; see Proposition \ref{prop5.5}. These estimates for $\xi$ then imply the estimates for $q$, so we finish the proof.

We finally remark that  the map \fm{\Omega\to[0,\infty)\times\R\times\mathbb{S}^2:\hspace{2em}(t,x)\mapsto (\eps\ln t-\delta,q(t,x),x/|x|):=(s,q,\omega)}
is an injective smooth function with a smooth inverse.  This is because $q_r>0$ everywhere in $\Omega$. Thus, a smooth function $F=F(t,x)$ induces a smooth function $F=F(s,q,\omega)$ and vice versa.

\subsection{The asymptotic equations and the scattering data}\label{sec1.5}
For each $(t,x)\in\Omega$, we define \fm{\mu(t,x):=(q_t-q_r)(t,x),\hspace{2em}U(t,x):=\eps^{-1}ru(t,x).}
We then obtain two smooth functions $\mu(s,q,\omega)$ and $U(s,q,\omega)$ as discussed at the end of the previous subsection.

To state the results in this subsection, we introduce a new notation $\mathfrak{R}_{s,p}$ for each $s,p\in\R$. For a function $F=F(t,x)$ defined in $\Omega\cap\{r-t<2R\}$, we write $F=\mathfrak{R}_{s,p}$ if for each integer $N\geq 1$ and for each $\eps\ll_N1$, we have
\fm{\sum_{|I|\leq N}|V^I(F)|\lesssim t^{s+C\eps}\lra{q}^{p},\hspace{2em}\forall(t,x)\in\Omega\cap\{r-t<2R\}.}
Here recall that $\{V_*\}$ is the weighted null frame.

By the chain rule, we have 
\fm{\partial_s=\eps^{-1}t(\partial_t-q_tq_r^{-1}\partial_r),\hspace{2em}\partial_q=q_r^{-1}\partial_r,\hspace{2em}\partial_{\omega_i}=r(\partial_i-q_iq_r^{-1}\partial_r).}
Then we can express $(\partial_s,\partial_q,\partial_\omega)$ in terms of the weighted null frame $\{V_*\}$. In fact, we have
\fm{\partial_s&=\sum_{a}\eps^{-1}\mathfrak{R}_{-1,0}V_a+(\eps^{-1}+\mathfrak{R}_{-1,0})V_4,\hspace{2em}
\partial_q=\sum_{k}\mathfrak{R}_{0,-1}V_k,\\
\partial_{\omega_i}&=\sum_{k\neq 3}\mathfrak{R}_{-1,0}V_k+\sum_ae_a^iV_a=\sum_{k\neq 3}\mathfrak{R}_{0,0}V_k.}
Meanwhile, from \eqref{qwe} and $e_4(e_3(q))=-\Gamma_{\alpha\beta}^0e_4^\alpha e_4^\beta e_3(q)$, we can show that
\fm{e_4(e_3(u))+r^{-1}e_3(u)=\eps\mathfrak{R}_{-3,0},\hspace{2em}e_4(e_3(q))=-\frac{1}{4}e_3(u)G(\omega)e_3(q)+\eps\mathfrak{R}_{-2,0}.}
Combine these estimates, and we obtain that 
\eq{\label{sec1.5f1}\left\{\begin{array}{l}\displaystyle
\partial_s\mu=\frac{1}{4}G(\omega)\mu^2 U_q+\eps^{-1}\mathfrak{R}_{-1,0},\\[1em]
\displaystyle\partial_sU_q=-\frac{1}{4}G(\omega)\mu U_q^2+\eps^{-1}\mathfrak{R}_{-1,0}.
\end{array}\right.}
That is, $(\mu,U)(s,q,\omega)$ is an approximate solution to the geometric reduced system \eqref{asypde2}.

Next, we note from \eqref{sec1.5f1} that $\partial_s(\mu U_q)=O(\eps^{-1}t^{-1+C\eps})$.  By integrating the remainder term $\eps^{-1}t^{-1+C\eps}$ (viewed as a function of $s$) with respect to $s$, we can show that $\{(\mu U_q)(s,q,\omega)\}_{s}$ is uniformly Cauchy for each $(q,\omega)\in\R\times\mathbb{S}^2$. Thus, the limit 
\fm{A(q,\omega):=-\frac{1}{2}\lim_{s\to\infty}(\mu U_q)(s,q,\omega)}exists and the convergence is uniform in $(q,\omega)$. This function $A$ is then the \emph{scattering data} in the asymptotic completeness problem.

Similarly, we can show that for each $m$ and $n$, the limit
\fm{A_{m,n}(q,\omega):=-\frac{1}{2}\lim_{s\to\infty}(\lra{q}\partial_q)^m\partial_\omega^n(\mu U_q)(s,q,\omega)}
exists and the convergence is uniform in $(q,\omega)$. The uniform convergences of these limits imply that the scattering data $A$ is smooth,
\fm{(\lra{q}\partial_q)^m\partial_\omega^nA(q,\omega)=A_{m,n}(q,\omega).} 

Following the same method, we can define 
\fm{A_1(q,\omega)&:=\lim_{s\to\infty}\exp(\frac{1}{2}G(\omega)A(q,\omega)s)\mu(s,q,\omega),\\
 A_2(q,\omega)&:=\lim_{s\to\infty}\exp(-\frac{1}{2}G(\omega)A(q,\omega)s)U_q(s,q,\omega).}
Both of these limits exist and have derivatives of any order with respect to $q$ and $\omega$, as long as $\eps$ is sufficiently small. It is clear that $A_1A_2\equiv -2A$, so we obtain an exact solution to the  geometric reduced system \eqref{asypde2}:
\eq{\label{sec1.5f2}\left\{\begin{array}{l}\displaystyle \wt{\mu}(s,q,\omega)=A_1(q,\omega)\exp(-\frac{1}{2}G(\omega)A(q,\omega)s),\\[1em] \displaystyle \wt{U}_q(s,q,\omega)=A_2(q,\omega)\exp(\frac{1}{2}G(\omega)A(q,\omega) s),\end{array}\right.}
By assuming $\lim_{q\to\infty}\wt{U}(s,q,\omega)=0$, we obtain a unique function $\wt{U}=\wt{U}(s,q,\omega)$.  By the definition of $(A,A_1,A_2)$, we expect the $(\mu-\wt{\mu},U-\wt{U})$, along with their derivatives with respect to $(s,q,\omega)$ of any order, decays faster than $\mu$ and $U$. We refer our readers to Proposition \ref{props5} for a complete list of estimates.

As defined, the scattering data $A$ depends on the initialization of the optical function $q$. In Section \ref{ss5.7}, see Proposition \ref{prop5.8gi}, we resolve this ambiguity and show a gauge independence result, which  states that the scattering data is independent of  the choice of $q$ in a precise  sense.

\subsection{An approximation}\label{sec1.6}
In the previous subsection, we have discussed how to obtain an exact solution \eqref{sec1.5f2} to the geometric reduced system \eqref{asypde2}. Our final objective is to show that this exact solution gives a good approximation of the exact solution $u$ to \eqref{qwe}.

We first solve
\fm{\wt{q}_t-\wt{q}_r=\wt{\mu}(\eps\ln(t)-\delta,\wt{q}(t,x),\omega)\quad\text{in }\Omega\cap\{r-t<2R\};\hspace{2em}\wt{q}=r-t\quad\text{when } r-t\geq 2R}
and set \fm{\wt{u}(t,x)=\eps r^{-1}\wt{U}(\eps\ln(t)-\delta,\wt{q}(t,x),\omega)\quad\text{in }\Omega\cap\{r-t<2R\}.}
Then, we can prove that $\wt{u}$ is an approximate solution to \eqref{qwe} in the following sense: for each integer $N\geq 1$ and $\eps\ll_N1$, we have
\eq{\label{sec1.6f1}\sum_{|I|\leq N}|Z^I(g^{\alpha\beta}(\wt{u})\partial_\alpha\partial_\beta\wt{u})|\lesssim \eps t^{-3+C\eps},\hspace{2em}\text{in }\Omega\cap\{r-t<2R\}.}
Here we denote by $Z$  any of the commuting vector fields: translations $\partial_\alpha$, scaling $t\partial_t+r\partial_r$, rotations $x_i\partial_j-x_j\partial_i$ and Lorentz boosts $x_i\partial_t+t\partial_i$. To make our proof simpler, we introduce a new function $F=F(q,\omega)$  such that $F_q=-2/A_1$. It can be shown that $q\mapsto F(q,\omega)$ has an inverse $q\mapsto \hat{F}(q,\omega)$. Now we define $\hat{A}(q,\omega):=A(\hat{F}(q,\omega),\omega)$ and define $(\hat{\mu},\hat{U}_q)(s,q,\omega)$ by replacing $(A_1,A_2,A)$ in \eqref{sec1.5f2} with $(-2,\hat{A},\hat{A})$. Then, $\hat{q}(t,x):=F(\wt{q}(t,x),\omega)$ is a solution to \fm{\hat{q}_t-\hat{q}_r=\hat{\mu}(\eps\ln t-\delta,\hat{q}(t,x),\omega)\quad\text{in }\Omega\cap\{r-t<2R\};\hspace{2em}\hat{q}=r-t\quad\text{when } r-t\geq 2R.}In addition, we have\fm{\wt{U}(\eps\ln(t)-\delta,\wt{q}(t,x),\omega)=\hat{U}(\eps\ln(t)-\delta,\hat{q}(t,x),\omega).} We can now follow the proof in Section 4 of \cite{yu2020} to prove \eqref{sec1.6f1}.

In order to estimate $u-\wt{u}$, we set $p(t,x):=F(q(t,x),\omega)-\hat{q}(t,x)$ in $\Omega$. We claim that, for each fixed $\gamma\in(0,1)$, an integer $N\geq 1$, and for each $\eps\ll_{\gamma,N}1$, whenever $(t,x)\in\Omega$ such that $|r-t|\lesssim t^{\gamma}$, we have $|Z^Ip(t,x)|\lesssim t^{-1+C\eps}\lra{r-t}$ for each $|I|\leq N$. To show this claim, we compute $p_t-p_r$ and apply a continuity argument. This claim then implies that, under the same assumptions on $\gamma$, $N$ and $\eps$, whenever $(t,x)\in\Omega$ such that $|r-t|\lesssim t^{\gamma}$, we have $|Z^I(u-\wt{u})(t,x)|\lesssim \eps t^{-2+C\eps}\lra{r-t}$ for each $|I|\leq N$. Recall from Lindblad \cite{lind} that we only have $Z^Iu=O(\eps t^{-1+C\eps})$, so $\wt{u}$ provides a good approximation of $u$.

\subsection{The main theorem}

We now state the main theorem which summarizes the outcome of the sequence of steps described in the  previous subsections. In this theorem, we say that a function $f=f(t,x)$ is smooth if for each large integer $N$, $f$ is $C^N$ whenever $\eps\ll_N 1$. See Section \ref{sec2.5} for details.

\thm{\label{mthm} Let $u$ be a smooth solution to the Cauchy problem \eqref{qwe} and \eqref{init}. Fix a constant $R>0$ such that $\text{supp }(u_0,u_1)\subset\{|x|\leq R\}$, so $u\equiv 0$ for $|x|\geq t+R$ by the finite speed of propagation. 
Set $T_0:=\exp(\delta/\eps)$ for a fixed constant $\delta>0$.  Then we have
\begin{enumerate}[\rm a)]
\item There exists a smooth solution to the eikonal equation 
\fm{ g^{\alpha\beta}(u)\partial_\alpha q\partial_\beta q=0\text{ in }\Omega;\hspace{1cm}q=|x|-t\text{ on }\partial\Omega.}
Here the region $\Omega\subset\R^{1+3}_{t,x}$ is defined by
\fm{ \Omega:=\{(t,x):\ t>T_0,\ |x|>(t+T_0)/2+2R\}.}
In $\Omega$, for each multiindex $I$ we have
\fm{|Z^Iq|\lesssim \lra{q}t^{C\eps},\hspace{2em}\sum_{1\leq i,j\leq 3}|Z^I\Omega_{ij}q|\lesssim t^{C\eps}.}
Moreover, the map \fm{\Omega\to[0,\infty)\times\R\times\mathbb{S}^2:\hspace{2em}(t,x)\mapsto (\eps\ln t-\delta,q(t,x),x/|x|)}
is an injective smooth function with a smooth inverse.  Thus, a smooth function $F=F(t,x)$ induces a smooth function $F=F(s,q,\omega)$ and vice versa.

\item In $\Omega$, we set $(\mu,U)(t,x):=(q_t-q_r,\eps^{-1}ru)(t,x)$ which induces a smooth function $(\mu,U)(s,q,\omega)$. Then, $(\mu,U)(s,q,\omega)$ is an approximate solution to the geometric reduced system \eqref{asypde2} in the sense that \fm{\left\{\begin{array}{l}\displaystyle
\partial_s\mu=\frac{1}{4}G(\omega)\mu^2 U_q+\eps^{-1}\mathfrak{R}_{-1,0},\\[1em]
\displaystyle\partial_sU_q=-\frac{1}{4}G(\omega)\mu U_q^2+\eps^{-1}\mathfrak{R}_{-1,0}.
\end{array}\right.}
Here the notation $\mathfrak{R}_{*,*}$ has been defined in Section \ref{sec1.5}. In addition, the following three limits exist for all $(q,\omega)\in\R\times\mathbb{S}^2$:
\fm{\left\{\begin{array}{l}
\displaystyle A(q,\omega):=-\frac{1}{2}\lim_{s\to\infty}(\mu U_q)(s,q,\omega),\\
\displaystyle A_1(q,\omega):=\lim_{s\to\infty}\exp(\frac{1}{2}G(\omega)A(q,\omega)s)\mu(s,q,\omega),\\
\displaystyle A_2(q,\omega):=\lim_{s\to\infty}\exp(-\frac{1}{2}G(\omega)A(q,\omega)s)U_q(s,q,\omega).
\end{array}\right.}
All of them are smooth functions of $(q,\omega)$ for $\eps\ll 1$, and we have $A_1A_2\equiv -2A$. By setting
\fm{\left\{\begin{array}{l}
\displaystyle\wt{\mu}(s,q,\omega):=A_1\exp(-\frac{1}{2}GAs),\\[1em]
\displaystyle \wt{U}_q(s,q,\omega):=A_2\exp(\frac{1}{2}GAs).
\end{array}\right.}
we obtain an exact solution to our reduced system \eqref{asypde2}.
\item We define $\wt{u}=\wt{u}(t,x)$ as in Section \ref{sec1.6}. Then the function $\wt{u}=\wt{u}(t,x)$ is an approximate solution to \eqref{qwe} in the following sense:
\fm{|Z^I(g^{\alpha\beta}(\wt{u})\partial_{\alpha}\partial_\beta\wt{u})(t,x)|\lesssim \eps t^{-3+C\eps},\hspace{2em}\forall (t,x)\in\Omega,\ \forall I.}
Moreover, we fix a constant $0<\gamma<1$ and a large integer $N$. Then, for $\eps\ll_{\gamma,N} 1$, at each $(t,x)\in\Omega$ such that $|r-t|\lesssim  t^{\gamma}$, we have \fm{|Z^I(u-\wt{u})|\lesssim_\gamma \eps t^{-2+C\eps}\lra{r-t},\qquad\forall|I|\leq N.}
\end{enumerate}}\rm

\rmkn{\rm Because of the special definition of smoothness in this paper, we emphasize that our main theorem only holds in the following sense: for each large integer $N$, there exists a sufficiently small constant $\eps_N>0$ depending on $N$ and the functions $u_0,u_1\in C_c^\infty(\R^3)$ given in \eqref{init}, such that the conclusions in Theorem \ref{mthm} hold for all $0<\eps<\eps_N$, with all ``smooth'' replaced by ``$C^N$'' in the statement.}

\rmkn{\label{mthm.1}\rm We expect the results above are gauge independent. That is, the scattering data $A=A(q,\omega)$ is independent of the choice of the optical function $q=q(t,x)$ in some suitable sense. In fact, we choose the region $\Omega$ in a way that $t\sim r$ in $\Omega\cap\{r-t<2R\}$, that $t\geq T_0=\exp(\delta/\eps)$ in $\overline{\Omega}$, and that $u\equiv 0$ in $\partial\Omega\cap\{t=T_0\}$. The proof in this paper is expected to work if we start with a different region $\Omega$ with these three properties hold. For example, we can replace the definition of $\Omega$ with \fm{\Omega=\Omega_{\kappa,\delta}:=\{(t,x):\ t>\exp(\delta/\eps),\ |x|-\exp(\delta/\eps)-2R>\kappa(t-\exp(\delta/\eps))\}}
for some fixed constants $\delta>0$ and $0<\kappa<1$. For a different choice of $(\kappa,\delta)$, we do not expect to get the same scattering data. However, Proposition \ref{prop5.8gi} states that the scattering data associated to different regions $\Omega_{\kappa,\delta}$ are in fact related to each other in some sense. }

\rmkn{\rm In our construction, we fix a parameter $\delta>0$ and solve the eikonal equation in a region contained in $\{t>\exp(\delta/\eps)\}$. In fact,  the proof in this paper is expected to work for each fixed $\delta>0$. However, we do not simply set $\delta=1$ here. Instead, we choose a sufficiently small $\delta>0$ which depends on the pair $(u_0,u_1)$, such that the nonlinear effects of \eqref{qwe} are negligible until we reach the time $\exp(\delta/\eps)$. For example, we can set $\delta$ to be the small constant $c$ in the almost global existence result.}

\rmkn{\rm The part c) of the main theorem is an approximation result near the light cone $r=t$.\footnote{We can also say that the main theorem is  an asymptotic result near the null infinity. In contrast, the result in \cite{cklind} is an asymptotic result near the timelike infinity. In that paper, the authors consider some limits of the form $\lim_{t\to\infty}tA_\mu(t,ty)$ for some $|y|<1$.} Outside the light cone, we have $u\equiv 0$ whenever  $r-t\geq R$ because of the finite speed of propagation. It is thus natural to ask whether we also have an approximation result in the interior region away from the light cone. For example, can we find an approximate solution $\wt{u}$ such that $u-\wt{u}$, along with its derivatives, decays faster than $\eps t^{-1+C\eps}$ whenever  $r<t- Ct^\gamma$, where $\gamma\in(0,1)$ is a fixed constant?

We first remark that $\wt{u}$ constructed in Section \ref{sec1.6} is not a good candidate in this case. One reason is that $\wt{u}$ is only defined in $\Omega$. Even in the region where it is defined, it does not give a good approximation near $\partial \Omega$. Note that  part c) of the main theorem implies that \fm{|Z^I(u-\wt{u})|\lesssim_\gamma \eps t^{-2+\gamma+C\eps},\qquad\text{whenever } (t,x)\in\Omega,\ r-t=- t^\gamma/4,\ \gamma\in(0,1).}
If we set $\gamma=1$ on the right hand side of this estimate, we get $\eps t^{-1+C\eps}$ which is the decay rate for the solution $u$ itself. Thus, $\wt{u}$ does not approximate $u$ very well away from the light cone in $\Omega$. Intuitively, this is because the geometric reduced system and the H\"{o}rmander's asymptotic PDE's are derived under the assumption $t\approx r\to\infty$.

By the pointwise estimates for $Z^Iu$ and Lemma \ref{l2.1} below, we have
$|\partial^k Z^Iu|\lesssim \eps t^{-1+C\eps}\lra{r-t}^{-k}$.
As a result, if $|r-t|\gtrsim t^{-\gamma}$ for some $\gamma>0$, we have \fm{|\partial^k Z^Iu|\lesssim \eps t^{-1-k\gamma+C\eps},\qquad \forall (k,I).}
So $\partial^kZ^Iu$ has a  decay rate better than $\eps t^{-1+C\eps}$ if $k>0$.  The case $k=0$ seems to be more complicated, since it is unclear what would be a good approximation for $Z^Iu$ in the interior region. We will not discuss this topic in this paper and we refer our readers to \cite{cklind} which is a paper on the asymptotic behavior of the Maxwell-Klein-Gordon system in this direction.
}

\rmkn{\rm We compare the results in this paper with those in Deng-Pusateri \cite{dengpusa}. First, the approximation result (i.e.\ part c) in Theorem \ref{mthm}) is better than that in \cite{dengpusa} (i.e.\ Theorem 2.3). This suggests that the geometric reduced system \eqref{asypde2} gives a more accurate description of the global solutions to \eqref{qwe} than the H\"{o}rmander's asymptotic PDE  \eqref{asypde1} does. Besides, the proof in this paper relies on the null geometry, and we only use estimates in  physical space. In contrast,  the authors in \cite{dengpusa} made use of the spacetime resonance method which  relies on estimates in frequency space.}
\rm

\subsection{Acknowledgement}
The author would like to thank his PhD advisor, Daniel Tataru, for suggesting this problem and for many helpful discussions. The author would like to thank Sung-Jin Oh for some helpful discussions on the optical function. The author is also grateful to the anonymous reviewers for their valuable comments and suggestions on this paper.

This research was partially supported by a James H. Simons Fellowship and by the NSF grant DMS-1800294. Most of the material here overlaps with the author's PhD dissertation~\cite{yuthesis}.

\section{Preliminaries}\label{sec2}

\subsection{Notations} We use $C$ to denote universal positive constants. We write $A\lesssim B$ or $A=O(B)$ if $|A|\leq CB$ for some  $C>0$. We write $A\sim B$ if $A\lesssim B$ and $B\lesssim A$. We  use $C_{v}$ or $\lesssim_v$ if we want to emphasize that the constant depends on a parameter $v$. We make an additional convention that the constants $C$ are always  independent of $\eps$; that is, we would never write $C_\eps$ or $\lesssim_\eps$ in this paper. The values of all constants in this paper may vary from line to line.

In this paper, we always assume that $\eps\ll 1$ which means $0<\eps<\eps_0$ for some sufficiently small constant $\eps_0<1$. Again,  we write $\eps\ll_v 1$ if we want to emphasize that $\eps_0$ depends on a parameter $v$.

Unless specified otherwise, we always assume that the Latin indices $i,j,l$ take values in $\{1,2,3\}$ and the Greek indices $\alpha,\beta$ take values in $\{0,1,2,3\}$.  We also assume $a,b\in \{1,2\}$ when we study the null frame introduced in Section \ref{s3}. We use subscript to denote partial derivatives, unless specified otherwise. For example, $u_{\alpha\beta}=\partial_\alpha\partial_\beta u$, $q_r=\partial_rq=\sum_i \omega_i\partial_iq$, $A_q=\partial_qA$ and etc. For a fixed integer $k\geq 0$, we  use $\partial^k$ to denote either a specific partial derivative  of order $k$, or the collection of partial derivatives  of order $k$.

To prevent confusion, we will only use $\partial_\omega$ to denote the angular derivatives  under the coordinate $(s,q,\omega)$, and will never use it under the coordinate $(t,r,\omega)$. For a fixed integer $k\geq 0$, we will use $\partial_\omega^k$ to denote  either a specific angular derivative  of order $k$, or the collection of all angular derivatives of order $k$.

\subsection{Commuting vector fields}

We denote by $Z$ any of the following vector fields:
\begin{equation}\label{vf} \partial_\alpha,\ \alpha=0,1,2,3;\ S=t\partial t+r\partial_r;\ \Omega_{ij}=x_i\partial_j-x_j\partial_i,\ 1\leq i<j\leq 3;\ \Omega_{0i}=x_i\partial_t+t\partial_i,\ i=1,2,3.\end{equation}We write these vector fields as $Z_1,Z_2,\dots,Z_{11}$, respectively. For any multiindex $I=(i_1,\dots,i_m)$ with length $m=|I|$ such that $1\leq i_*\leq 11$, we set $Z^I=Z_{i_1}Z_{i_2}\cdots Z_{i_m}$. Then we have the Leibniz's rule
\begin{equation}Z^I(fg)=\sum_{|J|+|K|=|I|}C^I_{JK}Z^JfZ^Kg,\hspace{1cm}\text{where $C_{JK}^I$ are constants.}\end{equation}

We have the following commutation properties.
\begin{equation}[S,\square]=-2\square,\hspace{1cm}[Z,\square]=0\text{ for other $Z$};\end{equation}
\begin{equation}[Z_1,Z_2]=\sum_{|I|=1} C_{Z_1,Z_2,I}Z^I,\hspace{1cm}\text{where $C_{Z_1,Z_2,I}$ are constants};\end{equation}
\begin{equation}\label{comf3}[Z,\partial_\alpha]=\sum_\beta C_{Z,\alpha\beta}\partial_\beta,\hspace{2em} \text{where $C_{Z,\alpha\beta}$ are constants}.\end{equation}

In Section \ref{sec7}, we need the following lemma related to the commuting vector fields. Here we use $f_0$ to denote an arbitrary  polynomial of $\{Z^I\omega\}$. It is then clear that $Z^If_0=f_0$ for each $I$.

\lem{\label{lemtrcom}For each multiindex $I$ and each function $F$,  we have
\eq{\label{lemtrcomf1}(\partial_t-\partial_r)Z^IF&=Z^I(F_t-F_r)+\sum_{|J|<|I|}[f_0Z^J(F_t-F_r)+\sum_if_0(\partial_i+\omega_i\partial_t)Z^JF].}
Besides, for each $1\leq k<k'\leq 3$, we have
\eq{\label{lemtrcomf2}(\partial_t-\partial_r)Z^I\Omega_{kk'}F&=Z^I\Omega_{kk'}(F_t-F_r)+\sum_{|J|<|I|}[f_0Z^J\Omega_{kk'}(F_t-F_r)+\sum_if_0(\partial_i+\omega_i\partial_t)Z^J\Omega_{kk'}F].}}
\begin{proof}First, note that $[\partial_t-\partial_r,Z]=f_0\cdot \partial$ and $\partial=f_0(\partial_t-\partial_r)+\sum_i f_0(\partial_i+\omega_i\partial_t)$. We now prove \eqref{lemtrcomf1} by induction on $|I|$. If $|I|=0$, there is nothing to prove. Now suppose we have proved \eqref{lemtrcomf1} for each $|I|<n$. Now we fix a multiindex $I$ with $|I|=n>0$. Then, by writing $Z^I=ZZ^{I'}$, we have \fm{&\hspace{1.5em}(\partial_t-\partial_r)Z^IF=[\partial_t-\partial_r,Z]Z^{I'}F+Z((\partial_t-\partial_r)Z^{I'}F)\\
&=f_0\cdot\partial Z^{I'}F+Z(Z^{I'}(F_t-F_r)+\sum_{|J|<n-1}[f_0Z^J(F_t-F_r)+\sum_if_0(\partial_i+\omega_i\partial_t)Z^JF]\\
&=f_0(f_0(\partial_t-\partial_r)+\sum_jf_0(\partial_j+\omega_j\partial_t))Z^{I'}F+Z^I(F_t-F_r)\\
&\hspace{1em}+\sum_{|J|<n-1}Z[f_0Z^J(F_t-F_r)+\sum_if_0(\partial_i+\omega_i\partial_t)Z^JF]\\
&=f_0(\partial_t-\partial_r)Z^{I'}F+\sum_jf_0(\partial_j+\omega_j\partial_t)Z^{I'}F+Z^I(F_t-F_r)\\
&\hspace{1em}+\sum_{|J|<n-1}[(Zf_0)Z^J(F_t-F_r)+\sum_i(Zf_0)(\partial_i+\omega_i\partial_t)Z^JF]\\&\hspace{1em}+\sum_{|J|<n-1}[f_0ZZ^J(F_t-F_r)+\sum_if_0Z(\partial_i+\omega_i\partial_t)Z^JF].}
In the second  equality, we can apply \eqref{lemtrcomf1} by the induction hypotheses. Moreover, we note that $[\partial_i+\omega_i\partial_t,Z]=f_0\cdot\partial$, so
\fm{Z(\partial_i+\omega_i\partial_t)Z^JF&=(\partial_i+\omega_i\partial_t)ZZ^JF+f_0\cdot\partial Z^JF\\
&=(\partial_i+\omega_i\partial_t)ZZ^JF+f_0(\partial_t-\partial_r) Z^JF+\sum_jf_0(\partial_j+\omega_j\partial_t)Z^JF.}
Now \eqref{lemtrcomf1} follows from the induction hypotheses and the computations above.

To prove \eqref{lemtrcomf2}, we replace $F$ with $\Omega_{kk'}F$ in \eqref{lemtrcomf1} and note that \fm{\ [\partial_t-\partial_r,\Omega_{kk'}]&=-\partial_r(x_k)\partial_{k'}+\partial_r(x_{k'})\partial_k+\sum_i\Omega_{kk'}(\omega_i)\partial_i\\
&=-\omega_k\partial_{k'}+\omega_{k'}\partial_k+\sum_i \omega_k(\delta_{ik'}-\omega_i\omega_{k'})\partial_i-\sum_i \omega_{k'}(\delta_{ik}-\omega_i\omega_{k})\partial_i=0.}
Now, \eqref{lemtrcomf2} is obvious.
\end{proof}\rm

\subsection{Several pointwise bounds}
We have the pointwise estimates for partial derivatives.
\begin{lem}\label{l2.1}
For any function $\phi$, we have 
\begin{equation}|\partial^k\phi|\leq C\lra{t-r}^{-k}\sum_{|I|\leq k}|Z^I\phi|,\hspace{1cm}\forall k\geq 1,\end{equation}
and
\begin{equation}|(\partial_t+\partial_r)\phi|+|(\partial_i-\omega_i\partial_r)\phi|\leq C\lra{t+r}^{-1}\sum_{|I|=1}|Z^I\phi|.\end{equation}
\end{lem}\rm

In addition, we have the Klainerman-Sobolev inequality.
\begin{prop} For $\phi\in C^\infty(\R^{1+3})$ which vanishes for large $|x|$, we have
\begin{equation}(1+t+|x|)(1+|t-|x||)^{1/2}|\phi(t,x)|\leq C\sum_{|I|\leq 2}\norm{Z^I\phi(t,\cdot)}_{L^2(\R^3)}.\end{equation}
\end{prop}

We also state the Gronwall's inequality.
\prop{Suppose $A,E,r$ are bounded functions from $[a,b]$ to $[0,\infty)$. Suppose that $E$ is increasing. If \fm{A(t)\leq E(t)+\int_a^br(s)A(s)\ ds,\hspace{1cm} \forall t\in[a,b],}then\fm{A(t)\leq E(t)\exp(\int_a^t r(s)\ ds),\hspace{1cm} \forall t\in[a,b].}}\rm
The proofs of these results are standard. See, for example, \cite{lind,sogg,horm} for the proofs.

\subsection{A key theorem and a convention}\label{sec2.5}
This paper is based on the following global existence result.

\thm[Lindblad  \cite{lind}]{\label{lindthm}Fix a large integer $N\gg1$. Then, for $\eps\ll_N1$, the Cauchy problem \eqref{qwe} with the initial data \eqref{init} has a global $C^N$ solution $u=u(t,x)$ for all $t\geq 0$. Moreover, we have  pointwise decays: $Z^Iu=O_I(\eps\lra{t}^{-1+C_I\eps})$ for each multiindex $I$ such that $|I|\leq N$. Moreover, we have $\partial u=O(\eps \lra{t}^{-1})$.}\rm
\bigskip

Most of the functions in this paper have similar properties. That is, they depend on a small parameter $\eps$, and they are $C^N$ for any large integer $N$ as long as $\eps\ll_N1$. For convenience, we make the following definition.

\defn{\rm Fix a function $f=f_\eps(t,x)$ which depends on a small parameter $\eps$. In this paper, we say that $f$ is \emph{smooth}, if for each large integer $N$, $f$ is $C^N$ whenever $\eps\ll_N 1$.

Following the same spirits, we say that all derivatives of a function satisfy some properties, if for each large integer $N$, all its derivatives of order $\leq N$ exist and satisfy such properties whenever $\eps\ll_N 1$.}\rm
\bigskip

We remark that under this definition, a smooth function does not need to be a $C^\infty$ function. It will be more convenient to work with this seemingly strange definition. 

Using such a convention, we can state Theorem \ref{lindthm} as follows: For $\eps\ll 1$, the Cauchy problem \eqref{qwe} with the initial data \eqref{init} has a global smooth solution $u=u(t,x)$ for all $t\geq 0$. Moreover, we have  pointwise decays: $Z^Iu=O_I(\eps\lra{t}^{-1+C_I\eps})$ for each multiindex $I$ and $\partial u=O(\eps\lra{t}^{-1})$. 

\subsection{The null condition of a matrix}
The definition and lemmas in this subsection will be used in Section \ref{sec4.2}. In this subsection, we assume that every matrix is in $\R^{4\times 4}$ and is a symmetric constant matrix.

\defn{\rm A matrix $m_0=(m_0^{\alpha\beta})_{\alpha,\beta=0,1,2,3}$ satisfies the \emph{null condition} if\fm{m_0^{\alpha\beta}\xi_\alpha\xi_\beta=0,\hspace{2em}\text{whenever }\xi\in\R^{1+3}\text{ and }|\xi_0|^2=|\xi_1|^2+|\xi_2|^2+|\xi_3|^2.}}\rm

We remark that a real symmetric constant matrix $m_0$ satisfies the null condition if and only if $m_0^{\alpha\beta}\xi_\alpha\eta_\beta$ is a linear combination of $-\xi_0\eta_0+\sum_{j=1}^3\xi_j\eta_j$ and $\xi_\alpha\eta_\beta-\xi_\beta\eta_\alpha$. 

We start with the following useful lemma. 
\lem{\label{lcommin} Suppose $m_0$ is a constant matrix satisfying the null condition. Then, for any two functions $\phi=\phi(t,x)$ and $\psi=\psi(t,x)$, we have
\fm{Z(m_0^{\alpha\beta}\phi_\alpha\psi_\beta)=m_0^{\alpha\beta}(\partial_\alpha Z\phi) \psi_\beta+m_0^{\alpha\beta}\phi_\alpha(\partial_\beta Z \psi)+m_1^{\alpha\beta}\phi_\alpha \psi_\beta.}Here $m_1$ is another symmetric constant matrix satisfying the null condition. Moreover, if $Z=\Omega_{ij}$ for $1\leq i,j\leq 3$ and if $(m_0^{\alpha\beta})=(m^{\alpha\beta})$ is the Minkowski metric, then $m_1=0$.}\rm

We refer our readers to Lemma 6.6.5 in \cite{horm} for the proof.

In addition, we have the following pointwise estimates related to the null condition.

\lem{\label{lnptwest} Suppose $m_0$ is a matrix satisfying the null condition. Then, for any two functions $\phi=\phi(t,x)$ and $\psi=\psi(t,x)$, if $t\sim r\gg1$, we have
\fm{|m_0^{\alpha\beta}\phi_\alpha \psi_\beta|\lesssim \lra{t}^{-1}(|Z\phi||\partial\psi|+|Z\psi||\partial \phi|).}
Here $|Zf|=\sum_{|J|=1}|Z^Jf|$ for a function $f=f(t,x)$.}\rm

We refer our readers to Lemma II.5.4 in \cite{sogg} for the proof.

\section{Construction of the optical function}\label{s3}
Let $u=u(t,x)$ be a global solution to \eqref{qwe} and \eqref{init} constructed in Theorem \ref{lindthm}. If we fix a constant $R>0$ such that $\text{supp }(u_0,u_1)\subset\{|x|\leq R\}$, then $u\equiv 0$ for $|x|\geq t+R$ by the finite speed of propagation. Our goal in this section is to construct an optical function, i.e.\ a solution to the eikonal equation 
\eq{\label{eik} g^{\alpha\beta}(u)\partial_\alpha q\partial_\beta q=0\text{ in }\Omega;\hspace{1cm}q=|x|-t\text{ on }\partial\Omega.}
The region $\Omega\subset\R^{1+3}_{t,x}$ is defined by
\eq{\label{defomega} \Omega:=\{(t,x):\ t>T_0,\ |x|>(t+T_0)/2+2R\}.}
Here $T_0:=\exp(\delta/\eps)$ for a fixed constant $\delta>0$. 

Our main result of this section is the following proposition.

\prop{\label{props3} The eikonal equation \eqref{eik} has a global $C^2$ solution in the region $\Omega$.}
\rm\bigskip

In Section \ref{s4}, we will show that this $C^2$ solution is in fact smooth (in the sense defined in Section \ref{sec2.5}).

Here we briefly explain how the optical function is constructed. In Section \ref{sec3.1}, we apply the method of characteristics and  solve the characteristic ODE's. Here the characteristics are in fact the null geodesics associated to the Lorentzian metric $(g_{\alpha\beta})$ which is the inverse of the coefficients $(g^{\alpha\beta}(u))$ in \eqref{eik}.  In Section \ref{sec3.2}, assuming that the optical function $q$ exists in some region, we prove several preliminary estiamtes for $q$ by studying the characteristic ODE's.

To finish the proof, we need to show that the characteristics, i.e.\ the geodesics, do not intersect with each other. This is related to the null geometry of the level sets of the optical function. In Section \ref{sec3.3} and \ref{sec3.4}, we construct a null frame  $\{e_k\}_{k=1}^4$ and then define several connection coefficients under the Lorentzian metric $(g_{\alpha\beta})$. Most importantly, we define the second fundamental form \fm{\chi_{ab}:=\lra{D_{e_a}e_4,e_b},\hspace{2em}a,b=1,2.}
Here $D$ is the Levi-Civita connection and $\lra{\cdot,\cdot}$ is the bilinear form, both associated to $(g_{\alpha\beta})$. One important quantity we need to estimate in our proof is the trace of $\chi$ which is sometimes called the \emph{null mean curvature}. We claim that it suffices to prove $\tr\chi>0$ everywhere. In fact,   for a $2$-sphere $S_{(t_0,q^0)}:=\{x:\ q(t_0,x)=q^0\}\subset\R^3$, we have \fm{\frac{d}{dt_0}|S_{(t_0,q^0)}|=\int_{S_{(t_0,q^0)}}\tr\chi \ dA.}
See, for example, Section 9.5 of \cite{alin}. If $\tr\chi>0$, then it implies that the $2$-sphere is expanding everywhere as the time increases. This  excludes the case when two distinct characteristics intersect with each other.

We now follow the idea in Christodoulou-Klainerman \cite{chriklai}. In Section \ref{sec3.5}, we derive an equation for $\chi$, called the Raychaudhuri equation. In Section \ref{sec3.6}, we use a continuity argument and the Raychaudhuri equation to prove that in the region where the optical function exists, we have\fm{\max_{a,b=1,2}|\chi_{ab}-\delta_{ab}r^{-1}|\lesssim t^{-2+C\eps}.}
We conclude that $\tr\chi>0$ everywhere. This implies that the characteristics will not intersect with each other, so we can extend the optical function to a slightly larger region. We thus finish the proof by making using of a continuity argument.

\subsection{The method of characteristics}\label{sec3.1}
Now we use the method of characteristics to solve \eqref{eik}.  We have the characteristic ODE's 
\eq{\label{charode}\left\{\begin{array}{l}\dot{x}^\alpha(s)=2g^{\alpha\beta}(x(s))p_\beta(s),\\
\dot{z}(s)=2g^{\alpha\beta}(x(s))p_\beta(s)p_\alpha(s)=0,\\
\dot{p}_\alpha(s)=-(\partial_\alpha g^{\mu\nu})(x(s))p_\mu(s)p_\nu(s).\end{array}\right.}
Here we write $g^{\alpha\beta}(t,x)=g^{\alpha\beta}(u(t,x))$ with an abuse of notation. We expect that $z(s)=q(x(s))$ and $p(s)=(\partial q)(x(s))$ for some optical function $q(t,x)$. By differentiating the first equation, we obtain the geodesic equation
\eq{\label{geoeqn}\ddot{x}^\alpha(s)+\Gamma^{\alpha}_{\mu\nu}\dot{x}^\mu(s)\dot{x}^\nu(s)=0.}
Here $\Gamma$ is the Christoffel symbol of the Levi-Civita connection $D$ of the Lorentzian metric $(g_{\alpha\beta})$. Thus, in this paper, the curve $x(s)$ is either called a \emph{characteristic curve}, or a \emph{geodesic}. 

To solve the eikonal equation \eqref{eik}, we only need to consider  the geodesics emanating from the surface \eq{\label{hdef}H&:=\{(t,x):\ t\geq T_0,\ r=(t+T_0)/2+2R\}\subset\partial\Omega.}
From these geodesics, later we will construct a solution $q(t,x)$ in the region $\Omega\cap\{r-t<2R\}$ such that $q=r-t$ in $\Omega\cap\{R<r-t< 2R\}$. Since $u\equiv 0$ in the region $r-t>R$, we can then extend our solution to the whole region $\Omega$ by defining $q=r-t$ when $r> t+R$.

To solve the characteristic ODE's \eqref{charode} and the geodesic equation \eqref{geoeqn}, we need to first determine $(\partial q)|_{H}$. Fix  $(t,x)\in H$ and recall that $q=r-t$ on $H$. Since $X_i:=\partial_i+2\omega_i\partial_t$ is tangent to $H$, we have $X_iq=X_i(r-t)=-\omega_i$ on $H$. Thus, for $(t,x)\in H$ we have $q_i=X_iq-2\omega_iq_t=-\omega_i-2\omega_iq_t$ and
\fm{0&=-q_t^2+2g^{0i}q_t(-\omega_i-2\omega_iq_t)+g^{ij}(-\omega_i-2\omega_iq_t)(-\omega_j-2\omega_jq_t)\\&=(-1-4g^{0i}\omega_i+4g^{ij}\omega_i\omega_j)q_t^2+(4g^{ij}\omega_i\omega_j-2g^{0i}\omega_i)q_t+g^{ij}\omega_i\omega_j.}
Since $g^{\alpha\beta}(u)=m^{\alpha\beta}+O(|u|)$, we have
\fm{0&=(-1+4m^{ij}\omega_i\omega_j+O(|u|))q_t^2+(4m^{ij}\omega_i\omega_j+O(|u|))q_t+(m^{ij}\omega_i\omega_j+O(|u|))\\
&=(3+O(|u|))q_t^2+(4+O(|u|))q_t+(1+O(|u|)).}
Since $|u|\ll 1$, by the root formula we can uniquely determine $q_t=-1+O(|u|)$ at $(t,x)$ (the other root $q_t=-1/3+O(|u|)$ is discarded since we expect $q$ to behave like $r-t$). We also have $q_i=-\omega_i-2\omega_iq_t=\omega_i+O(|u|)$ and $q_r=\omega_iq_i=1+O(|u|)$. If moreover $t< T_0+2R$, then $r=(t+T_0)/2+2R> t+R$ and thus $g^{\alpha\beta}\equiv m^{\alpha\beta}$. Thus, we have $q_t=-1$ and $q_i=\omega_i$ for $(t,x)\in H$ such that $t< T_0+2R$.

Now fix $x(0)\in H$. We set \fm{z(0)=r(x(0))-x^0(0),\hspace{2em} p_\alpha(0)=(\partial_\alpha q)(x(0))}where we set  \fm{r(V):=\Big(\sum_{i=1}^3 (V^i)^2\Big)^{1/2},\hspace{2em}\text{for a vector }V=(V^\alpha)_{\alpha=0}^3.} We have the following lemma.

\lem{\label{l3.1}Fix $x(0)\in H$ and construct $z(0),p(0)$  as above. Then the system \eqref{charode} along with the initial data $(x(0),z(0),p(0))$ has a unique solution $(x(s),z(s),p(s))$ on $[0,\infty)$. In addition, we have $\dot{x}^0(s)>0$ for all $s\geq 0$, and $x^0(s)\to\infty$ as $s\to \infty$.

If moreover we have $x(0)\in H\cap\{t<T_0+2R\}$, then $x(s)=(2s,2s\omega)+x(0)$. In other words, the geodesics emanating from $H\cap\{t<T_0+2R\}$ are straight lines. Thus $q=r-t$ whenever $r>t+R$.}
\begin{proof}
We apply the Picard existence and uniqueness theorem, e.g.\ Theorem 1.17 in \cite{tao}, to \eqref{charode}. From the theorem, we obtain a unique solution $(x(s),z(s),p(s))$ for all $0\leq s<s_{\max}$. By the blowup criterion in the theorem, either we have $s_{\max}<\infty$ and $|x(s)|+|z(s)|+|p(s)|\to\infty$ as $s\to s_{\max}$, or we have $s_{\max}=\infty$. Here $|x(s)|+|z(s)|+|p(s)|\to\infty$ is equivalent to $|x(s)|+|\dot{x}(s)|\to \infty$ due to $z(s)=z(0)$ and  the first equation in \eqref{charode}.

We claim that, along each geodesic, for all $s\geq 0$ we have \eq{\label{l3.1f01}4g^{\alpha\beta}(x(s))p_\alpha(s)p_\beta(s)=2\dot{x}^\alpha(s)p_\alpha(s)=g_{\alpha\beta}(x(s))\dot{x}^\alpha(s)\dot{x}^\beta(s)=0.}
In other words, the geodesics $x(s)$ are null curves. The first two equations follow from the first equation in \eqref{charode}, so here we only prove the last one. Note that the equality holds for $s=0$ by the construction of $(\partial q)|_H$. In addition,
\fm{\frac{d}{ds}(g^{\alpha\beta}(x(s))p_\alpha(s)p_\beta(s))&=2g^{\alpha\beta}(x(s))\dot{p}_\alpha(s)p_\beta(s)+(\partial_\mu g^{\alpha\beta})(x(s))\dot{x}^\mu(s)p_\alpha(s)p_\beta(s)\\
&=\dot{x}^\alpha(s)\dot{p}_\alpha(s)-\dot{p}_\mu(s)\dot{x}^\mu(s)=0.}
In the last line we use the third equation in \eqref{charode}. This ends the proof of \eqref{l3.1f01}.

Next we claim that $\dot{x}^0(s)>0$ for all $s$. Since $g^{\alpha\beta}(u)=m^{\alpha\beta}+O(|u|)$ for $|u|\ll1$,  its inverse $(g_{\alpha\beta}(u))$ is also a small pertubation of the Minkowski metric, i.e.\ $g_{\alpha\beta}=m_{\alpha\beta}+O(|u|)$. Thus, \eqref{l3.1f01} implies \fm{0=g_{00}(\dot{x}^0)^2+2g_{0i}\dot{x}^0\dot{x}^i+g_{ij}\dot{x}^i\dot{x}^j=-(\dot{x}^0(s))^2+\sum_i(\dot{x}^i(s))^2+O(|u(x(s))||\dot{x}|^2).}
We first show that $\dot{x}^0(s)\neq 0$ for all $s$. If $\dot{x}^0(s_0)=0$ for some $s_0>0$, then we have $g_{ij}\dot{x}^i\dot{x}^j=0$ at $s=s_0$. Since $g_{ij}=\delta_{ij}+O(|u|)$, the symmetric matrix $(g_{ij})$ is  positive definite. Then $\dot{x}(s_0)=0$. However, recall that $x(s)$ is a geodesic, and the only geodesic passing through $x(s_0)$ with $\dot{x}(s_0)=0$ is the constant curve $x(s)=x(s_0)$. This leads to a contradiction. In addition, since $q_t=-1+O(|u|)$ on $H$ and $\dot{x}^0(0)=2g^{0\beta}p_\beta(0)$, we have $\dot{x}^0(0)=2+O(|u|)$. Thus $\dot{x}^0(s)>0$ for all $s$.

Moreover, since $u=O(\eps\lra{t}^{-1+C\eps})$, we have
\fm{|-(\dot{x}^0(s))^2+\sum_i(\dot{x}^i(s))^2|&\leq C\eps \lra{x^0(s)}^{-1+C\eps}(|\dot{x}^0(s)|^2+\sum_i(\dot{x}^i(s))^2).}
By choosing $\eps\ll 1$, we can make $C\eps\leq 1/2$. Thus, for $\eps\ll 1$, we have
\fm{\sum_i(\dot{x}^i(s))^2&\leq (\dot{x}^0(s))^2+\frac{1}{2} (|\dot{x}^0(s)|^2+\sum_i(\dot{x}^i(s))^2)\Longrightarrow \sum_i(\dot{x}^i(s))^2\lesssim (\dot{x}^0(s))^2.}
Thus, for each $i$ we have
\fm{|x^i(s)|&=|x^i(0)+\int_0^s \dot{x}^i(\tau)\ d\tau|\leq |x^i(0)|+C\int_0^s \dot{x}^0(\tau)\ d\tau=|x^i(0)|+Cx^0(s).}
In conclusion, if $|x(s)|+|\dot{x}(s)|\to\infty$, then we must have $x^0(s)+\dot{x}^0(s)\to\infty$. 

If we differentiate the first equation in \eqref{charode} and use the third one, we obtain
\fm{|\ddot{x}^0(s)|&\leq|2g^{0\beta}\dot{p}_\beta|+|2(\partial_\mu g^{0\beta})\dot{x}^\mu p_\beta| \lesssim|\partial u(x(s))||\dot{x}(s)|^2\lesssim\eps\lra{x^0(s)}^{-1}(\dot{x}^0(s))^2.}
The last inequality follows since $|\dot{x}^i(s)|\lesssim \dot{x}^0(s)$ and since $\partial u=O(\eps \lra{t}^{-1})$. Since $\dot{x}^0>0$, we then have \fm{|\frac{d}{ds}\ln \dot{x}^0|=\frac{|\ddot{x}^0|}{\dot{x}^0}\leq C\eps\frac{\dot{x}^0}{x^0}=C\eps\frac{d}{ds}\ln x^0,}
which implies that
\fm{|\ln\dot{x}^0(s)-\ln\dot{x}^0(0)|&\lesssim\eps(\ln x^0(s)-\ln x^0(0)).}
The last inequality is equivalent to
\fm{\dot{x}^0(0)(\frac{x^0(s)}{x^0(0)})^{-C\eps}\leq\dot{x}^0(s)\leq \dot{x}^0(0)(\frac{x^0(s)}{x^0(0)})^{C\eps}.}
It follows that
\fm{\frac{d}{ds}((x^0(s))^{1-C\eps})&=(1-C\eps)(x^0(s))^{-C\eps}\dot{x}^0(s)\leq \dot{x}^0(0)(x^0(0))^{-C\eps},\\
\frac{d}{ds}((x^0(s))^{1+C\eps})&=(1+C\eps)(x^0(s))^{C\eps}\dot{x}^0(s)\geq \dot{x}^0(0)(x^0(0))^{C\eps}>0,}
and thus
\eq{\label{l3.1f1}(x^0(s))^{1-C\eps}&\leq (x^0(0))^{1-C\eps}+\dot{x}^0(0)s(x^0(0))^{-C\eps},}
\eq{\label{l3.1f2}(x^0(s))^{1+C\eps}&\geq (x^0(0))^{1+C\eps}+\dot{x}^0(0)s(x^0(0))^{C\eps}.}
If $s_{\max}<\infty$, then  $x^0(s)\to\infty$ as $s\to s_{\max}$ fails because of \eqref{l3.1f1}. On the other hand, if $s_{\max}<\infty$, then $x^0(s)+\dot{x}^0(s)\to\infty$ as discussed above. But since $\dot{x}^0(s)\leq \dot{x}^0(0)(x^0(s)/x^0(0))^{C\eps}$, we must have $x^0(s)\to\infty$ as $s\to s_{\max}$, which is a contradiction. Thus, $s_{\max}=\infty$. We thus conclude $x^0(s)\to\infty$ as $s\to\infty$ by \eqref{l3.1f2}.

The proof of the second half of this lemma is easy. We simply use the fact that $g^{\alpha\beta}(u)=m^{\alpha\beta}$ when $r\geq t+R$.
\end{proof}
\rmk{\rm \label{rmk3.1.1}We let $\mcl{A}$ denote the set of all the geodesics constructed in this lemma.}

\rm

\subsection{Estimates for the optical function}\label{sec3.2}
Fix a time $T> T_0=\exp(\delta/\eps)$ and we set $\Omega_T=\Omega\cap\{t\leq T,\ r-t\leq 2R\}$. Note that $r\sim t$ in $\Omega_T$. From now on, we assume that the optical function $q=q(t,x)$ exists in $\Omega_T$, that $q$ is $C^2$  and that $q_t<0$ everywhere. We remark that the assumptions are true for $T=T_0+2R$ since $g^{\alpha\beta}\equiv m^{\alpha\beta}$ in $\Omega_{T_0+2R}$. Our goal is to derive some estimates which allow us to extend the optical function to $\Omega_{T+\epsilon}$ for some $\epsilon>0$.

First of all, we claim that each point in $\Omega_T$ lies on exactly one geodesic in $\mcl{A}$ (which is defined in Remark  \ref{rmk3.1.1}). A direct corollary is that to define a function $F(t,x)$ in $\Omega_T$, we can define $F(x(s))$ along each geodesic in $\mcl{A}$. To prove this claim, we define a vector field $L=L^\alpha\partial_\alpha$ by  $L^\alpha:=2g^{\alpha\beta}q_\beta$. Note that $L^0> 0$ everywhere. In fact, we have \fm{g_{\alpha\beta}L^\alpha L^\beta=4g_{\alpha\beta}g^{\alpha\alpha'}g^{\beta\beta'}q_{\alpha'}q_{\beta'}=4g^{\alpha'\beta'}q_{\alpha'}q_{\beta'}=0.} If $L^0=0$, then $g_{ij}L^iL^j=0$. But $g_{ij}=\delta_{ij}+O(|u|)$, so $(g_{ij})$ is positive definite for $\eps\ll1$. Thus,  $L^\alpha=0$ and $q_t=\frac{1}{2}g_{0\beta}L^\beta=0$. This contradicts with the assumption that $q_t<0$. And since $L^0=-2q_t+O(|u\partial q|)=2+O(|u|)>0$ on $\partial\Omega$, we have $L^0>0$ in $\Omega_T$. Moreover, because of the characteristic ODE's \eqref{charode}, a curve in $\Omega_T$ is a geodesic in $\mcl{A}$  if and only if it is an integral curve of $L$ emanating from $H$. By the existence and uniqueness of integral curves, we finish the proof of the claim. 

We also claim that each geodesic emanating from $H\cap\partial\Omega_T$ must stay in $\Omega_T$ until it intersects with $\{t=T\}$. This claim simply follows from the fact that the optical function remains constant along each geodesic and that the optical function is injective when restricted to $(\partial\Omega_T)\setminus\{t=T\}$.

Here a useful lemma which follows directly from the chain rule and the pointwise estimates in Theorem \ref{lindthm} (also see Proposition 6.1 in Lindblad \cite{lind}).

\lem{\label{l3enf1} For each $k\geq 0$ and $\eps\ll_k1$, we have\fm{\sum_{|I|\leq k}(|Z^I(g^{\alpha\beta}-m^{\alpha\beta})|+|Z^I(g_{\alpha\beta}-m_{\alpha\beta})|)\lesssim_k \sum_{|I|\leq k}|Z^Iu|\lesssim_k \eps\lra{t}^{-1+C_k\eps}.} Moreover, \fm{|\partial g^{\alpha\beta}|+|\partial g_{\alpha\beta}|+|\Gamma_{\mu\nu}^\alpha|\lesssim |\partial u|\lesssim\eps \lra{t}^{-1}.}}
\rm

Now we can prove several useful estimates for $q$ in $\Omega_T$.

\lem{\label{l3enf2} In $\Omega_T$, we  have $|Sq|+\sum_i|\Omega_{0i}q|\lesssim |q|+t^{C\eps}$, $|\partial q|+
\sum_{i,j}|\Omega_{ij}q|\lesssim t^{C\eps}$ and $\sum_i|q_i-\omega_iq_r|\lesssim t^{-1+C\eps}$.}
\begin{proof}
If we apply a vector field $Z$ defined by \eqref{vf} to the eikonal equation, we obtain 
\fm{0=(Zg^{\alpha\beta})q_\alpha q_\beta+2g^{\alpha\beta} q_\alpha Zq_\beta=(Zg^{\alpha\beta})q_\alpha q_\beta+2g^{\alpha\beta} q_\alpha \partial_\beta Zq+2g^{\alpha\beta} q_\alpha  [Z,\partial_\beta]q.}
It is easy to check that $2m^{\alpha\beta} q_\alpha[Z,\partial_\beta]q=0$ if $Z\neq S$ and $[S,\partial_\beta]=-\partial_\beta$.  Thus, for some geodesic $x(s)$, we have
\fm{|\frac{d}{ds}(Zq(x(s)))|&\lesssim (|Zg^{\alpha\beta}|+|g^{\alpha\beta}-m^{\alpha\beta}|)|p(s)|^2\lesssim \eps(x^0(s))^{-1+C\eps}|\dot{x}(s)|^2\lesssim \eps(x^0(s))^{-1+C\eps}\dot{x}^0(s).}
Recall that $p(s)=(\partial q)(x(s))$ and that  we have $|\dot{x}^i(s)|\lesssim \dot{x}^0(s)\lesssim (x^0(s))^{C\eps}$ from the proof of Lemma \ref{l3.1}.  Since $\partial q=(-1,\omega)+O(|u|)$ on $H$, we have $|Sq|+|\Omega_{0j}q|=O(|q|+\eps t^{C\eps})$ and $|\Omega_{ij}q|=O(\eps t^{C\eps})$ on $H$. By integrating the inequality, we have \fm{|Zq(x(s))-Zq(x(0))|\lesssim\int_0^s \eps(x^0(\tau))^{-1+C\eps}\dot{x}^0(\tau)\ d\tau \lesssim (x^0(s))^{C\eps},}
so we have
\fm{|Zq(x(s))|\lesssim |Zq(x(0))|+(x^0(s))^{C\eps}\lesssim 1+|q(x(0))|+ (x^0(s))^{C\eps}=1+|q(x(s))|+ (x^0(s))^{C\eps}.}
In conclusion, we have $|Zq|=O(|q|+t^{C\eps})$ in $\Omega_T$. For $Z=\partial_\alpha$ or $\Omega_{ij}$ we have better bounds $|\Omega_{ij}q|+|\partial q|=O(t^{C\eps})$, since the estimates for $\partial q|_H$ and $\Omega_{ij}q|_H$ are better. In addition, we have $|q_i-\omega_iq_r|=r^{-1}|\sum_j\omega_j\Omega_{ij}q|\lesssim t^{-1+C\eps}$.
\end{proof}

\lem{\label{l3enf3} For each $(t,x)\in \Omega_T$, we have $q_r\geq C^{-1}t^{-C\eps}$, $-q_t\geq C^{-1}t^{-C\eps}$ and $|q_t+q_r|\lesssim\eps t^{-1+C\eps}$.}
\begin{proof}
Recall that from the proof of Lemma \eqref{l3.1}, we have $|\dot{x}^i(s)|\lesssim \dot{x}^0(s)$ and \fm{(x^0(s))^{-C\eps}\leq\dot{x}^0(0)(\frac{x^0(s)}{x^0(0)})^{-C\eps}\leq \dot{x}^0(s)\leq \dot{x}^0(0)(\frac{x^0(s)}{x^0(0)})^{C\eps}\leq (x^0(s))^{C\eps}} along each geodesic $x(s)$ in $\mcl{A}$.  At $(t_0,x_0)=x(s_0)$ for some geodesic $x(s)$ in $\mcl{A}$, we have
\eq{\label{l3enf3.2}q_t&=\frac{1}{2}g_{0\alpha}\dot{x}^\alpha(s_0)=-\frac{1}{2}\dot{x}^0(s_0)+O(|u(x(s_0))||\dot{x}(s_0)|)\leq -\frac{1}{2}t_0^{-C\eps}+C\eps t_0^{-1+C\eps}\leq -\frac{1}{4}t_0^{-C\eps}.}
Here we take $\eps\ll 1$ as usual.

To prove the estimate for $q_r$, we first prove that $q_r>0$ in $\Omega_T$. Assume  $q_r=0$ at some $(t_0,x_0)\in\Omega_T$. By the eikonal equation \eqref{eik} and the previous lemma, at $(t_0,x_0)$ we have
\eq{\label{l3enf3.1}0&=g^{00}q_t^2+2g^{0i}q_t(q_i-q_r\omega_i)+g^{ij}(q_i-\omega_iq_r)(q_j-\omega_j q_r)\\&=-q_t^2+O(|u||q_t|\sum_i|q_i-q_r\omega_i|)+O((\sum_{i}|q_i-\omega_iq_r|)^2)\\&=-q_t^2+O(t_0^{-2+C\eps}).}
Plug \eqref{l3enf3.2} into \eqref{l3enf3.1}, and we conclude that $t_0^{-2C\eps}\lesssim q_t^2\lesssim t_0^{-2+C\eps}$ and $t_0^{2-3C\eps}\lesssim 1$. This is impossible, since $t_0^{2-3C\eps}\geq t_0\geq T_0=\exp(\delta/\eps)\gg 1$ for $\eps\ll 1$. So we have $q_r\neq 0$ everywhere in $\Omega_T$. Since $q_r=1+O(|u|)>0$ on $H$, we have $q_r>0$ everywhere in $\Omega_T$. By \eqref{l3enf3.2}, we  have $-q_t+q_r\geq -q_t\geq \frac{1}{4}t^{-C\eps}$. Then since
\fm{0&=-q_t^2+\sum_i q_i^2+O(|u||\partial q|^2)=(q_t+q_r)(-q_t+q_r)+\sum_i(q_i-q_r\omega_i)^2+O(\eps t^{-1+C\eps}|\partial q|^2)\\&=(q_t+q_r)(-q_t+q_r)+O(t^{-2+2C\eps}+\eps t^{-1+C\eps})} and since $t^{-1}\leq T_0^{-1}\ll \eps$, we have
\fm{|q_t+q_r|=(-q_t+q_r)^{-1}O(\eps t^{-1+C\eps})\lesssim t^{C\eps}\cdot \eps t^{-1+C\eps}\lesssim \eps t^{-1+C\eps}.}
Then we have $q_r=-q_t+(q_t+q_r)\geq C^{-1}t^{-C\eps}-C\eps t^{-1+C\eps}\geq C^{-1}t^{-C\eps}$.
\end{proof}

\rm

\subsection{A null frame}\label{sec3.3}
We construct a null frame $\{e_1,e_2,e_3,e_4\}$ in $\Omega_T$ as follows. Define two vector fields $e_3,e_4$ by
\fm{e_4:=(L^0)^{-1}L,\hspace{2em} e_3:=e_4+2g^{0\alpha}\partial_\alpha.}
Since $g^{00}\equiv -1$, we have $e_4^0\equiv 1$ and $e_3^0\equiv-1$. Moreover,  we have
\eq{\label{nfcond}\lra{e_4,e_4}&=(L^0)^{-2}\lra{L,L}=(L^0)^{-2}g_{\alpha\beta}L^\alpha L^\beta=0,\\
\lra{e_4,e_3}=\lra{e_3,e_4}&=\lra{2g^{0\alpha}\partial_\alpha,e_4}=2g_{\alpha\beta}g^{0\alpha}e_4^\beta=2e_4^0=2,\\
\lra{e_3,e_3}&=\lra{e_4,e_3}+\lra{2g^{0\alpha}\partial_\alpha,e_3}=2+2g_{\alpha\beta}g^{0\alpha}e_3^\beta=2+2e_3^0=0.}
Here $\lra{\cdot,\cdot}$ is the bilinear form defined by the Lorentzian metric  $(g_{\alpha\beta})=(g^{\alpha\beta})^{-1}$.

Next we define $\{e_a\}_{a=1,2}$.  When restricted to the $2$-sphere $H\cap\{t=T'\}$ for some $T'\geq T_0$, the metric $(g_{\alpha\beta})$ is positive definite. Thus, we can choose a smooth orthonormal basis $\{E_a\}_{a=1,2}$ locally on this $2$-sphere. Here we make our choice such that $E_a|_H$ depends only on $\omega$ and not on $t$. Note that $E_a$ is tangent to $H\cap\{t=T'\}$, that $E_a^0=0$ and that $\lra{E_a,E_b}=\delta_{ab}$. Then we take the parallel transport of $E_a$ along the geodesics. That is, we consider the equations $D_4E_a=0$ for $a=1,2$. Here $D$ is the Levi-Civita connection of the Lorentzian metric, and $D_4:=D_{e_4}$. Since $e_4$ is tangent to the geodesic, equivalently we need to solve the ODE's
\eq{\label{patrode}\frac{d}{ds}E_a^\alpha(x(s))+\dot{x}^\mu(s) E_a^\nu(x(s))\Gamma^\alpha_{\mu\nu}(x(s))=0.}
By the existence and uniqueness for linear ODE's (e.g.\ Theorem 4.12 in \cite{leebook}), these ODE's admit a unique solution for all $0\leq s\leq s_0$.  Finally, we define 
\fm{e_a&:=E_a-E_a^0e_4,\hspace{1cm}a=1,2.}
Thus $e_a^0=0$. Unlike $e_3,e_4$, the vector fields $e_1,e_2$ cannot be defined globally in $\Omega_T$. This is because there is no global orthonormal basis on a $2$-sphere. In the rest of this paper, when we state a property of $e_a$ on $\Omega_{T}$, we mean that any locally defined $e_a$ satisfies this property.

We conclude that  $\{e_k\}_{k=1,2,3,4}$ is a null frame by \eqref{nfcond} and the following lemma. 

\lem{\label{l3.2} In $\Omega_T$ we have $\lra{e_a,e_b}=\delta_{ab}$ and $\lra{e_4,e_a}=\lra{e_3,e_a}=0$ for each $a,b=1,2$.}
\begin{proof}
We first prove that $\lra{E_a,E_b}=\delta_{ab}$ and $\lra{e_4,E_a}=0$ on $H$. The first equality follows directly from the construction of $\{E_a\}$. To prove the second one, we recall that $q_i=q_r\omega_i$ on $H$; see the computations right above Lemma \ref{l3.1}. Moreover, note that $\sum_ix^i(0)E_a^i=0$ since $E_a$ is tangent to the sphere on $H$. Thus, on $H$, we have
\fm{\lra{L,E_a}&=g_{\alpha\beta}L^\alpha E_a^\beta=2q_\beta E_a^\beta=2q_i E_a^i=2q_r\omega_iE_a^i=0.}
And since $e_4=(L^0)^{-1}L$, we have $\lra{e_4,E_a}=0$ at $x(0)$.

Along each geodesic $x(s)$ in $\mcl{A}$, we have
\fm{e_4\lra{E_a,E_b}&=\lra{D_4E_a,E_b}+\lra{E_a,D_4E_b}=0,\\
e_4\lra{L,E_a}&=\lra{D_4L,E_a}+\lra{L,D_4E_a}=0.}
Because of the equalities at $s=0$, we conclude that $\lra{E_a,E_b}=\delta_{ab}$ and $\lra{L,E_a}=0$ (and thus $\lra{e_4,E_a}=0$) along each geodesic.

Finally, note that\fm{\lra{e_a,e_b}&=\lra{E_a,E_b}-E_a^0\lra{e_4,E_b}-E_b^0\lra{E_a,e_4}+E_a^0 E_b^0\lra{e_4,e_4}=\delta_{ab},\\
\lra{e_4,e_a}&=\lra{e_4,E_a}-E_a^0\lra{e_4,e_4}=0,\\
\lra{e_3,e_a}&=\lra{2g^{0\alpha}\partial_\alpha,e_a}+\lra{e_4,e_a}=2g_{\alpha\beta}g^{0\alpha}e_a^\beta=2e_a^0=0.}
This finishes the proof.
\end{proof}\rm

Before we move on to the next lemma, we summarize some important properties of a null frame. First, any vector field $X$  can be uniquely expressed as a linear combination of the null frame: 
\eq{\label{xvfformula}X=\sum_{a=1,2}\lra{X,e_a}e_a+\frac{1}{2}\lra{X,e_4}e_3+\frac{1}{2}\lra{X,e_3}e_4.}

In addition,  for each $k=1,2,3,4$ we have
\fm{\lra{g^{\alpha\beta}\partial_\beta,e_k}=g^{\alpha\beta} g_{\beta\mu}e_k^\mu=e_k^\alpha,}
so we obtain\eq{\label{gpaformula}g^{\alpha\beta}\partial_\beta=\sum_{a=1,2}e_a^\alpha e_a+\frac{1}{2}e_4^\alpha e_3+\frac{1}{2}e_3^\alpha e_4\Longrightarrow g^{\alpha\beta}=\sum_{a=1,2}e_a^\alpha e_a^\beta+\frac{1}{2}e_4^\alpha e_3^\beta+\frac{1}{2}e_3^\alpha e_4^\beta.}

Finally, we have $e_1(q)=e_2(q)=e_4(q)=0$ and $e_3(q)=L^0$ in $\Omega_T$. In fact, since $q_\alpha=\frac{1}{2}g_{\alpha\beta}L^\beta$, we have $X q=\frac{1}{2}\lra{X,L}=\frac{1}{2}L^0\lra{e_4,X}$ for each vector field $X$. Then we use the properties of a null frame. The equality $e_1(q)=e_2(q)=e_4(q)=0$ implies that $e_1,e_2,e_4$ are tangent to the level set of $q$, so $e_1,e_2,e_4$ are sometimes called the tangential derivatives.

The next lemma shows several better estimates for the tangential derivatives.

\lem{\label{l3betterest} In $\Omega_T$, we have $e_4=\partial_t+\partial_r+O(t^{-1+C\eps})\partial$, $e_3=e_4+2g^{0\alpha}\partial_\alpha=-\partial_t+\partial_r+O(t^{-1+C\eps})\partial$ and $e_a=O(1)\partial$. Then, for all $I,s,l$, we have\fm{\sum_{k=1,2,4}(|e_k(\partial^sZ^Iu)|+|e_k(\partial^sZ^Ig^{\alpha\beta})|+|e_k(\partial^sZ^Ig_{\alpha\beta})|)\lesssim\eps t^{-2+C\eps}\lra{r-t}^{-s}.}
Here we use the convention given in Section \ref{sec2.5}.
Moreover, we have \fm{|e_1(\partial_\alpha g_{\mu\nu})e_2^\alpha|+|e_2(\partial_\alpha g_{\mu\nu})e_1^\alpha|+|e_1(\partial_\alpha g_{\mu\nu})e_1^\alpha-e_2(\partial_\alpha g_{\mu\nu})e_2^\alpha|\lesssim \eps t^{-3+C\eps}.}}
\begin{proof}
By the lemmas in Section \ref{sec3.2}, we have \fm{e_4^i-\omega_i&=\frac{L^i-L^0\omega_i}{L^0}=\frac{2q_i+2q_t\omega_i+O(|u||\partial q|)}{-2q_t+O(|u||\partial q|)}=\frac{2(q_i-q_r\omega_i)+2(q_r+q_t)\omega_i+O(|u||\partial q|)}{-2q_t+O(|u||\partial q|)}.}By Lemma \ref{l3enf2} and Lemma \ref{l3enf3}, the denominator has a lower bound $C^{-1}t^{-C\eps}-C\eps t^{-1+C\eps}\geq (2C)^{-1}t^{-C\eps}$ and the numerator is $O(t^{-1+C\eps})$. In conclusion, $e_4=\partial_t+\partial_r+O(t^{-1+C\eps})\partial$. It follows that for each $I$, \fm{|e_4(\partial^sZ^Iu)|&\lesssim |(\partial_t+\partial_r)\partial^sZ^Iu|+t^{-1+C\eps}|\partial \partial^sZ^Iu|
\\
&\lesssim \lra{t+r}^{-1}\sum_{|J|=1}|Z^J\partial^sZ^Iu|+t^{-1+C\eps}\lra{r-t}^{-s-1}\sum_{|J|\leq s+1}|Z^JZ^Iu|
\\&\lesssim \lra{t+r}^{-1}\sum_{|J|\leq 1}|\partial^sZ^JZ^Iu|+t^{-1+C\eps}\lra{r-t}^{-s-1}\cdot \eps t^{-1+C\eps}\\
&\lesssim \lra{t+r}^{-1}\cdot \eps t^{-1+C\eps}\lra{r-t}^{-s}+\eps t^{-2+C\eps}\lra{r-t}^{-s-1}\\&\lesssim\eps t^{-2+C\eps}\lra{r-t}^{-s}.}
Here we apply Lemma \ref{l2.1}, the pointwise decays in Theorem \ref{lindthm}, and \eqref{comf3}.
By the chain rule and Leibniz's rule, we can express $e_4(\partial^sZ^I(g^{\alpha\beta},g_{\alpha\beta}))$ as a linear combination of terms of the form\fm{\frac{d^m}{du^m}(g^{\alpha\beta},g_{\alpha\beta})(u)\cdot (\partial^{s_1}Z^{I_1}u)\cdots(\partial^{s_{m-1}}Z^{I_{m-1}}u)\cdot e_4(\partial^{s_m}Z^{I_m}u)}
where $\sum s_*=s$, $\sum |I_*|=|I|$ and $m>0$. These terms have an upper bound
\fm{\eps t^{-1+C\eps}\lra{r-t}^{-s_1}\cdots \eps t^{-1+C\eps}\lra{r-t}^{-s_{m-1}}\cdot \eps t^{-2+C\eps}\lra{r-t}^{-s_m}\lesssim\eps t^{-2+C\eps}\lra{r-t}^{-s}.}
We thus have  $e_4(\partial^sZ^I(g^{\alpha\beta},g_{\alpha\beta}))= O(\eps t^{-2+C\eps}\lra{r-t}^{-s})$.

Next we fix $(t_0,x_0)\in\Omega_T$. Without loss of generality, we assume $|q_3|=\max\{|q_j|:\ j=1,2,3\}$ at $(t_0,x_0)$. For $i=1,2$, we define \fm{Y_i:=q_i\partial_3-q_3\partial_i=r^{-1}q_r\Omega_{i3}+(q_i-\omega_iq_r)\partial_3-(q_3-\omega_3q_r)\partial_i=r^{-1}q_r\Omega_{i3}+O(t^{-1+C\eps})\partial.}
Here $\{Y_1,Y_2\}$ is a basis of the tangent space of the $2$-sphere $\Sigma_{(t_0,x_0)}=\{t=t_0,q=q(t_0,x_0)\}$  at $(t_0,x_0)$. Since $e_a$ lies in the tangent space (as $e_a^0=0$ and $e_a(q)=0$), we can write $e_a=\sum_{i=1,2}c_{ai}Y_{i}$ in a unique way. Since \fm{\lra{Y_i,Y_j}=q_iq_jg_{33}+q_3^2g_{ij}-q_iq_3g_{3j}-q_jq_3g_{3i}=q_iq_j+q_3^2\delta_{ij}+O(|u|q_3^2),\hspace{2em} i,j=1, 2,} we have 
\fm{1&=\lra{e_a,e_a}=\sum_{i,j} c_{ai}c_{aj}\lra{Y_i,Y_j}
=(\sum_ic_{ai}q_i)^2+(1+O(|u|))q_3^2\sum_{i} c_{ai}^2.}
Then, for $\eps\ll 1$ we have
\fm{1\geq 0+(1+O(\eps t^{-1+C\eps}))q_3^2\sum_{i} c_{ai}^2\geq \frac{1}{2}q_3^2\sum_{i} c_{ai}^2.}
Thus, we have $|q_3c_{ai}|\lesssim 1$ for each $a,i$ and thus $e_a^\alpha=\sum_i c_{ai}Y_i^\alpha=O(|c_{ai}q_3|)=O(1)$. And since $C^{-1}t^{-C\eps}\leq |q_r|=|\sum_i\omega_iq_i|\leq\sum_i|q_i|\leq 3|q_3|$, for each multiindex $I$, we have
\fm{|e_{a}(\partial^sZ^Iu)|&\leq \sum_{i}|c_{ai}Y_i(\partial^sZ^Iu)|\lesssim \sum_i|c_{ai}|(r^{-1}|q_r||\Omega \partial^sZ^Iu|+t^{-1+C\eps}|\partial\partial^sZ^Iu|)\\&\lesssim \eps t^{-2+C\eps}\lra{r-t}^{-s}.}
By the chain rule and Leibniz's rule, we finish the proof of the first estimate.

In addition,
\eq{\label{l3betterestf1}0&=\lra{e_1,e_1}-\lra{e_2,e_2}\\&=(\sum_ic_{1i}q_i)^2-(\sum_ic_{2i}q_i)^2+q_3^2\sum_i (c_{1i}^2-c_{2i}^2)+O(|u|q_3^2\sum_{a,i} c_{ai}^2)\\&=(\sum_ic_{1i}q_i)^2-(\sum_ic_{2i}q_i)^2+q_3^2\sum_i (c_{1i}^2-c_{2i}^2)+O(|u|)\\
&=\sum_{i,j}(c_{1i}c_{1j}-c_{2i}c_{2j})q_iq_j-(\sum_ic_{2i}q_i)^2+q_3^2\sum_i (c_{1i}^2-c_{2i}^2)+O(|u|),}
\eq{\label{l3betterestf2}0&=\lra{e_1,e_2}\\&=\sum_{i,j}c_{1i}c_{2j}\lra{Y_i,Y_j}=\sum_{i,j}c_{1i}c_{2j}q_iq_j+\sum_i c_{1i}c_{2i}q_3^2+O(|u|q_3^2\sum_{i,j}|c_{1i}c_{2j}|)\\
&=\sum_{i,j}c_{1i}c_{2j}q_iq_j+\sum_i c_{1i}c_{2i}q_3^2+O(|u|).}
Then, we have
\fm{Y_i(Zg)&=r^{-1}q_r\Omega_{i3}g+O(t^{-1+C\eps}|\partial g|)=O(\eps t^{-2+C\eps}),}
\fm{Y_i(\partial_\alpha g)Y_j^\alpha&=(r^{-1}q_r\Omega_{i3}(\partial_\alpha g)+(q_i-\omega_iq_r)\partial_3\partial_\alpha g-(q_3-\omega_3q_r)\partial_i\partial_\alpha g)Y_j^\alpha\\
&=r^{-1}q_r(Y_j^\alpha[\Omega_{i3},\partial_\alpha] g+Y_j\Omega_{i3}g)+(q_i-\omega_iq_r)Y_j(\partial_3 g)-(q_3-\omega_3q_r)Y_j(\partial_i g)\\&=r^{-1}q_r(-Y_j^i\partial_3g+Y_j^3\partial_ig) +r^{-1}q_rY_j\Omega_{i3}g+O(t^{-1+C\eps}|Y_j(\partial g)|)\\&=r^{-1}q_r(\delta_{ij}q_3\partial_3g+q_j\partial_ig) +O(\eps t^{-3+C\eps}),}
\fm{e_a(\partial_\alpha g)e_b^\alpha&=\sum_{i,j} c_{ai}Y_i(\partial_\alpha g)c_{bj}Y_j^\alpha=\sum_{i,j}c_{ai}c_{bj}(r^{-1}q_r(\delta_{ij}q_3\partial_3g+q_j\partial_ig) +O(\eps t^{-3+C\eps}))\\
&=\sum_{i}r^{-1}c_{ai}c_{bi}q_rq_3\partial_3g+\sum_{i,j}r^{-1}c_{ai}c_{bj}q_r q_j\partial_ig+O(\sum_{i,j}|c_{ai}c_{bj}||q_3| \eps t^{-3+C\eps})\\
&=\sum_{i}r^{-1}c_{ai}c_{bi}q_rq_3\partial_3g+\sum_{i,j}r^{-1}c_{ai}c_{bj}q_r q_j\partial_ig+O(\eps t^{-3+C\eps}).}
When $a\neq b$,  by \eqref{l3betterestf2} we have 
\fm{e_a(\partial_\alpha g)e_b^\alpha&=r^{-1}q_rq_3^{-1}(-\sum_{i,j}c_{ai}c_{bj}q_iq_j+O(|u|))\partial_3g+\sum_{i,j}r^{-1}c_{ai}c_{bj}q_r q_j\partial_ig+O(\eps t^{-3+C\eps})\\
&=r^{-1}q_rq_3^{-1}\sum_{i,j}c_{ai}c_{bj}q_j(-q_i\partial_3g+q_3 \partial_ig)+O(r^{-1}|q_rq_3^{-1}||u||\partial g|)+O(\eps t^{-3+C\eps})\\
&=r^{-1}q_rq_3^{-1}\sum_{i,j}c_{ai}c_{bj}q_j(-Y_ig)+O(\eps t^{-3+C\eps})=O(\eps t^{-3+C\eps}).}
By \eqref{l3betterestf1} we have
\fm{&\hspace{1em}e_1(\partial_\alpha g)e_1^\alpha-e_2(\partial_\alpha g)e_2^\alpha\\
&=\sum_{i}r^{-1}(c_{1i}^2-c_{2i}^2)q_rq_3\partial_3g+\sum_{i,j}r^{-1}(c_{1i}c_{1j}-c_{2i}c_{2j})q_r q_j\partial_ig+O(\eps t^{-3+C\eps})\\
&=r^{-1}q_rq_3^{-1}(-\sum_{i,j}(c_{1i}c_{1j}-c_{2i}c_{2j})q_iq_j)\partial_3g+\sum_{i,j}r^{-1}(c_{1i}c_{1j}-c_{2i}c_{2j})q_r q_j\partial_ig+O(\eps t^{-3+C\eps})\\
&=\sum_{i,j}r^{-1}q_rq_3^{-1}q_j(c_{1i}c_{1j}-c_{2i}c_{2j})(-Y_ig)+O(\eps t^{-3+C\eps})=O(\eps t^{-3+C\eps}).}

It is clear that our proof would still work if we assume $|q_1|=\max\{|q_j|:\ j=1,2,3\}$ or $|q_2|=\max\{|q_j|:\ j=1,2,3\}$. This ends the proof.
 
\end{proof}

\lem{\label{l3qrt} In $\Omega_T$, we have $|q-(r-t)|\lesssim t^{C\eps}$.}
\begin{proof}
By the previous lemma and Lemma \ref{l3enf3}, we have
\fm{e_4^i-\omega_i&=\frac{2(q_i-q_r\omega_i)+2(q_r+q_t)\omega_i+O(|u||\partial q|)}{L^0}=2(L^0)^{-1}(q_i-q_r\omega_i)+O(\eps t^{-1+C\eps}).}
Thus, \fm{e_4(q-r+t)=(\partial_t+\partial_r)(-r+t)-2(L^0)^{-1}\sum_i(q_i-q_r\omega_i)\omega_i+O(\eps t^{-1+C\eps})=O(\eps t^{-1+C\eps}).}
Suppose $(t,x)\in\Omega_T$ lies on a geodesic $x(s)$ in $\Omega_T$. Since $q-r+t=0$ on $H$, by integrating $e_4(q-r+t)$ along this geodesic, we have 
\fm{|q-r+t|\lesssim \int_{x^0(0)}^t\eps \tau^{-1+C\eps}\ d\tau\lesssim t^{C\eps}.}
\end{proof}

\rm

\subsection{The connection coefficients}\label{sec3.4}
From now on, we write $D_k=D_{e_k}$ for $k=1,2,3,4$ for simplicity. 

\lem{\label{l3.3} In $\Omega_T$, we have
\fm{D_4e_k&=(\Gamma_{\alpha\beta}^0e_4^\alpha e_k^\beta) e_4,\hspace{2em}k=1,2,4.}
As a result, we have $e_4(e_k^\alpha)=O(\eps t^{-2+C\eps})$ for each $k=1,2,3,4$.}
\begin{proof}
Since a geodesic in $\mcl{A}$ is an integral curve of $L$, we have $L^\alpha=\dot{x}^\alpha(s)$ at $x(s)$. Then, the geodesic equation \eqref{geoeqn}  implies \fm{L(L^0)=\dot{x}^\alpha(s)(\partial_{\alpha}L^0)=\frac{d}{ds}L^0(x(s))=\ddot{x}^0(s)=-\Gamma_{\mu\nu}^0 L^\mu L^\nu,\hspace{2em}\text{at }x(s).}
Divide both sides by $L^0$, and we conclude $e_4(L^0)=-\Gamma_{\mu\nu}^0 e_4^\mu L^\nu$ in $\Omega_T$ and thus $e_4(\ln L^0)=-\Gamma_{\mu\nu}^0 e_4^\mu e_4^\nu$. Similarly, from \eqref{patrode} we obtain $e_4(E_a^0)=-\Gamma_{\mu\nu}^0 e_4^\mu E_a^\nu$. Thus, we have
\fm{D_4e_4&=D_4((L^0)^{-1}L)=-(L^0)^{-2}e_4(L^0)L+(L^0)^{-1}D_4L=-(L^0)^{-1}e_4(L^0)e_4=(\Gamma_{\mu\nu}^0 e_4^\mu e_4^\nu) e_4.}

For $a=1,2$, since $D_4E_a=0$, we have
\fm{D_4e_a&=D_4(E_a-E_a^0e_4)=-D_4(E_a^0e_4)=-e_4(E_a^0)e_4-E_a^0D_4e_4\\
&=(\Gamma_{\mu\nu}^0 e_4^\mu E_a^\nu) e_4-(E_a^0\Gamma_{\mu\nu}^0 e_4^\mu e_4^\nu) e_4=\Gamma_{\mu\nu}^0e_4^\mu(  E_a^\nu-E_a^0  e_4^\nu) e_4\\&=(\Gamma_{\mu\nu}^0e_4^\mu e_a^\nu) e_4.}
In addition, $D_4e_k=e_4(e_k^\alpha)\partial_\alpha+\Gamma_{\mu\nu}^\alpha e_4^\mu e_k^\nu\partial_\alpha$. If we consider the coefficients of $\partial_\alpha$ in $D_4e_k$  for $k=1,2,4$, we have $e_4(e_k^\alpha)=\Gamma_{\mu\nu}^0e_4^\mu e_k^\nu e_4^\alpha-\Gamma_{\mu\nu}^\alpha e_4^\mu e_k^\nu$. By Lemma \ref{l3betterest}, we have
\eq{\label{l3.3fchrit}\Gamma^{\alpha}_{\mu\nu}&=\frac{1}{2}g^{\alpha\beta}(\partial_\mu g_{\nu\beta}+\partial_\nu g_{\mu\beta}-\partial_\beta g_{\mu\nu})\\&=\frac{1}{2}g^{\alpha\beta}(\partial_\mu g_{\nu\beta}+\partial_\nu g_{\mu\beta})-\frac{1}{2}(\sum_a e_a^\alpha e_a(g_{\mu\nu})+\frac{1}{2}(e_3^\alpha e_4(g_{\mu\nu})+e_4^\alpha e_3(g_{\mu\nu})))\\
&=\frac{1}{2}g^{\alpha\beta}(\partial_\mu g_{\nu\beta}+\partial_\nu g_{\mu\beta})-\frac{1}{4}e_4^\alpha e_3(g_{\mu\nu})+O(\eps t^{-2+C\eps}).}
Then, since $e_4^0=1$, for $k=1,2,4$ we have
\fm{e_4(e_k^\alpha)&=(\frac{1}{2}g^{0\beta}(\partial_\mu g_{\nu\beta}+\partial_\nu g_{\mu\beta})-\frac{1}{4}e_4^0 e_3(g_{\mu\nu})+O(\eps t^{-2+C\eps}))e_4^\mu e_k^\nu e_4^\alpha\\&\hspace{1em}-(\frac{1}{2}g^{\alpha\beta}(\partial_\mu g_{\nu\beta}+\partial_\nu g_{\mu\beta})-\frac{1}{4}e_4^\alpha e_3(g_{\mu\nu})+O(\eps t^{-2+C\eps})) e_4^\mu e_k^\nu\\
&=\frac{1}{2}g^{0\beta}(e_4(g_{\nu\beta})e_k^\nu e_4^\alpha+e_k( g_{\mu\beta})e_4^\mu e_4^\alpha)+\frac{1}{2}g^{\alpha\beta}(e_4( g_{\nu\beta})e_k^\nu+e_k( g_{\mu\beta})e_4^\mu)\\&\hspace{1em}-\frac{1}{4}e_3(g_{\mu\nu})(e_4^\mu e_k^\nu e_4^\alpha e_4^0-e_4^\mu e_k^\nu e_4^\alpha)+O(\eps t^{-2+C\eps})\\
&=O(\eps t^{-2+C\eps}).}
It follows that $e_4(e_3^\alpha)=e_4(e_4^\alpha)+e_4(2g^{0\alpha})=O(\eps t^{-2+C\eps})$. This finishes the proof.
\end{proof}
\rmk{\rm\label{rmk3.3.1} Since $e_3(q)=L^0$, we have \fm{e_4(e_3(q))=e_4(L^0)=-\Gamma^0_{\alpha\beta}e_4^\alpha L^\beta=-\Gamma^0_{\alpha\beta}e_4^\alpha e_4^\beta e_3(q).}This equality is useful in the rest of this paper.}
\rm

Next, we set $\chi_{ab}:=\lra{D_ae_4,e_b}$ for $a,b=1,2$. 
\lem{\label{l3.4} In $\Omega_T$, we have

{\rm (a)}  $\chi_{12}=\chi_{21}$.

{\rm (b)} $\tr\chi:=\chi_{11}+\chi_{22}$ is independent of the choice of $e_1$ and $e_2$.

{\rm (c)} \fm{\ [e_4,e_a]&=-\sum_b\chi_{ab}e_b,\ 
D_ae_4 =\sum_b\chi_{ab}e_b+(e_4^\mu e_a^\nu\Gamma_{\mu\nu}^0)e_4,\ 
e_a(e_4^\alpha)&=\sum_b\chi_{ab}e_b^\alpha+O(\eps t^{-2+C\eps}).}}

\begin{proof}
(a) Since $e_a(q)=0$, we have\fm{\lra{e_4,[e_1,e_2]}=(L^0)^{-1}\lra{L,[e_1,e_2]}=2(L^0)^{-1}[e_1,e_2]q=2(L^0)^{-1}(e_1(e_2(q))-e_2(e_1(q)))=0.}
And since \fm{\lra{D_ke_l,e_m}=e_k(\lra{e_l,e_m})-\lra{e_l,D_ke_m}=-\lra{e_l,D_ke_m},\hspace{2em}k,l,m=1,2,3,4,}
we have
\fm{\chi_{12}-\chi_{21}&=\lra{D_1e_4,e_2}-\lra{D_2e_4,e_1}=\lra{e_4,-D_1e_2+D_2e_1}=-\lra{e_4,[e_1,e_2]}=0.}

(b) Suppose that $\{e_k'\}$ is another null frame with $e_3=e_3'$ and $e_4=e_4'$. Then we have $e_a'=\sum_b\lra{e_a',e_b}e_b$, $e_a=\sum_b\lra{e_a,e_b'}e_b'$ and thus \fm{e_{a}=\sum_b\lra{e_a,e_b'}e_b'=\sum_{b,c}\lra{e_a,e_b'}\lra{e_b',e_c}e_c\Longrightarrow\sum_{b,c}\lra{e_a,e_b'}\lra{e_b',e_c}=\delta_{ac}.} Then,
\fm{\chi_{11}'+\chi_{22}'&=\sum_a\lra{D_{e_a'}e_4,e_a'}=\sum_{a}\sum_{b,c}\lra{e_a',e_b}\lra{e_a',e_c}\lra{D_be_4,e_c}\\
&=\sum_{b,c}\sum_{a}\lra{e_a',e_b}\lra{e_a',e_c}\chi_{bc}=\sum_{b,c}\delta_{bc}\chi_{bc}=\chi_{11}+\chi_{22}.}

(c) Since $D_4e_k=(\Gamma_{\alpha\beta}^0e_4^\alpha e_k^\beta) e_4$ for $k=1,2,4$, we have $\lra{D_4e_k,e_a}=0$ for $k=1,2,4$ and thus
\fm{\lra{e_4,[e_4,e_a]}&=\lra{e_4,D_4e_a-D_ae_4}=-\lra{D_4e_4,e_a}-\frac{1}{2}e_a\lra{e_4,e_4}=0,\\
\lra{e_b,[e_4,e_a]}&=\lra{e_b,D_4e_a-D_ae_4}=\lra{e_b,D_4e_a}-\chi_{ab}=-\chi_{ab}.}
Since $e_4^0=1$ and $e_a^0=0$, we have $[e_4,e_a]^0=0$ (where $[e_4,e_a]=[e_4,e_a]^\alpha\partial_\alpha$) and thus \fm{\lra{e_3,[e_4,e_a]}=\lra{e_4,[e_4,e_a]}+2g^{0\alpha}g_{\alpha\beta}[e_4,e_a]^\beta=0+2[e_4,e_a]^0=0.}
By \eqref{xvfformula} we conclude that $[e_4,e_a]=-\sum_{b=1,2}\chi_{ab}e_b$. The second equality follows from $D_ae_4=[e_a,e_4]+D_4e_a$. The third one follows from $e_a(e_4^\alpha)-e_4(e_a^\alpha)=[e_a,e_4]^\alpha$ and the previous lemma.
\end{proof}
\subsection{The Raychaudhuri equation}\label{sec3.5}
\rm It turns out that the estimates for $\chi_{ab}$ are crucial in the proof of the global existence of the optical function.  To obtain such estimates, we need the Raychaudhuri equation
\eq{\label{rayeqn} e_4(\chi_{ab})&=-\sum_c\chi_{ac}\chi_{cb}+\Gamma_{\alpha\beta}^0e_4^\alpha e_4^\beta\chi_{ab}+\lra{R(e_4,e_a)e_4,e_b}.}
Here $\lra{R(X,Y)Z,W}:=\lra{D_XD_YZ-D_YD_XZ-D_{[X,Y]}Z,W}$ is the curvature tensor. In fact, since $2\lra{D_ae_4,e_4}=e_a\lra{e_4,e_4}=0$, we have
\fm{e_4(\chi_{ab})&=e_4\lra{D_ae_4,e_b}=\lra{D_4D_ae_4,e_b}+\lra{D_ae_4,D_4e_b}\\&=\lra{D_aD_4e_4,e_b}+\lra{D_{[e_4,e_a]}e_4,e_b}+\lra{R(e_4,e_a)e_4,e_b}+\Gamma_{\alpha\beta}^0e_4^\alpha e_b^\beta\lra{D_ae_4,e_4}\\
&=\lra{D_a(\Gamma_{\alpha\beta}^0e_4^\alpha e_4^\beta e_4),e_b}-\sum_c\chi_{ac}\lra{D_{c}e_4,e_b}+\lra{R(e_4,e_a)e_4,e_b}\\
&=e_a(\Gamma_{\alpha\beta}^0e_4^\alpha e_4^\beta)\lra{ e_4,e_b}+\Gamma_{\alpha\beta}^0e_4^\alpha e_4^\beta\chi_{ab}-\sum_c\chi_{ac}\chi_{cb}+\lra{R(e_4,e_a)e_4,e_b}\\&=\Gamma_{\alpha\beta}^0e_4^\alpha e_4^\beta\chi_{ab}-\sum_c\chi_{ac}\chi_{cb}+\lra{R(e_4,e_a)e_4,e_b}.}

From \eqref{rayeqn}, we can compute $e_4(\chi_{11}-\chi_{22})$, $e_4(\chi_{12})$ and $e_4(\tr\chi)$. Note that 
\fm{\sum_c\chi_{1c}\chi_{c1}-\sum_c\chi_{2c}\chi_{c2}&=\chi_{11}^2-\chi_{22}^2=\tr\chi(\chi_{11}-\chi_{22}),}
\fm{\sum_c\chi_{1c}\chi_{c2}=\sum_c\chi_{2c}\chi_{c1}&=\chi_{11}\chi_{12}+\chi_{12}\chi_{22}=\chi_{12}\tr\chi,}
\fm{\sum_c\chi_{1c}\chi_{c1}+\sum_c\chi_{2c}\chi_{c2}&=\chi_{11}^2+\chi_{22}^2+2\chi_{12}^2=\frac{1}{2}(\tr\chi)^2+\frac{1}{2}(\chi_{11}-\chi_{22})^2+2\chi_{12}^2.}
As for the curvature tensor, we have the following lemma.
\lem{\label{l3curest} In $\Omega_T$, we have 
\fm{\lra{R(e_4,e_a)e_4,e_b}&=e_4(f_{ab})+\frac{1}{2}e_4^\alpha e_a^\beta e_4^\mu e_b^\nu\partial_\beta\partial_\nu g_{\alpha\mu}+O(\eps^2 t^{-3+C\eps})}
where\fm{f_{ab}:=\frac{1}{2}(e_a^\beta  e_b^\nu e_4(g_{\beta\nu})-e_a^\beta e_4^\mu e_b( g_{\beta\mu})) -\frac{1}{2}e_4^\alpha e_a(g_{\alpha\nu}) e_b^\nu=O(\eps t^{-2+C\eps}).}
Moreover,
\fm{\lra{R(e_4,e_1)e_4,e_1}-\lra{R(e_4,e_2)e_4,e_2}&=e_4(f_{11}-f_{22})+O(\eps t^{-3+C\eps}),}
\fm{\lra{R(e_4,e_1)e_4,e_2}&=e_4(f_{12})+O(\eps t^{-3+C\eps}),}
\fm{\lra{R(e_4,e_1)e_4,e_1}+\lra{R(e_4,e_2)e_4,e_2}&=e_4(\tr f-\frac{1}{2}e_4^\alpha e_4^\mu e_3(g_{\alpha\mu}))+O(\eps^2t^{-3+C\eps}).}
}
\begin{proof}
We have $\lra{R(e_4,e_a)e_4,e_b}=e_4^\alpha e_a^\beta e_4^\mu e_b^\nu R_{\alpha\beta\mu\nu}$ where $R_{\alpha\beta\mu\nu}$ is  given by \fm{R_{\alpha\beta\mu\nu}&:=\lra{R(\partial_\alpha,\partial_\beta)\partial_\mu,\partial_\nu}=g_{\sigma\nu}(\partial_\alpha\Gamma_{\beta\mu}^\sigma-\partial_\beta\Gamma_{\alpha\mu}^\sigma+\Gamma_{\beta\mu}^\delta\Gamma_{\alpha\delta}^\sigma-\Gamma_{\alpha\mu}^\delta\Gamma_{\beta\delta}^\sigma)\\
&=\partial_\alpha\Gamma_{\nu\beta\mu}-\partial_\beta\Gamma_{\nu\alpha\mu}-\Gamma_{\beta\mu}^\sigma\partial_\alpha g_{\sigma\nu}+\Gamma_{\alpha\mu}^\sigma\partial_\beta  g_{\sigma\nu}+\Gamma_{\beta\mu}^\delta\Gamma_{\nu\alpha\delta}-\Gamma_{\alpha\mu}^\delta\Gamma_{\nu\beta\delta}\\
&=\partial_\alpha\Gamma_{\nu\beta\mu}-\partial_\beta\Gamma_{\nu\alpha\mu}-\Gamma_{\beta\mu}^\delta\Gamma_{\delta\nu\alpha}+\Gamma_{\alpha\mu}^\delta\Gamma_{\delta\nu\beta}\\
&=\frac{1}{2}(\partial_\alpha\partial_\mu g_{\beta\nu}-\partial_\alpha\partial_\nu g_{\beta\mu}-\partial_\beta\partial_\mu g_{\alpha\nu}+\partial_\beta\partial_\nu g_{\alpha\mu})-\Gamma_{\beta\mu}^\delta\Gamma_{\delta\nu\alpha}+\Gamma_{\alpha\mu}^\delta\Gamma_{\delta\nu\beta}.}
Here for simplicity we set $\Gamma_{\alpha\mu\nu}:=g_{\alpha\beta} \Gamma^{\beta}_{\mu\nu}=\frac{1}{2}(\partial_\mu g_{\alpha\nu}+\partial_\nu g_{\alpha\mu}-\partial_\alpha g_{\mu\nu})$. Then
\fm{&\hspace{1.5em}\frac{1}{2}e_4^\alpha e_a^\beta e_4^\mu e_b^\nu(\partial_\alpha\partial_\mu g_{\beta\nu}-\partial_\alpha\partial_\nu g_{\beta\mu}-\partial_\beta\partial_\mu g_{\alpha\nu}+\partial_\beta\partial_\nu g_{\alpha\mu})\\
&=\frac{1}{2}e_4(\partial_\mu g_{\beta\nu}-\partial_\nu g_{\beta\mu}) e_a^\beta e_4^\mu e_b^\nu-\frac{1}{2}e_4^\alpha e_a^\beta e_4(\partial_\beta g_{\alpha\nu}) e_b^\nu+\frac{1}{2}e_4^\alpha e_a^\beta e_4^\mu e_b^\nu\partial_\beta\partial_\nu g_{\alpha\mu}\\
&=e_4(\frac{1}{2}(\partial_\mu g_{\beta\nu}-\partial_\nu g_{\beta\mu}) e_a^\beta e_4^\mu e_b^\nu-\frac{1}{2}e_4^\alpha e_a^\beta (\partial_\beta g_{\alpha\nu}) e_b^\nu)+\frac{1}{2}e_4^\alpha e_a^\beta e_4^\mu e_b^\nu\partial_\beta\partial_\nu g_{\alpha\mu}\\&\hspace{1em}+O(|\partial g|\sum_{k=1,2,4}|e_4(e_k^\alpha)|)\\
&=e_4(f_{ab})+\frac{1}{2}e_4^\alpha e_a^\beta e_4^\mu e_b^\nu\partial_\beta\partial_\nu g_{\alpha\mu}+O(\eps^2t^{-3+C\eps}).}
To finish the proof of the first part, we note that
\fm{\Gamma_{\beta\mu}^\delta\Gamma_{\delta\nu\alpha}&=g^{\sigma\delta}\Gamma_{\sigma\beta\mu}\Gamma_{\delta\nu\alpha}=\frac{1}{4}g^{\sigma\delta}(\partial_\beta g_{\sigma\mu}+\partial_\mu g_{\beta\sigma}-\partial_\sigma g_{\beta\mu})(\partial_\alpha g_{\delta\nu}+\partial_\nu g_{\alpha\delta}-\partial_\delta g_{\alpha\nu}).}
By \eqref{gpaformula}, we have
\fm{e_4^\alpha e_a^\beta e_4^\mu e_b^\nu\Gamma^\delta_{\beta\mu}\Gamma_{\delta\nu\alpha}&=\frac{1}{4}g^{\sigma\delta}\partial_\sigma g\partial_\delta g+\sum_{k=1,2,4} O(1) e_k(g)\partial g\\
&=\frac{1}{4}\sum_{c=1,2} e_c(g) e_c(g)+\frac{1}{8}e_3(g)e_4(g)+\frac{1}{8}e_4(g)e_3(g)+O(\sum_{k=1,2,4}|e_k(g)||\partial g|)\\&=O(\eps t^{-2+C\eps}\cdot \eps t^{-1+C\eps})=O(\eps^2t^{-3+C\eps}).}
Similarly, we have $e_4^\alpha e_a^\beta e_4^\mu e_b^\nu\Gamma_{\alpha\mu}^\delta\Gamma_{\delta\nu\beta}=O(\eps^2t^{-3+C\eps})$.

To prove the second half, we only need to consider the term $\frac{1}{2}e_4^\alpha e_a^\beta e_4^\mu e_b^\nu\partial_\beta\partial_\nu g_{\alpha\mu}$. By Lemma \ref{l3betterest}, we have 
\fm{\frac{1}{2}e_4^\alpha e_1^\beta e_4^\mu e_2^\nu\partial_\beta\partial_\nu g_{\alpha\mu}&=\frac{1}{2}e_4^\alpha  e_4^\mu e_1^\beta e_2(\partial_\beta  g_{\alpha\mu})=O(\eps t^{-3+C\eps}),}
\fm{\frac{1}{2}e_4^\alpha e_1^\beta e_4^\mu e_1^\nu\partial_\beta\partial_\nu g_{\alpha\mu}-\frac{1}{2}e_4^\alpha e_2^\beta e_4^\mu e_2^\nu\partial_\beta\partial_\nu g_{\alpha\mu}&=\frac{1}{2}e_4^\alpha  e_4^\mu (e_1^\beta e_1(\partial_\beta g_{\alpha\mu})- e_2^\beta e_2 (\partial_\beta  g_{\alpha\mu}))=O(\eps t^{-3+C\eps}).}
Finally, note that 
\fm{\sum_ae_4^\alpha e_a^\beta e_4^\mu e_a^\nu\partial_\beta\partial_\nu g_{\alpha\mu}&=\frac{1}{2}e_4^\alpha e_4^\mu(g^{\beta\nu}-\frac{1}{2}e_3^\beta e_4^\nu-\frac{1}{2}e_4^\beta e_3^\nu) \partial_\beta\partial_\nu g_{\alpha\mu}\\
&=\frac{1}{2}e_4^\alpha e_4^\mu g^{\beta\nu} \partial_\beta\partial_\nu g_{\alpha\mu}-\frac{1}{2}e_4^\alpha e_4^\mu e_3^\beta e_4( \partial_\beta  g_{\alpha\mu})\\
&=-e_4( \frac{1}{2}e_4^\alpha e_4^\mu e_3^\beta\partial_\beta  g_{\alpha\mu})+O(\eps^2 t^{-3+C\eps}).}
We briefly explain how we obtain the third estimate here. If $F=F(u)$ is a function of $u$ which is a solution to \eqref{qwe}, then  by \eqref{gpaformula}
\fm{g^{\beta\nu}\partial_\beta\partial_\nu(F(u))&=F'(u)g^{\beta\nu}u_{\beta\nu}+F''(u)g^{\beta\nu} u_\beta u_\nu=0+F''(u)(\sum_ce_c(u)e_c(u)+e_3(u)e_4(u))\\&=O(\eps t^{-3+C\eps}).} We thus have $e_4^\alpha e_4^\mu g^{\beta\nu} \partial_\beta\partial_\nu g_{\alpha\mu}=O(\eps t^{-3+C\eps})$. To handle the other term, we note that \fm{e_4( \frac{1}{2}e_4^\alpha e_4^\mu e_3^\beta\partial_\beta  g_{\alpha\mu})-\frac{1}{2}e_4^\alpha e_4^\mu e_3^\beta e_4( \partial_\beta  g_{\alpha\mu})&=\frac{1}{2}e_4(e_4^\alpha e_4^\mu e_3^\beta)\partial_\beta  g_{\alpha\mu}=O(\eps^2 t^{-3+C\eps}).} 
\end{proof}
\rm \bigskip

Thus, it follows from \eqref{rayeqn} that
\eq{\label{rayeqn1} 
\left\{\begin{array}{l}\displaystyle e_4(\chi_{11}-\chi_{22})=-\tr\chi(\chi_{11}-\chi_{22})+\Gamma_{\alpha\beta}^0e_4^\alpha e_4^\beta(\chi_{11}-\chi_{22})+e_4(f_{11}-f_{22})+O(\eps t^{-3+C\eps}),\\[1em]
\displaystyle e_4(\chi_{12})=-\chi_{12}\tr\chi+\Gamma_{\alpha\beta}^0e_4^\alpha e_4^\beta\chi_{12}+e_4(f_{12})+O(\eps t^{-3+C\eps}),\\[1em]
\displaystyle e_4(\tr\chi)=-\frac{1}{2}(\tr\chi)^2-\frac{1}{2}(\chi_{11}-\chi_{22})^2-2\chi_{12}^2+\Gamma_{\alpha\beta}^0e_4^\alpha e_4^\beta\tr\chi\\[.5em]\displaystyle\hspace{4em}+e_4(\tr f-\frac{1}{2}e_4^\alpha e_4^\mu e_3(g_{\alpha\mu}))+O(\eps^2t^{-3+C\eps}).\end{array}\right.}
It turns out to be more convenient to work with \eqref{rayeqn1} instead of \eqref{rayeqn}.

\subsection{Continuity argument}\label{sec3.6}
\rm
Fix a geodesic $x(s)$ in $\mcl{A}$ with $x^0(0)\in H\cap\{t<T\}$. Since $\dot{x}^0(s)>0$ for all $s\geq 0$ and $\lim_{s\to \infty} x^0(s)=\infty$, there exists a unique $0<s_0<\infty$ such that $x^0(s_0)=T$. Also fix some $s_1\in[0,s_0]$. Our assumption is that for all $s\in[0,s_1]$, at $(t,x)=x(s)\in\Omega_T$ we have 
\eq{\label{ca} \max_{a,b=1,2}|\chi_{ab}-\delta_{ab}r^{-1}|\leq A t^{-2+B\eps}.}Here $A$ and $B$ are large constants which are independent of $T,\eps,s_1,s_0$ and the geodesic $x(s)$. In the derivation below, we always assume that the constants $C$ in the inequalities are given before we choose $A,B$, and that the constants $C$ are also independent of $T,\eps,s_1,s_0$ and  $x(s)$. Note that for $A,B\gg 1$, we have \eqref{ca} for $s_1=0$ by the next lemma.

\lem{\label{l3init} On $H$, we have $|\partial^2q|\lesssim t^{-1}$ and $\max_{a,b=1,2}|\chi_{ab}-\delta_{ab}r^{-1}|\lesssim t^{-2+C\eps}$.}
\begin{proof}
Recall from Section \ref{sec3.1} that  on $H$ we have
\fm{(-1-4g^{0i}\omega_i+4g^{ij}\omega_i\omega_j)q_t^2+(4g^{ij}\omega_i\omega_j-2g^{0i}\omega_i)q_t+g^{ij}\omega_i\omega_j=0.}
To compute $X_iq_t$ where $X_i=\partial_i+2\omega_i\partial_t$, we apply $X_i$ to the equation and then solve for $X_iq_t$. Then,
\fm{X_iq_t&=-\frac{q_t^2X_i(-1-4g^{0i}\omega_i+4g^{ij}\omega_i\omega_j)+q_tX_i(4g^{ij}\omega_i\omega_j-2g^{0i}\omega_i)+X_i(g^{ij}\omega_i\omega_j)}{2q_t(-1-4g^{0i}\omega_i+4g^{ij}\omega_i\omega_j)+4g^{ij}\omega_i\omega_j-2g^{0i}\omega_i}.}
Note that every term on the right hand side is known. The denominator is equal to $-2+O(|u|)$ on $H$, so it is nonzero for $\eps\ll 1$. In addition, we have $X_i\omega_j=O(r^{-1})=O(t^{-1})$ and $X_iu=O(|\partial u|)=O(\eps t^{-1})$, so $X_iq_t=O(t^{-1})$. Next, we have \fm{X_iq_j=X_i(-\omega_j-2\omega_jq_t)=-(\partial_i\omega_j) (1+2q_t)-\omega_iX_iq_t=O(t^{-1}).} By applying $\partial_t$ to the eikonal equation, we have
\fm{0&=2g^{\alpha\beta}q_\beta q_{t\alpha}+(\partial_tg^{\alpha\beta})q_\alpha q_\beta=2g^{0\beta}q_\beta q_{tt}+2g^{i\beta}q_\beta (X_iq_{t}-2\omega_iq_{tt})+(\partial_tg^{\alpha\beta})q_\alpha q_\beta.}And since $(q_t,q_i)=(-1,\omega)+O(|u|)$ on $H$, we have
\fm{q_{tt}&=-\frac{2g^{i\beta}q_\beta X_iq_{t}+(\partial_tg^{\alpha\beta})q_\alpha q_\beta}{2g^{0\beta} q_\beta-4g^{i\beta}\omega_iq_\beta}=-\frac{O(|\partial q|t^{-1}+\eps t^{-1}|\partial q|^2)}{-2q_t-4q_r+O(|u||\partial q|)}=O(t^{-1}).}
Finally we note that $q_{it}=X_iq_t-2\omega_iq_{tt}=O(t^{-1})$ and $q_{ij}=X_iq_j-2\omega_iq_{jt}=O(t^{-1})$.

We  move on to the estimates for $\chi$. By definition, we have \fm{\chi_{ab}&=\lra{D_ae_4,e_b}=(e_a(e_4^\alpha)+e_a^\mu e_4^\nu\Gamma_{\mu\nu}^\alpha)e_b^\beta g_{\alpha\beta}.}
As computed in Lemma \ref{l3betterest}, we have
\fm{e_a^\mu e_4^\nu\Gamma_{\mu\nu}^\alpha e_b^\beta g_{\alpha\beta}&=(\frac{1}{2}g^{\alpha\gamma}(\partial_\mu g_{\nu\gamma}+\partial_\nu g_{\mu\gamma})-\frac{1}{4}e_4^\alpha e_3(g_{\mu\nu})+O(\eps t^{-2+C\eps}))e_a^\mu e_4^\nu e_b^\beta g_{\alpha\beta}\\
&=\frac{1}{2}(e_a( g_{\nu\beta})e_4^\nu e_b^\beta g_{\alpha\beta}+e_4(g_{\mu\beta})e_a^\mu e_b^\beta g_{\alpha\beta})-\frac{1}{4} e_3(g_{\mu\nu})e_a^\mu e_4^\nu \lra{e_4,e_b}+O(\eps t^{-2+C\eps})\\&=O(\eps t^{-2+C\eps}).}
In addition, recall from Section \ref{sec3.1} that $q_i=\omega_iq_r$ on $H$. Since $e_a$ is tangent to $H$, on $H$ we have \fm{e_a(q_i)=e_a(\omega_iq_r)=e_a^jr^{-1}(\delta_{ij}-\omega_i\omega_j)q_r+\omega_ie_a(q_r)=e_a^ir^{-1}-\omega_iq_rr^{-1}e_a(r)+\omega_ie_a(q_r).}  Since $e_a$ is tangent to the $2$-sphere $\{t=t_0,q=q(t_0,x_0)\}=\{t=t_0,|x|=|x_0|\}$ at $(t_0,x_0)\in H$, we have $e_a(r)=e_a^i\omega_i=0$ on $H$. Thus, on $H$ we have  \fm{e_b^\gamma e_a(q_\gamma)&=e_b^ie_a(q_i)=\sum_ie_b^i(e_a^ir^{-1}-0+\omega_ie_a(q_r))\\&=r^{-1}g_{ij}e_a^ie_b^j-r^{-1}(g_{ij}-\delta_{ij})e_a^ie_b^j+0=r^{-1}\delta_{ab}+O(\eps t^{-2+C\eps}).}
It follows that
\fm{e_a(e_4^\alpha)&=e_a(\frac{L^\alpha}{L^0})=\frac{L^0e_a(2g^{\alpha\gamma}q_\gamma)-L^\alpha e_a(2g^{0\gamma}q_\gamma)}{(L^0)^2}=\frac{2(g^{\alpha\gamma}-e_4^\alpha g^{0\gamma})e_a(q_\gamma)}{L^0}+O(\eps t^{-2+C\eps}),}
\fm{e_a(e_4^\alpha)e_b^\beta g_{\alpha\beta}&=\frac{2(e_b^\gamma -\lra{e_4,e_b} g^{0\gamma})e_a(q_\gamma)}{-2q_t+O(|u||\partial q|)}+O(\eps t^{-2+C\eps})=\frac{2e_b^\gamma e_a(q_\gamma)}{2+O(|u|)}+O(\eps t^{-2+C\eps})\\
&=r^{-1}\delta_{ab}+O(\eps t^{-2+C\eps}).}
This finishes the proof.
\end{proof}
\rm \bigskip

To complete the continuity argument, we need to prove \eqref{ca} with $A$ replaced by $A/2$. We start with $\chi_{12}$ and $\chi_{11}-\chi_{22}$. By \eqref{rayeqn1}, we have
\fm{e_4(r^2(\chi_{12}-f_{12}))&=2re_4(r)(\chi_{12}-f_{12})+r^2e_4(\chi_{12}-f_{12})\\
&=2re_4(r)(\chi_{12}-f_{12})+r^2((-\tr\chi+\Gamma_{\alpha\beta}^0e_4^\alpha e_4^\beta)\chi_{12}+O(\eps t^{-3+C\eps}))\\
&=r(2e_4(r)-r\tr\chi+r\Gamma_{\alpha\beta}^0e_4^\alpha e_4^\beta)\chi_{12}-2re_4(r)f_{12}+O(\eps t^{-1+C\eps}).}
Recall that $e_4(r)=1+O(t^{-1+C\eps})$, $f_{12}=O(\eps t^{-2+C\eps})$ and $r\Gamma_{\alpha\beta}^0e_4^\alpha e_4^\beta=O(r|\partial g|)=O(\eps)$. By \eqref{ca}, we have $|2-r\tr\chi|\leq 2Art^{-2+B\eps}$. In conclusion, 
\fm{|e_4(r^2(\chi_{12}-f_{12}))|&\leq r(2Art^{-2+B\eps}+C\eps+Ct^{-1+C\eps})\cdot At^{-2+B\eps}+C\eps t^{-1+C\eps}\\
&\leq  CA^2t^{-2+2B\eps}+CA\eps t^{-1+B\eps}+CAt^{-2+(B+C)\eps}+C\eps t^{-1+C\eps}\\
&\leq CA^2t^{-2+2B\eps}+CA\eps t^{-1+B\eps}.}
By choosing $A,B\gg C$, we obtain the last inequality. On $H$, we have $|r^2(\chi_{12}-f_{12})|\leq Ct^{C\eps}$ by the previous lemma. Thus, by integrating $e_4(r^2(\chi_{12}-f_{12}))$ along the geodesic, we have
\fm{|r^2(\chi_{12}-f_{12})|&\leq C(x^0(0))^{C\eps}+CA^2(x^0(0))^{-1+2B\eps}+CA B^{-1} t^{B\eps}\\
&\leq Ct^{C\eps}+CA^2T_0^{-1+2B\eps}+CA B^{-1} t^{B\eps}.}
Since $T_0\gg \eps^{-1}$, we have $A^2T_0^{-1+2B\eps}\leq 1$ for $\eps\ll 1$. In addition, by choosing $B\geq A$, we have  \fm{|\chi_{12}|\leq r^{-2}(|f_{12}|+Ct^{C\eps}+C+Ct^{B\eps})\leq Ct^{-2+B\eps}.}
Here $C$ is independent of $A$ and $B$, so if we choose $A\geq 4C$, we obtain with $|\chi_{12}|\leq \frac{1}{4}At^{-2+B\eps}$. The proof for $|\chi_{11}-\chi_{22}|\leq \frac{1}{4}At^{-2+B\eps}$ is essentially the same.

To finish the continuity argument, we need to prove that $|\tr\chi-2r^{-1}|\leq \frac{1}{4}At^{-2+B\eps}$. For $h=\tr\chi-\tr f=\tr\chi+O(\eps t^{-2+C\eps})$, by \eqref{ca} we have $h=2r^{-1}+O(At^{-2+B\eps})\sim 2r^{-1}$. Then, for $\eps \ll 1$, by the last equation in \eqref{rayeqn1} we have
\fm{&\hspace{1.5em}e_4(h^{-1})=-h^{-2}e_4(h)\\&=-h^{-2}(-\frac{1}{2}(\tr\chi)^2+\Gamma_{\alpha\beta}^0 e_4^\alpha e_4^\beta\tr\chi-\frac{1}{2}e_4(e_4^\alpha e_4^\beta e_3(g_{\alpha\beta}))+O(\eps^2 t^{-3+C\eps}+(\chi_{11}-\chi_{22})^2+\chi_{12}^2))\\
&=-h^{-2}(-\frac{1}{2}h^2+\Gamma_{\alpha\beta}^0 e_4^\alpha e_4^\beta h-\frac{1}{2}e_4(e_4^\alpha e_4^\beta e_3(g_{\alpha\beta}))+O(\eps t^{-3+C\eps}+\eps^2 t^{-3+C\eps}+A^2t^{-4+2B\eps}))\\
&=\frac{1}{2}-\Gamma_{\alpha\beta}^0 e_4^\alpha e_4^\beta h^{-1}+\frac{1}{2}h^{-2}e_4^\alpha e_4^\beta e_4(e_3(g_{\alpha\beta}))+O(\eps t^{-1+C\eps}).}
In the last line we use the product rule and the estimate $e_4(e_4^\alpha)=O(\eps t^{-2+C\eps})$. In addition, we have
\fm{|h^{-1}-r/2|&=\frac{|2-r(\tr\chi-\tr f)|}{2h}\lesssim r(|2-r\tr\chi|+|r\tr f|)\lesssim At^{B\eps};}
by \eqref{l3.3fchrit}, we have \fm{\Gamma_{\alpha\beta}^0 e_4^\alpha e_4^\beta&=\frac{1}{2}g^{0\gamma}(e_4^\beta e_4(g_{\beta \gamma})+e_4^\alpha e_4( g_{\alpha\gamma}))-\frac{1}{4}e_4^0e_3(g_{\alpha\beta})e_4^\alpha e_4^\beta+O(\eps t^{-2+C\eps})\\&=-\frac{1}{4}e_3(g_{\alpha\beta})e_4^\alpha e_4^\beta+O(\eps t^{-2+C\eps}).}
Thus, we have 
\eq{\label{negheqn}e_4(h^{-1})&=\frac{1}{2}+\frac{1}{4}e_4^\alpha e_4^\beta e_3(g_{\alpha\beta})h^{-1}+\frac{1}{4}rh^{-1}e_4^\alpha e_4^\beta e_4(e_3(g_{\alpha\beta}))\\&\hspace{1em}+O(\eps t^{-1+C\eps}+h^{-1}\eps t^{-2+C\eps}+At^{B\eps}h^{-1}|e_4(e_3(g))|)\\
&=\frac{1}{2}+\frac{1}{4} h^{-1}e_4^\alpha e_4^\beta (e_3(g_{\alpha\beta})+re_4(e_3(g_{\alpha\beta})))+O(At^{1+B\eps}|e_4(e_3(g_{\alpha\beta}))|+\eps t^{-1+C\eps}).}
The next three lemmas are necessary for us to control $e_3(g_{\alpha\beta})+re_4(e_3(g_{\alpha\beta}))$ and $e_4(e_3(g_{\alpha\beta}))$.
\lem{\label{l3q2} Under the assumption \eqref{ca}, in $\Omega_T$ we have $|e_a(e_3(q))|+|e_a(\partial q)|\lesssim t^{-1+C\eps}$, $|e_a(\Omega_{ij}q)|\lesssim At^{-1+B\eps}|e_3(q)|+t^{-1+C\eps}$ and $|\partial^2q|\lesssim t^{C\eps}$.}
\begin{proof}We have (assuming $\{a,a'\}=\{1,2\}$)
\fm{e_4(e_a(e_3(q)))&=[e_4,e_a]e_3(q)+e_a(e_4(e_3(q)))=-\sum_b\chi_{ab}e_b(e_3(q))-e_a(\Gamma_{\mu\nu}^0 e_4^\mu e_4^\nu e_3(q))\\
&=-\sum_b\chi_{ab}e_b(e_3(q))-2\Gamma_{\mu\nu}^0 (\sum_b \chi_{ab}e_b^\mu+O(\eps t^{-2+C\eps})) e_4^\nu e_3(q)\\&\hspace{1em}-\Gamma_{\mu\nu}^0 e_4^\mu e_4^\nu e_a(e_3(q))-e_a(\Gamma_{\mu\nu}^0) e_4^\mu e_4^\nu e_3(q)\\
&=-(\chi_{aa}+\Gamma_{\mu\nu}^0 e_4^\mu e_4^\nu)e_a(e_3(q))-\chi_{12}e_{a'}(e_3(q))\\&\hspace{1em}-(2\Gamma_{\mu\nu}^0\sum_b \chi_{ab}e_b^\mu e_4^\nu+e_a(\Gamma_{\mu\nu}^0) e_4^\mu e_4^\nu+O(\eps t^{-2+C\eps}|\Gamma|))e_3(q).}
Since $\chi_{ab}=r^{-1}\delta_{ab}+O(At^{-2+B\eps})\sim r^{-1}$ for $\eps\ll_{A,B}1$, the last term is $O(\eps t^{-2+C\eps}|e_3(q)|)=O(\eps t^{-2+C\eps})$.  Then, 
\fm{&\hspace{1.5em}|e_4(re_a(e_3(q)))|=|e_4(r)e_a(e_3(q))+re_4(e_a(e_3(q)))|\\&\leq |(1+O(t^{-1+C\eps}))e_a(e_3(q))-r(\chi_{aa}+\Gamma_{\mu\nu}^0 e_4^\mu e_4^\nu)e_a(e_3(q))-r\chi_{12}e_{a'}(e_3(q))|+C\eps t^{-1+C\eps}\\
&\leq(|r^{-1}-\chi_{aa}|+|\Gamma_{\mu\nu}^0 e_4^\mu e_4^\nu|+O(t^{-2+C\eps}))|re_a(e_3(q))|+|r\chi_{12}e_{a'}(e_3(q))|+C\eps t^{-1+C\eps}\\
&\leq (At^{-2+B\eps}+C\eps t^{-1}+Ct^{-2+C\eps})|re_a(e_3(q))|+CA t^{-2+B\eps}|re_{a'}(e_3(q))|+C\eps t^{-1+C\eps}\\
&\leq C\eps t^{-1}\sum_b|re_b(e_3(q))|+C\eps t^{-1+C\eps}.}
In the last line, we choose $\eps\ll 1$ so that $C\eps t^{-1}\geq At^{-2+B\eps}+t^{-2+C\eps}$ for $t\geq T_0= \exp(\delta/\eps)$. Since $e_a$ is tangent to $H$, on $H$ we have $e_a(e_3(q))=e_a(2g^{0\alpha}q_\alpha)=O(|\partial^2q|+|e_a(g)\partial q|)=O(t^{-1})$ by Lemma \ref{l3init}. In conclusion,  if $(t,x)\in \Omega_T$ lies on a geodesic $x(s)$ in $\mcl{A}$, at $(t,x)$ we have
\fm{\sum_a|re_a(e_3(q))|&\leq \sum_a|re_a(e_3(q))|(x(0))+\int_{x^0(0)}^t C\eps \tau^{-1}\sum_a|re_a(e_3(q))|(\tau,\widetilde{x}(\tau))\ d\tau+Ct^{C\eps}\\
&\leq C+Ct^{C\eps}+\int_{x^0(0)}^t C\eps \tau^{-1}\sum_a|re_a(e_3(q))|(\tau,\widetilde{x}(\tau))\ d\tau.}
Here $(\tau,\widetilde{x}(\tau))$ is a reparametrization  of the geodesic $x(s)$. We conclude that $\sum_a|re_a(e_3(q))|\lesssim C t^{C\eps}$ by the Gronwall's inequality. In addition, in $\Omega_T$ we have
\fm{e_a(q_\alpha)&=e_a(\frac{1}{2}\lra{\partial_\alpha,e_4}e_3(q))=e_a(\frac{1}{2}e_4^\beta g_{\alpha\beta}e_3(q))\\
&=\frac{1}{2}e_a(e_4^\beta) g_{\alpha\beta}e_3(q)+\frac{1}{2}e_4^\beta e_a(g_{\alpha\beta})e_3(q)+\frac{1}{2}e_4^\beta g_{\alpha\beta}e_a(e_3(q))=O(t^{-1+C\eps}).}

Next we compute $e_a(\Omega_{ij}q)$. Note that \fm{\Omega_{ij}q&=\frac{1}{2}\lra{\Omega_{ij},e_4}e_3(q)=\frac{1}{2}(x_ig_{j\beta}-x_jg_{i\beta})e_4^\beta e_3(q)=\frac{1}{2}r(\omega_ig_{j\beta}e_4^\beta-\omega_jg_{i\beta}e_4^\beta) e_3(q).}
We have \fm{\omega_ig_{j\beta}e_4^\beta-\omega_jg_{i\beta}e_4^\beta=\omega_ie_4^j-\omega_je_4^i+O(|u|)=O(\sum_j|e_4^j-\omega_j|)+O(|u|)=O(t^{-1+C\eps}),}
so $r(\omega_ig_{j\beta}e_4^\beta-\omega_jg_{i\beta}e_4^\beta)e_a(e_3(q))=O(t^{-1+C\eps})$. In addition,
\fm{&\hspace{1.5em}e_a((x_ig_{j\beta}-x_jg_{i\beta})e_4^\beta)\\&=
(e_a^ig_{j\beta}-e_a^jg_{i\beta})e_4^\beta+(x_ig_{j\beta}-x_jg_{i\beta})e_a(e_4^\beta)+O(|e_a(g)|)\\
&=e_a^ie_4^j-e_a^je_4^i+(x_ig_{j\beta}-x_jg_{i\beta})\sum_b (\chi_{ab}e_b^\beta+O(\eps t^{-2+C\eps}))+O(|e_a(g)|+|u|)\\
&=e_a^ie_4^j-e_a^je_4^i+\sum_b \chi_{ab}(x_ie_b^j-x_je_b^i+O(r|u|))+O(\eps t^{-1+C\eps})\\
&=e_a^ie_4^j-e_a^je_4^i+r^{-1}(x_ie_a^j-x_je_a^i)+O(r(|\chi_{aa}-r^{-1}|+|\chi_{12}|))+O(\eps t^{-1+C\eps})\\
&=e_a^i(e_4^j-\omega_j)-e_a^j(e_4^i-\omega_i)+O(At^{-1+B\eps})+O(\eps t^{-1+C\eps})=O(At^{-1+B\eps}).}
By the product rule we obtain the second estimate.

Finally, we consider $\partial^2q$. Recall that $e_4^\alpha=L^\alpha/L^0$ and that $|\partial q|\sim|q_r|\sim|q_t|\sim e_3(q)$. By the characteristic ODE's, we have
\fm{e_4(q_\alpha)&=\frac{-(\partial_\alpha g^{\mu\nu})q_\mu q_\nu}{e_3(q)}=O(\eps t^{-1})e_3(q)}
and thus
\fm{\partial_\alpha(e_4(q_\beta))&=\frac{-\partial_\alpha((\partial_\beta g^{\mu\nu})q_\mu q_\nu)e_3(q)+(\partial_\beta g^{\mu\nu})q_\mu q_\nu\cdot 2\partial_\alpha(g^{0\gamma}q_\gamma)}{(e_3(q))^2}\\
&=\frac{-2(\partial_\beta g^{\mu\nu})q_\mu q_{\alpha\nu}e_3(q)+(\partial_\beta g^{\mu\nu})q_\mu q_\nu\cdot 2g^{0\gamma}q_{\alpha\gamma}}{(e_3(q))^2}+O(\eps t^{-1+C\eps})\\
&=O(|\partial g||\partial^2q|)+O(\eps t^{-1+C\eps})=O(\eps t^{-1}|\partial^2q|)+O(\eps t^{-1+C\eps}).}
In the second line, we take out those terms without $\partial^2q$ and control them using the estimates for $g$ and $\partial q$. In the last line, we use the estimate $|\partial q|\sim e_3(q)$.
Besides, we have
\fm{\partial_\alpha e_4^\beta&=\frac{\partial_\alpha (L^\beta)L^0-L^\beta\partial_\alpha(L^0)}{(L^0)^2}=\frac{2\partial_\alpha (g^{\beta\nu}q_\nu)-2e_4^\beta\partial_\alpha (g^{0\nu}q_\nu)}{e_3(q)}\\
&=\frac{2(g^{\beta\nu}-e_4^\beta  g^{0\nu})q_{\alpha\nu}}{e_3(q)}+O(|\partial g||\partial q|(e_3(q))^{-1})\\
&=\frac{(\sum_a 2e_a^\beta e_a^\nu+e_3^\beta e_4^\nu+e_4^\beta e_4^\nu)q_{\alpha\nu}}{e_3(q)}+O(\eps t^{-1})\\
&=\frac{2\sum_ae_a^\beta e_a(q_\alpha)+(e_3^\beta+e_4^\beta)e_4(q_\alpha)}{e_3(q)}+O(\eps t^{-1})=\frac{2\sum_ae_a^\beta e_a(q_\alpha)}{e_3(q)}+O(\eps t^{-1}).}
Thus, we have
\fm{e_4(q_{\alpha\beta})&=[e_4,\partial_\alpha]q_\beta+\partial_\alpha (e_4(q_\beta))=-\partial_\alpha (e_4^\mu)\partial_\mu(q_\beta)+\partial_\alpha (e_4(q_\beta))\\
&=O((e_3(q))^{-1}\sum_a|e_a(q_\beta)e_a(q_\alpha)|)+O(\eps t^{-1}|\partial^2q|)+O(\eps t^{-1+C\eps})\\
&=O(\eps t^{-1}|\partial^2q|)+O(\eps t^{-1+C\eps}+t^{-2+C\eps}).}
In the last line we use the estimate $e_3(q)\geq C^{-1}t^{-C\eps}$.
Since $\partial^2 q=O(t^{-1})$ on $H$, we conclude $\partial^2q=O(t^{C\eps})$ by the Gronwall's inequality.

\end{proof}

\lem{\label{l3ca1} Set $h_i:=r(\partial_i(ru)-q_iq_r^{-1}\partial_r(ru))$. Under the assumption \eqref{ca}, in $\Omega_T$ we have $|h_i|\lesssim \eps t^{C\eps}$, $|e_a(h_i)|\lesssim A\eps t^{-1+B\eps}$ and $e_a(ru)=\sum_i e_a(\omega_i)h_i$.}
\begin{proof}We have
\fm{h_i&=r(\omega_iu+ru_i-q_iq_r^{-1}u-q_iq_r^{-1}ru_r)=ruq_r^{-1}(q_r\omega_i-q_i)+r^2(u_i-q_iq_r^{-1}u_r)\\
&=(ru+r^2u_r)q_r^{-1}(q_r\omega_i-q_i)+r^2(u_i-\omega_iu_r)=(u+ru_r)q_r^{-1}\sum_j\omega_j\Omega_{ij}q+\sum_jx_j\Omega_{ji}u.}
Since $|u|+|u_r|\lesssim \eps t^{-1+C\eps}$, $|q_i-\omega_iq_r|\lesssim t^{-1+C\eps}$ and $|u_i-\omega_iu_r|\lesssim\eps t^{-2+C\eps}$, we obtain $|h_i|\lesssim\eps t^{C\eps}$. Moreover, \fm{e_a(x_j\Omega_{ij}u)=e_a^j\Omega_{ij}u+x_je_a(\Omega_{ij}u)=O(\eps t^{-1+C\eps}),}
\fm{e_a((u+ru_r)q_r^{-1}\omega_j\Omega_{ij}q)&=e_a(u+ru_r)q_r^{-1}\omega_j\Omega_{ij}q-(u+ru_r)q_r^{-2}e_a(q_r)\omega_j\Omega_{ij}q\\&\hspace{1em}+(u+ru_r)q_r^{-1}e_a(\omega_j)\Omega_{ij}q+(u+ru_r)q_r^{-1}\omega_je_a(\Omega_{ij}q)\\
&=O(\eps t^{-1+C\eps})+O(\eps |q_r|^{-1}|e_a(\Omega q)|)\\
&=O(\eps t^{-1+C\eps})+O(A\eps t^{-1+B\eps}\frac{e_3(q)}{q_r})=O(A\eps t^{-1+B\eps}).}
Here we apply many estimates such as $e_a(r)=O(1)$, $e_a(\omega_i)=O(r^{-1})$, $\Omega q=O(t^{C\eps})$, $q_r\geq C^{-1}t^{-C\eps}$ and etc. In particular, we apply $e_a(\Omega q)=O(At^{-1+B\eps}e_3(q)+t^{-1+C\eps})$ from the previous lemma. Thus, we have $e_a(h_i)=O(A\eps t^{-1+B\eps})$.

Finally, we have
\fm{\sum_i e_a(\omega_i)h_i&=\sum_{i,j} e_a^jr^{-1}(\delta_{ij}-\omega_i\omega_j)h_i\\&=\sum_{i} e_a^i(\partial_i(ru)-q_iq_r^{-1}\partial_r(ru))-\sum_{i,j} e_a^j\omega_i\omega_j (\partial_i(ru)-q_iq_r^{-1}\partial_r(ru))\\
&=e_a(ru)-e_a(q)q_r^{-1}\partial_r(ru)-\sum_{j} e_a^j\omega_j \sum_i(\omega_i\partial_i(ru)-\omega_iq_iq_r^{-1}\partial_r(ru))\\
&=e_a(ru).}
\end{proof}

\lem{\label{l3ca}Under the assumption \eqref{ca}, in $\Omega_T$ we have $|r^{-1}e_3(u)+e_4(e_3(u))|\lesssim \eps At^{-3+B\eps}$ and  $|e_4(e_3(u))|\lesssim \eps t^{-2}$.}
\begin{proof} The second inequality follows directly from the first one. To prove the first one, we note that for each function $F=F(t,x)$, we have
\fm{g^{\alpha\beta}\partial_{\alpha}\partial_{\beta} F&=(\sum_a e_a^\alpha e_a^\beta +\frac{1}{2}e_4^\alpha e_3^\beta +\frac{1}{2}e_3^\alpha e_4^\beta)\partial_{\alpha}\partial_{\beta} F\\&= \sum_a (e_a(e_a(F))-e_a(e_a^\alpha)F_\alpha)+e_4(e_3(F))-e_4(e_3^\alpha)F_\alpha\\
&=\sum_a (e_a(e_a(F))-(D_ae_a)F+e_a^\mu e_a^\nu \Gamma_{\mu\nu}^\alpha F_\alpha)+e_4 (e_3(F))-(D_4e_3)F+e_4^\mu e_3^\nu\Gamma_{\mu\nu}^\alpha F_\alpha.}
By \eqref{l3.3fchrit}, we have
\fm{e_a^\mu e_a^\nu \Gamma_{\mu\nu}^\alpha F_\alpha&=\frac{1}{2}g^{\alpha\beta} F_\alpha(e_a^\nu e_a( g_{\nu\beta})+e_a^\mu e_a(g_{\mu\beta}))-\frac{1}{4} e_3(g_{\mu\nu})e_a^\mu e_a^\nu  e_4(F)+O(\eps t^{-2+C\eps}  |\partial F|)\\
&=O(\eps t^{-2+C\eps}|\partial F|+\eps t^{-1}|e_4(F)|),}
\fm{e_4^\mu e_3^\nu\Gamma_{\mu\nu}^\alpha F_\alpha&=\frac{1}{2}g^{\alpha\beta}F_\alpha(e_3^\nu e_4( g_{\nu\beta})+e_4^\mu e_3( g_{\mu\beta}))-\frac{1}{4}e_4^\mu e_3^\nu  e_3(g_{\mu\nu})e_4(F)+O(\eps t^{-2+C\eps}|\partial F|)\\
&=\frac{1}{2}(\sum_a e_a^\beta e_a(F)+\frac{1}{2}e_3^\beta e_4(F)+\frac{1}{2}e_4^\beta e_3(F)) e_4^\mu e_3(g_{\mu\beta})+O(\eps t^{-2+C\eps}|\partial F|+\eps t^{-1}|e_4(F)|)\\
&=\frac{1}{4} e_3(F) e_4^\beta e_4^\mu e_3(g_{\mu\beta})+O(\eps t^{-2+C\eps}|\partial F|+\eps t^{-1}\sum_{k=1,2,4}|e_k(F)|).}
Moreover, since \fm{D_ae_a&=\lra{D_ae_a,e_{a'}}e_{a'}+\frac{1}{2}\lra{D_ae_a,e_4}e_3+\frac{1}{2}\lra{D_ae_a,e_3}e_4\\&=\lra{D_ae_a,e_{a'}}e_{a'}+\frac{1}{2}(-\chi_{aa})e_3+(-\frac{1}{2}\chi_{aa}+e_a^\mu e_a^\nu\Gamma_{\mu\nu}^0)e_4,\hspace{2em}a\neq a'\\
D_4e_3&=\sum_b\lra{D_4e_3,e_b}e_b+\frac{1}{2}\lra{D_4e_3,e_4}e_3+\frac{1}{2}\lra{D_4e_3,e_3}e_4\\
&=-2\sum_b\Gamma_{\mu\nu}^0e_4^\mu e_b^\nu e_b-\Gamma_{\mu\nu}^0e_4^\mu e_4^\nu e_3,}
we have
\fm{\sum_a(D_ae_a)F&=\lra{D_1e_1,e_2}e_2(F)+\lra{D_2e_2,e_1}e_1(F)-\frac{1}{2}(\tr\chi)(e_3(F)+e_4(F))+\sum_ae_a^\mu e_a^\nu\Gamma_{\mu\nu}^0 e_4(F)\\
&=\lra{D_1e_1,e_2}e_2(F)+\lra{D_2e_2,e_1}e_1(F)-\frac{1}{2}(\tr\chi)e_3(F)+O(t^{-1}|e_4(F)|)\\
&=\lra{D_1e_1,e_2}e_2(F)+\lra{D_2e_2,e_1}e_1(F)-r^{-1}e_3(F)+O(t^{-1}|e_4(F)|+At^{-2+B\eps}|e_3(F)|),\\(D_4e_3)F&=-2\sum_b\Gamma_{\mu\nu}^0e_4^\mu e_b^\nu e_b(F)-\Gamma_{\mu\nu}^0e_4^\mu e_4^\nu e_3(F)\\&=\frac{1}{4}e_3(g_{\alpha\beta})e_4^\alpha e_4^\beta e_3(F)+O(\eps t^{-1}\sum_{b}|e_b(F)|+\eps t^{-2+C\eps}|e_3(F)|).}
Here we use the assumption \eqref{ca} and $|e_3(u)|\lesssim|\partial u|\lesssim\eps t^{-1}$. In conclusion, we have
\fm{g^{\alpha\beta}\partial_{\alpha}\partial_{\beta} F&=\sum_a e_a(e_a(F))-\lra{D_1e_1,e_2}e_2(F)-\lra{D_2e_2,e_1}e_1(F)+e_4 (e_3(F))+r^{-1}e_3(F)\\&\hspace{1em}+O(t^{-1}|e_4(F)|+At^{-2+B\eps}|e_3(F)|)+O(\eps t^{-2+C\eps}|\partial F|+\eps t^{-1}\sum_{k=1,2,4}|e_k(F)|).}
By taking $F=u$, we obtain
\eq{\label{l3caf1}0=g^{\alpha\beta}\partial_\alpha\partial_\beta u&=\sum_a e_a(e_a(u))-\lra{D_1e_1,e_2}e_2(u)-\lra{D_2e_2,e_1}e_1(u)\\&\hspace{1em}+r^{-1}e_3(u)+e_4 (e_3(u))+O(A\eps t^{-3+B\eps}).}
In addition, note that 
\fm{e_4(e_3(F))+r^{-1}e_3(F)&=e_4(2g^{0\alpha}+e_4^\alpha) F_\alpha+(2g^{0\alpha}+e_4^\alpha)e_4(F_\alpha)+r^{-1}e_3(F)\\
&=O((|e_4(g^{0\alpha})|+|e_4(e_4^\alpha)|)|\partial F|+|e_4(F_\alpha)|+r^{-1}|e_3(F)|)\\
&=O(\eps t^{-2+C\eps}|\partial F|+|e_4(\partial F)|+r^{-1}|e_3(F)|).}
Thus, we have
\fm{&\hspace{1.5em}|\sum_a e_a(e_a(F))-\lra{D_1e_1,e_2}e_2(F)-\lra{D_2e_2,e_1}e_1(F)|\\&\lesssim |\partial^2F|+\eps t^{-2+C\eps}|\partial F|+r^{-1}|e_3(F)|+t^{-1}|e_4(F)|+At^{-2+B\eps}|e_3(F)|)+\eps t^{-1}\sum_{k=1,2,4}|e_k(F)|.}
When $F=r^{-1}$, the right hand side has an upper bound $Ct^{-3+C\eps}$. When $F=\omega_i$, the right hand side has an upper bound $Ct^{-2+C\eps}$. Here we choose $\eps\ll_{A,B}1$ so that  $At^{-2+B\eps}|e_3(r^{-1})|\lesssim At^{-4+B\eps}\lesssim t^{-3}$ and $At^{-2+B\eps}|e_3(\omega_i)|\lesssim At^{-3+B\eps}\lesssim t^{-2}$.

We set $U(t,x)=ru(t,x)$. Then, by the previous lemma,
\fm{e_a(u)&=e_a(r^{-1}U)=e_a(r^{-1})U+r^{-1}e_a(U)=e_a(r^{-1})U+r^{-1}\sum_i e_a(\omega_i)h_i,\\
e_a(e_a(u))&=e_a(e_a(r^{-1}))U+2e_a(r^{-1})\sum_i e_a(\omega_i)h_i+r^{-1}\sum_i e_a(e_a(\omega_i))h_i+r^{-1}\sum_i e_a(\omega_i)e_a(h_i)\\
&=e_a(e_a(r^{-1}))U+r^{-1}\sum_i e_a(e_a(\omega_i))h_i+O(A\eps t^{-3+B\eps}+\eps t^{-3+C\eps}).}
Thus, we have
\fm{&\hspace{1.5em}\sum_a e_a(e_a(u))-\lra{D_1e_1,e_2}e_2(u)-\lra{D_2e_2,e_1}e_1(u)\\
&=(\sum_ae_a(e_a(r^{-1}))-\lra{D_1e_1,e_2}e_2(r^{-1})-\lra{D_2e_2,e_1}e_1(r^{-1}))U\\
&\hspace{1em}+r^{-1}\sum_i (\sum_ae_a(e_a(\omega_i))-\lra{D_1e_1,e_2}e_2(\omega_i)-\lra{D_2e_2,e_1}e_1(\omega_i))h_i+O(A\eps t^{-3+B\eps}+\eps t^{-3+C\eps})\\
&=O(t^{-3+C\eps}|ru|+t^{-2+C\eps}r^{-1}|h_i|+A\eps t^{-3+B\eps}+\eps t^{-3+C\eps})=O(A\eps t^{-3+B\eps}).}
We finish the proof by this estimate and \eqref{l3caf1}.
\end{proof}\rm

We now finish the continuity argument. By writing $g'_{\alpha\beta}:=\frac{d}{du}|_{u=0}g^{\alpha\beta}(u)$, we have
\fm{e_3(g_{\alpha\beta})&=g'_{\alpha\beta}(u)e_3(u),\\e_4(e_3(g_{\alpha\beta}))&=g'_{\alpha\beta}(u)e_4(e_3(u))+g''_{\alpha\beta}(u)e_4(u)e_3(u)\\&=O(\eps t^{-2}+\eps t^{-2+C\eps}\cdot\eps t^{-1})=O(\eps t^{-2}),}and thus
\fm{e_3(g_{\alpha\beta})+re_4(e_3(g_{\alpha\beta}))&=g'_{\alpha\beta}(u)(e_3(u)+re_4(e_3(u)))+g''_{\alpha\beta}(u)e_4(u)e_3(u)\\
&=O(rA\eps t^{-3+B\eps}+r\eps t^{-2+C\eps}\cdot \eps t^{-1})=O(A\eps t^{-2+B\eps}).}
Thus, by \eqref{negheqn},
\fm{|e_4(h^{-1})-\frac{1}{2}|&\lesssim t\cdot A\eps t^{-2+B\eps}+At^{1+B\eps}\cdot\eps t^{-2}+\eps t^{-1+C\eps}\lesssim A\eps t^{-1+B\eps}.}
By the initial condition, on $H$ we have \fm{|h^{-1}-r/2|&=\frac{|2-r(\tr\chi-\tr f)|}{2h}\lesssim r(|2-r\tr\chi|+|r\tr f|)\lesssim t^{C\eps}}where the constants are known before we choose $A,B$. Now, suppose that $(t,x)\in\Omega_T$ lies on a geodesic $x(s)$ in $\mcl{A}$. At $x(0)$, we have $h^{-1}|_{x(0)}=r(x(0))/2+O((x^0(0))^{C\eps})$. Thus,
\fm{|h^{-1}|_{(t,x)}-\frac{1}{2}r(x(0))-\frac{1}{2}(t-x^0(0))|&\leq  |h^{-1}|_{(t,x)}-h^{-1}|_{x(0)}-\frac{1}{2}(t-x^0(0))|+Ct^{C\eps}\\&\lesssim\int_{x^0(0)}^t A \eps \tau^{-1+B\eps}\ d\tau+t^{C\eps}\lesssim B^{-1}At^{B\eps}+t^{C\eps}.} 
Also note that $r(x(0))-x^0(0)+t=q(t,x)+t=r+O(t^{C\eps})$ by Lemma \ref{l3qrt}. In conclusion, $|h^{-1}-r/2|\lesssim t^{C\eps}+B^{-1}At^{B\eps}$ at $(t,x)$. This  implies that $h^{-1}\sim r$ and 
\fm{|\tr\chi-\frac{2}{r}|&\leq|h-\frac{2}{r}|+C\eps t^{-2+C\eps}\lesssim |\frac{r-2h^{-1}}{rh^{-1}}|+C\eps t^{-2+C\eps}\\
&\leq Cr^{-2}(Ct^{C\eps}+CB^{-1}At^{B\eps})+C\eps t^{-2+C\eps}\leq Ct^{-2+C\eps}+CB^{-1}At^{-2+B\eps}.}
By choosing $B\geq A\gg_C1$, we conclude that $|\tr\chi-2/r|\leq \frac{1}{4}At^{-2+B\eps}$. This finishes the continuity argument as we have proved that \eqref{ca} holds with $A$ replaced by $A/4$.
\rm

\section{Derivatives of the optical function}\label{s4}\rm
In this section, we aim to prove that $q$ is smooth in $\Omega$, where smoothness is defined in Section \ref{sec2.5}. Our main result is the following proposition.
\prop{\label{props4}The optical function $q=q(t,x)$ constructed in Proposition \ref{props3} is a smooth function in $\Omega$. Moreover, in $\Omega$, we have $Z^Iq=O(\lra{q}t^{C\eps})$ and $Z^I\Omega_{ij} q=O(t^{C\eps})$ for each multiindex $I$ and $1\leq i<j\leq 3$.}\rm
\bigskip

In Section \ref{s5.1}, we define the commutator coefficients $\xi_{**}^*$ with respect to the null frame $\{e_k\}$, and derive several differential equations for $\xi$ and their derivatives. Note that the estimates for these $\xi$ would imply the estimates for $q$ in Proposition \ref{props4}. We also define a weighted null frame $\{V_k\}$ which will be used in the rest of this paper. In Section \ref{sec4.2}, we focus on the estimates for $q$ on the surface $H$ where the initial data of $q$ are assigned. In Section \ref{sec4.3new}, we prove Proposition \ref{prop5.5} which gives several important estimates for $\xi$. Here we make use of  the differential equations and the  estimates  on $H$ proved in  the first two subsections. Finally, in Section \ref{sec5.4}, we conclude the proof of Proposition \ref{props4} by applying Proposition \ref{prop5.5}.

To end this section, in Section \ref{sec5.5} we  derive two equations \eqref{secend} and \eqref{secend2} for $e_3(u)$ and $e_3(q)$, respectively. In these two equations, we have estimates for all derivatives of the remainder terms. While they are not related to the proof of Proposition \ref{props4}, they will be very useful  in the next section.

\subsection{Setup}\label{s5.1}
As a convention, we use $k,l$ to denote a number in $\{1,2,3,4\}$, and we use $a,b,c$ to denote a number in $\{1,2\}$. For a finite sequence of indices $K=(k_1,\dots,k_m)$, we set $|K|=m$, $n_{K,k}=\{j:\ k_j=k\}$ and $e_K=e_{k_1}e_{k_2}\cdots e_{k_m}$. 

\subsubsection{Commutator coefficients}
We define \fm{\xi_{kl}^a=\lra{[e_{k},e_l],e_a},\ a=1,2;\ \xi_{kl}^3=\frac{1}{2}\lra{[e_{k},e_l],e_4},\ \xi_{kl}^4=\frac{1}{2}\lra{[e_{k},e_l],e_3}.}
By \eqref{xvfformula} we have $[e_{k_1},e_{k_2}]=\xi_{k_1k_2}^le_l$. Thus these $\xi_{**}^*$'s are also called \emph{commutator coefficients} in this paper.

We now derive several equations for $\xi$. Note that $\xi_{k_1k_2}^l=-\xi_{k_2k_1}^l$ (so $\xi_{kk}^l=0$) and that $\xi_{kl}^3=\xi_{kl}^4$ since $[e_k,e_l]$ never contains $\partial_t$. Thus, we only need to study those $\xi_{k_1k_2}^l$'s with $k_1<k_2$ and $l\leq 3$.

We start with $[e_3,e_4]$. By Lemma \ref{l3.3} we have 
\fm{\lra{[e_3,e_4],e_4}&=\lra{D_3e_4-D_4e_3,e_4}=-\lra{D_4e_3,e_4}=\lra{e_3,D_4e_4}=2\Gamma_{\alpha\beta}^0e_4^\alpha e_4^\beta,} so $\xi_{34}^3=\Gamma_{\alpha\beta}^0e_4^\alpha e_4^\beta$. For $\xi_{34}^a$, we have the following equation
\fm{e_4(\xi_{34}^a)&=e_4(\lra{D_3e_4-D_4e_3,e_a})=e_4(\lra{D_3e_4,e_a})+e_4(\lra{e_3,D_4e_a})\\&=\lra{D_4D_3e_4,e_a}+\lra{D_3e_4,D_4e_a}+2e_4(\Gamma_{\alpha\beta}^0e_4^\alpha e_a^\beta)\\
&=\lra{D_3D_4e_4,e_a}+\lra{D_{[e_4,e_3]}e_4,e_a}+\lra{R(e_4,e_3)e_4,e_a}+\lra{D_3e_4,(\dots)e_4}+2e_4(\Gamma_{\alpha\beta}^0e_4^\alpha e_a^\beta)\\
&=\lra{D_3((\Gamma_{\alpha\beta}^0e_4^\alpha e_4^\beta)e_4),e_a}-\xi_{34}^l\lra{D_le_4,e_a}+\lra{R(e_4,e_3)e_4,e_a}+2e_4(\Gamma_{\alpha\beta}^0e_4^\alpha e_a^\beta)\\
&=-\chi_{ba}\xi_{34}^b+\lra{R(e_4,e_3)e_4,e_a}+2e_4(\Gamma_{\alpha\beta}^0e_4^\alpha e_a^\beta).}

Next we consider $[e_a,e_4]$. From Lemma \ref{l3.4}, we have $\xi_{a4}^b=\chi_{ab}$ and $\xi_{a4}^3=0$. Thus we have the Raychaudhuri equation \fm{e_4(\chi_{ab})&=\Gamma_{\alpha\beta}^0e_4^\alpha e_4^\beta\chi_{ab}-\sum_c\chi_{ac}\chi_{cb}+\lra{R(e_4,e_a)e_4,e_b}.}

Next we consider $[e_1,e_2]$. Note that $\xi_{12}^3=0$ as $\lra{[e_1,e_2],e_4}=0$. For $\xi_{12}^a$, we have $\xi_{12}^1=\lra{D_1e_2-D_2e_1,e_1}=\lra{D_1e_2,e_1}$ and $\xi_{12}^2=\lra{D_1e_2-D_2e_1,e_2}=-\lra{D_2e_1,e_2}=\lra{D_2e_2,e_1}$. So, $\xi_{12}^a=\lra{D_ae_2,e_1}$ and
\fm{e_4(\xi_{12}^a)&=e_4(\lra{D_ae_2,e_1})=\lra{D_4D_ae_2,e_1}+\lra{D_ae_2,D_4e_1}\\&=\lra{D_aD_4e_2,e_1}+\lra{D_{[e_4,e_a]}e_2,e_1}+\lra{R(e_4,e_a)e_2,e_1}+\Gamma^0_{\alpha\beta}e_4^\alpha e_1^\beta\lra{D_ae_2,e_4}\\
&=\Gamma_{\alpha\beta}^0e_4^\alpha e_2^\beta\chi_{a1}-\Gamma^0_{\alpha\beta}e_4^\alpha e_1^\beta\chi_{a2}-\chi_{ac}\xi_{12}^c+\lra{R(e_4,e_a)e_2,e_1}.}

We end with $[e_a,e_3]$. Note that
\fm{\xi_{a3}^3&=\frac{1}{2}\lra{D_ae_3-D_3e_a,e_4}=-\frac{1}{2}\lra{e_3,D_ae_4}+\frac{1}{2}\lra{e_a,D_3e_4}\\
&=-\frac{1}{2}\xi_{a4}^4-\frac{1}{2}\lra{e_3,D_4e_a}+\frac{1}{2}\xi_{34}^a+\frac{1}{2}\lra{e_a,D_4e_3}=-\lra{e_3,D_4e_a}+\frac{1}{2}\xi_{34}^a\\&=-2\Gamma_{\alpha\beta}^0e_4^\alpha e_a^\beta+\frac{1}{2}\xi_{34}^a,}
\fm{\xi_{a3}^a&=\lra{D_ae_3-D_3e_a,e_a}=\lra{D_ae_3,e_a}=\chi_{aa}+\lra{D_a(2g^{0\alpha}\partial_\alpha),e_a}\\
&=\chi_{aa}+2e_a(g^{0\alpha})g_{\alpha\beta}e_a^\beta+2g^{0\alpha}e_a^\beta\Gamma_{\beta\alpha}^\mu g_{\mu\nu}e_a^{\nu}.}
For $\xi_{a3}^b$ where $a\neq b$, we have
\fm{e_4(\xi_{a3}^b)&=e_4(\lra{D_ae_3-D_3e_a,e_b})=e_4(\chi_{ab}+\lra{D_a(2g^{0\alpha}\partial_\alpha),e_b}-\lra{D_3e_a,e_b})\\
&=e_4(\chi_{ab}+2e_a(g^{0\alpha})g_{\alpha\beta}e_b^\beta+2g^{0\alpha}e_a^\beta\Gamma_{\beta\alpha}^\mu g_{\mu\nu}e_b^{\nu})-\lra{D_4D_3e_a,e_b}-\lra{D_3e_a,D_4e_b}\\
&=e_4(\chi_{ab}+2e_a(g^{0\alpha})g_{\alpha\beta}e_b^\beta+2g^{0\alpha}e_a^\beta\Gamma_{\beta\alpha}^\mu g_{\mu\nu}e_b^{\nu})-\lra{D_3D_4e_a,e_b}-\lra{D_{[e_4,e_3]}e_a,e_b}\\
&\hspace{1em}-\lra{R(e_4,e_3)e_a,e_b}-\Gamma_{\alpha\beta}^0e_4^\alpha e_b^\beta\lra{D_3e_a,e_4}\\
&=(e_4+\Gamma_{\mu\nu}^0e_4^\mu e_4^\nu)(\chi_{ab}+2e_a(g^{0\alpha})g_{\alpha\beta}e_b^\beta+2g^{0\alpha}e_a^\beta\Gamma_{\beta\alpha}^\mu g_{\mu\nu}e_b^{\nu})-\Gamma_{\mu\nu}^0e_4^\mu e_4^\nu\xi_{a3}^b-\sum_c\xi_{34}^c\xi_{ab}^c\\
&\hspace{1em}-\lra{R(e_4,e_3)e_a,e_b}-\Gamma_{\alpha\beta}^0e_4^\alpha e_a^\beta\xi_{34}^b+\Gamma_{\alpha\beta}^0e_4^\alpha e_b^\beta\xi_{34}^a.}

Given $\xi$, we can express $e_{k_1}(e_{k_2}^\alpha)$ in terms of $e_*^*$ and $\xi_{**}^*$. In fact, the formulas for $e_4(e_k^\alpha)$ follow from Lemma \ref{l3.3}. Besides,
\fm{e_k(e_4^\alpha)&=[e_k,e_4]^\alpha+e_4(e_k^\alpha)=\xi_{k4}^le_l^\alpha+e_4(e_k^\alpha),\\
e_k(e_3^\alpha)&=e_k(e_4^\alpha)+2e_k(g^{0\alpha}),\\
e_3(e_k^\alpha)&=[e_3,e_k]^\alpha+e_k(e_3^\alpha)=\xi_{3k}^l e_l^\alpha+e_k(e_3^\alpha),\\
e_a(e_b^\alpha)&=(D_ae_b)^\alpha-e_a^\mu e_b^\nu\Gamma_{\mu\nu}^\alpha\\&=\sum_c\lra{D_ae_b,e_c}e_c^\alpha+\frac{1}{2}\lra{D_ae_b,e_3}e_4^\alpha+\frac{1}{2}\lra{D_ae_b,e_4}e_3^\alpha-e_a^\mu e_b^\nu\Gamma_{\mu\nu}^\alpha\\
&=-\sum_c\xi_{bc}^ae_c^\alpha-\frac{1}{2}\chi_{ab}(e_4^\alpha+e_3^\alpha)-\lra{e_b,D_a(g^{0\beta}\partial_\beta)}e_4^\alpha-e_a^\mu e_b^\nu\Gamma_{\mu\nu}^\alpha\\
&=-\sum_c\xi_{bc}^ae_c^\alpha-\frac{1}{2}\chi_{ab}(e_4^\alpha+e_3^\alpha)-(e_b^\mu g_{\mu\beta}e_a(g^{0\beta})+e_b^\mu g_{\mu\nu} g^{0\beta}e_a^\sigma\Gamma_{\sigma\beta}^\nu)e_4^\alpha-e_a^\mu e_b^\nu\Gamma_{\mu\nu}^\alpha. }
\bigskip

\subsubsection{A weighted null frame}A new frame $\{V_k\}$ defined below turns out to be very useful in this section.
\defn{\rm We define a new frame $\{V_k\}_{k=1}^4$ by $V_a=re_a$ for $a=1,2$ and $V_3=(3R-r+t)e_3$ and $V_4=te_4$. We call $\{V_k\}_{k=1}^4$ a \emph{weighted null frame}, since  $V_k$ is  a multiple of $e_k$ for each $k$.

As usual, for each multiindex $K=(k_1,\dots,k_m)$ with $k_*\in\{1,2,3,4\}$, we define $V^I=V_{k_1}\cdots V_{k_m}$ as the product of $|I|$ vector fields. }\rm

It is easy to see that
\eq{\label{derqvz}
\left\{\begin{array}{l}
V_4=t(t+r)^{-1}S+(t+r)^{-1}t\omega_j\Omega_{0j}+t(e_4^i-\omega_i)\partial_i,\\
V_3=(3R-r+t)r^{-1}V_4+2g^{0\alpha}(3R-r+t)\partial_\alpha,\\
V_a=V_a(r)\omega_i\partial_i+e_a^i\omega_j\Omega_{ji};
\end{array}\right.}
\eq{\label{derqzv}
Z=r^{-1}\sum_a\lra{Z,e_a}V_a+\frac{1}{2}t^{-1}\lra{Z,e_3}V_4+\frac{1}{2}(3R-r+t)^{-1}\lra{Z,e_4}V_3.}
These formulas illustrate the connection between the weighted null frame and the commuting vector fields. 

Here we briefly explain why we work with $\{V_k\}$. First,  we note that  \fm{Z\approx \sum_{k\neq 3}O(t)e_k+O(\lra{r-t})e_3\approx\sum_k O(1)V_k.} If we work with a usual null frame, then in order to prove $Z^Iq=O(\lra{q}t^{C\eps})$, we might need to prove
\eq{\label{sec4.1fe1}|e_I(q)|\lesssim \lra{r-t}^{1-n_{I,3}}t^{-n_{I,1}-n_{I,2}-n_{I,4}+C\eps}}where $e_I$ and $n_{I,*}$ are defined at the beginning of Section \ref{s5.1}. In contrast, if we work with a weighted null frame, then we can prove \eq{\label{sec4.1fe2}|V^Iq|\lesssim \lra{r-t}t^{C\eps}.}
Since \eqref{sec4.1fe1} is much more complicated than \eqref{sec4.1fe2}, we expect the proof to be much simpler if we choose to work with the new weighted null frame.

Next, to prove an estimate for $V^Iq$, we need to compute \fm{e_4(V^Iq)=t^{-1}\sum_{I=(J,j,J')}V^J[V_4,V_j]V^{J'}q.}
Since $V_k$ is a multiple of $e_k$ for each $k$, we expect $[V_4,V_k]$ to be relatively simple. If we choose to work with the commuting vector fields defined in \eqref{vf}, then we need to compute either $[e_4,Z]$ or $[V_4,Z]$. Neither of these two terms has a simple form.

\subsection{Estimates on $H$}\label{sec4.2}
We start with the estimates on the surface $H$. Recall that  the vector fields $X_i=\partial_i+2\omega_i\partial_t$ are tangent to $H$ for $i=1,2,3$. For a multiindex $I=(i_1,\dots,i_m)$ where $i_j\in\{1,2,3\}$, we write $X^I=X_{i_1}\cdots X_{i_m}$ and $|I|=m$.

In this subsection, we keep using the convention stated in Section \ref{sec2.5}. 

We have the following pointwise estimate.  We ask our readers to compare this lemma with Lemma \ref{l2.1}.
\lem{\label{derqlem1} Suppose that $F=F(t,x)$ is a smooth function whose domain is contained in $\{(t,x)\in\R^{1+3}:\ r\sim t\gtrsim 1\}$. Then, for nonnegative integers $m,n$, we have
\fm{\sum_{|I|=m,\ |J|=n}|Z^IX^JF|\lesssim\lra{r-t}^{-n}\sum_{|I|\leq m+n}|Z^IF|.} }
\begin{proof}We induct first on $m+n$ and then on $n$. There is nothing to prove when $n=0$. If $m=0$ and $n=1$, we  simply apply Lemma \ref{l2.1}. In general, we fix multiindicies $I,J$ such that $|I|=m$ and $|J|=n$, such that $m+n>1$ and $n>0$. We can write $X^J=X^{J'}X_j$. Then, by our induction hypotheses, we have \fm{|Z^IX^JF|&\leq |Z^IX^{J'}\partial_jF|+|Z^IX^{J'}(\omega_j\partial_tF)|\\
&\lesssim \lra{r-t}^{1-n}\sum_{ |K|\leq n+m-1}(|Z^K\partial F|+|Z^K(\omega_j\partial_t F)|).}
Since $Z^K\omega=O(1)$ for each $|K|\geq 0$, by the Leibniz's rule we have 
\fm{|Z^IX^JF|&\lesssim \lra{r-t}^{1-n}\sum_{ |K|\leq n+m-1}|Z^K\partial F|\lesssim \lra{r-t}^{1-n}\sum_{ |K|\leq n+m-1}|\partial Z^K F|\\&\lesssim\lra{r-t}^{-n}\sum_{ |K|\leq n+m}| Z^K F|.}
In the second inequality here we use the commutation property $[Z,\partial]=C\partial$. 
\end{proof}
\rm

The next lemma is a variant of Lemma \ref{lcommin} with $Z$ replaced by $X$. Note that we do not need to  assume that $(m_{0}^{\alpha\beta})$ satisfies the null condition defined in Section \ref{sec2}.
\lem{\label{derqlem2} Fix two functions $\phi(t,x)$ and $\psi(t,x)$. Let $(m_0^{\alpha\beta})$ be a constant matrix. Then,  \fm{X_i(m_0^{\alpha\beta}\phi_\alpha\psi_\beta)=m_0^{\alpha\beta}(\partial_\alpha X_i\phi)\psi_\beta+m_0^{\alpha\beta}\phi_\alpha(\partial_\beta X_i\psi)+r^{-1}\sum_{\alpha,\beta}f_0\phi_\alpha \psi_\beta.}
Here $f_0$ denotes a polynomial of $\omega$; we allow $f_0$ to vary from line to line.}
\begin{proof}
We have $[X_i,\partial_\alpha]=-2(\partial_\alpha\omega_i)\partial_t$. By the Leibniz's rule, we have
\fm{X_i(m_0^{\alpha\beta}\phi_\alpha\psi_\beta)&=m_0^{\alpha\beta}(\partial_\alpha X_i\phi)\psi_\beta+m_0^{\alpha\beta}\phi_\alpha(\partial_\beta X_i\psi)-2m_0^{\alpha\beta}(\partial_\alpha\omega_i)\phi_t\psi_\beta-2m_0^{\alpha\beta}(\partial_\beta\omega_i)\psi_t\phi_\alpha\\
&=m_0^{\alpha\beta}(\partial_\alpha X_i\phi)\psi_\beta+m_0^{\alpha\beta}\phi_\alpha(\partial_\beta X_i\psi)\\&\hspace{1em}-2r^{-1}[m_0^{j\beta}(\delta_{ji}-\omega_j\omega_i)\phi_t\psi_{\beta}+m_0^{\alpha j}(\delta_{ji}-\omega_{j}\omega_i)\psi_t\phi_\alpha]\\
&=m_0^{\alpha\beta}(\partial_\alpha X_i\phi)\psi_\beta+m_0^{\alpha\beta}\phi_\alpha(\partial_\beta X_i\psi)+r^{-1}\sum_{\alpha,\beta}f_0\phi_\alpha \psi_\beta.}
\end{proof}
\rm
Using the previous two lemmas, we can now prove the estimates for $Z^Iq$ on $H$. In the next two lemmas, $\Omega^I$ denotes the product of $|I|$ vector fields in $\{\Omega_{12},\Omega_{23},\Omega_{13}\}$. In the rest of Section \ref{sec4.2}, we would use $\Omega$ to denote any vector field in $\{\Omega_{12},\Omega_{23},\Omega_{13}\}$ instead of the region. There should be no confusion as we focus on estimates on $H$.

\lem{\label{derqlh} On $H$, for all multiindices $I$, we have $Z^{I}q=O(\lra{q}t^{C\eps})$ and $Z^I\Omega q=O(t^{C\eps})$. }
\begin{proof}For convenience, we set
\fm{\mcl{O}_{m,n,p}&=\mcl{O}_{m,n,p}(t,x):=\sum_{|I|=m,\ |J|=n,\ |K|=p}|Z^IX^J\Omega^Kq|.}
On $H$, we claim that \fm{\mcl{O}_{m,n,0}\lesssim \lra{q}^{1-n}t^{C\eps},\ \forall m,n\geq 0;\hspace{2em}
\mcl{O}_{m,n,p}\lesssim \lra{q}^{-n}t^{C\eps},\ \forall m,n\geq 0,\ p>0.}

We first assume $m=0$. Since $\Omega$ and $X$ are tangent to $H$ and since $q|_H=r-t$, we have $X^J\Omega^Kq=X^J\Omega^K(r-t)$  for all multiindices $J,K$. If $|K|>0$, we have $X^J\Omega^K(r-t)=0$; if  $|J|>0$, we have $X^J(r-t)=O(r^{1-|J|})=O(\lra{q}^{1-|J|})$. Then, on $H$ we have $\mcl{O}_{0,0,0}=|q|$, $\mcl{O}_{0,n,p}=0$ for $p>0$, and $\mcl{O}_{0,n,0}=O(\lra{q}^{1-n})$ for  $n>0$. So the claim is true for $m=0$.

In general, we fix $(m,n,p)$ with $m>0$. Suppose we have proved \eq{\label{derqlhih}\mcl{O}_{m',n',0}\lesssim \lra{q}^{1-n'}t^{C\eps},\quad &\forall m',n'\geq 0 \text{ such that }m'+n'<m+n+p\\&\text{or }m'+n'=m+n+p,\ m'<m;\\
\mcl{O}_{m',n',p'}\lesssim \lra{q}^{-n'}t^{C\eps},\quad &\forall m',n'\geq 0,\ p'>0 \text{ such that }m'+n'+p'<m+n+p\\&\text{or }m'+n'+p'=m+n+p,\ m'<m.} From now on, we fix three multiindices $I,J,K$ such that $|I|=m$, $|J|=n$, and $|K|=p$.

We write $Z^I=ZZ^{I'}$ and apply $Z^{I'}X^{J}\Omega^K$ to the eikonal equation.  We have
\fm{0&=2g^{\alpha\beta}q_\beta (\partial_\alpha Z^{I'}X^J \Omega^{K}q)+\mcl{R}_1+\mcl{R}_2+\mathcal{R}_3}
where the remainders are given by
\fm{\mcl{R}_1&=Z^{I'}X^J \Omega^{K}(m^{\alpha\beta}q_\alpha q_\beta)-2m^{\alpha\beta} (\partial_\alpha Z^{I'}X^J \Omega^{K} q)q_\beta,\\
\mcl{R}_2&=Z^{I'}X^J \Omega^{K}((g^{\alpha\beta}-m^{\alpha\beta})q_\alpha q_\beta)-2(g^{\alpha\beta}-m^{\alpha\beta})q_\beta ( Z^{I'}X^J \Omega^{K}q_\alpha),\\
\mcl{R}_3&=2(g^{\alpha\beta}-m^{\alpha\beta})q_\beta( Z^{I'}X^J \Omega^{K}q_\alpha-\partial_\alpha Z^{I'}X^J \Omega^{K}q)}

We start with $\mcl{R}_3$. Recall that $g-m=O(\eps t^{-1+C\eps})$ and $q_\beta=O(1)$ on $H$. Besides, $Z^{I'}X^J \Omega^{K}q_\alpha-\partial_\alpha Z^{I'}X^J \Omega^{K}q$ is a linear combination of terms of the following forms
\fm{Z^{I_1}[Z,\partial_\alpha]Z^{I_2}X^J \Omega^{K}q=CZ^{I_1}\partial Z^{I_2}X^J \Omega^{K}q,&\hspace{2em}Z^{I_1}ZZ^{I_2}=Z^{I'};\\
Z^{I'}X^{J_1}[X,\partial_\alpha]X^{J_2}\Omega^Kq=CZ^{I'}X^{J_1}((\partial_\alpha\omega)\partial_tX^{J_2}\Omega^Kq),&\hspace{2em}X^{J_1}XX^{J_2}=X^{J};\\
Z^{I'}X^J\Omega^{K_1}[\Omega,\partial_\alpha]\Omega^{K_2}q=CZ^{I'}X^J\Omega^{K_1}\partial \Omega^{K_2}q,&\hspace{2em}\Omega^{K_1}\Omega\Omega^{K_2}=\Omega^K.}
The first  row has an upper bound
\fm{\sum_{|K'|\leq |I_1|+|I_2|}|\partial Z^{K'}X^J\Omega^Kq|&\lesssim\lra{r-t}^{-1}\sum_{|K'|\leq m-1}|Z^{K'}X^J\Omega^Kq|=\lra{q}^{-1}\sum_{m'\leq m-1}\mcl{O}_{m',n,p}\\
&\lesssim\lra{q}^{-1}\cdot\lra{q}^{1-n}t^{C\eps}\lesssim \lra{q}^{-n}t^{C\eps}.}
We can use the induction hypotheses \eqref{derqlhih} to control the sum $\sum_{m'\leq m-1}\mcl{O}_{m',n,p}$, since $m'+n+p\leq m-1+n+p<m+n+p$. 
The second row has an upper bound
\fm{&\hspace{1.5em}\sum_{|I_1|+|I_2|=m-1\atop |J'_1|+|J'_2|=|J_1|}|Z^{I_1}X^{J'_1}\partial\omega|\cdot|Z^{I_2}X^{J'_2}\partial X^{J_2} \Omega^{K}q|\\
&\lesssim\sum_{ |J_1'|+|J_2'|=|J_1|}r^{-1-|J_1'|}\cdot\lra{r-t}^{-|J_2'|-1-|J_2|}\sum_{|K'|\leq |I_2|+|J_2'|+1+|J_2|}|Z^{K'}\Omega^Kq|\\
&\lesssim \lra{q}^{-1-n}\sum_{m'\leq m-1+n}\mcl{O}_{m',0,p}\lesssim \lra{q}^{-n}t^{C\eps}.}
In the first inequality we apply Lemma \ref{l2.1} and Lemma \ref{derqlem1}.
In the second line, we apply \eqref{derqlhih}. The third row has an upper bound\fm{\lra{r-t}^{-n}\sum_{|K'|\leq m-1+n}|Z^{K'}\Omega^{K_1}\partial\Omega^{K_2}q|&\lesssim\lra{r-t}^{-1-n}\sum_{|K'|\leq m-1+n+|K_1|+1}|Z^{K'}\Omega^{K_2}q|\\
&\lesssim\lra{q}^{-1-n}\sum_{m'\leq m-1+n+p}\mcl{O}_{m',0,0}\lesssim\lra{q}^{-n}t^{C\eps}.} In conclusion, $\mcl{R}_3=O(\eps t^{-1+C\eps}\lra{q}^{-n})$.

We move on to $\mcl{R}_2$. By the Leibniz's rule, we can express $\mcl{R}_2$ as a linear combination of terms of the form\fm{Z^{I_1}X^{J_1}\Omega^{K_1}(g^{\alpha\beta}-m^{\alpha\beta})\cdot Z^{I_2}X^{J_2}\Omega^{K_2}q_\alpha\cdot Z^{I_3}X^{J_3}\Omega^{K_3}q_\beta,}where $\sum|I_*|=m-1$, $\sum|J_*|=n$, $\sum|K_*|=p$, $\max_{l=2,3}\{|I_l|+|J_l|+|K_l|\}<m+n+p-1$. On $H$, by Lemma \ref{derqlem1} and \eqref{derqlhih} we have 
\fm{|Z^{I_2}X^{J_2} \Omega^{K_2}q_\alpha|\lesssim\lra{q}^{-|J_2|}\sum_{|K'|\leq |I_2|+|J_2|+|K_2|}|Z^{K'}q_\alpha|\lesssim\lra{q}^{-|J_2|-1}\sum_{|K'|< m+n+p}|Z^{K'}q|\lesssim \lra{q}^{-|J_2|}t^{C\eps}.}
We can estimate $Z^{I_3}X^{J_3} \Omega^{K_3}q_\beta$ in the same way. And since 
$Z^{I_1}X^{J_1}\Omega^{K_1}(g^{\alpha\beta}-m^{\alpha\beta})=O(\eps\lra{q}^{-|J_1|}t^{-1+C\eps})$ by Lemma \ref{derqlem1}, we conclude that $\mcl{R}_2=O(\eps\lra{q}^{-n}t^{-1+C\eps})$ on $H$.

We move on to $\mcl{R}_1$.  By Lemma \ref{lcommin}, we can write $\Omega^K(m^{\alpha\beta}q_\alpha q_\beta)$ as a linear combination (with real constant coefficients) of terms of the form
\eq{\label{derqlht1}m^{\alpha\beta}(\partial_\alpha \Omega^{K_1}q)(\partial_\beta \Omega^{K_2}q),\hspace{2em}\min\{1,p\}\leq|K_1|+|K_2|\leq p.}
Here $(m^{\alpha\beta})$ is the usual Minkowski metric. In fact, if $p=0$, then \eqref{derqlht1} is  $m^{\alpha\beta}q_\alpha q_\beta$ so there is nothing to prove; if $p>0$, then we guarantee that $|K_1|+|K_2|>0$ in \eqref{derqlht1} since  
\fm{\Omega^K(m^{\alpha\beta}q_\alpha q_\beta)=\Omega^{K'}(m^{\alpha\beta}(\partial_\alpha\Omega  q)q_\beta+m^{\alpha\beta}q_\alpha (\partial_\beta\Omega q)),\hspace{2em}\Omega^K=\Omega^{K'}\Omega.}
Next we consider $X^J\Omega^K(m^{\alpha\beta}q_\alpha q_\beta)$, so we apply $X^J$ to \eqref{derqlht1}. By Lemma \ref{derqlem2}, we can write $X^J\Omega^K(m^{\alpha\beta}q_\alpha q_\beta)$ as a linear combination (with real constant coefficients) of terms of the form
\fm{\left\{\begin{array}{ll}m^{\alpha\beta}(\partial_\alpha X^{J_1}\Omega^{K_1}q)(\partial_\beta X^{J_2}\Omega^{K_2}q),&|J_1|+|J_2|=n,\\& \min\{1,p\}\leq |K_1|+|K_2|\leq p;\\[1em]
X^{J_1}(r^{-1}f_0)\cdot(X^{J_2}\partial X^{J_2'}\Omega^{K_1} q)(X^{J_3}\partial X^{J_3'} \Omega^{K_2} q),& \sum|J_*|+|J_*'|=n-1,\\& \min\{1,p\}\leq |K_1|+|K_2|\leq p.\end{array}\right.}
Again $(m^{\alpha\beta})$ is the Minkowski metric. We finally apply $Z^{I'}$ to each of these terms. By Lemma \ref{lcommin} and the Leibniz's rule, we can write $\mcl{R}_1$ as a linear combination (with real constant coefficients) of terms of the form
\eq{\label{derqlht2}\left\{\begin{array}{l}m_0^{\alpha\beta}(\partial_\alpha Z^{I_1}X^{J_1}\Omega^{K_1}q)(\partial_\beta Z^{I_2}X^{J_2}\Omega^{K_2}q),\\\hspace{4em}|I_1|+|I_2|\leq m-1,\ |J_1|+|J_2|= n,\ \min\{1,p\}\leq|K_1|+|K_2|\leq p
\\\hspace{4em}|I_1|+|J_1|+|K_1|,\ |I_2|+|J_2|+|K_2|<m-1+n+p;
\\[1em]
Z^{I_3}X^{J_3}(r^{-1}f_0)\cdot(Z^{I_1}X^{J_1}\partial X^{J_1'}\Omega^{K_1} q)(Z^{I_2}X^{J_2}\partial X^{J_2'} \Omega^{K_2} q),\\\hspace{4em} \sum|I_*|=m-1,\ \sum|J_*|+|J_*'|=n-1,\ \min\{1,p\}\leq|K_1|+|K_2|\leq p.\end{array}\right.}
Here $(m_0^{\alpha\beta})$ is some constant matrix satisfying the null condition defined in Section \ref{sec2}.  It follows from Lemma \ref{lnptwest} that on $H$ the terms of the first type in \eqref{derqlht2} has an upper bound \fm{&\hspace{1.5em} \lra{t}^{-1}\sum_{|L|=1}(|Z^LZ^{I_1}X^{J_1}\Omega^{K_1}q||\partial Z^{I_2}X^{J_2}\Omega^{K_2}q|+|\partial Z^{I_1}X^{J_1}\Omega^{K_1}q||Z^L Z^{I_2}X^{J_2}\Omega^{K_2}q|)\\
&\lesssim t^{-1}\lra{q}^{-1}\sum_{|L_1|=|L_2|=1}|Z^{L_1}Z^{I_1}X^{J_1}\Omega^{K_1}q||Z^{L_2} Z^{I_2}X^{J_2}\Omega^{K_2}q|\lesssim t^{-1}\lra{q}^{-1}\mcl{O}_{1+|I_1|,|J_1|,|K_1|}\mcl{O}_{1+|I_2|,|J_2|,|K_2|}.}
Since $\min_{l=1,2}\{|I_l|+|J_l|+|K_l|+1\}< m+n+p$ and since $|J_1|+|J_2|=n$, we can apply \eqref{derqlhih} to conclude that on $H$ \fm{|m_0^{\alpha\beta}(\partial_\alpha Z^{I_1}X^{J_1}\Omega^{K_1}q)(\partial_\beta Z^{I_2}X^{J_2}\Omega^{K_2}q)|\lesssim t^{-1+C\eps}\lra{q}^{1-n},&\hspace{2em}\text{if }p=0;\\|m_0^{\alpha\beta}(\partial_\alpha Z^{I_1}X^{J_1}\Omega^{K_1}q)(\partial_\beta Z^{I_2}X^{J_2}\Omega^{K_2}q)|\lesssim t^{-1+C\eps}\lra{q}^{-n},&\hspace{2em}\text{if }p>0.}
Meanwhile, by Lemma \ref{derqlem1} and \eqref{derqlhih}, on $H$ we have \fm{|Z^{I_3}X^{J_3}(r^{-1}f_0)|&\lesssim t^{-1+C\eps}\lra{q}^{-|J_3|},\\
 |Z^{I_1}X^{J_1}\partial X^{J_1'}\Omega^{K_1} q|&\lesssim \lra{q}^{-1-|J_1|-|J'_1|}\sum_{m'\leq |I_1|+1+|J_1|+|J_1'|}\mcl{O}_{m',0,|K_1|},\\
 |Z^{I_2}X^{J_2}\partial X^{J_2'}\Omega^{K_2} q|&\lesssim \lra{q}^{-1-|J_2|-|J'_2|}\sum_{m'\leq |I_2|+1+|J_2|+|J_2'|}\mcl{O}_{m',0,|K_2|}.} Here we can  apply \eqref{derqlhih}  as $\max_{l=1,2}\{|I_l|+|J_l|+|J_l'|+|K_l|+1\}<m+n+p$. Thus, the product of these terms is $O(t^{-1+C\eps}\lra{q}^{1-n})$ if $p=0$, or $O(t^{-1+C\eps}\lra{q}^{-n})$ if $p>0$. Thus, on $H$ we have $\mcl{R}_1=O(t^{-1+C\eps}\lra{q}^{1-n})$ if $p=0$, and $\mcl{R}_1=O(t^{-1+C\eps}\lra{q}^{-n})$ if $p>0$.  In conclusion, we have
\fm{2g^{\alpha\beta}q_\beta(\partial_\alpha Z^{I'}X^J\Omega^Kq)=O(t^{-1+C\eps}\lra{q}^{1-n}),&\hspace{2em}\text{if }p=0;\\
2g^{\alpha\beta}q_\beta(\partial_\alpha Z^{I'}X^J\Omega^Kq)=O(t^{-1+C\eps}\lra{q}^{-n}),&\hspace{2em}\text{if }p>0.}

Next, we note that
\fm{X_jZ^{I'}X^J\Omega^Kq&=Z^{I'}X_jX^J\Omega^Kq+\sum_{I'=(I_1,i,I_2)}Z^{I_1}[X_j,Z_i]Z^{I_2}X^J\Omega^Kq,\\
\Omega_{kk'}Z^{I'}X^J\Omega^Kq&=Z^{I'}X^J\Omega_{kk'}\Omega^Kq+\sum_{I'=(I_1,i,I_2)}Z^{I_1}[\Omega_{kk'},Z_i]Z^{I_2}X^J\Omega^Kq\\&\hspace{2em}+\sum_{J=(J_1,j,J_2)}Z^{I'}X^{J_1}[\Omega_{kk'},X_j]X^{J_2}\Omega^Kq.}
Recall that $[\Omega,Z]=\sum f_0Z$ and $[X,Z]=\sum f_0\partial$ where $f_0$ denotes any function such that $Z^{K'}f_0=O(1)$ for all $K'$. By Lemma \ref{l2.1} we have
\fm{|X_jZ^{I'}X^J\Omega^Kq|&\lesssim \mcl{O}_{m-1,n+1,p}+\sum_{I'=(I_1,i,I_2)}|Z^{I_1}(f_0\partial Z^{I_2}X^J\Omega^Kq)|\\&\lesssim  \mcl{O}_{m-1,n+1,p}+\lra{q}^{-1}\sum_{m'\leq m-1}\mcl{O}_{m',n,p},}
\fm{&\hspace{1.5em}|\Omega_{kk'}Z^{I'}X^J\Omega^Kq|\\
&\lesssim \mcl{O}_{m-1,n,p+1}+\sum_{I'=(I_1,i,I_2)}|Z^{I_1}(f_0ZZ^{I_2}X^J\Omega^Kq)|+\sum_{J=(J_1,j,J_2)}|Z^{I'}X^{J_1}(f_0\partial X^{J_2}\Omega^Kq)|\\
&\lesssim \mcl{O}_{m-1,n,p+1}+\sum_{m'\leq m-1}\mcl{O}_{m',n,p}+\sum_{|J_1|+|J_2|=n-1}\lra{q}^{-|J_1|}|Z^{I'}Z^{J_1}(f_0\partial X^{J_2}\Omega^Kq)|\\
&\lesssim \mcl{O}_{m-1,n,p+1}+\sum_{m'\leq m-1}\mcl{O}_{m',n,p}+\lra{q}^{-n}\sum_{m'\leq m+n-1}\mcl{O}_{m',0,p}.}
In conclusion,  on $H$ we have
\fm{|XZ^{I'}X^J\Omega^Kq|\lesssim \lra{q}^{-n}t^{C\eps},&\hspace{2em}\text{if }p=0;\hspace{2em}
|XZ^{I'}X^J\Omega^Kq|\lesssim \lra{q}^{-1-n}t^{C\eps},&\hspace{2em}\text{if }p>0;\\
|\Omega Z^{I'}X^J\Omega^Kq|\lesssim \lra{q}^{1-n}t^{C\eps},&\hspace{2em}\text{if }p=0;\hspace{2em}
|\Omega Z^{I'}X^J\Omega^Kq|\lesssim \lra{q}^{-n}t^{C\eps},&\hspace{2em}\text{if }p>0.}

We now end the proof. By setting $L^\alpha=2g^{\alpha\beta}q_\beta$ and $L=L^\alpha\partial_\alpha$, we have
\fm{\partial_t&=\frac{L-L^iX_i}{L^0-2\omega_iL^i}=-\frac{1}{2}L+\sum_i\omega_iX_i+O(|u|)L+\sum_iO(|u|)X_i,\\
\partial_j&=X_j-2\omega_j\partial_t=\omega_jL+X_j-2\omega_j\sum_i\omega_iX_i+O(|u|)L+\sum_iO(|u|)X_i.}
Note that $L^0=2+O(|u|)$ and $L^i=2\omega_i+O(|u|)$ on $H$. Then, we have
\fm{S&=(-\frac{1}{2}t+r)L+(t-r)\sum_i\omega_iX_i+O((r+t)|u|)L+\sum_iO((r+t)|u|)X_i\\
&=O(t+\eps t^{C\eps})L+\sum_iO(\lra{q}+\eps t^{C\eps})X_i.}
And since $\Omega_{kk'}=x_{k}X_{k'}-x_{k'}X_k$, we have $\sum_kr^{-1}\omega_k\Omega_{kk'}=X_{k'}-\sum_k\omega_{k'}\omega_kX_k$. Thus,
\fm{\Omega_{0j}&=(-\frac{1}{2}x_j+t\omega_j)L+tX_j+(x_j-2t\omega_j)\sum_i\omega_iX_i+O((r+t)|u|)L+\sum_iO((r+t)|u|)X_i\\
&=t(X_j-\omega_j\omega_iX_i)+O(t+\eps t^{C\eps})L+\sum_iO(\lra{q}+\eps t^{C\eps})X_i\\
&=tr^{-1}\sum_i \omega_{i}\Omega_{ij}+O(t+\eps t^{C\eps})L+\sum_iO(\lra{q}+\eps t^{C\eps})X_i.}
In conclusion, for each $Z\in\{\partial_\alpha,S,\Omega_{0j}\}$, we have
\fm{|ZZ^{I'}X^J\Omega^Kq|&\lesssim \sum_{1\leq i<j\leq 3}|\Omega_{ij}Z^{I'}X^J\Omega^Kq|+t|LZ^{I'}X^J\Omega^Kq|+(\lra{q}+t^{C\eps})\sum_i|X_iZ^{I'}X^J\Omega^Kq|. }
If $p=0$, the right hand side has an upper  bound $\lra{q}^{1-n}t^{C\eps}$; if $p>0$, the right hand side has an upper  bound $\lra{q}^{-n}t^{C\eps}$. We finish the proof by induction.
\end{proof}\rm

\lem{\label{lemexipr} On $H$, we have $Z^I(q_i-\omega_iq_r)=O(t^{-1+C\eps})$ and $Z^I(q_t+q_r)=O(\eps t^{-1+C\eps})$ for each $I$. As a result, $Z^I(q_i+\omega_iq_t)=O(t^{-1+C\eps})$.}
\begin{proof}
Recall that $q_i-\omega_iq_r=\sum_jr^{-1}\omega_j\Omega_{ji}q$. By Lemma \ref{derqlh} and the Leibniz's rule, for each $I$ we have \fm{|Z^{I}(r^{-1}\omega_j\Omega_{ji}q)|\lesssim\sum_{|I_1|+|I_2|=|I|}|Z^{I_1}(r^{-1}\omega_j)|\cdot|Z^{I_2}\Omega q|\lesssim t^{-1+C\eps}.}
So $Z^I(q_i-\omega_iq_r)=O(t^{-1+C\eps})$. Moreover, by the eikonal equation we have
\fm{-(q_t+q_r)(q_t-q_r)+\sum_i(q_i-\omega_iq_r)^2+(g^{\alpha}(u)-m^{\alpha\beta})q_\alpha q_\beta=0,}
so
\fm{q_t+q_r=\frac{\sum_i(q_i-\omega_iq_r)^2+(g^{\alpha}(u)-m^{\alpha\beta})q_\alpha q_\beta}{q_t-q_r}.}
Thus, $Z^I(q_t+q_r)$ is a linear combination of terms of the form
\fm{(q_t-q_r)^{-1-s}\cdot Z^{I_1}(q_t-q_r)\cdots Z^{I_s}(q_t-q_r)\cdot Z^{I_0}(\sum_i(q_i-\omega_iq_r)^2+(g^{\alpha}(u)-m^{\alpha\beta})q_\alpha q_\beta)}
where $\sum|I_*|=|I|$. It is clear that $Z^{I_*}(q_t-q_r)=O(t^{C\eps})$ and that $q_t-q_r=-2+O(\eps t^{-1+C\eps})\leq -1$ on $H$. Moreover, since $Z^I(r^{-1}\Omega q)=O(t^{-1+C\eps})$ for each $I$, we have $Z^{I_0}((q_i-\omega_iq_r)^2)=O(t^{-2+C\eps})$. Finally, for each $I$ we have \fm{|Z^{I}((g^{\alpha\beta}-m^{\alpha\beta})q_\alpha q_\beta)|&\lesssim \sum_{|I_1|+|I_2|+|I_3|=|I|}|Z^{I_1}(g-m)||Z^{I_2}\partial q||Z^{I_3}\partial q|\lesssim \eps t^{-1+C\eps}.}
In conclusion, $Z^I(q_t+q_r)=O(t^{-2+C\eps}+\eps t^{-1+C\eps})=O(\eps t^{-1+C\eps})$, as $t\geq T_0=\exp(\delta/\eps)$. Since $q_i+\omega_iq_t=q_i-\omega_iq_r+\omega_i(q_t+q_r)$, we can easily show $Z^I(q_i+\omega_iq_t)=O(t^{-1+C\eps})$ by the Leibniz's rule.
\end{proof}\rm

We move on to estimates for $e_*^*$ and $\xi_{**}^*$ on $H$.

\lem{\label{lemexi}On $H$, we have $Z^Ie_k^\alpha=O(t^{C\eps})$ and $Z^I(e_3^i-\omega_i,e_4^i-\omega_i)=O(t^{-1+C\eps})$ for each~$I$.}
\begin{proof}Since $e_4^0=1$, $e_3^0=-1$ and $e_a^0=0$, we can ignore the case $\alpha=0$. We write 
\fm{e_4^i-\omega_i&=(g^{0\mu}q_\mu)^{-1}(g^{i\beta}q_\beta-\omega_ig^{0\beta}q_\beta)\\&=(g^{0\mu}q_\mu)^{-1}(q_i+\omega_iq_t+(g^{i\beta}-m^{i\beta})q_\beta-\omega_i(g^{0\beta}-m^{0\beta})q_\beta)\\
&=:(g^{0\mu}q_\mu)^{-1}\mcl{Q}.}
By Lemma \ref{derqlh}, Lemma \ref{lemexipr} and the Leibniz's rule,  we have \fm{Z^I\mcl{Q}=O(t^{-1+C\eps}),\hspace{2em}Z^I(g^{0\mu}q_\mu)=O(t^{C\eps}),\hspace{2em}g^{0\mu}q_\mu=1+O(\eps t^{-1+C\eps})\geq 1/2.}
Besides,  $Z^I(e_4^i-\omega_i)$ is a linear combination of terms of the form \fm{(g^{0\mu}q_\mu)^{-1-s}Z^{I_1}(g^{0\mu}q_\mu)\cdots Z^{I_s}(g^{0\mu}q_\mu)Z^{I_0}\mcl{Q},\hspace{2em}\sum|I_*|=|I|,\ |I_j|>0\text{ for }j\neq 0.}
We conclude that $Z^I(e_4^i-\omega_i)=O(t^{-1+C\eps})$. Since $Z^I\omega=O(1)$ on $H$, we conclude that $Z^Ie_4^i=O(t^{C\eps})$. And since $Z^I(e_3^i-e_4^i)=2Z^Ig^{0i}=O(\eps t^{-1+C\eps})$, we conclude that $Z^I(e_3^i-\omega_i)=O(t^{-1+C\eps})$ and $Z^Ie_3^i=O(t^{C\eps})$ on $H$  for each $I$. The proofs of these estimates do not rely on the estimates for $Z^Ie_a^*$, so we can use them freely in the following proof.

Next, we claim that $Z^IX^J \Omega^Ke_a^i=O(\lra{q}^{-|J|}t^{C\eps})$ on $H$ for all $I,J,K$ and $a=1,2$. Recall that $\Omega^K$ is the product of $|K|$ vector fields in $\{\Omega_{12},\Omega_{23},\Omega_{13}\}$. We induct first on $|I|+|J|+|K|$ and then on $|I|$. When $|I|+|J|+|K|=0$, there is nothing to prove. When $|I|=0$ and $|J|+|K|>0$, we have $X^J\Omega^Ke_a^i=O(r^{-|K|})$ on $H$, since $e_a^i|_H$ is a locally defined  function of $\omega$ and it is independent of $t$.

In general, we fix $I,J,K$ such that $|I|>0$. Suppose we have proved the claim for all $(I',J',K')$ such that $|I'|+|J'|+|K'|<|I|+|J|+|K|$, or $|I'|+|J'|+|K'|=|I|+|J|+|K|$ and $|I'|<|I|$. We write $Z^I=ZZ^{I'}$.  For $a=1,2$ we have
\fm{Z^{I'}X^J\Omega^Ke_4(e_a^i)&=Z^{I'}X^J\Omega^K(e_4^\alpha e_a^\beta\Gamma_{\alpha\beta}^0e_4^i-e_4^\alpha e_a^\beta \Gamma_{\alpha\beta}^i).}
Since we can write $\Gamma=g\cdot \partial g$, for each $K'$, we have $Z^{K'}\Gamma=O(\eps t^{-1+C\eps}\lra{q}^{-1})$ on $H$. By induction hypotheses, Lemma \ref{derqlem1} and the Leibniz's rule, we conclude that \fm{Z^{I'}X^J\Omega^Ke_4(e_a^i)=O(\eps t^{-1+C\eps}\lra{q}^{-1-|J|}).}

Moreover, $Z^{I'}X^J\Omega^Ke_4(e_a^i)$ is equal to the sum of $e_4(Z^{I'}X^J\Omega^Ke_a^i)$ and a linear combination of terms of the form
\eq{\label{lemexif1}&Z^{I_1}[e_4,Z^{I_2}]Z^{I_3}X^J\Omega^Ke_a^i,&(I_1,I_2,I_3)=I',\ |I_2|=1;\\
&Z^{I'}X^{J_1}[e_4,X^{J_2}]X^{J_3}\Omega^Ke_a^i,&(J_1,J_2,J_3)=J,\ |J_2|=1;\\
&Z^{I'}X^J\Omega^{K_1}[e_4,\Omega^{K_2}]\Omega^{K_3}e_a^i,&(K_1,K_2,K_3)=K,\ |K_2|=1.}
Note that \fm{\ [e_4,Z]&=e_4(z^\nu)\partial_\nu-Z(e_4^\nu)\partial_\nu=e_4(z^\nu)\partial_\nu-Z(\omega_j)\partial_j-Z(e_4^j-\omega_j)\partial_j,\\
[e_4,X_l]&=e_4(2\omega_l)\partial_t-X_l(e_4^j)\partial_j=2r^{-1}(e_4^l-\omega_l-(\omega_j-e_4^j)\omega_j\omega_l)\partial_t-(\partial_l\omega_j)\partial_j-X_l(e_4^j-\omega_j)\partial_j} where we write $Z=z^\nu(t,x)\partial_\nu$. We have
\fm{e_4(z^\nu)\partial_\nu-Z(\omega_j)\partial_j=\left\{\begin{array}{ll}
-\partial(\omega_j)\partial_j,&Z=\partial;\\
(r+t)^{-1}S+(r+t)^{-1}\omega_l\Omega_{0l}+(e_4^j-\omega_j)\partial_j,&Z=S;\\
r^{-1}\Omega_{ij}+(e_4^i-\omega_{i})\partial_j-(e_4^j-\omega_{j})\partial_i-r^{-1}\Omega_{ij},&Z=\Omega_{ij};\\
r^{-1}\Omega_{0i}+r^{-1}(t-r)\partial_i+(e_4^i-\omega_i)\partial_t-tr^{-2}\omega_l\Omega_{li},&Z=\Omega_{0i}.
\end{array}\right.}
In conclusion, \fm{\ [e_4,Z]=f_1\cdot Z,\hspace{2em} [e_4,X]=f_1\cdot \partial} where $f_1$ denotes any function satisfying $Z^{J'}f_1=O(t^{-1+C\eps})$ for each $J'$ on $H$. Thus, the first row  in \eqref{lemexif1} has an upper bound
\fm{|Z^{I_1}(f_1ZZ^{I_3}X^{J}\Omega^Ke_a^i)|\lesssim\sum_{|J'|\leq |I_1|}t^{-1+C\eps}|Z^{J'}ZZ^{I_3}X^{J}\Omega^Ke_a^i|\lesssim t^{-1+C\eps}\lra{q}^{-|J|}.}
We can use the induction hypotheses here as \fm{|J'|+1+|I_3|+|J|+|K|\leq |I_1|+1+|I_3|+|J|+|K|=|I'|+|J|+|K|<|I|+|J|+|K|.}
The second row in \eqref{lemexif1} has an upper bound
\fm{|Z^{I'}X^{J_1}(f_1\partial X^{J_3}\Omega^Ke_a^i)|&\lesssim\sum_{|J'|\leq |I'|+|J_1|}\lra{q}^{-|J_1|}|Z^{J'}(f_1\partial X^{J_3}\Omega^{K}e_a^i)|\\
&\lesssim \sum_{|J'|\leq |I'|+|J_1|}\lra{q}^{-|J_1|}t^{-1+C\eps}|Z^{J'}\partial X^{J_3}\Omega^{K}e_a^i|\\
&\lesssim \sum_{|J'|\leq |I'|+|J_1|+1}\lra{q}^{-|J_1|-1}t^{-1+C\eps}|Z^{J'} X^{J_3}\Omega^{K}e_a^i|\lesssim \lra{q}^{-|J|}t^{-1+C\eps}.}
We can use the induction hypotheses here as \fm{|J'|+|J_3|+|K|\leq |I'|+|J_1|+1+|J_3|+|K|=|I'|+|J|+|K|<|I|+|J|+|K|.}
The third row in \eqref{lemexif1} has an upper bound
\fm{|Z^{I'}X^J\Omega^{K_1}(f_1Z\Omega^{K_3}e_a^i)|&\lesssim\sum_{|J'|\leq |I'|+|J|} \lra{q}^{-|J|}|Z^{J'}\Omega^{K_1}(f_1Z\Omega^{K_3}e_a^i)|\\&\lesssim\sum_{|J'|\leq |I'|+|J|+|K_1|} \lra{q}^{-|J|}t^{-1+C\eps}|Z^{J'}Z\Omega^{K_3}e_a^i|\lesssim  \lra{q}^{-|J|}t^{-1+C\eps}.}
We can use the induction hypotheses here as \fm{|J'|+|K_3|+1\leq |I'|+|J|+|K_1|+1+|K_3|=|I'|+|J|+|K|<|I|+|J|+|K|.}
In conclusion, on $H$ we have
\fm{e_4(Z^{I'}X^J\Omega^Ke_a^i)=Z^{I'}X^J\Omega^K e_4(e_a^i)+O(t^{-1+C\eps}\lra{q}^{-|J|})=O(t^{-1+C\eps}\lra{q}^{-|J|}).}

We recall from the proof of Lemma \ref{derqlem1} that $[Z,\Omega]=C\cdot Z$ and $[Z,X]=f_0\cdot \partial$ where $f_0$ denotes any function such that $Z^{K'}f_0=O(t^{C\eps})$ on $H$ for each $K'$. If we keep commuting $\Omega$ with each vector field in $Z^{I'}X^J$  and applying the Leibniz's rule, we get $\Omega Z^{I'}X^{J}\Omega^Ke_a^i=O(t^{C\eps}\lra{q}^{-|J|})$. If we keep commuting $X_l$ with each vector field in $Z^{I'}$ and applying the Leibniz's rule, we get $X_lZ^{I'}X^J\Omega^Ke_a^i=O(t^{C\eps}\lra{q}^{-1-|J|})$. Finally, we recall from the proof of Lemma \ref{derqlem1} that we can write \fm{(\partial,S,\Omega_{0j})=O(t)L+O(1)\cdot\Omega +O(\lra{q}+\eps t^{C\eps})\cdot X}where $L=2g^{\alpha\beta}q_\beta\partial_\alpha=O(1)e_4$ on $H$.
In conclusion, when $Z=\partial,S,\Omega_{0j}$, we have
\fm{|ZZ^{I'}X^J\Omega^Ke_a^i|&\lesssim t|e_4(Z^{I'}X^J\Omega^Ke_a^i)|+|\Omega Z^{I'}X^J\Omega^Ke_a^i|+\lra{q}t^{C\eps}|XZ^{I'}X^J\Omega^Ke_a^i|\lesssim t^{C\eps}\lra{q}^{-|J|}.}
We finish the proof by induction.
\end{proof}

\rm We  now prove the following lemma which illustrates the connection between the weighted null frame and the commuting vector fields.
\lem{\label{derqvzlem} Let $F=F(t,x)$ be a smooth function defined near $H$. Then, on $H$ we have \fm{|V^IF|\lesssim \sum_{|J|\leq |I|}t^{C\eps}|Z^IF|.}}
\begin{proof}We induct on $|I|$. When $|I|=0$, there is nothing to prove. Suppose we have proved the estimate for each function $F$ and for each multiindex $I'$ such that $|I'|<|I|$. Then, by writing $V^I=V^{I'}V_k$ and applying the induction hypotheses, we have
\fm{|V^IF|&\lesssim \sum_{|J|\leq |I|-1}t^{C\eps}|Z^J(V_kF)|.}
We then apply \eqref{derqvz}. When $k=4$, we have $V_4F=f_0\cdot ZF$. Here $f_0$ denotes any function such that $Z^{J'}f_0=O(t^{C\eps})$ on $H$ for each $J'$. In particular, since $Z^{J'}(e_4^i-\omega_i)=O(t^{-1+C\eps})$ for each $J'
$ by Lemma \ref{lemexi}, we have $Z^{J'}(t(e_4^i-\omega_i))=O(t^{-1+C\eps})$ and thus $t(e_4^i-\omega_i)=f_0$.  By the Leibniz's rule, we have
\fm{|V^IF|\lesssim \sum_{|J|\leq |I|-1}t^{C\eps}|Z^J(f_0\cdot ZF)|\lesssim \sum_{|J|\leq |I|-1}t^{C\eps}|Z^JZF|\lesssim \sum_{|J|\leq |I|}t^{C\eps}|Z^JF|.}
The proof for $k=3$ follows from the case $k=4$ and the estimate $Z^{J'}(r-t)=O(\lra{r-t})$  for all $J'$. Finally, when $k=a\in\{1,2\}$, we note that \fm{V_a(r)=re_a^j \omega_j=re_a^\alpha (-g^{\alpha\beta}+m^{\alpha\beta})e_4^\beta+re_a^j m^{j l}(-e_4^l+\omega_l).}
By Lemma \ref{lemexi}, we have $Z^{J'}(\omega_* ,e_*^*)=O(t^{C\eps})$ and thus $Z^{J'}(V_a(r))=O(t^{C\eps})$ on $H$ for each $|J'|$. Thus, for all $|J|\leq |I|-1$, we have
\fm{|Z^{J}(V_aF)|&\lesssim |Z^{J}(V_a(r)\omega_i\partial_iF)|+|Z^{J}(e_a^i\omega_j\Omega_{ji}F)|\\&\lesssim t^{C\eps}\sum_{|K|\leq |J|}|Z^K\partial F|+t^{C\eps}\sum_{|K|\leq |J|}|Z^KF|\lesssim t^{C\eps}\sum_{|K|\leq |I|}|Z^K F|.}
This finishes the proof.
\end{proof}
\rmk{\rm\label{derqvzrmk1} With the help of this lemma, we conclude immediately that \fm{&V^I(g-m)=O(\eps t^{-1+C\eps}),\hspace{1em} V^I((3R-r+t)^{-1})=O(\lra{q}^{-1}t^{C\eps}),\hspace{1em}V^I(r^{-1},t^{-1})=O(t^{-1+C\eps}),}
\fm{V^I(q)=\lra{q}t^{C\eps},\hspace{2em} V^Ie_k^\alpha=O(t^{C\eps}),\hspace{2em} V^I(e_3^i-\omega_i,e_4^i-\omega_i)=O(t^{-1+C\eps})} on $H$ for each $I$. }\rm

\lem{\label{derqxi}For each $I$, on $H$ we have $V^I(\xi_{13}^2,\xi_{23}^1)=O(\lra{q}^{-1}t^{C\eps})$, $V^I(\xi_{34}^a)=O(t^{-1+C\eps}\lra{q}^{-1})$ and  $V^I(\xi_{k_1k_2}^a)=O(t^{-1+C\eps})$ for all other $k_1<k_2$ and $a\in\{1,2\}$; $V^{I}(\xi_{k_1k_2}^3)=O(t^{-1+C\eps}\lra{q}^{-1})$ for all $k_1<k_2$;  $V^I(\chi_{ab}-r^{-1}\delta_{ab})=O(t^{-2+C\eps})$.}
\begin{proof} First, for any function $F=F(t,x)$ and for each $1\leq k\leq 4$, on $H$ we have \eq{\label{derqxif1}|V^I(e_k(F))|\lesssim \lra{q}^{-1}t^{C\eps}\sum_{|J|\leq |I|+1}|V^J(F)|.}
This inequality easily follows from the Leibniz's rule, Remark \ref{derqvzrmk1} and the estimate $\lra{r-t}\lesssim t$ on $H$.

Since $e_l(\lra{e_{k_1},e_{k_2}})=0$ for each $k_1,k_2,l$, we have
\fm{2\xi_{k_1k_2}^3&=\lra{[e_{k_1},e_{k_2}],e_4}=e_{k_1}(e_{k_2}^\alpha)g_{\alpha\beta}e_4^\beta-e_{k_2}(e_{k_1}^\alpha)g_{\alpha\beta}e_4^\beta\\
&=-e_{k_2}^\alpha e_{k_1}(g_{\alpha\beta})e_4^\beta-e_{k_2}^\alpha g_{\alpha\beta}e_{k_1}(e_4^\beta)+e_{k_1}^\alpha e_{k_2}(g_{\alpha\beta})e_4^\beta+e_{k_1}^\alpha g_{\alpha\beta}e_{k_2}(e_4^\beta).}
We assume $k_1\neq k_2$ as  $\xi_{k_1k_1}^*\equiv 0$. By \eqref{derqxif1} and the Leibniz's rule, on $H$ for each $I$ we have
\fm{|V^I(-e_{k_2}^\alpha e_{k_1}(g_{\alpha\beta})e_4^\beta+e_{k_1}^\alpha e_{k_2}(g_{\alpha\beta})e_4^\beta)|\lesssim \eps t^{-1+C\eps}\lra{q}^{-1}.}
Moreover, since $e_4^0\equiv 1$, we have
\fm{e_{k_2}^\alpha g_{\alpha\beta}e_{k_1}(e_4^\beta)&=e_{k_2}^\alpha g_{\alpha j}e_{k_1}(e_4^j-\omega_j)+e_{k_2}^\alpha g_{\alpha j}e_{k_1}(\omega_j)\\
&=e_{k_2}^\alpha g_{\alpha j}e_{k_1}(e_4^j-\omega_j)+r^{-1}e_{k_2}^\alpha g_{\alpha j}(e_{k_1}^j-e_{k_1}^l\omega_l\omega_j).}
Again, by \eqref{derqxif1} and the Leibniz's rule, on $H$ for each $I$ we have
\fm{|V^I(e_{k_2}^\alpha g_{\alpha j}e_{k_1}(e_4^j-\omega_j))|\lesssim t^{-1+C\eps}\lra{q}^{-1}.}
If $k_1=3$ or $4$, then since \fm{e_{k_1}^j-e_{k_1}^l\omega_l\omega_j=e_{k_1}^j-\omega_j+(1-e_{k_1}^l\omega_l)\omega_j=e_{k_1}^j-\omega_j+\sum_l(\omega_l-e_{k_l})\omega_l\omega_j,} by the Leibniz's rule and the estimate $V^I(e_3^i-\omega_i,e_4^i-\omega_i)=O(t^{-1+C\eps})$ for each $I$, we conclude that 
\fm{|V^I(r^{-1}e_{k_2}^\alpha g_{\alpha j}(e_{k_1}^j-\omega_j+(1-e_{k_1}^l\omega_l)\omega_j))|\lesssim t^{-2+C\eps},\hspace{2em}k_1\geq 3.}
If $k_1=1$ or $2$, then $e_{k_1}^0=0$. 
\fm{r^{-1}e_{k_2}^\alpha g_{\alpha j}(e_{k_1}^j-e_{k_1}^l\omega_l\omega_j)&=r^{-1}\lra{e_{k_2},e_{k_1}}-r^{-1}e_{k_2}^\alpha g_{\alpha j}e_{k_1}^l\omega_l\omega_j=-r^{-1}e_{k_2}^\alpha g_{\alpha j}e_{k_1}^l\omega_l\omega_j.}
Note that
\fm{e_{k_1}^l\omega_l&=e_{k_1}^l\delta_{ll'}e_4^{l'}+e_{k_1}^l\delta_{ll'}(\omega_{l'}-e_4^{l'})=e_{k_1}^\mu g_{\mu\nu}e_4^{\nu}-e_{k_1}^\mu (g_{\mu\nu}-m_{\mu\nu})e_4^{\nu}+e_{k_1}^l\delta_{ll'}(\omega_{l'}-e_4^{l'})\\
&=-e_{k_1}^\mu (g_{\mu\nu}-m_{\mu\nu})e_4^{\nu}+e_{k_1}^l\delta_{ll'}(\omega_{l'}-e_4^{l'}).}
Thus, by the Leibniz's rule, we have $V^I(e_{k_1}^l\omega_l)=O(t^{-1+C\eps})$ and thus \fm{|V^I(r^{-1}e_{k_2}^\alpha g_{\alpha j}(e_{k_1}^j-e_{k_1}^l\omega_l\omega_j))|\lesssim t^{-2+C\eps},\hspace{2em}k_1\leq 2.}
In conclusion, for each $I$, on $H$ we have
\fm{|V^I(\xi_{k_1k_2}^3)|\lesssim t^{-1+C\eps}\lra{q}^{-1}+t^{-2+C\eps}\lesssim t^{-1+C\eps}\lra{q}^{-1}.}

Next, we have
\fm{\xi_{k_1k_2}^c&=\lra{[e_{k_1},e_{k_2}],e_c}=e_{k_1}(e_{k_2}^\alpha)g_{\alpha\beta}e_c^\beta-e_{k_2}(e_{k_1}^\alpha)g_{\alpha\beta}e_c^\beta.} We first prove some estimates for  $e_{k_1}(e_{k_2}^\alpha)g_{\alpha\beta}e_c^\beta$ with $k_1\neq k_2$. If $k_1=a\in\{1,2\}$ and $k_2=b\in\{1,2\}$, we have $e_a=r^{-1}V_{a}$ and thus $V^I(e_a(e_b^\alpha)g_{\alpha\beta}e_c^\beta)=O(t^{-1+C\eps})$ on $H$. If $k_2=3$ and $k_1=a\in\{1,2\}$, then
\fm{e_{a}(e_3^\alpha)g_{\alpha\beta}e_c^\beta&=e_{a}(\omega_i)g_{i\beta}e_c^\beta+e_{a}(e_3^i-\omega_i)g_{i\beta}e_c^\beta=r^{-1}(e_a^i-e_a^l\omega_l\omega_i)g_{ij}e_c^j+e_{a}(e_3^i-\omega_i)g_{i\beta}e_c^\beta\\
&=r^{-1}\delta_{ac}-r^{-1}(e_a^l\omega_l)\omega_i g_{ij}e_c^j+r^{-1}V_{a}(e_3^i-\omega_i)g_{i\beta}e_c^\beta.}
Recall that $V^I(e_a^l\omega_l)=O(t^{-1+C\eps})$ on $H$. By Remark  \ref{derqvzrmk1}, we have   $V^I(e_{a}(e_3^\alpha)g_{\alpha\beta}e_c^\beta-r^{-1}\delta_{ac})=O(t^{-2+C\eps})$ on $H$. Following the same proof, we can show that $V^I(e_{a}(e_4^\alpha)g_{\alpha\beta}e_c^\beta-r^{-1}\delta_{ac})=O(t^{-2+C\eps})$  on $H$.
  Next, for $k\neq 3$ we have
\fm{e_4(e_k^\alpha)g_{\alpha\beta}e_c^\beta&=e_4^\mu e_k^\nu (\Gamma_{\mu\nu}^0e_4^\alpha-\Gamma_{\mu\nu}^\alpha)g_{\alpha\beta}e_c^\beta=-e_4^\mu e_k^\nu \Gamma_{\mu\nu}^\alpha g_{\alpha\beta}e_c^\beta\\
&=-\frac{1}{2}e_4^\mu e_k^\nu e_c^\beta(\partial_\mu g_{\beta\nu}+\partial_\nu g_{\beta\mu}-\partial_\beta g_{\mu\nu})\\
&=-\frac{1}{2}(t^{-1} e_k^\nu e_c^\beta V_4( g_{\beta\nu})+e_4^\mu  e_c^\beta (t^{-1},r^{-1})V_k (g_{\beta\mu})-r^{-1}e_4^\mu e_k^\nu V_c( g_{\mu\nu})).}
\fm{e_4(e_3^\alpha)g_{\alpha\beta}e_c^\beta&=e_4(2g^{0\alpha})g_{\alpha\beta}e_c^\beta+e_4(e_4^\alpha)g_{\alpha\beta}e_c^\beta=t^{-1}V_4(2g^{0\alpha})g_{\alpha\beta}e_c^\beta+e_4(e_4^\alpha)g_{\alpha\beta}e_c^\beta.}
Then,  on $H$ we have $V^I(e_4(e_k^\alpha)g_{\alpha\beta}e_c^\beta)=O(\eps t^{-2+C\eps})$. Next, we  have
\fm{e_3(e_4^\alpha)g_{\alpha\beta}e_c^\beta&=e_3(\omega_j)g_{j\beta}e_c^\beta+(3R-r+t)^{-1}V_3(e_4^j-\omega_j)g_{j\beta}e_c^\beta\\
&=r^{-1}(e_3^j-\omega_j+(1-\sum_le_3^l\omega_l)\omega_j)g_{j\beta}e_c^\beta+(3R-r+t)^{-1}V_3(e_4^j-\omega_j)g_{j\beta}e_c^\beta.}
Then,  on $H$ we have $V^I(e_3(e_4^\alpha)g_{\alpha\beta}e_c^\beta)=O(t^{-1+C\eps}\lra{q}^{-1})$. Besides, we have
\fm{e_3(e_c^\alpha)g_{\alpha\beta}e_c^\beta&=-e_c^\alpha g_{\alpha\beta}e_3(e_c^\beta)-e_c^\alpha e_3(g_{\alpha\beta})e_c^\beta\Longrightarrow e_3(e_c^\alpha)g_{\alpha\beta}e_c^\beta=-\frac{1}{2}(3R-r+t)^{-1}e_c^\alpha V_3(g_{\alpha\beta})e_c^\beta,}
so we have $V^I(e_3(e_c^\alpha)g_{\alpha\beta}e_c^\beta)=O(\eps t^{-1+C\eps}\lra{q}^{-1})$ on $H$. If $c'\neq c$, then
\fm{e_3(e_{c'}^\alpha)g_{\alpha\beta}e_c^\beta&=(3R-r+t)^{-1}V_3(e_{c'}^\alpha)g_{\alpha\beta}e_c^\beta,}
so we have $V^I(e_3(e_{c'}^\alpha)g_{\alpha\beta}e_c^\beta)=O(\lra{q}^{-1}t^{C\eps})$ on $H$ if $c\neq c'$. All these estimates imply that on $H$, we have
\fm{V^I(\xi_{ab}^c,\xi_{a4}^c,\xi_{c3}^c)=O(t^{-1+C\eps});\hspace{1.5em}V^I(\xi_{c'3}^c)=O(\lra{q}^{-1}t^{C\eps}),\ c\neq c';\hspace{1.5em}V^I(\xi_{34}^c)=O(t^{-1+C\eps}\lra{q}^{-1}).}
Moreover,
\fm{|V^I(\chi_{ab}-r^{-1}\delta_{ab})|&\leq|V^I(e_a(e_4^\alpha)g^{\alpha\beta}e_b^\beta-r^{-1}\delta_{ab})|+|V^I(e_a(e_4^\alpha)g^{\alpha\beta}e_b^\beta)|\lesssim t^{-2+C\eps}. }
\end{proof}
\rm

\subsection{Estimates in $\Omega$}\label{sec4.3new}\rm
Recall that we defined a weighted null frame $\{V_k\}_{k=1}^4$ in Section \ref{s5.1}. Our goal in this section is to prove the following proposition. Note that the estimates here are the same as those in Lemma \ref{derqxi}.

\prop{\label{prop5.5} In $\Omega\cap\{r-t<2R\}$, for each $I$ we have the following estimates:
\eq{\label{prop5.5c1}|V^I(\xi_{13}^{2})|+|V^I(\xi_{23}^{1})|\lesssim \lra{q}^{-1}t^{C\eps};}
\eq{\label{prop5.5c5}|V^I(\xi_{34}^{a})|\lesssim \lra{q}^{-1}t^{-1+C\eps};}
for all other $(k_1,k_2,a)$ such that $k_1<k_2$ and $a=1,2$, we have 
\eq{\label{prop5.5c4}|V^I(\xi_{k_1k_2}^a)|\lesssim t^{-1+C\eps};}
for all $k_1<k_2$, we have
\eq{\label{prop5.5c2}|V^I(\xi_{k_1k_2}^3)|\lesssim t^{-1+C\eps}\lra{q}^{-1};}
for $\xi_{a4}^b=\chi_{ab}$, we have
\eq{\label{prop5.5c3}|V^I(\chi_{ab}-r^{-1}\delta_{ab})|\lesssim t^{-2+C\eps}.}
In this proposition we use the convention given in Section \ref{sec2.5}. That is, for each fixed integer $N>0$, we can choose $\eps\ll_N1$, such that  the  estimates in this proposition hold for all multiindices $I$ with $|I|\leq N$.}
\rm\bigskip

Since it is known that $q=r-t$ for $r-t>R$, we only care about the region where $r-t<2R$ in this subsection.  Recall that every point in $\Omega\cap\{r-t<2R\}$ lies on exactly one geodesic in $\mcl{A}$ emanating from $H$. The  following lemma would be the key lemma in the proof of Proposition \ref{prop5.5}.

\lem{\label{l5.8}Fix $0<\eps\ll 1$. Let $Q_1,\dots,Q_m$ be $m$ functions defined in $\Omega\cap\{r-t<2R\}$. For each $i=1,\dots,m$, suppose  in $ \Omega\cap\{r-t<2R\}$ we have
\eq{\label{l5.8f1}e_4(Q_i)=(-n_0r^{-1}+n_1e_4(\ln (3R-r+t)))Q_i+O(\eps t^{-1}\sum_j|Q_j|)+O(f(t)).}
Here $n_0,n_1\geq 0$  are two fixed real numbers which do not depend on $i$. Moreover, for some fixed $s\geq 1$, we suppose that $Q_i|_H=O(h(t))$ for each $i$. Then, in $\Omega\cap\{r-t<2R\}$ we have  \eq{\label{l5.8c}\sum_i|Q_i|\lesssim t^{-n_0+C\eps}((x^0(0))^{n_0}h(x^0(0))+\int_{x^0(0)}^{t} \tau^{n_0+C\eps}f(\tau)\ d\tau).}
Here we suppose that $(t,x)$ lies on the geodesic $x(s)$ in $\mcl{A}$ and that the integral is taken along the geodeisc $x(s)$.}
\begin{proof}Recall that $e_4(r)=1+O(t^{-1+C\eps})$. If we define $Q'_i=(3R-r+t)^{-n_1}r^{n_0}Q_i$, then by \eqref{l5.8f1}, we have
\fm{e_4(Q'_i)&=-n_1(3R-r+t)^{-n_1-1}e_4(3R-r+t)r^{n_0}Q_i+n_0(3R-r+t)^{-n_1}r^{n_0-1}e_4(r)Q_i\\&\hspace{1em}+(3R-r+t)^{-n_1}r^{n_0}e_4(Q_i)\\
&=n_0r^{-1}(e_4(r)-1)Q'_i+O(\eps t^{-1}\sum_j|Q'_j|+(3R-r+t)^{-n_1}r^{n_0}f(t))\\
&=O(\eps t^{-1}\sum_j|Q'_j|+(3R-r+t)^{-n_1}r^{n_0}f(t)).}
To get the last equality, we note that $r^{-1}(e_4(r)-1)=O(t^{-2+C\eps})=O(\eps t^{-1})$ as $t\geq \exp(\delta/\eps)$.

In addition, we have $\lra{q}/\lra{r-t}=t^{O(\eps)}$. In fact, by Lemma \ref{l3qrt}, we have $|q-(r-t)|\lesssim t^{C\eps}$ and thus \fm{1+|q|\lesssim 1+|r-t|+t^{C\eps}\lesssim t^{C\eps}\lra{r-t}&\Longrightarrow \lra{r-t}^{-1}\lesssim \lra{q}^{-1}t^{C\eps}\\1+|r-t|\lesssim 1+|q|+t^{C\eps}\lesssim t^{C\eps}\lra{q}&\Longrightarrow \lra{q}^{-1}\lesssim \lra{r-t}^{-1}t^{C\eps}.}
Thus, in $\Omega\cap\{r-t<2R\}$ we have \fm{(3R-r+t)^{-n_1}r^{n_0}f(t)\lesssim \lra{q}^{-n_1}t^{n_0+C\eps}f(t).}

Fix a point $(t_0,x_0)$ in $\Omega\cap\{r-t<2R\}$, and let $x(s)$ be the unique geodesic in $\mcl{A}$ passing through $(t_0,x_0)$. Note that $t_0\geq x^0(0)\geq T_0$ and that $q$ remains constant along each geodesic in $\mcl{A}$. Then by integrating $e_4(Q'_i)$, we have \fm{\sum_i|Q'_i(t_0,x_0)|&\lesssim \sum_i|Q'_i(x(0))|+ \int_{x^0(0)}^{t_0} \eps \tau^{-1}\sum_j|Q'_j(\tau,y(\tau))|+ \lra{q}^{-n_1}\tau^{n_0+C\eps}f(\tau)\ d\tau\\
&\lesssim \lra{q}^{-n_1}(x^0(0))^{n_0}h(x^0(0))+ \int_{x^0(0)}^{t_0} \eps \tau^{-1}\sum_j|Q'_j(\tau,y(\tau))|+ \lra{q}^{-n_1}\tau^{n_0+C\eps}f(\tau)\ d\tau.}
Here $(\tau,y(\tau))$ is a reparameterization of $x(s)$ such that $y(t_0)=x_0$. By the Gronwall's inequality, we conclude that \fm{\sum_i|Q'_i(t_0,x_0)|\lesssim t_0^{C\eps}\lra{q}^{-n_1}((x^0(0))^{n_0}h(x^0(0))+\int_{x^0(0)}^{t_0} \tau^{n_0+C\eps}f(\tau)\ d\tau).} To end the proof, we multiply both sides by $r^{-n_0}(3R-r+t)^{n_1}$, and recall that $t\sim r$ in $\Omega\cap\{r-t<2R\}$. 
\end{proof}
\rm\bigskip

To prove Proposition \ref{prop5.5}, we induct  on $|I|$.

\subsubsection{The base case $I=0$.} From Section \ref{s5.1}, in $\Omega\cap\{r-t<2R\}$ we already have the following estimates: $\xi_{34}^3=O(|\Gamma|)=O(\min\{\eps t^{-1},\eps t^{-1+C\eps}\lra{r-t}^{-1}\})$, $\xi_{a4}^b=\chi_{ab}=\delta_{ab}r^{-1}+O(t^{-2+C\eps})=O(t^{-1})$, $\xi_{a3}^a=\chi_{aa}+O(\eps t^{-1})=O(t^{-1})$, $\xi_{a4}^3=\xi_{12}^3=0$. To control the rest $\xi$, we recall that
\eq{\label{p5.5f1}&\hspace{1.5em}\lra{R(e_k,e_l)e_r,e_s}=e_k^\alpha e_l^\beta e_r^\mu e_s^\nu R_{\alpha\beta\mu\nu}\\
&=e_k^\alpha e_l^\beta e_r^\mu e_s^\nu (\frac{1}{2}(\partial_\alpha\partial_\mu g_{\beta\nu}-\partial_\alpha\partial_\nu g_{\beta\mu}-\partial_\beta\partial_\mu g_{\alpha\nu}+\partial_\beta\partial_\nu g_{\alpha\mu})-\Gamma_{\beta\mu}^\delta\Gamma_{\delta\nu\alpha}+\Gamma_{\alpha\mu}^\delta\Gamma_{\delta\nu\beta}).}
If at most one of $k,l,r,s$ is equal to $3$, then we have $\lra{R(e_k,e_l)e_r,e_s}=O(\eps t^{-2+C\eps}\lra{r-t}^{-1})$ by Lemma \ref{l3betterest}. From the equations in Section \ref{s5.1} we have
\fm{|e_4(\xi_{34}^a)+r^{-1}\xi_{34}^a|&\lesssim t^{-2+C\eps}\sum_b|\xi_{34}^b|+\eps t^{-2+C\eps}\lra{q}^{-1},\\
|e_4(\xi_{12}^a)+r^{-1}\xi_{12}^a|&\lesssim  t^{-2+C\eps}\sum_b|\xi_{12}^b|+\eps t^{-2+C\eps}\lra{q}^{-1}.}
By Lemma \ref{l5.8} with $n_0=1$, $n_1=0$ and $f(t)=\eps t^{-2+C\eps}\lra{q}^{-1}$, we have \fm{|\xi_{34}^a|&\lesssim t^{-1+C\eps}(\lra{q}^{-1}(x^0(0))^{C\eps}+t^{C\eps}\lra{q}^{-1})\lesssim t^{-1+C\eps}\lra{q}^{-1},\\|\xi_{12}^a|&\lesssim t^{-1+C\eps}((x^0(0))^{C\eps}+t^{C\eps}\lra{q}^{-1})\lesssim t^{-1+C\eps}.}Here we get different estimates for $\xi_{34}^a$ and $\xi_{12}^a$ because their estimates on $H$ are different; see Lemma \ref{derqxi}.  

It  follows from Section \ref{s5.1} that $\xi_{a3}^3=\frac{1}{2}\xi_{34}^a+O(\eps t^{-1})=O(t^{-1+C\eps})$. It remains to estimate $\xi_{a3}^{a'}$ where $a\neq a'$. Note that
\fm{e_4(\xi_{a3}^{a'})&=(e_4+\Gamma_{\mu\nu}^0e_4^\mu e_4^\nu)(\chi_{aa'}+2e_a(g^{0\alpha})g_{\alpha\beta}e_{a'}^\beta+2g^{0\alpha}e_a^\beta\Gamma_{\beta\alpha}^\mu g_{\mu\nu}e_{a'}^{\nu})-\Gamma_{\mu\nu}^0e_4^\mu e_4^\nu\xi_{a3}^{a'}-\sum_c\xi_{34}^c\xi_{aa'}^c\\
&\hspace{1em}-\lra{R(e_4,e_3)e_a,e_{a'}}-\Gamma_{\alpha\beta}^0e_4^\alpha e_a^\beta\xi_{34}^{a'}+\Gamma_{\alpha\beta}^0e_4^\alpha e_{a'}^\beta\xi_{34}^a\\
&=e_4(\chi_{aa'})-\Gamma_{\mu\nu}^0e_4^\mu e_4^\nu\xi_{a3}^{a'}-\sum_c\xi_{34}^c\xi_{aa'}^c+O(\eps t^{-2+C\eps})=O(\eps t^{-1}|\xi_{a3}^{a'}|)+O(t^{-2+C\eps}).}
By Lemma \ref{l5.8} with $n_0=n_1=0$ and $f(t)=t^{-2+C\eps}$ we have $|\xi_{a3}^{a'}|\lesssim (x^0(0))^{-1+C\eps}t^{C\eps}$. Here note that if $(t,x)$ lies on a geodesic $x(s)$ in $\mathcal{A}$, then 
\fm{q(t,x)&=q(x(0))=r(x(0))-x^0(0)=\frac{T_0-x^0(0)}{2}+2R\Longrightarrow x^0(0)=T_0-2(q-2R).}
And since we only care about the region where $q< 2R$, we have $t\geq x^0(0)\sim(T_0+\lra{q})\geq \lra{q}$.
In conclusion, we prove Proposition \ref{prop5.5} in the case $I=0$.

\subsubsection{The general case.} Fix $m>0$. Suppose we have proved Proposition \ref{prop5.5} for all $|I|<m$. Our goal is to prove Proposition \ref{prop5.5} for $|I|= m$.

Under the induction hypotheses, we can prove a key lemma which is Lemma \ref{l5.6} below.  For  convenience, we introduce the following notation.

\defn{\rm Let $F=F(t,x)$ be a function with domain $\Omega\cap\{r-t<2R\}$.  For any integer $m\geq 0$ and any real numbers $s,p$, we write $F=\mathfrak{R}_{s,p}^m$ if for $\eps\ll_{s,p,m}1$ we have
\fm{\sum_{|I|\leq m}|V^I(F)|\lesssim t^{s+C\eps}\lra{q}^{p}\hspace{2em}\text{in }\Omega\cap\{r-t<2R\}.}

By the Leibniz's rule, we can easily prove that $\mathfrak{R}_{s_1,p_1}^{m_1}\cdot \mathfrak{R}_{s_2,p_2}^{m_2}=\mathfrak{R}_{s_1+s_2,p_1+p_2}^{\min\{m_1,m_2\}}$. In addition, under the induction hypotheses, we have
\eq{\label{prop5.5hypo}\begin{array}{c}(\xi_{13}^2,\xi_{23}^1)=\mathfrak{R}_{0,-1}^{m-1};\ \xi_{34}^a=\mathfrak{R}_{-1,-1}^{m-1};\ \xi_{k_1k_2}^a=\mathfrak{R}_{-1,0}^{m-1} \text{ for all other }k_1<k_2\text{ and }a=1,2;\\[.5em]
\xi_{k_1k_2}^3=\mathfrak{R}_{-1,-1}^{m-1}\text{ for all }k_1<k_2;\ \chi_{ab}-r^{-1}\delta_{ab}=\mathfrak{R}_{-2,0}^{m-1}.\end{array}} }\rm

\lem{\label{l5.6} For $\eps\ll_{m}1$, we have 
\eq{\label{l5.6f1}e_k^\alpha=\mathfrak{R}_{0,0}^m;}
\eq{\label{l5.6f11}(e_4^i-\omega_i,e_3^i-\omega_i)=\mathfrak{R}_{-1,0}^m;}
\eq{\label{l5.6f2}(g^{\alpha\beta}-m^{\alpha\beta}, g_{\alpha\beta}-m_{\alpha\beta})=\eps\mathfrak{R}_{-1,0}^{m+1},\ \Gamma_{\mu\nu}^\alpha=\eps\mathfrak{R}_{-1,-1}^{m+1};}
for each fixed $s\in\R$, we have
\eq{\label{l5.6f3}\omega_i=\mathfrak{R}_{0,0}^{m+1},\ (t^s,r^s)=\mathfrak{R}_{s,0}^{m+1},\ (3R-r+t)^s=\mathfrak{R}_{0,s}^{m+1}.}
}
\begin{proof}
We prove by induction. First, since $e_*^*=O(1)$, we have $e_k^\alpha=\mathfrak{R}_{0,0}^0$; by Lemma \ref{l3betterest}, we have $(e_4^i-\omega_i,e_3^i-\omega_i)=\mathfrak{R}_{-1,0}^0$. Besides, $(g_{**}-m_{**},g^{**}-m^{**})=O(\eps t^{-1+C\eps})$ and
\fm{|\Gamma|\lesssim |g||\partial g|\lesssim\eps t^{-1+C\eps}\lra{r-t}^{-1}\lesssim\eps t^{-1+C\eps}\lra{q}^{-1}.} Here we use the estimate $\lra{r-t}/\lra{q}=t^{O(\eps)}$. Besides, \fm{\sum_k|V_k(g)|\lesssim\sum_{k\neq 3}(t+r)|e_k(g)|+\lra{r-t}|\partial g|\lesssim \eps t^{-1+C\eps}.} Since $\Gamma$ is a linear combination of terms of the form $g\cdot\partial g$ with constant real coefficients, by Lemma \ref{l3betterest} we have \fm{\sum_k|V_k(\Gamma)|&\lesssim\sum_{k} (|V_k(g)||\partial g|+|g|\cdot|V_k(\partial g)|)\\
&\lesssim \eps t^{-1+C\eps}\cdot \eps\lra{r-t}^{-1}t^{-1+C\eps}+\sum_{k\neq 3}(t+r)|e_k(\partial g)|+\lra{r-t}|\partial^2g|\\
&\lesssim \eps  \lra{q}^{-1} t^{-1+C\eps}.} We thus obtain \eqref{l5.6f2} with $m=0$. Since $3R-r+t\sim\lra{r-t}$ in $\Omega\cap\{r-t<2R\}$,  \eqref{l5.6f3} with $m+1$ replaced by $0$ is obvious. In addition,  by writing $Vf:=(V_1f,V_2f,V_3f,V_4f)$, we have
\eq{\label{vformula}\left\{\begin{array}{l}V(t)=(0,0,-(3R-r+t),t);\\
V(r)=(re_1(r),re_2(r),(3R-r+t)(e_3^i\omega_i),te_4^i\omega_i);\\
V(\omega_i)=(e_1^i-\omega_ie_1(r),e_2^i-\omega_ie_2(r),r^{-1}(3R-r+t)(e_3^i-\omega_ie_3^j\omega_j),r^{-1}t(e_4^i-\omega_ie_4^j\omega_j));\\
V(3R-r+t)=(-re_1(r),-re_2(r),(3R-r+t)(-1-e_3^i\omega_i),t(1-e_4^i\omega_i))\end{array}\right.}
Since $e_3,e_4=\pm\partial_t+\partial_r+O(t^{-1+C\eps})\partial$, we have\eq{\label{a4formula}e_a(r)&=e_a^i\omega_i=\sum_i e_a^ie_4^i+\sum_ie_a^i(\omega_i-e_4^i)\\&=\lra{e_a,e_4}-(g^{\alpha\beta}-m^{\alpha\beta})e_4^\alpha e_a^\beta+\sum_ie_a^i(\omega_i-e_4^i)=O(t^{-1+C\eps}),\\
1-e_4^i\omega_i&=-\sum_i(e_4^i-\omega_i)\omega_i=O(t^{-1+C\eps}).}
Also note that for each fixed $s\in\R$ and for each funtion $\phi(t,x)$, $V(\phi^s)=s\phi^{s-1}V(\phi)$. Then, we have $V(\omega)=O(t^{C\eps})$, $V(t^s,r^s)=O(t^{s+C\eps})$, $V((3R-r+t)^s)=O(\lra{r-t}^st^{C\eps})$. We thus obtain \eqref{l5.6f3} with $m=0$. This finishes the proof in the base case.

In general, we assume that we have proved \eqref{l5.6f1}-\eqref{l5.6f3} with $m$ replaced by $n$ where $0\leq n<m$. We first prove \eqref{l5.6f1} with $m$ replaced by $n+1$. Fix a multiindex $I$ such that $|I|=n+1$. If $I=(I',4)$, note that $te_4(e_k^\alpha)$ is a linear combination (with constant real coefficients) of terms of the form $t\Gamma_{**}^* (e_*^*)(e_*^*)(e_*^*)$, $-t\Gamma_{**}^* (e_*^*)(e_*^*)$ and $V_4(g^{0\alpha})$. By the induction hypotheses,  we notice that
\fm{t\Gamma_{**}^* (e_*^*)(e_*^*)(e_*^*)=\mathfrak{R}_{1,0}^{n+1}\cdot\eps\mathfrak{R}_{-1,-1}^{n+1}\cdot \mathfrak{R}_{0,0}^{n}\cdot \mathfrak{R}_{0,0}^{n}\cdot \mathfrak{R}_{0,0}^{n}=\eps \mathfrak{R}_{0,-1}^{n}}and similarly \fm{t\Gamma_{**}^* (e_*^*)(e_*^*)=\eps\mathfrak{R}_{0,-1}^n.}
Besides, \fm{g^{0\alpha}-m^{0\alpha}=\eps\mathfrak{R}_{-1,0}^{n+1}\Longrightarrow V_k(g^{0\alpha})=\eps \mathfrak{R}_{-1,0}^n.}
So in conclusion, \fm{V_4(e_k^\alpha)=\eps\mathfrak{R}_{0,-1}^n\Longrightarrow V^I(e_k^\alpha)=O(\eps\lra{q}^{-1}t^{C\eps}).} 

If $I=(I',k')$ where $k'\neq 4$, then by the formulas at the end of Section \ref{s5.1}, we have
\fm{V_{k'}(e_4^\alpha)&=r\xi_{a4}^l e_l^\alpha+rt^{-1}V_4(e_{a}^\alpha)\\&=\mathfrak{R}_{1,0}^{n+1}\cdot\mathfrak{R}_{-1,0}^{m-1}\cdot \mathfrak{R}_{0,0}^{n}+\mathfrak{R}_{1,0}^{n+1}\cdot\mathfrak{R}_{-1,0}^{n+1}\cdot \eps\mathfrak{R}_{0,-1}^{n}=\mathfrak{R}_{0,0}^{n},\ \hspace{2em}k'=a=1,2;\\
V_3(e_4^\alpha)&=(3R-r+t)\xi_{34}^l e_l^\alpha+t^{-1}(3R-r+t)V_4(e_{3}^\alpha)\\&=\mathfrak{R}_{0,1}^{n+1}\cdot \mathfrak{R}_{-1,-1}^{m-1}\cdot \mathfrak{R}_{0,0}^{n}+\mathfrak{R}_{-1,0}^{n+1}\cdot \mathfrak{R}_{0,1}^{n+1}\cdot \eps\mathfrak{R}_{0,-1}^{n}=\mathfrak{R}_{-1,0}^{n}.}
In addition, note that $e_3^\alpha=e_4^\alpha+2g^{0\alpha}$, so  
\fm{V_{k'}(e_4^\alpha,e_3^\alpha)=\mathfrak{R}_{0,0}^{n}\Longrightarrow V^I(e_4^\alpha,e_3^\alpha)=O(t^{C\eps}).}

If $I=(I',3)$, we have
\fm{V_3(e_a^\alpha)&=(3R-r+t)\xi_{a3}^le_l^\alpha+r^{-1}(3R-r+t)V_a(e_3^\alpha)\\&=\mathfrak{R}_{0,1}^{n+1}\cdot\mathfrak{R}_{0,-1}^{m-1}\cdot \mathfrak{R}_{0,0}^{n}+\mathfrak{R}_{-1,0}^{n+1}\cdot \mathfrak{R}_{0,1}^{n+1}\cdot \mathfrak{R}_{0,0}^{n}=\mathfrak{R}_{0,0}^{n}.}
Here we recall that $t\gtrsim x^0(0)\sim\lra{q}+T_0$, so $\mathfrak{R}_{-s,s}^n=\mathfrak{R}_{0,0}^n$ for each $s>0$. Thus,\fm{V^I(e_a^\alpha)=O(t^{C\eps}).}

If $I=(I',a)$, then 
\fm{V_a(e_b^\alpha)&=-\sum_cr\xi_{bc}^ae_c^\alpha-\frac{1}{2}r\chi_{ab}(e_4^\alpha+e_3^\alpha)-(e_b^\mu g_{\mu\beta}V_a(g^{0\beta})+re_b^\mu g_{\mu\nu} g^{0\beta}e_a^\sigma\Gamma_{\sigma\beta}^\nu)e_4^\alpha-re_a^\mu e_b^\nu\Gamma_{\mu\nu}^\alpha.}
Again, by our induction hypotheses, we conclude that \fm{V_a(e_b^\alpha)=\mathfrak{R}_{0,0}^n\Longrightarrow V^I(e_b^\alpha)=O(t^{C\eps}).}
Summarize all the results above and we conclude that $e_*^*=\mathfrak{R}_{0,0}^{n+1}$. Note that the computations above work as long as $n\leq m-1$.

Next we prove \eqref{l5.6f11} with $m$ replaced by $n+1$. It suffices to consider $e_4^i-\omega_i$ as $e_3^i-e_4^i=2g^{0i}=\eps\mathfrak{R}_{-1,0}^{n+1}$. Fix a multiindex $I$ with $|I|=n+1$. Note that\fm{V_{a}(e_4^i-\omega_i)&=re_a(e_4^i-\omega_i)=r(\xi_{a4}^le_l^i+e_4(e_a^i)-r^{-1}(e_a^i-\omega_ie_a(r)))\\
&=r(\chi_{ab}-\delta_{ab}r^{-1})e_b^i+re_4(e_a^i)+r^{-1}\omega_iV_a(r)\\
&=\mathfrak{R}_{1,0}^{n+1}\cdot\mathfrak{R}_{-2,0}^{m-1}\cdot \mathfrak{R}_{0,0}^{n}+re_4(e_a^i)+\mathfrak{R}_{-1,0}^n=re_4(e_a^i)+\mathfrak{R}_{-1,0}^n,\\
V_4(e_4^i-\omega_i)&=te_4(e_4^i-\omega_i)=t(e_4(e_4^i)-(e_4^j-\omega_j)\partial_j\omega_i)\\
&=te_4(e_4^i)-tr^{-1}(e_4^i-\omega_i-\omega_i\omega_j(e_4^j-\omega_j))\\
&=te_4(e_4^i)+\mathfrak{R}_{0,0}^{n+1}\cdot(\mathfrak{R}_{-1,0}^{n}+\mathfrak{R}_{0,0}^{n+1}\cdot\mathfrak{R}_{-1,0}^{n})=te_4(e_4^i)+\mathfrak{R}_{-1,0}^{n},\\
V_3(e_4^i-\omega_i)&=(3R-r+t)e_3(e_4^i-\omega_i)=(3R-r+t)(\xi_{34}^le_l^i+e_4(e_3^i)-(e_3^j-\omega_j)\partial_j\omega_i)\\
&=(3R-r+t)(\xi_{34}^le_l^i+e_4(e_4^i)+2t^{-1}V_4(g^{0i})-r^{-1}(e_3^i-\omega_i-(e_3^j-\omega_j)\omega_i\omega_j))\\
&=(3R-r+t)e_4(e_4^i)+\mathfrak{R}_{0,1}^{n+1}\cdot (\mathfrak{R}_{-1,-1}^{m-1}\cdot\mathfrak{R}_{0,0}^{n}+\eps\mathfrak{R}_{-2,0}^{n}+\mathfrak{R}_{-1,0}^{n+1}\cdot \mathfrak{R}_{-1,0}^n )\\
&=(3R-r+t)e_4(e_4^i)+ \mathfrak{R}_{-1,0}^n .}
Here we use \eqref{prop5.5hypo}. To finish the proof, we note that for $k\neq 3$,
\fm{2e_4(e_k^i)&=2e_4^\alpha e_k^\beta(\Gamma_{\alpha\beta}^0e_4^i-\Gamma_{\alpha\beta}^i)=e_4^\alpha e_k^\beta (g^{0\delta}e_4^i-g^{i\delta})(\partial_\alpha g_{\delta\beta}+\partial_\beta g_{\delta\alpha}-\partial_\delta g_{\alpha\beta})\\
&=(g^{0\delta}e_4^i-g^{i\delta})(e_4(g_{\delta\beta})e_k^\beta+e_k(g_{\delta\alpha})e_4^\alpha)+e_4^\alpha e_k^\beta (-\frac{1}{2}e_4(g_{\alpha\beta}) (e_4^i+e_3^i)-\sum_b e_b^ie_b(g_{\alpha\beta}))\\
&=\mathfrak{R}_{0,0}^{n+1} t^{-1}V_4(g)+\mathfrak{R}_{0,0}^{n+1} r^{-1}V_a(g)=\eps\mathfrak{R}_{-2,0}^{n+1}.}
Also note that $e_4(g)=t^{-1}V_4(g)=\eps \mathfrak{R}_{-2,0}^{n+1}$ and that $e_k^0$ is a constant, so we have $e_4(e_k^\alpha)=\eps\mathfrak{R}_{-2,0}^{n+1}$ for each $k,\alpha$. Thus, \fm{V(e_4^i-\omega_i)=\mathfrak{R}_{-1,0}^{n+1}\Longrightarrow e_4^i-\omega_i=\mathfrak{R}_{-1,0}^{n+1}.}

Finally, we prove \eqref{l5.6f2} and \eqref{l5.6f3} with $m+1$ replaced by $n+2$. Fix a multiindex $I$ such that $|I|=n+2$.  Note that \fm{(3R+t-r)\partial_t&=3R\partial_t+\frac{tS-x_i\Omega_{0i}}{r+t}=\mathfrak{R}_{0,0}^{n+1}\cdot Z,\\
(3R+t-r)\partial_r&=3R\partial_r+\frac{t\omega_i\Omega_{0i}-rS}{r+t}=\mathfrak{R}_{0,0}^{n+1}\cdot Z,\\
(3R+t-r)\partial_i&=3R\partial_i+(t-r)\omega_i\partial_r+(t-r)r^{-1}\omega_j\Omega_{ji}=\mathfrak{R}_{0,0}^{n+1}\cdot Z.}Thus, $\partial=(3R+t-r)^{-1}\mathfrak{R}_{0,0}^{n+1}\cdot Z=\mathfrak{R}_{0,-1}^{n+1}\cdot Z$. Since we have just proved $e_*^*=\mathfrak{R}_{0,0}^{n+1}$ and $e_4^i-\omega_i=\mathfrak{R}_{-1,0}^{n+1}$, by \eqref{a4formula} we have $e_a(r)=\mathfrak{R}_{-1,0}^{n+1}$. In conclusion, by \eqref{derqvz} we have
\fm{V_4&=t(t+r)^{-1}S+(t+r)^{-1}t\omega_j\Omega_{0j}+t(e_4^i-\omega_i)\partial_i=\mathfrak{R}_{0,0}^{n+1}\cdot Z,\\
V_3&=(3R-r+t)r^{-1}V_4+2g^{0\alpha}(3R-r+t)\partial_\alpha=\mathfrak{R}_{0,0}^{n+1}\cdot Z,\\
V_a&=re_a(r)\omega_i\partial_i+e_a^i\omega_j\Omega_{ji}=\mathfrak{R}_{-1,0}^{n+1}\cdot \mathfrak{R}_{0,-1}^{n+1}\cdot Z+\mathfrak{R}_{0,0}^{n+1}\cdot Z=\mathfrak{R}_{0,0}^{n+1}\cdot Z.}
Now, given a function $F=F(t,x)$, if $|I|=n+2$, we can write $V^IF$ as a linear combination of terms of the form
\eq{\label{sfformula}V^{I_1}(\mathfrak{R}_{0,0}^{n+1})\cdots V^{I_s}(\mathfrak{R}_{0,0}^{n+1})Z^sF,\hspace{2em}\sum|I_*|+s=n+2,\ s>0.}
Since $|I_j|<n+2$ for each $j$, we have $V^{I_j}(\mathfrak{R}_{0,0}^{n+1})=O(t^{C\eps})$. Note that for each $J$ with $|J|>0$, we have $Z^Jg=O(\eps t^{-1+C\eps})$,  $Z^J\omega=O(1)$, $Z^J(t^s,r^s)=O(t^{s})$, $Z^J((3R-r+t)^s)=O(\lra{r-t}^{s})$ and $Z^J(\Gamma)=O(\eps t^{-1+C\eps}\lra{q}^{-1})$. The last one is true because $Z^J\Gamma$ is a linear combination (with constant real coefficients) of terms of the form $(Z^{J_1} g)\cdot (Z^{J_2}\partial g)=O(\eps t^{-1+C\eps}\lra{r-t}^{-1})$. By plugging these estimates into \eqref{sfformula}, we conclude \eqref{l5.6f2} and \eqref{l5.6f3} with $m+1$ replaced by $n+2$.
\end{proof}
\rmk{\rm \label{rmk5.6.1} We have $Z^I\partial^k g=\eps\mathfrak{R}_{-1,-k}^{m+1}$ for each $I$ and $k$, as long as $\eps\ll_{I,k}1$. This follows directly from \eqref{sfformula}, Lemma \ref{l2.1} and $[Z,\partial]=C\cdot \partial$.

From the proof, we note that $e_4(e_k^\alpha)=\eps\mathfrak{R}_{-2,0}^m$ and $e_a(r)=\mathfrak{R}_{-1,0}^m$. These estimates are better than what we can get from \eqref{l5.6f1} and \eqref{l5.6f3}.}
\rm 
\bigskip

By Lemma \ref{l5.6}, we have $e_4^i\omega_i-1=(e_4^i-\omega_i)\omega_i=\mathfrak{R}_{-1,0}^m$. This result can be improved as shown in the next lemma.
\lem{\label{derqle4r}For $\eps\ll_m1$, we have $e_4^i\omega_i-1=\eps\mathfrak{R}_{-1,0}^m$.}
\begin{proof}

By Lemma \ref{l5.6}, we have \fm{e_a^j\omega_j=-(g^{\alpha\beta}-m^{\alpha\beta})e_4^\alpha e_a^\beta+\sum_i e_a^i(\omega_i-e_4^i)=\mathfrak{R}_{-1,0}^{m}.}
Recall that \fm{g^{\alpha\beta}&=\sum_a e_a^\alpha e_a^\beta+\frac{1}{2}(e_4^\alpha e_3^\beta+e_3^\alpha e_4^\beta).}
Then,
\fm{g^{\alpha\beta}(\partial_\alpha (r-t))(\partial_\beta (r-t))&=\sum_a (e_a^i\omega_i) (e_a^j\omega_j)+(e_4^i\omega_i-1) (e_3^j\omega_j+1)\\&=\mathfrak{R}_{-2,0}^m+(e_4^i\omega_i-1) (2+(e_3^j-\omega_j)\omega_j).}
Meanwhile,  we have
\fm{g^{\alpha\beta}(\partial_\alpha (r-t))(\partial_\beta (r-t))&=g^{00}-2g^{0i}\omega_i+g^{ij}\omega_i\omega_j\\
&=-2g^{0i}\omega_i+(g^{ij}-m^{ij})\omega_i\omega_j=\eps\mathfrak{R}_{-1,0}^{m+1}.}
Thus,
\fm{e_4^i\omega_i-1&=(2+(e_3^j-\omega_j)\omega_j)^{-1}(\eps\mathfrak{R}_{-1,0}^m+\mathfrak{R}_{-2,0}^m)=(2+(e_3^j-\omega_j)\omega_j)^{-1}\cdot\eps\mathfrak{R}_{-1,0}^m.}
Here we note that $\mathfrak{R}_{-2,0}^m=\eps\mathfrak{R}_{-1,0}^m$ as $t\geq \exp(\delta/\eps)$.

Fix a multiindx $I$ with $|I|\leq m$. Then, $V^I(e_4^i\omega_i-1)$ is a linear combination of terms of the form
\fm{(2+(e_3^j-\omega_j)\omega_j)^{-s-1}V^{I_0}(\eps\mathfrak{R}_{-1,0}^m)V^{I_2}(2+(e_3^j-\omega_j)\omega_j)\cdots V^{I_s}(2+(e_3^j-\omega_j)\omega_j)}
where $\sum|I_*|=|I|\leq m$ such that $|I_k|>0$ for each $k>0$. Thus, we can replace $V^{I_*}(2+(e_3^j-\omega_j)\omega_j)$ with $V^{I_*}((e_3^j-\omega_j)\omega_j)$ in the product. By Lemma \ref{l5.6} we have $(e_3^j-\omega_j)\omega_j=\mathfrak{R}_{-1,0}^m$. Since $e_3^j-\omega_j=O(t^{-1+C\eps})$, we have
$2+(e_3^j-\omega_j)\omega_j\geq 1$ for $\eps\ll 1$. In conclusion, we have
\fm{|V^I(e_4^i\omega_i-1)|&\lesssim \eps t^{-1+C\eps}\cdot \max_{0\leq s\leq m}\{(t^{-1+C\eps})^s\}\lesssim \eps t^{-1+C\eps}.}
Thus, $e_4^i\omega_i-1=\eps \mathfrak{R}_{-1,0}^m$.
\end{proof}

\rm\bigskip

We can now control the curvature tensor terms.

\lem{\label{l5.7}We have  $\lra{R(e_4,e_k)e_l,e_p})=\eps 
\mathfrak{R}_{-2,-1}^m$ if $l,p\neq 3$. }
\begin{proof} 
By \eqref{p5.5f1}, we can express $e_4^\alpha e_k^\beta e_l^\mu e_p^\nu R_{\alpha\beta\mu\nu}$ as a linear combination of terms of the form
\fm{e_4(\partial_\mu g_{\beta\nu}-\partial_\nu g_{\beta\mu})e_k^\beta e_l^\mu e_p^\nu,\ e_l(\partial_{\beta}g_{\alpha\nu})e_4^\alpha e_k^\beta  e_p^\nu,\ e_p(\partial_\beta g_{\alpha\mu})e_4^\alpha e_k^\beta e_l^\mu,\ e_4^\alpha e_k^\beta e_l^\mu e_p^\nu\cdot \Gamma\cdot (g\cdot \Gamma).}
By Lemma \ref{l5.6} and Remark \ref{rmk5.6.1}, we have 
\fm{e_4(\partial_\mu g_{\beta\nu}-\partial_\nu g_{\beta\mu})e_k^\beta e_l^\mu e_p^\nu&=t^{-1}V_4(\partial g)\cdot\mathfrak{R}_{0,0}^{m}=\mathfrak{R}_{-1,0}^{m}\cdot Z(\partial g)=\eps\mathfrak{R}_{-2,-1}^m.}
Since $l\neq 3$, we either have $e_l=t^{-1}V_l$ or $e_l=r^{-1}V_l$. In both cases, we can follow the same proof as above to conclude that \fm{e_l(\partial_{\beta}g_{\alpha\nu})e_4^\alpha e_k^\beta  e_p^\nu=\eps\mathfrak{R}_{-2,-1}^m.}
Similarly, we also have
\fm{e_p(\partial_\beta g_{\alpha\mu})e_4^\alpha e_k^\beta e_l^\mu&=\eps\mathfrak{R}_{-2,-1}^m.}
Finally, note that 
\fm{e_4^\alpha e_k^\beta e_l^\mu e_p^\nu\cdot \Gamma\cdot (g\cdot \Gamma)&=(\eps\mathfrak{R}_{-1,-1}^{m+1})^2\cdot \mathfrak{R}_{0,0}^m=\eps^2\mathfrak{R}_{-2,-2}^m.}
Thus, $\lra{R(e_4,e_k)e_l,e_p}=\eps\mathfrak{R}_{-2,-1}^m$.
\end{proof}
\rm\bigskip

Lemma \ref{l5.7} can be improved in a special case.

\lem{\label{l410ct}{\rm (a)} We have
\fm{\lra{R(e_4,e_a)e_4,e_b}=e_4(f_{ab})+\frac{1}{4}e_4^\alpha  e_4^\mu r^{-1}\delta_{ab}e_3(g_{\alpha\mu})+\eps\mathfrak{R}_{-3,0}^m.}
Here we set
\fm{f_{ab}=\frac{1}{2}(e_a^\beta  e_b^\nu e_4(g_{\beta\nu})-e_a^\beta e_4^\mu e_b( g_{\beta\mu})) -\frac{1}{2}e_4^\alpha e_a(g_{\alpha\nu}) e_b^\nu=\eps\mathfrak{R}_{-2,0}^m.}

{\rm (b)} Assume that $\chi_{ab}=\mathfrak{R}_{-1,0}^m$. Then we have \fm{\Gamma_{\alpha\beta}^0 e_4^\alpha e_4^\beta\chi_{ab}+\frac{1}{4}e_4^\alpha e_4^\beta e_3(g_{\alpha\beta})\chi_{ab}=\eps\mathfrak{R}_{-3,0}^m.}

}
\begin{proof} (a) Recall that $\lra{R(e_4,e_a)e_4,e_b}=e_4^\alpha e_a^\beta e_4^\mu e_b^\nu R_{\alpha\beta\mu\nu}$ where $R_{\alpha\beta\mu\nu}$ is  given by \fm{R_{\alpha\beta\mu\nu}&=\frac{1}{2}(\partial_\alpha\partial_\mu g_{\beta\nu}-\partial_\alpha\partial_\nu g_{\beta\mu}-\partial_\beta\partial_\mu g_{\alpha\nu}+\partial_\beta\partial_\nu g_{\alpha\mu})-\Gamma_{\beta\mu}^\delta\Gamma_{\delta\nu\alpha}+\Gamma_{\alpha\mu}^\delta\Gamma_{\delta\nu\beta}.}
Note that (for simplicity we take the sum over all the indices without writing the summation)
\fm{&\hspace{1.5em}\frac{1}{2}e_4^\alpha e_a^\beta e_4^\mu e_b^\nu\partial_\beta\partial_\nu g_{\alpha\mu}\\&=\frac{1}{2}e_4^\alpha e_a^\beta e_4^\mu e_b^i(\omega_i\partial_r)(\partial_\beta g_{\alpha\mu})+\frac{1}{2}e_4^\alpha e_a^\beta e_4^\mu e_b^i(\partial_i-\omega_i\partial_r)(\partial_\beta g_{\alpha\mu})\\
&=\frac{1}{2}e_4^\alpha  e_4^\mu e_b(r)e_a^\beta\partial_r(\partial_\beta g_{\alpha\mu})+\frac{1}{2}e_4^\alpha e_a^\beta e_4^\mu e_b^ir^{-1}\omega_j\Omega_{ji}(\partial_\beta g_{\alpha\mu})\\
&=\frac{1}{2}e_4^\alpha  e_4^\mu e_b(r)\omega_je_a(\partial_j g_{\alpha\mu})+\frac{1}{2}e_4^\alpha e_a^\beta e_4^\mu e_b^i r^{-1}\omega_j[\Omega_{ji},\partial_\beta] (g_{\alpha\mu})+\frac{1}{2}e_4^\alpha  e_4^\mu e_b^i r^{-1}\omega_je_a(\Omega_{ji} g_{\alpha\mu})\\
&=\frac{1}{2}e_4^\alpha  e_4^\mu e_b(r)\omega_je_a(\partial_j g_{\alpha\mu})+\frac{1}{2}e_4^\alpha  e_4^\mu r^{-1}(-e_a(r)e_b(g_{\alpha\mu})+e_a^ie_b^i\partial_r(g_{\alpha\mu})) \\&\hspace{1em}+\frac{1}{2}e_4^\alpha  e_4^\mu e_b^i r^{-1}\omega_je_a(\Omega_{ji} g_{\alpha\mu})\\
&=\frac{1}{2}e_4^\alpha  e_4^\mu e_b(r)\omega_je_a(\partial_j g_{\alpha\mu})+\frac{1}{2}e_4^\alpha  e_4^\mu r^{-1}(-e_a(r)e_b(g_{\alpha\mu})+(\delta_{ab}-e_a^\beta(g_{\beta\nu}-m_{\beta\nu})e_b^\nu)\partial_r(g_{\alpha\mu}))\\&\hspace{1em} +\frac{1}{2}e_4^\alpha  e_4^\mu e_b^i r^{-1}\omega_je_a(\Omega_{ji} g_{\alpha\mu})\\
&=\frac{1}{2}e_4^\alpha  e_4^\mu r^{-1}e_b(r)\omega_jV_a(\partial_j g_{\alpha\mu})+\frac{1}{2}e_4^\alpha  e_4^\mu r^{-1}(-r^{-1}e_a(r)V_b(g_{\alpha\mu})+(\delta_{ab}-e_a^\beta(g_{\beta\nu}-m_{\beta\nu})e_b^\nu)\partial_r(g_{\alpha\mu}))\\&\hspace{1em} +\frac{1}{2}e_4^\alpha  e_4^\mu e_b^i r^{-2}\omega_jV_a(\Omega_{ji} g_{\alpha\mu}).}
Recall that in Lemma \ref{l5.6}, we have proved that $e_a(r)=\mathfrak{R}_{-1,0}^m$. Thus, we have
\fm{&\hspace{1.5em}\frac{1}{2}e_4^\alpha e_a^\beta e_4^\mu e_b^\nu\partial_\beta\partial_\nu g_{\alpha\mu} =\frac{1}{2}e_4^\alpha e_4^\mu r^{-1}\delta_{ab}(\partial_rg_{\alpha\mu})+\eps\mathfrak{R}_{-3,0}^m\\&=\frac{1}{2}e_4^\alpha e_4^\mu r^{-1}\delta_{ab}(\omega_j-\frac{1}{2}e_3^j-\frac{1}{2}e_4^j)\partial_jg_{\alpha\mu}+\frac{1}{4}e_4^\alpha e_4^\mu r^{-1}\delta_{ab}(e_3(g_{\alpha\mu})+e_4(g_{\alpha\mu}))+\eps\mathfrak{R}_{-3,0}^m\\&=\frac{1}{4}e_4^\alpha e_4^\mu r^{-1}\delta_{ab}e_3(g_{\alpha\mu})+\eps\mathfrak{R}_{-3,0}^m.}

Next, we note that
\fm{&\hspace{1.5em}\frac{1}{2}e_4^\alpha e_a^\beta e_4^\mu e_b^\nu(\partial_\alpha\partial_\mu g_{\beta\nu}-\partial_\alpha\partial_\nu g_{\beta\mu}-\partial_\beta\partial_\mu g_{\alpha\nu})\\
&=\frac{1}{2} e_a^\beta e_4^\mu e_b^\nu e_4(\partial_\mu g_{\beta\nu}-\partial_\nu g_{\beta\mu})-\frac{1}{2}e_4^\alpha e_a^\beta e_b^\nu e_4(\partial_\beta g_{\alpha\nu})\\
&=e_4(f_{ab})-\frac{1}{2} e_4(e_a^\beta e_4^\mu e_b^\nu) (\partial_\mu g_{\beta\nu}-\partial_\nu g_{\beta\mu})-\frac{1}{2}e_4(e_4^\alpha e_a^\beta e_b^\nu )(\partial_\beta g_{\alpha\nu}).}
In Lemma \ref{l5.6}, we have proved that $e_4(e_k^\alpha)=\eps\mathfrak{R}_{-2,0}^m$. By Lemma \ref{l5.6}, we can easily prove that $f_{ab}=\eps\mathfrak{R}_{-2,0}^m$. This implies that 
\fm{\frac{1}{2}e_4^\alpha e_a^\beta e_4^\mu e_b^\nu(\partial_\alpha\partial_\mu g_{\beta\nu}-\partial_\alpha\partial_\nu g_{\beta\mu}-\partial_\beta\partial_\mu g_{\alpha\nu})=e_4(f_{ab})+\eps\mathfrak{R}_{-3,0}^m.}

Finally, we note that
\fm{&\hspace{1.5em}e_4^\alpha e_a^\beta e_4^\mu e_b^\nu (-\Gamma_{\beta\mu}^\delta\Gamma_{\delta\nu\alpha}+\Gamma_{\alpha\mu}^\delta\Gamma_{\delta\nu\beta})\\
&=-\frac{1}{2}e_4^\alpha e_a^\beta e_4^\mu e_b^\nu \Gamma_{\beta\mu}^\delta\Gamma_{\delta\nu\alpha}+\frac{1}{2}e_4^\alpha e_a^\beta e_4^\mu e_b^\nu \Gamma_{\alpha\mu}^\delta\Gamma_{\delta\nu\beta}\\
&=-\frac{1}{2}e_4^\alpha e_a^\beta e_4^\mu e_b^\nu g^{\delta\sigma}(\partial_\beta g_{\mu\sigma}+\partial_\mu g_{\beta\sigma}-\partial_\sigma g_{\beta\mu})(\partial_\alpha g_{\nu\delta}+\partial_\nu g_{\alpha\delta}-\partial_\delta g_{\alpha\nu})\\
&\hspace{1em}+\frac{1}{2}e_4^\alpha e_a^\beta e_4^\mu e_b^\nu g^{\delta\sigma}(\partial_\alpha g_{\mu\sigma}+\partial_\mu g_{\alpha\sigma}-\partial_\sigma g_{\alpha\mu})(\partial_\beta g_{\nu\delta}+\partial_\nu g_{\beta\delta}-\partial_\delta g_{\beta\nu}).}
Note that in the expansion of the right hand side, each term contains a product $e_k(g)\cdot e_l(g)$ where $l\neq 3$, except
\fm{I:=-\frac{1}{2}e_4^\alpha e_a^\beta e_4^\mu e_b^\nu g^{\delta\sigma}\partial_\sigma g_{\beta\mu}\partial_\delta g_{\alpha\nu}+\frac{1}{2}e_4^\alpha e_a^\beta e_4^\mu e_b^\nu g^{\delta\sigma}\partial_\sigma g_{\alpha\mu}\partial_\delta g_{\beta\nu}.}
Now we apply $g^{\delta\sigma}=\sum_ae_a^\delta e_a^\sigma+\frac{1}{2}(e_3^\delta e_4^\sigma+e_3^\sigma e_4^\delta)$. Then, we can also write $I$ as a sum of several terms containing $e_k(g)\cdot e_l(g)$ where $l\neq 3$. Since $e_l(g)=V_l(g)\cdot\mathfrak{R}_{-1,0}^{m+1}$, the whole sum is $\eps^2\mathfrak{R}_{-3,0}^m$. Combine all the disccussion above and we finish the proof.

(b) We have
\fm{\Gamma_{\alpha\beta}^0 e_4^\alpha e_4^\beta\chi_{ab}&=\frac{1}{2}g^{0\mu}(e_4^\beta e_4(g_{\beta\mu})+e_4^\alpha e_4(g_{\alpha\mu})-e_4^\alpha e_4^\beta\partial_\mu g_{\alpha\beta})\chi_{ab}\\
&=-\frac{1}{2}g^{0\mu}e_4^\alpha e_4^\beta\partial_\mu g_{\alpha\beta}\chi_{ab}+\mathcal{R}=-\frac{1}{4}e_4^\alpha e_4^\beta(e_3(g_{\alpha\beta})-e_4(g_{\alpha\beta}))\chi_{ab}+\mathcal{R}\\
&=-\frac{1}{4}e_4^\alpha e_4^\beta e_3(g_{\alpha\beta})\chi_{ab}+\mathcal{R}.}
Here the remainder $\mathcal{R}$ is a linear combination of $g\cdot (e_*^*)\cdot e_4(g)\cdot \chi$ or $(e_*^*)\cdot (e_*^*)\cdot e_4(g)\cdot \chi$. Since $e_4(g)=t^{-1}V_4(g)=\eps\mathfrak{R}_{-2,0}^m$ and $(g,e_*^*)=\mathfrak{R}_{0,0}^m$, under our assumption on $\chi$, it follows from the Leibniz's rule that $\mathcal{R} =\eps \mathfrak{R}_{-3,0}^m$.
\end{proof}
\rmk{\rm \label{rmk410ct.1} Note  we only have $\chi=\mathfrak{R}_{-1,0}^{m-1}$ from our induction hypotheses, so we cannot apply (b) directly assuming \eqref{prop5.5hypo} only.}
\rm\bigskip

We now prove Proposition \ref{prop5.5} for $|I|= m$. Fix a multiindex $I$ such that $|I|=m$. We have
\fm{\ [V_4,V_4]&=0,\\
[V_4,V_a]&=t(e_4^i-\omega_i)\omega_ie_a-t(r\chi_{ab}-\delta_{ab})e_b,\\
[V_4,V_3]&=-t(e_4^i-\omega_i)\omega_i e_3+(3R-r+t)e_4-t(3R-r+t)\xi_{34}^l e_l.}
We write $[V_4,V_k]:=\eta_k^lV_l$. Then by Lemma \ref{l5.6},   Lemma \ref{derqle4r} and the inudction hypotheses \eqref{prop5.5hypo}, we have
\eq{\label{derqetaformula}\left\{\begin{array}{l}
\eta_a^a=(e_4^i-\omega_i)\omega_itr^{-1}-t(\chi_{aa}-r^{-1})=\mathfrak{R}_{-1,0}^{m-1};\\[.5em]
\eta_{a}^{a'}=-t\chi_{12}=\mathfrak{R}_{-1,0}^{m-1},\hspace{1em}a\neq a'\\[.5em]
\eta_3^3=-t(e_4^i-\omega_i)\omega_i(3R-r+t)^{-1}-t\xi_{34}^3=\eps\mathfrak{R}_{0,-1}^m;\\[.5em]
\eta_3^4=(3R-r+t)t^{-1}-(3R-r+t)\xi_{34}^4=\mathfrak{R}_{-1,1}^m\\[.5em]
\eta_{3}^a=-(3R-r+t)\xi_{34}^atr^{-1}=\mathfrak{R}_{-1,0}^{m-1};\\[.5em]
\eta_*^*\equiv 0\text{ in all other cases.}\end{array}\right.}
In summary we have $\eta^*_*=\mathfrak{R}_{-1,1}^{m-1}$. Here we briefly explain why $\eta_{3}^3=\eps\mathfrak{R}_{0,-1}^m$, since all other estimates are clear. Note that $(e_4^i-\omega_i)\omega_i=\eps\mathfrak{R}_{-1,0}^m$ by Lemma  \ref{derqle4r}. Also note that $\xi_{34}^4=\xi_{34}^3=e_4^\alpha e_4^\beta\Gamma_{\alpha\beta}^0=\eps\mathfrak{R}_{-1,-1}^m$. Thus,
\fm{\eta_3^3&=-t(e_4^i-\omega_i)\omega_i(3R-r+t)^{-1}-t\xi_{34}^3=\mathfrak{R}_{1,0}^{m+1}\cdot \eps \mathfrak{R}_{-1,0}^m\cdot\mathfrak{R}_{0,-1}^{m+1}+\mathfrak{R}_{1,0}^{m+1}\cdot \eps\mathfrak{R}_{-1,-1}^m=\eps\mathfrak{R}_{0,-1}^m.}
In addition, since $\Gamma=O(\eps t^{-1})$, we have \fm{\eta^{3}_3&=(3R-r+t)^{-1}te_4(3R-r+t)-t\xi_{34}^3=V_4(\ln(3R-r+t))+O(\eps).}

Next, we note that
\eq{\label{p5.5f2}&\hspace{1.5em}V_4(V^I(\xi_{k_1k_2}^{l_1}))\\&=\sum_{(J,k,J')=I}V^J[V_4,V_k]V^{J'}(\xi_{k_1k_2}^{l_1})+V^I(V_4(\xi_{k_1k_2}^{l_1}))\\
&=\sum_{(J,k,J')=I}V^J(\eta_k^lV_l(V^{J'}(\xi_{k_1k_2}^{l_1})))+V^I(V_4(\xi_{k_1k_2}^{l_1}))\\
&=\sum_{(J,k,J')=I} \eta_k^lV^{(J,l,J')}(\xi_{k_1k_2}^{l_1})+\sum_{|J_1|+|J_2|=m\atop 0<|J_1|<m}C_{J_1,J_2}V^{J_1}(\eta_k^l)V^{J_2}(\xi_{k_1k_2}^{l_1})+V^I(V_4(\xi_{k_1k_2}^{l_1}))\\
&=:Q_1+Q_2+Q_3.}

In $Q_1$, we note that if $\eta_k^l\not\equiv 0$, then we must have $n_{(J,l,J'),3}\leq n_{(J,k,J'),3}$. Recall that $n_{J,3}$ denotes the number of $V_3$ in the product $V^J$. This is because $\eta_{k}^3\equiv 0$ for $k\neq 3$. In addition, we note that $n_{(J,l,J'),3}<n_{(J,k,J'),3}$ if $k=3$ and $l\neq 3$. Then,
\eq{\label{q1formula}Q_1&=(n_{I,3}\eta_{3}^3-\sum_an_{I,a}\eta_{a}^a)V^{I}(\xi_{k_1k_2}^{l_1})+O((|\eta_{1}^2|+|\eta_2^1|)\sum_{|J|=m\atop n_{J,3}=n_{I,3}}|V^J(\xi_{k_1k_2}^{l_1})|)\\&\hspace{1em}+O(
\sum_{l\neq 3}|\eta_{3}^l|\sum_{(J_1,3,J_2)=I}|V^{(J_1,l,J_2)}(\xi_{k_1k_2}^{l_1})|)\\
&=n_{I,3}V_4(\ln(3R-r+t))V^{I}(\xi_{k_1k_2}^{l_1})+O((\eps+t^{-1+C\eps})\sum_{|J|=m,\atop n_{J,3}=n_{I,3}}|V^J(\xi_{k_1k_2}^{l_1})|)\\&\hspace{1em}+O(\lra{q}t^{-1+C\eps}\sum_{|J|=m,\atop n_{J,3}<n_{I,3}}|V^J(\xi_{k_1k_2}^{l_1})|).}

In $Q_2$, we have $|J_1|,|J_2|<m$.  Since $\eta_*^*=\mathfrak{R}_{-1,1}^{m-1}$, we have \eq{\label{q2formula}|Q_2|&\lesssim \sum_{|J_1|+|J_2|=m\atop 0<|J_1|<m}|V^{J_1}(\mathfrak{R}_{-1,1}^{m-1})V^{J_2}(\xi_{k_1k_2}^{l_1})|\lesssim t^{-1+C\eps}\lra{q}\sum_{0<|J|<m}|V^{J}(\xi_{k_1k_2}^{l_1})|.}

Now we combine  \eqref{p5.5f2} with  Section \ref{s5.1}. First, note that $\xi_{34}^3=\Gamma^{0}_{\alpha\beta}e_4^\alpha e_4^\beta=\eps\mathfrak{R}_{-1,-1}^{m}$ by Lemma \ref{l5.6}, so $|V^I(\xi_{34}^3)|\lesssim \eps t^{-1+C\eps}\lra{q}^{-1}$ whenever $|I|\leq m$. There is no need to apply \eqref{p5.5f2}.

Next, we consider $\chi_{ab}=\xi_{a4}^b$. 

\prop{\label{prop5.5part1} Under our induction hypotheses \eqref{prop5.5hypo}, for $|I|=m$ we have \fm{|V^I(\chi_{ab})|\lesssim t^{-1+C\eps} ,\hspace{2em}|V^I(\chi_{ab}-r^{-1}\delta_{ab})|\lesssim t^{-2+C\eps}.}So $\chi_{ab}=\mathfrak{R}_{-1,0}^m$ and $\chi_{ab}-r^{-1}\delta_{ab}=\mathfrak{R}_{-2,0}^m$. }
\begin{proof} We first prove that $V^I(\chi_{ab})=O(t^{-1+C\eps})$ whenever $|I|=m$. Fix $I$ such that $|I|=m$ and $n_{I,3}=n\leq m$. Recall from \eqref{prop5.5hypo} that  $\chi_{ab}=\mathfrak{R}_{-1,0}^{m-1}$ and $\chi_{ab}-r^{-1}\delta_{ab}=\mathfrak{R}_{-2,0}^{m-1}$. Suppose that we have proved $V^J(\chi_{ab})=O(t^{-1+C\eps})$  for all $J$ such that $|J|=m$ and $n_{J,3}<n$. Note that
\fm{\chi_{ac}\chi_{cb}&=\delta_{ab}r^{-2}+2(\chi_{ab}-\delta_{ab}r^{-1})r^{-1}+(\chi_{ac}-\delta_{ac}r^{-1})(\chi_{cb}-\delta_{cb}r^{-1}).} By Lemma \ref{l5.6}, we have  $r^{-1}=\mathfrak{R}_{-1,0}^{m+1}$ and $t=\mathfrak{R}_{1,0}^{m+1}$. Also note that $V(tr^{-1})=V((t-r)r^{-1})=\mathfrak{R}_{-1,1}^{m}$. Thus, \fm{&\hspace{1.5em}|\sum_cV^I(t\chi_{ac}\chi_{cb})-2tr^{-1}V^I(\chi_{ab}-\delta_{ab}r^{-1})-V^I(\delta_{ab}r^{-2}t)|
\\&\lesssim\sum_{|J_1|+|J_2|=m\atop |J_1|>0}|V^{J_1}(tr^{-1})V^{J_2}(\chi_{ab}-r^{-1}\delta_{ab})|+t|\chi_{**}-\delta_{**}r^{-1}||V^I(\chi_{**}-\delta_{**}r^{-1})|\\&\hspace{1em}+\sum_{|J_1|+|J_2|+|J_3|=m\atop |J_2|<m,\ |J_3|<m}|V^{J_1}(t)V^{J_2}(\chi_{**}-\delta_{**}r^{-1})V^{J_3}(\chi_{**}-\delta_{**}r^{-1})|
\\&\lesssim \lra{q}t^{-3+C\eps}+t^{-1+C\eps}|V^I(\chi_{**}-\delta_{**}r^{-1})|.}
By the Raychaudhuri equation, we have \fm{&\hspace{1.5em}V^I(V_4(\chi_{ab}))=V^I(t\Gamma_{\alpha\beta}^0e_4^\alpha e_4^\beta\chi_{ab})-\sum_cV^I(t\chi_{ac}\chi_{cb})+V^I(t\lra{R(e_4,e_a)e_4,e_b})\\
&=t\Gamma_{\alpha\beta}^0e_4^\alpha e_4^\beta V^I(\chi_{ab})+O(\sum_{|J_1|+|J_2|=m\atop |J_2|<m}|V^{J_1}(\eps \mathfrak{R}_{0,-1}^m)V^{J_2}(\chi_{ab})|)\\&\hspace{1em}-2tr^{-1}V^I(\chi_{ab}-\delta_{ab}r^{-1})-V^I(\delta_{ab}r^{-2}t)+O(\lra{q}t^{-3+C\eps}+t^{-1+C\eps}|V^I(\chi_{**}-\delta_{**}r^{-1})|)\\&\hspace{1em}+V^I(\eps t\mathfrak{R}_{-2,-1}^m)\\
&=-2tr^{-1}V^I(\chi_{ab})+O((\eps+t^{-1+C\eps})|V^I(\chi_{**})|)+O(t^{-1+C\eps}).}
Besides, by \eqref{q1formula} and our induction hypotheses, we have
\fm{|Q_1-nV_4(\ln(3R-r+t))V^I(\chi_{ab})|&\lesssim \eps\sum_{|J|=m,\atop n_{J,3}=n}|V^J(\chi_{ab})|+\lra{q}t^{-1+C\eps}\sum_{|J|=m,\atop n_{J,3}<n }|V^J(\chi_{ab})|\\
&\lesssim\eps\sum_{|J|=m\atop n_{J,3}=n}|V^J(\chi_{ab})|+\lra{q}t^{-2+C\eps}.}
By \eqref{q2formula}  and our induction hypotheses, we have
\fm{|Q_2|&\lesssim t^{-1+C\eps}\lra{q}\sum_{|J|<m}|V^{J}(\chi_{ab})|\lesssim t^{-2+C\eps}\lra{q}.}
In conclusion, by \eqref{p5.5f2} we have
\fm{&\hspace{1.5em}|e_4(V^I(\chi_{ab}))+(-ne_4(\ln(3R-r+t))+2r^{-1})V^I(\chi_{ab})|\\
&\lesssim t^{-1}(|Q_1-nV_4(\ln(3R-r+t))V^I(\chi_{ab})|+|Q_2|+|V^I(V_4(\chi_{ab}))+2tr^{-1}V^I(\chi_{ab})|)\\&\lesssim \eps t^{-1}\sum_{c,c'}\sum_{|J|=m\atop n_{J,3}=n}|V^J(\chi_{cc'})|+t^{-2+C\eps}+\lra{q}t^{-3+C\eps}\lesssim \eps t^{-1}\sum_{c,c'}\sum_{|J|=m\atop n_{J,3}=n}|V^J(\chi_{cc'})|+t^{-2+C\eps}.}
The last inequality holds as $\lra{q}\lesssim t$. By Lemma \ref{l5.8} with $n_0=2$, $n_1=n$ and Lemma \ref{derqxi}, we conclude that \fm{\sum_{a,b}\sum_{|I|=m\atop n_{I,3}=n}|V^I(\chi_{ab})|&\lesssim t^{-2+C\eps}(x^0(0)^{2}\cdot x^{0}(0)^{-1+C\eps}+\int_{x^0(0)}^t \tau^{2+C\eps}\cdot\tau^{-2+C\eps}\ d\tau)\\&\lesssim t^{-2+C\eps}\cdot t^{1+C\eps}\lesssim t^{-1+C\eps}.}
By induction we obtain  $\chi_{ab}=\mathfrak{R}_{-1,0}^m$.

Next we prove $V^I(\chi_{ab}-r^{-1}\delta_{ab})=O(t^{-2+C\eps})$ whenever $|I|=m$. Again  fix $I$ such that $|I|=m$ and $n_{I,3}=n\leq m$. Suppose we have proved that $V^J(\chi_{ab}-r^{-1}\delta_{ab})=O(t^{-2+C\eps})$  for $|J|=m$ and $n_{J,3}<n$. Now we can apply Lemma \ref{l410ct}. We have
\fm{&\hspace{1.5em}V^I(V_4(\chi_{ab}))=V^I(t\Gamma_{\alpha\beta}^0e_4^\alpha e_4^\beta\chi_{ab})-\sum_cV^I(t\chi_{ac}\chi_{cb})+V^I(t\lra{R(e_4,e_a)e_4,e_b})\\
&=V^I(-\frac{1}{4}e_4^\alpha e_4^\beta e_3(g_{\alpha\beta})t\chi_{ab}+t\eps \mathfrak{R}_{-3,0}^m)+V^I(V_4(f_{ab})+\frac{1}{4}e_4^\alpha e_4^\beta tr^{-1}\delta_{ab}e_3(g_{\alpha\beta})+t\eps \mathfrak{R}_{-3,0}^m)\\&\hspace{1em}-2tr^{-1}V^I(\chi_{ab}-\delta_{ab}r^{-1})-V^I(\delta_{ab}r^{-2}t)+O(t^{-3+C\eps}\lra{q}+t^{-1+C\eps}|V^I(\chi_{**}-r^{-1}\delta_{**})|)\\
&=V^I(-\frac{1}{4}e_4^\alpha e_4^\beta e_3(g_{\alpha\beta})t(\chi_{ab}-r^{-1}\delta_{ab}))+V^I(V_4(f_{ab}))+O(\eps t^{-2+C\eps})\\&\hspace{1em}-2tr^{-1}V^I(\chi_{ab}-\delta_{ab}r^{-1})-V^I(\delta_{ab}r^{-2}t)+O(t^{-3+C\eps}\lra{q}+t^{-1+C\eps}|V^I(\chi_{**}-r^{-1}\delta_{**})|).}
Also note that
\fm{V^I(V_4(r^{-1}))&=V^I(te_4(r^{-1}))=V^I(-tr^{-2}e_4(r))}
and that $e_4(r)-1=\eps\mathfrak{R}_{-1,0}^m$ by Lemma \ref{derqle4r}.
In conclusion,
\fm{&\hspace{1.5em}V^I(V_4(\chi_{ab}-r^{-1}\delta_{ab}-f_{ab}))\\&=V^I(-\frac{1}{4}e_4^\alpha e_4^\beta e_3(g_{\alpha\beta})t(\chi_{ab}-r^{-1}\delta_{ab}))-2tr^{-1}V^I(\chi_{ab}-\delta_{ab}r^{-1})+V^I(\delta_{ab}r^{-2}t(e_4(r)-1))\\&\hspace{1em}+O(t^{-3+C\eps}\lra{q}+\eps t^{-2+C\eps}+t^{-1+C\eps}|V^I(\chi_{**}-r^{-1}\delta_{**})|)\\
&=V^I(-\frac{1}{4}e_4^\alpha e_4^\beta e_3(g_{\alpha\beta})t(\chi_{ab}-r^{-1}\delta_{ab}))-2tr^{-1}V^I(\chi_{ab}-\delta_{ab}r^{-1})\\&\hspace{1em}+O(t^{-3+C\eps}\lra{q}+\eps t^{-2+C\eps}+t^{-1+C\eps}|V^I(\chi_{**}-r^{-1}\delta_{**})|).}
 Besides, we note that
\fm{V^I(-\frac{1}{4}e_4^\alpha e_4^\beta e_3(g_{\alpha\beta})t(\chi_{ab}-r^{-1}\delta_{ab}))+\frac{1}{4}e_4^\alpha e_4^\beta e_3(g_{\alpha\beta})tV^I(\chi_{ab}-r^{-1}\delta_{ab})} is a linear combination of terms of the form
\fm{V^{I_1}(e_4^\alpha e_4^\beta t(3R-r+t)^{-1}V_3(g_{\alpha\beta}))V^{I_2}(\chi_{ab}-r^{-1}\delta_{ab})}
where $|I_1|+|I_2|=|I|=m$ and $|I_2|<m$. By the induction hypotheses and since 
\fm{e_4^\alpha e_4^\beta t(3R-r+t)^{-1}V_3(g_{\alpha\beta})=\mathfrak{R}_{1,-1}^m\cdot\eps \mathfrak{R}_{-1,0}^{m}=\eps \mathfrak{R}_{0,-1}^m} by Lemma \ref{l5.6}, we conclude that \fm{V^I(-\frac{1}{4}e_4^\alpha e_4^\beta e_3(g_{\alpha\beta})t(\chi_{ab}-r^{-1}\delta_{ab}))+\frac{1}{4}e_4^\alpha e_4^\beta e_3(g_{\alpha\beta})tV^I(\chi_{ab}-r^{-1}\delta_{ab})=O(\eps t^{-2+C\eps}\lra{q}^{-1}).}
Thus, by setting $F_{ab}=\chi_{ab}-r^{-1}\delta_{ab}-f_{ab}=\mathfrak{R}_{-2,0}^{m-1}$ and noting that $f_{ab}=\eps\mathfrak{R}_{-2,0}^m$, we have
\fm{V^I(V_4(F_{ab}))&=-2tr^{-1}V^I(F_{ab}+f_{ab})+O(\eps |V^I(F_{ab}+f_{ab})|)\\&\hspace{1em}+O(\eps t^{-2+C\eps}+t^{-3+C\eps}\lra{q}+t^{-1+C\eps}|V^I(F_{**}+f_{**})|)\\
&=-2tr^{-1}V^I(F_{ab})+O(\eps |V^I(F_{ab})|+\eps t^{-2+C\eps}+t^{-3+C\eps}\lra{q}+t^{-1+C\eps}|V^I(F_{**})|).}
In \eqref{p5.5f2}, \eqref{q1formula} and \eqref{q2formula}, we can replace $\xi_{k_1k_2}^{l_1}$ with $F_{ab}$. Thus, we have $V_4(V^I(F_{ab}))=Q_1+Q_2+V^I(V_4(F_{ab}))$, where by the induction hypotheses we have
\fm{Q_1&=nV_4(\ln(3R-r+t))V^I(F_{ab})+O(\eps \sum_{|J|=m\atop n_{J,3}=n}|V^J(F_{ab})|)+O(\lra{q}t^{-1+C\eps}\sum_{|J|=m\atop n_{J,3}<n}|V^J(F_{ab})|)\\
&=nV_4(\ln(3R-r+t))V^I(F_{ab})+O(\eps \sum_{|J|=m\atop n_{J,3}=n}|V^J(F_{ab})|)+O(\lra{q}t^{-3+C\eps}),}
\fm{|Q_2|&\lesssim \lra{q}t^{-1+C\eps}\sum_{0<|J|<m}|V^J(F_{ab})|\lesssim \lra{q}t^{-3+C\eps}.}
Thus,
\fm{&\hspace{1.5em}|e_4(V^I(F_{ab}))-ne_4(\ln(3R-r+t))V^I(F_{ab})+2r^{-1}V^I(F_{ab})|\\&\lesssim \eps t^{-1}\sum_{|J|=m\atop n_{J,3}=n}|V^J(F_{ab})|+t^{-2+C\eps}|V^I(F_{**})|+t^{-4+C\eps}\lra{q}+\eps t^{-3+C\eps}.}
By Lemma \ref{l5.8} with $n_0=2$, $n_1=n$ and Lemma \ref{derqxi}, we have
\fm{\sum_{a,b}\sum_{|I|=m\atop n_{I,3}=n}|V^I(F_{ab})|&\lesssim t^{-2+C\eps}(x^0(0)^{C\eps}+\int_{x^0(0)}^t \lra{q}\tau^{-2+C\eps}+\eps\tau^{-1+C\eps}\ d\tau)\\
&\lesssim t^{-2+C\eps}(x^0(0)^{C\eps}+ \lra{q}(x^0(0))^{-1+C\eps}+t^{C\eps})\lesssim t^{-2+C\eps}.}
Here we recall that $t\geq x^0(0)\sim T_0+\lra{q}$. We then finish the proof by induction.
\end{proof}\rm

Next, we consider $\xi_{12}^a$.
\prop{\label{prop5.5part2} Under our induction hypotheses \eqref{prop5.5hypo}, for  $|I|=m$, we have 
\fm{|V^I(\xi_{12}^a)|\lesssim t^{-1+C\eps}.}So $\xi_{12}^a=\mathfrak{R}_{-1,0}^m$.}
\begin{proof}
Fix $I$ such that $|I|=m$ and $n_{I,3}=n\leq m$. Recall from \eqref{prop5.5hypo} that $\xi_{12}^a=\mathfrak{R}_{-1,0}^{m-1}$. Suppose that $V^J(\xi_{12}^a)=O(t^{-1+C\eps})$ for $|J|=m$ and $n_{J,3}<n$.
By the equation in Section \ref{s5.1} we have
\eq{\label{prop5.5part2f1}V^I(V_4(\xi_{12}^a))&=V^I(t\Gamma_{\alpha\beta}^0e_4^\alpha e_2^\beta\chi_{a1}-t\Gamma^0_{\alpha\beta}e_4^\alpha e_1^\beta\chi_{a2})-V^I(t\chi_{ac}\xi_{12}^c)+V^I(t\lra{R(e_4,e_a)e_2,e_1}).}
By Lemma \ref{l5.7}, the last term is $O(\eps\lra{q}^{-1}t^{-1+C\eps})$. By Lemma \ref{l5.6} and Proposition \ref{prop5.5part1}, we note that \fm{t\Gamma_{\alpha\beta}^0e_4^\alpha e_2^\beta\chi_{a1}-t\Gamma^0_{\alpha\beta}e_4^\alpha e_1^\beta\chi_{a2}=\mathfrak{R}_{1,0}^{m+1}\cdot \eps\mathfrak{R}_{-1,-1}^{m+1}\cdot\mathfrak{R}_{0,0}^m\cdot \mathfrak{R}_{0,0}^m\cdot\mathfrak{R}_{-1,0}^m=\eps\mathfrak{R}_{-1,-1}^m.}
Thus, the first term in \eqref{prop5.5part2f1} is also $O(\eps \lra{q}^{-1}t^{-1+C\eps})$. Next, by the Leibniz's rule we have \fm{|V^I(t\chi_{ac}\xi_{12}^c)-t\chi_{ac}V^I(\xi_{12}^c)|&\lesssim\sum_{|J_1|+|J_2|=m\atop |J_1|>0}|V^{J_1}(t\chi_{ac})V^{J_2}(\xi_{12}^c)|\\
&\lesssim \sum_{|J_1|+|J_2|=m\atop |J_1|>0}(|V^{J_1}(t(\chi_{ac}-\delta_{ac}r^{-1}))V^{J_2}(\xi_{12}^c)|+|V^{J_1}(tr^{-1})V^{J_2}(\xi_{12}^a)|).}
By Proposition \ref{prop5.5part1} we have $t(\chi_{ac}-\delta_{ac}r^{-1})=\mathfrak{R}_{-1,0}^m$. Also recall that  $V(tr^{-1})=V((t-r)r^{-1})=\mathfrak{R}_{-1,1}^{m}$. Thus,
\fm{|V^I(t\chi_{ac}\xi_{12}^c)-tr^{-1}V^I(\xi_{12}^a)|&\lesssim |V^I(t\chi_{ac}\xi_{12}^c)-t\chi_{ac}V^I(\xi_{12}^c)|+|t(\chi_{ac}-r^{-1}\delta_{ac})V^I(\xi_{12}^c)|\\
&\lesssim t^{-2+C\eps}\lra{q}+t^{-1+C\eps}|V^I(\xi_{12}^*)|.}
In conclusion, we have
\fm{V^I(V_4(\xi_{12}^a))&=-tr^{-1}V^{I}(\xi_{12}^a)+O(t^{-1+C\eps}|V^I(\xi_{12}^*)|+t^{-2+C\eps}\lra{q}+\eps\lra{q}^{-1}t^{-1+C\eps}).}

Moreover, by \eqref{q1formula}, we have
\fm{|Q_1-nV_4(\ln(3R-r+t))V^I(\xi_{12}^a)|&\lesssim \eps \sum_{|J|=m\atop n_{J,3}=n}|V^J(\xi_{12}^{a})|+\lra{q}t^{-1+C\eps}\sum_{|J|=m\atop n_{J,3}<n}|V^J(\xi_{12}^{a})|\\
&\lesssim \eps \sum_{|J|=m\atop n_{J,3}=n}|V^J(\xi_{12}^{a})|+\lra{q}t^{-2+C\eps}.}
By \eqref{q2formula}, we have
\fm{|Q_2|&\lesssim t^{-1+C\eps}\lra{q}\sum_{0<|J|<m}|V^J(\xi_{12}^a)|\lesssim t^{-2+C\eps}\lra{q}.}
Thus,
\fm{&\hspace{1.5em}|e_4(V^I(\xi_{12}^a))+(-ne_4(\ln(3R-r+t))+r^{-1})V^I(\xi_{12}^a)|\\&\lesssim \eps t^{-1}\sum_{|J|=m\atop n_{J,3}=n}|V^J(\xi_{12}^{a})|+t^{-2+C\eps}|V^I(\xi_{12}^*)|+t^{-3+C\eps}\lra{q}+\eps \lra{q}^{-1}t^{-2+C\eps}.}
We now apply Lemma \ref{l5.8} with $n_0=1$, $n_1=n$ and Lemma \ref{derqxi}. Then,
\fm{\sum_a\sum_{|I|=m\atop n_{I,3}=n}|V^I(\xi_{12}^a)|
&\lesssim t^{-1+C\eps}(x^0(0)^{C\eps}+\int_{x^0(0)}^t\tau^{-2+C\eps}\lra{q}+\eps \lra{q}^{-1}\tau^{-1+C\eps}\ d\tau)\\
&\lesssim t^{-1+C\eps}(x^0(0)^{C\eps}+x^0(0)^{-1+C\eps}\lra{q}+ \lra{q}^{-1}t^{C\eps})\lesssim t^{-1+C\eps}.}
Again recall that $t\geq x^0(0)\sim \lra{q}+T_0$. We finish the proof by induction.
\end{proof}
\rm

Next we study $\xi_{34}^a$. The proof of the following proposition is very similar to that of the previous one.
\prop{\label{prop5.5part3} Under our induction hypotheses \eqref{prop5.5hypo}, for  $|I|=m$, we have 
\fm{|V^I(\xi_{34}^a)|\lesssim t^{-1+C\eps}\lra{q}^{-1}.} So $\xi_{34}^a=\mathfrak{R}_{-1,-1}^m$.}
\begin{proof}
Fix $I$ such that $|I|=m$ and $n_{I,3}=n\leq m$. Recall from \eqref{prop5.5hypo} that $\xi_{34}^a=\mathfrak{R}_{-1,-1}^{m-1}$. Suppose that $V^J(\xi_{34}^a)=O(t^{-1+C\eps}\lra{q}^{-1})$ for $|J|=m$ and $n_{J,3}<n$. By the equation in Section \ref{s5.1} we have
\fm{V^I(V_4(\xi_{34}^a))&=-V^I(t\chi_{ba}\xi_{34}^b)+V^I(t\lra{R(e_4,e_3)e_4,e_a})+2V^{I}(V_4(\Gamma_{\alpha\beta}^0 e_4^\alpha e_a^\beta)).}
By Lemma \ref{l5.7}, the second term is $O(\eps t^{-1+C\eps}\lra{q}^{-1})$. In the third term, we note that
\fm{V_4(\Gamma_{\alpha\beta}^0 e_4^\alpha e_a^\beta)&=V_4(\Gamma_{\alpha\beta}^0) e_4^\alpha e_a^\beta+\Gamma_{\alpha\beta}^0 V_4(e_4^\alpha) e_a^\beta+\Gamma_{\alpha\beta}^0 e_4^\alpha  V_4(e_a^\beta)\\
&=\eps\mathfrak{R}_{-1,-1}^m+\eps\mathfrak{R}_{-1,-1}^m \cdot \eps\mathfrak{R}_{-1,0} +\eps\mathfrak{R}_{-1,-1}^m   \cdot \eps\mathfrak{R}_{-1,0}=\eps\mathfrak{R}_{-1,-1}^m.}
We recall from Remark \ref{rmk5.6.1} that $e_4(e_*^*)=\eps\mathfrak{R}_{-2,0}^m$. Thus, $V^I(V_4(\Gamma_{\alpha\beta}^0))=O(\eps \lra{q}^{-1}t^{-1+C\eps})$. Following the computation in Proposition \ref{prop5.5part2}, we can prove that \fm{&\hspace{1.5em}|V^{I}(t\chi_{ba}\xi_{34}^b)-tr^{-1}V^I(\xi_{34}^a)|\lesssim |V^I(t\chi_{ab}\xi_{34}^b)-t\chi_{ab}V^I(\xi_{34}^b)|+|t(\chi_{ab}-r^{-1}\delta_{ab})V^I(\xi_{34}^b)|\\
&\lesssim \sum_{|J_1|+|J_2|=m\atop |J_1|>0}(|V^{J_1}(t(\chi_{ab}-\delta_{ab}r^{-1}))V^{J_2}(\xi_{34}^b)|+|V^{J_1}(tr^{-1})V^{J_2}(\xi_{34}^a)|)+t^{-1+C\eps}|V^I(\xi_{34}^b)|\\
&\lesssim t^{-2+C\eps}+t^{-1+C\eps}|V^I(\xi_{34}^*)|.}
Moreover, by \eqref{q1formula} we have
\fm{|Q_1-nV_4(\ln(3R-r+t))V^I(\xi_{34}^a)|&\lesssim \eps \sum_{|J|=m\atop n_{J,3}=n}|V^J(\xi_{34}^{a})|+\lra{q}t^{-1+C\eps}\sum_{|J|=m\atop n_{J,3}<n}|V^J(\xi_{34}^{a})|\\
&\lesssim \eps \sum_{|J|=m\atop n_{J,3}=n}|V^J(\xi_{34}^{a})|+t^{-2+C\eps}.}
By \eqref{q2formula}, we have
\fm{|Q_2|&\lesssim t^{-1+C\eps}\lra{q}\sum_{0<|J|<m}|V^J(\xi_{34}^a)|\lesssim t^{-2+C\eps}.}
Thus,
\fm{&\hspace{1.5em}|e_4(V^I(\xi_{34}^a))+(-ne_4(\ln(3R-r+t))+r^{-1})V^I(\xi_{34}^a)|\\&\lesssim \eps t^{-1}\sum_{|J|=m\atop n_{J,3}=n}|V^J(\xi_{34}^{a})|+t^{-2+C\eps}|V^I(\xi_{34}^*)|+t^{-3+C\eps}+\eps \lra{q}^{-1}t^{-2+C\eps}.}
We now apply Lemma \ref{l5.8} with $n_0=1$, $n_1=n$ and Lemma \ref{derqxi}. Then,
\fm{\sum_a\sum_{|I|=m\atop n_{I,3}=n}|V^I(\xi_{12}^a)|
&\lesssim t^{-1+C\eps}(x^0(0)^{C\eps}\lra{q}^{-1}+\int_{x^0(0)}^t\tau^{-2+C\eps}+\eps \lra{q}^{-1}\tau^{-1+C\eps}\ d\tau)\\
&\lesssim t^{-1+C\eps}(x^0(0)^{C\eps}\lra{q}^{-1}+x^0(0)^{-1+C\eps}+ \lra{q}^{-1}t^{C\eps})\lesssim t^{-1+C\eps}\lra{q}^{-1}.}
Again recall that $t\geq x^0(0)\sim \lra{q}+T_0$. We finish the proof by induction.
\end{proof}\rm

\bigskip

Finally, we consider $\xi_{a3}^l$. The case when $l\in\{a,3\}$ is easy.
\prop{\label{prop5.5part4} Under our induction hypotheses \eqref{prop5.5hypo}, for  $|I|=m$, we have 
\fm{\lra{q}|V^I(\xi_{a3}^3)|+|V^I(\xi_{a3}^a)|\lesssim t^{-1+C\eps}.}So $\xi_{a3}^3=\mathfrak{R}_{-1,-1}^m$ and $\xi_{a3}^a=\mathfrak{R}_{-1,0}^m$.}
\begin{proof}
Recall from Section \ref{s5.1} that 
\fm{\xi_{a3}^3&=-2\Gamma_{\alpha\beta}^0e_4^\alpha e_a^\beta+\frac{1}{2}\xi_{34}^a,\hspace{2em}
\xi_{a3}^a&=\chi_{aa}+2e_a(g^{0\alpha})g_{\alpha\beta}e_a^\beta+2g^{0\alpha}e_a^\beta\Gamma_{\beta\alpha}^\mu g_{\mu\nu}e_a^{\nu}.}
Now we apply Lemma \ref{l5.6}. Since $\Gamma=\eps\mathfrak{R}_{-1,-1}^{m+1}$ and $(g,e_*^*)=\mathfrak{R}_{0,0}^m$, we have $\Gamma_{\alpha\beta}^0e_4^\alpha e_a^\beta =\eps \mathfrak{R}_{-1,-1}^m$ and $g^{0\alpha}e_a^\beta\Gamma_{\beta\alpha}^\mu g_{\mu\nu}e_a^{\nu}=\eps \mathfrak{R}_{-1,-1}^m$. Since $e_4(g^{0\alpha})=t^{-1}V_4(g)=\eps\mathfrak{R}_{-2,0}^m$ and $e_a(g^{0\alpha})=r^{-1}V_a(g)=\eps\mathfrak{R}_{-2,0}^m$, we have $e_a(g^{0\alpha})g_{\alpha\beta}e_a^\beta=\eps\mathfrak{R}_{-2,0}^m$. We thus conclude that \fm{(\xi_{a3}^3,\xi_{a3}^a)=(\frac{1}{2}\xi_{34}^a,\chi_{aa})+\eps\mathfrak{R}_{-1,-1}^m.}
We finally apply Proposition \ref{prop5.5part1}, Proposition \ref{prop5.5part2} and Proposition \ref{prop5.5part3} to conclude that $\xi_{a3}^3=\mathfrak{R}_{-1,-1}^m$ and $\xi_{a3}^a=\mathfrak{R}_{-1,0}^m$.
\end{proof}\rm

The case $l=a'$ where $\{a,a'\}=\{1,2\}$ is harder.

\prop{\label{prop5.5part5} Under our induction hypotheses \eqref{prop5.5hypo}, for $|I|=m$, we have 
\fm{|V^I(\xi_{a3}^{a'})|\lesssim \lra{q}^{-1}t^{C\eps}.}So $\xi_{a3}^{a'}=\mathfrak{R}_{0,-1}^m$.}
\begin{proof}Fix $I$ such that $|I|=m$ and $n_{I,3}=n\leq m$. Recall from \eqref{prop5.5hypo} that $\xi_{a3}^{a'}=\mathfrak{R}_{0,-1}^{m-1}$. Suppose that $V^J(\xi_{a3}^{a'})=O(\lra{q}^{-1}t^{C\eps})$ for $|J|=m$ and $n_{J,3}<n$. By the equation in Section \ref{s5.1} we have
\fm{V^I(V_4(\xi_{a3}^{a'}))&=V^I((V_4+t\Gamma_{\mu\nu}^0e_4^\mu e_4^\nu)(\chi_{aa'}+2e_a(g^{0\alpha})g_{\alpha\beta}e_{a'}^\beta+2g^{0\alpha}e_a^\beta\Gamma_{\beta\alpha}^\mu g_{\mu\nu}e_{a'}^{\nu}))-V^I(t\Gamma_{\mu\nu}^0e_4^\mu e_4^\nu\xi_{a3}^{a'})\\
&\hspace{1em}-\sum_cV^I(t\xi_{34}^c\xi_{aa'}^c)-V^I(t\lra{R(e_4,e_3)e_a,e_{a'}})-V^I(t\Gamma_{\alpha\beta}^0e_4^\alpha e_a^\beta\xi_{34}^{a'}+t\Gamma_{\alpha\beta}^0e_4^\alpha e_{a'}^\beta\xi_{34}^a).}
By the Leibniz's rule and all the previous results, we conclude that the second line has an upper bound
\fm{t^{-1+C\eps}\lra{q}^{-1}+\eps \lra{q}^{-1}t^{-1+C\eps}\lesssim t^{-1+C\eps}\lra{q}^{-1}.}
In the first line, we note that 
\fm{t\Gamma_{\mu\nu}^0e_4^\mu e_4^\nu(2e_a(g^{0\alpha})g_{\alpha\beta}e_{a'}^\beta+2g^{0\alpha}e_a^\beta\Gamma_{\beta\alpha}^\mu g_{\mu\nu}e_{a'}^{\nu})=\eps\mathfrak{R}_{0,-1}^m\cdot(\eps\mathfrak{R}_{-2,0}^m+\eps\mathfrak{R}_{-1,-1}^m)=\eps^2\mathfrak{R}_{-1,-2}^m.}
Besides, since $\chi_{aa'}=\mathfrak{R}_{-2,0}^m$ and since $\sum_c\chi_{ac}\chi_{ca'}=\chi_{12}\tr\chi$, we have
\fm{&\hspace{1.5em}|V^I(V_4(\chi_{aa'})+t\Gamma_{\mu\nu}^0e_4^\mu e_4^\nu\chi_{aa'})|\\&\lesssim|V^I(2t\Gamma_{\alpha\beta}^0e_4^\alpha e_4^\beta\chi_{aa'})|+|V^I(t\chi_{12}(\chi_{11}+\chi_{22}))|+|V^I(t\lra{R(e_4,e_a)e_4,e_{a'}})|\\
&\lesssim |V^I(\eps\mathfrak{R}_{-2,-1}^m)|+|V^I(\mathfrak{R}_{1,0}^{m+1}\cdot \mathfrak{R}_{-2,0}^{m}\cdot \mathfrak{R}_{-1,0}^{m})|+|V^I(\eps\mathfrak{R}_{-1,-1}^m)|\lesssim t^{-2+C\eps}+\eps t^{-1+C\eps}\lra{q}^{-1}\lesssim t^{-1+C\eps}\lra{q}^{-1}.}
Moreover, recall that $V_4(e_*^*)=\eps\mathfrak{R}_{-1,0}^m$. We also have $\partial g=\eps\mathfrak{R}_{-1,-1}^{m+1}$ by Remark \ref{rmk5.6.1}. Thus, we have
\fm{&\hspace{1.5em}V_{4}(2e_a(g^{0\alpha})g_{\alpha\beta}e_{a'}^\beta+2g^{0\alpha}e_a^\beta\Gamma_{\beta\alpha}^\mu g_{\mu\nu}e_{a'}^{\nu})=2V_4(e_a(g^{0\alpha}))g_{\alpha\beta}e_{a'}^\beta+\eps \mathfrak{R}_{-1,-1}^m\\
&=2e_a^\sigma V_4(\partial_\sigma g^{0\alpha})g_{\alpha\beta}e_{a'}^\beta+2V_4(e_a^\sigma)(\partial_\sigma g^{0\alpha})g_{\alpha\beta}e_{a'}^\beta+\eps \mathfrak{R}_{-1,-1}^m=\eps \mathfrak{R}_{-1,-1}^m.}
In conclusion,
\fm{|V^I(V_4(\xi_{a3}^{a'})|&\lesssim |V^I(t\Gamma^{0}_{\mu\nu} e_4^\mu e_4^\nu \xi_{a3}^{a'})|+t^{-1+C\eps}\lra{q}^{-1}\\
&\lesssim |t\Gamma^{0}_{\mu\nu} e_4^\mu e_4^\nu V^I(\xi_{a3}^{a'})|+\sum_{|J_1|+|J_2|=m\atop |J_2|<m}|V^{J_1}(\eps\mathfrak{R}_{0,-1}^m) V^{J_2}(\xi_{a3}^{a'})|+t^{-1+C\eps}\lra{q}^{-1}\\
&\lesssim \eps|V^I(\xi_{a3}^{a'})|+\eps\lra{q}^{-2}t^{C\eps}+t^{-1+C\eps}\lra{q}^{-1}.}

Next, by \eqref{q1formula}, we have
\fm{|Q_1-nV_4(\ln(3R-r+t))V^I(\xi_{a3}^{a'})|&\lesssim \eps\sum_{|J|=m\atop n_{J,3}=n}|V^J(\xi_{a3}^{a'})|+\lra{q}t^{-1+C\eps}\sum_{|J|=m\atop n_{J,3}<n}|V^J(\xi_{a3}^{a'})|\\
&\lesssim \eps\sum_{|J|=m\atop n_{J,3}=n}|V^J(\xi_{a3}^{a'})|+t^{-1+C\eps}.}
By \eqref{q2formula}, we have
\fm{|Q_2|&\lesssim t^{-1+C\eps}\lra{q}\sum_{0<|J|<m}|V^J(\xi_{a3}^{a'})|\lesssim t^{-1+C\eps}.}
Thus,
\fm{|e_4(V^I(\xi_{a3}^{a'}))-ne_4(\ln(3R-r+t))V^I(\xi_{a3}^{a'})| \lesssim\eps t^{-1}\sum_{|J|=m\atop n_{J,3}=n}|V^J(\xi_{a3}^{a'})|+\eps\lra{q}^{-2}t^{-1+C\eps}+t^{-2+C\eps}.}
By Lemma \ref{l5.8} with $n_0=0$, $n_1=n$ and Lemma \ref{derqxi}, we have
\fm{\sum_{a,a'}\sum_{|I|=m\atop n_{I,3}=n}|V^I(\xi_{a3}^{a'})|
&\lesssim t^{C\eps}(\lra{q}^{-1}x^0(0)^{C\eps}+\int_{x^0(0)}^t\eps\lra{q}^{-2}\tau^{-1+C\eps}+\tau^{-2+C\eps}\ d\tau)\\
&\lesssim t^{C\eps}(\lra{q}^{-1}t^{C\eps}+\lra{q}^{-2}t^{C\eps}+(x^0(0))^{-1+C\eps})\lesssim \lra{q}^{-1}t^{C\eps}.}
We finish the proof by induction.
\end{proof}\rm

\bigskip

Combining Proposition \ref{prop5.5part1}-\ref{prop5.5part5}, we finish the proof of Proposition \ref{prop5.5} by induction.

\subsection{Estimates for higher derivatives of $q$}\label{sec5.4}
Now we can prove the estimates for higher derivatives of $q$. We first note that \eqref{derqetaformula} holds for each $m\geq 1$, as long as $\eps\ll_m1$. This is because \eqref{derqetaformula} is a result of \eqref{prop5.5hypo} which then results from  Proposition \ref{prop5.5}.

\lem{\label{l4.18}In $\Omega\cap\{r-t<2R\}$,  we have $V^Iq=O(\lra{q}t^{C\eps})$ for each multiindex $I$. }
\begin{proof}
We induct on $|I|$. If $|I|=0$, there is nothing to prove. If $|I|=1$, the estimates are clear since $V_1(q)=V_2(q)=V_4(q)=0$ and $V_3(q)=O((3R-r+t)|\partial q|)=O(\lra{q}t^{C\eps})$.

In general, we fix an integer $m>1$. By choosing $\eps\ll_m1$, we can assume that Proposition \ref{prop5.5} holds for all $|I|\leq m$. Suppose we have proved the estimates for $|I|<m$, so $q=\mathfrak{R}_{0,1}^{m-1}$. Fix  a multiindex $I$ such that $|I|=m$. If $n_{I,4}>0$, we can write $I=(J',4,J)$. Here we can assume $|J|>0$ since otherwise we have $V^I(q)=V_{J'}(V_4(q))=0$. By \eqref{derqetaformula}, we have
\fm{V^I(q)&=V^{J'}(V_4(V^J(q))=\sum_{J=(J_1,k,J_2)}V^{(J',J_1)}[V_4,V_k]V^{J_2}(q)\\&=\sum_{J=(J_1,k,J_2)}V^{(J',J_1)}(\eta_{k}^lV^{(l,J_2)}(q))=\sum_{J=(J_1,k,J_2)}V^{(J',J_1)}(\mathfrak{R}_{-1,1}^{m-1}\cdot\mathfrak{R}_{0,1}^{m-1-(1+|J_2|)})\\&=\sum_{J=(J_1,k,J_2)}V^{(J',J_1)}(\mathfrak{R}_{-1,2}^{m-1-(1+|J_2|)})=O(\lra{q}^2t^{-1+C\eps})=O(\lra{q}t^{C\eps}).}
Here we  note that $|J_2|+1=|J|-|J_1|= m-1-|J'|-|J_1|$, so we are able to apply the definition of $\mathfrak{R}_{**}^*$ here.

Next suppose $n_{I,3}<m$ and $n_{I,4}=0$. Thus we can write $I=(J',a,J)$ where $n_{J,3}=|J|$. Here we can assume $|J|>0$  since $V_a(q)=0$.  Then
\fm{V^I(q)&=V^{J'}V_a(V^J(q))=\sum_{J=(J_1,3,J_2)}V^{(J',J_1)}[V_a,V_3]V^{J_2}(q).}
Note that
\fm{\ [V_a,V_3]F&=V_a((3R-r+t)e_3(F))-V_3(re_a(F))\\
&=V_a(3R-r+t)e_3(F)-V_3(r)e_a(F)+(3R-r+t)r[e_a,e_3](F)\\
&=-(3R-r+t)^{-1}V_a(r)V_3(F)-r^{-1}V_3(r)V_a(F)\\&\hspace{1em}+(3R-r+t)\xi_{a3}^bV_b(F)+r\xi_{a3}^3V_3(F)+(3R-r+t)rt^{-1}\xi_{a3}^4V_4(F).}
By Lemma \ref{l5.6} and Remark \ref{rmk5.6.1}, we have $V_a(r)=\mathfrak{R}_{0,0}^m$, $V_3(r)=(3R-r+t)e_3^i\omega_i=\mathfrak{R}_{0,1}^m$. By Proposition \ref{prop5.5}, we have
\fm{\ [V_a,V_3]&=\sum_{k=1}^4\mathfrak{R}_{0,0}^m\cdot V_k=\mathfrak{R}_{0,0}^m\cdot V.}
Thus,
\fm{V^I(q)&=\sum_{J=(J_1,3,J_2)}V^{(J',J_1)}(\mathfrak{R}_{0,0}^m\cdot V(V^{J_2}(q)))\\&=\sum_{J=(J_1,3,J_2)}V^{(J',J_1)}(\mathfrak{R}_{0,0}^m\cdot \mathfrak{R}_{0,1}^{m-1-(1+|J_2|)})=O(t^{C\eps}\lra{q}).}
Again, we have $m-1=1+|J_2|+|J_1|+|J'|$.

Finally, suppose $n_{I,3}=|I|$. We have
\fm{V_4(V^I(q))=\sum_{I=(J_1,3,J_2)\atop n_{J_1,3}=|J_1|,\ n_{J_2,3}=|J_2|}V^{J_1}[V_4,V_3]V^{J_2}(q)=\sum_{I=(J_1,3,J_2)\atop n_{J_1,3}=|J_1|,\ n_{J_2,3}=|J_2|}V^{J_1}(\eta_{3}^lV^{(l,J_2)}(q)).}
By the Leibniz's rule, we can express $V^{J_1}(\eta_{3}^lV^{(l,J_2)}(q))$ as a linear combination of terms of the form $V^{K_1}(\eta_{3}^l)V^{K_2}(q)$, where $|K_1|+|K_2|=m$, $K_2$ contains $l$, and $(K_1,K_2)$ is an rearrangement of $(J_1,l,J_2)$. Now recall from \eqref{derqetaformula} that $\eta_{3}^l=\mathfrak{R}_{-1,1}^{m-1}+\eps\mathfrak{R}_{0,-1}^m$. Since $V^J(q)=O(\lra{q}t^{C\eps})$ for $|J|=m$ and $n_{J,3}<|J|$, we have
\fm{&\hspace{1.5em}V^{J_1}(\eta_3^lV^{(l,J_2)}(q))\\
&=\eta_3^3V^{I}(q)+O(\sum_{|K_1|+|K_2|=m,\ 0<|K_1|<m\atop n_{|K_1|,3}=|K_1|,\ n_{|K_2|,3}=|K_2|}|V^{K_1}(\eta_3^3)V^{K_2}(q)|)\\
&\hspace{1em}+O(\sum_{l\neq 3}\sum_{|K_1|+|K_2|=m,\ |K_2|>0\atop n_{K_1,3}=|K_1|,\ n_{K_2,3}=|K_2|-1}|V^{K_1}(\eta_{3}^l)V^{K_2}(q)|)\\
&=(te_4(\ln(3R-r+t))+O(\eps))V^{I}(q)+O(\sum_{0<|K_1|<m}|V^{K_1}(\eps\mathfrak{R}_{0,-1}^m+\mathfrak{R}_{-1,1}^{m-1})|\cdot t^{C\eps}\lra{q})\\
&\hspace{1em}+O(\sum_{|K_1|<m}|V^{K_1}(\eps\mathfrak{R}_{0,-1}^m+\mathfrak{R}_{-1,1}^{m-1})\cdot\lra{q}t^{C\eps}|\\
&=te_4(\ln(3R-r+t))V^I(q)+O(\eps|V^I(q)|)+O(\eps t^{C\eps}+t^{-1+C\eps}\lra{q}^2).}
Thus,
\fm{|e_4(V^I(q))-me_4(\ln(3R-r+t))V^I(q)|&\lesssim \eps t^{-1}|V^I(q)|+\eps t^{-1+C\eps}+t^{-2+C\eps}\lra{q}^2.}
Recall from Remark \ref{derqvzrmk1} that $V^I(q)=O(t^{C\eps}\lra{q})$ on $H$. Then, by Lemma \ref{l5.8} with $n_0=0$ and $n_1=|I|$, we have
\fm{|V^I(q)|&\lesssim t^{C\eps}(\lra{q}x^0(0)^{C\eps}+\int_{x^0(0)}^t \eps\tau^{-1+C\eps}+\tau^{-2+C\eps}\lra{q}^2\ d\tau)\\
&\lesssim t^{C\eps}(\lra{q}t^{C\eps}+ t^{C\eps}+(x^0(0))^{-1+C\eps}\lra{q}^2)\lesssim\lra{q}t^{C\eps}.}
\end{proof}\rm

We have the following important corollary.
\cor{\label{cor4.21}The function  $q(t,x)$ is a smooth function (in the sense defined in Section \ref{sec2.5}) in $\Omega$. Moreover, we have $Z^Iq=O(\lra{q}t^{C\eps})$ and $Z^I\Omega_{ij}q=O(t^{C\eps})$ for each multiindex $I$ and $1\leq i<j\leq 3$.}
\begin{proof} Fix an integer $m>1$. We seek to prove that for $\eps\ll_m1$, $q$ is a $C^m$ function and $Z^Iq=O(\lra{q}t^{C\eps})$ for  $|I|\leq m$. By writing $Z=z^\nu(t,x)\partial_\nu$, we have
\fm{r^{-1}\lra{Z,e_a}&=r^{-1}z^\alpha e_a^\beta g_{\alpha\beta}=\mathfrak{R}_{0,0}^{m},\hspace{2em}
t^{-1}\lra{Z,e_3}=t^{-1}z^\alpha e_3^\beta g_{\alpha\beta}=\mathfrak{R}_{0,0}^{m}.}
Moreover, 
\fm{\lra{Z,e_4}&=z^\alpha e_4^\beta g_{\alpha\beta}=z^\alpha e_4^\beta (g_{\alpha\beta}-m_{\alpha\beta})+z^\alpha e_4^\beta m_{\alpha\beta}\\
&=\eps\mathfrak{R}_{0,0}^m-z^0+z^i(e_4^i-\omega_i)+z^i\omega_i=\mathfrak{R}_{0,1}^m+Z(r-t).}
We can easily check that $Z(r-t)=\mathfrak{R}_{0,1}^m$, so $(3R-r+t)^{-1}\lra{Z,e_4}=\mathfrak{R}_{0,0}^m$. Then, by \eqref{derqzv}, $Z=\mathfrak{R}_{0,0}^m\cdot V$, so $Z^Iq$ is a linear combination of terms of the form
\fm{Z^{I_1}(\mathfrak{R}_{0,0}^m)\cdots Z^{I_s}(\mathfrak{R}_{0,0}^m)V^s(q),\hspace{2em}\sum|I_*|+s=|I|,\ s>0.}
Each of such terms is $O(t^{C\eps}\lra{q})$ if $|I|\leq m$, so we have $Z^Iq=O(t^{C\eps}\lra{q})$ for $|I|\leq m$.

Moreover, for each $m>1$, as long as $\eps\ll_m1$, we have $q=\mathfrak{R}_{0,1}^{m+1}$ by Lemma \ref{l4.18}. Then we have
\fm{\Omega_{ij}q&=\frac{1}{2}\lra{\Omega_{ij},e_4}e_3(q)=\frac{1}{2}(x_ig_{j\beta}-x_jg_{i\beta})e_4^\beta e_3(q)\\
&=\frac{1}{2}(x_i m_{jk}-x_j m_{ik})\omega_k e_3(q)+\frac{1}{2}(x_i(g_{jk}-m_{jk})-x_j(g_{ik}-m_{ik}))\omega_k e_3(q)\\&\hspace{1em}+\frac{1}{2} (x_ig_{jk}-x_jg_{ik})(e_4^k-\omega_k) e_3(q)\\
&=0+\eps\mathfrak{R}_{0,0}^m+\mathfrak{R}_{0,0}^m=\mathfrak{R}_{0,0}^m.}
Again, for each multiindex $I$ with $|I|\leq m$, we can write $Z^I\Omega_{ij}q$ as a linear combination of terms of the form
\fm{Z^{I_1}(\mathfrak{R}_{0,0}^m)\cdots Z^{I_s}(\mathfrak{R}_{0,0}^m)V^s\Omega_{ij}(q),\hspace{2em}\sum|I_*|+s=m,\ s>0.}
Each of such terms is $O(t^{C\eps})$, so we have $Z^I\Omega_{ij}q=O(t^{C\eps})$ for $|I|\leq m$.
\end{proof}
\rm

\subsection{More estimates}\label{sec5.5}
We end this section with some estimates derived from our original wave equation \eqref{qwe}. We first introduce a new definition.
\defn{\rm Let $F=F(t,x)$ be a function with domain $\Omega\cap\{r-t<2R\}$.  For any integer $m\geq 0$ and any real numbers $s,p$, we have defined $F=\mathfrak{R}_{s,p}^m$ in Section \ref{sec4.3new} prior to Lemma \ref{l5.6}. We now define $F=\mathfrak{R}_{s,p}$, if $F=\mathfrak{R}_{s,p}^m$ for each $m\geq 0$.

Again, by the Leibniz's rule, we have $V^I(\mathfrak{R}_{s,p})=\mathfrak{R}_{s,p}$ and  $\mathfrak{R}_{s_1,p_1}\cdot\mathfrak{R}_{s_2,p_2}=\mathfrak{R}_{s_1+s_2,p_1+p_2}$. In addition, by Proposition \ref{prop5.5}, we have
\fm{\begin{array}{c}(\xi_{13}^2,\xi_{23}^1)=\mathfrak{R}_{0,-1};\ \xi_{34}^a=\mathfrak{R}_{-1,-1};\ \xi_{k_1k_2}^a=\mathfrak{R}_{-1,0} \text{ for all other }k_1<k_2\text{ and }a=1,2;\\[.5em]
\xi_{k_1k_2}^3=\mathfrak{R}_{-1,-1}\text{ for all }k_1<k_2;\ \chi_{ab}-r^{-1}\delta_{ab}=\mathfrak{R}_{-2,0}.\end{array}} There are many other estimates in Section \ref{sec4.3new} invovling $\mathfrak{R}_{*,*}^*$. They would still hold if all the superscripts are removed, because they all rely on Proposition \ref{prop5.5}. For example, by Lemma \ref{l5.6} we have
\fm{\begin{array}{c}e_*^*=\mathfrak{R}_{0,0},\ (e_4^i-\omega_i,e_3^i-\omega_i)=\mathfrak{R}_{-1,0};\ \partial^sZ^I(g-m)=\eps\mathfrak{R}_{-1,-s},\ \Gamma_{**}^*=\eps\mathfrak{R}_{-1,-1};\\[.5em]
\omega=\mathfrak{R}_{0,0},\ (t^s,r^s)=\mathfrak{R}_{s,0},\ (3R-r+t)^s=\mathfrak{R}_{0,s}.\end{array}}

We remark that this definition follows the spirits of the convention in Section \ref{sec2.5}. In the defintion of $\mathfrak{R}_{s,p}^m$, we require some estimates to hold for all $\eps\ll_{s,p,m}1$. The dependence  on $m$ here should be emphasized.
}\rm

\bigskip

Our  goal in this subsection is to prove that  \eq{\label{secend}e_4(e_3(u))+r^{-1}e_3(u)=\eps\mathfrak{R}_{-3,0},\hspace{2em}e_4(e_3(u))=\eps\mathfrak{R}_{-2,0};}
\eq{\label{secend2}e_4(e_3(q))=-\frac{1}{4}e_3(u)G(\omega)e_3(q)+\eps\mathfrak{R}_{-2,0}.}
We start our proof with the following lemma.

\lem{\label{l414}We have the following estimates.
\begin{enumerate}[\rm (a)]
\item $q_\alpha=\mathfrak{R}_{0,0}$, $q_r^{-1}=\mathfrak{R}_{0,0}$; $e_k(q_r)=\mathfrak{R}_{-1,-1}$, $e_k(q_r^{-1})=\mathfrak{R}_{-1,-1}$ for $k\neq 3$. 
\item $q_i+\omega_iq_t=\mathfrak{R}_{-1,0}$, $u_i+\omega_iu_t=\eps\mathfrak{R}_{-2,0}$. 
\item $e_k(q_i+\omega_iq_t)=\mathfrak{R}_{-2,0}$, $e_k(u_i+\omega_iu_t)=\eps\mathfrak{R}_{-3,0}$, for $k\neq 3$.
\item In {\rm (b)} and {\rm (c)} we can replace $q_i+\omega_iq_t$ with $q_t+q_r$ or $q_i-\omega_iq_r$, and replace $u_i+\omega_iu_t$ with $u_t+u_r$ or $u_i-\omega_iu_r$. The results are the same.
\end{enumerate}
\begin{proof}
(a) By Lemma \ref{l4.18}, we have  $V_3(q)=\mathfrak{R}_{0,1}$ and $e_3(q)=V_3(q)=\mathfrak{R}_{0,0}$. Then,
\fm{q_\alpha&=\frac{1}{2}g_{\alpha\beta}e_4^\beta e_3(q)=\mathfrak{R}_{0,0}\cdot \mathfrak{R}_{0,0}\cdot \mathfrak{R}_{0,0}=\mathfrak{R}_{0,0}.}
Since $\omega_i=\mathfrak{R}_{0,0}$, we have $q_r=\mathfrak{R}_{0,0}$. Since $q_r\geq C^{-1}t^{-C\eps}$ and since $V^I(q_r^{-1})$ is a linear combination of terms of the form
\eq{\label{l414f1}q_r^{-s-1}V^{I_1}(q_r)\cdots V^{I_s}(q_r),\hspace{1em}\text{where }\sum|I_j|=|I|,\ |I_j|>0,}
we conclude that $V^I(q_r^{-1})=O(t^{C\eps})$ for each $I$ and thus $q_r^{-1}=\mathfrak{R}_{0,0}$. Besides,  we have \fm{e_k(e_3(q))&=[e_k,e_3]q=\xi_{k3}^3e_3(q),\hspace{2em}k=1,2,3,4;\\
2\omega_ig_{i\beta}e_4^\beta&=\lra{e_3+e_4,e_4}+(2\omega_i-e_4^i-e_3^i)g_{i\beta}e_4^\beta=2+\mathfrak{R}_{-1,0}.} Thus, for $k\neq 3$, \fm{e_k(q_r)&=e_k(\frac{1}{2}\omega_ig_{i\beta}e_4^\beta e_3(q))=e_k(\frac{1}{2}\omega_ig_{i\beta}e_4^\beta) e_3(q)+\frac{1}{2}\omega_ig_{i\beta}e_4^\beta e_k(e_3(q))\\
&=e_k(\frac{1}{2}+\mathfrak{R}_{-1,0}) e_3(q)+(\frac{1}{2}+\mathfrak{R}_{-1,0}) \xi_{k3}^3e_3(q)\\
&=\mathfrak{R}_{-1,0}\cdot V_k(\mathfrak{R}_{-1,0})\cdot \mathfrak{R}_{0,0}+\mathfrak{R}_{-1,-1}=\mathfrak{R}_{-1,-1}.}
Now if we expand $V^I(e_k(q_r^{-1}))$, each term is still of the form \eqref{l414f1} with $s>0$ and $V^{I_s}(q_r)$ replaced by $V^{I_s}(e_k(q_r))$. We thus conclude that $e_k(q_r^{-1})=\mathfrak{R}_{-1,-1}$ for $k\neq 3$.

(b) We have \fm{q_i+\omega_iq_t&=\frac{1}{2}(g_{i\beta}+\omega_ig_{0\beta})e_4^\beta e_3(q)}and
\fm{u_i+\omega_iu_t&=\frac{1}{2}(g_{i\beta}+\omega_ig_{0\beta})e_4^\beta e_3(u)+\frac{1}{2}(g_{i\beta}+\omega_ig_{0\beta})e_3^\beta e_4(u)+\sum_a(g_{i\beta}+\omega_ig_{0\beta})e_a^\beta e_a(u)\\
&=\frac{1}{2}(g_{i\beta}+\omega_ig_{0\beta})e_4^\beta (3R-r+t)^{-1}V_3(u)+\eps\mathfrak{R}_{-2,0}.}
Here we have \fm{(g_{i\beta}+\omega_ig_{0\beta})e_4^\beta=e_4^i-\omega_i+((g_{i\beta}-m_{i\beta})+\omega_i(g_{0\beta}-m_{0\beta}))e_4^\beta=\mathfrak{R}_{-1,0}.}
We thus conclude that $q_i+\omega_iq_t=\mathfrak{R}_{-1,0}$ and $u_i+\omega_iu_t=\eps\mathfrak{R}_{-2,0}$.

(c) Recall that $e_a(r)=\mathfrak{R}_{-1,0}$, $e_4(\omega_i)=r^{-1}(e_4^i-\omega_i+(1-e_4^j\omega_j)\omega_i)=\mathfrak{R}_{-2,0}$ and $e_4(e_k^\alpha)=\eps\mathfrak{R}_{-2,0}$ by Lemma \ref{l5.6} and Lemma \ref{derqle4r}. Besides, note that
\fm{e_a(\omega_i)&=r^{-1}(e_a^i-e_a(r)\omega_i)=r^{-1}e_a^i+\mathfrak{R}_{-2,0},\\
e_4(\omega_i)&=(e_4^j-\omega_j)\partial_j\omega_i=r^{-1}(e_4^i-\omega_i-(e_4^j-\omega_j)\omega_j\omega_i)=\mathfrak{R}_{-2,0}.}
Thus we have
\fm{&\hspace{1.5em}e_a((g_{i\beta}+\omega_ig_{0\beta})e_4^\beta)=e_a(g_{i\beta}+\omega_ig_{0\beta})e_4^\beta+(g_{i\beta}+\omega_ig_{0\beta})e_a(e_4^\beta)\\
&=(e_a(g_{i\beta})+\omega_ie_a(g_{0\beta})+e_a(\omega_i)g_{0\beta})e_4^\beta+(g_{i\beta}+\omega_ig_{0\beta})(\xi_{a4}^le_l^\beta+e_4(e_a^\beta))\\
&=(\eps\mathfrak{R}_{-2,0}+(r^{-1}e_a^i+\mathfrak{R}_{-2,0})g_{0\beta})e_4^\beta+(g_{i\beta}+\omega_ig_{0\beta})(\xi_{a4}^be_b^\beta+\eps\mathfrak{R}_{-2,0})
\\
&=r^{-1}e_a^ig_{0\beta}e_4^\beta+r^{-1}(g_{i\beta}+\omega_ig_{0\beta})e_a^\beta+(g_{i\beta}+\omega_ig_{0\beta})(\chi_{ab}-\delta_{ab}r^{-1})e_b^\beta+\mathfrak{R}_{-2,0}\\
&=r^{-1}(-e_a^i+e_a^i(g_{0\beta}-m_{0\beta})e_4^\beta+e_a^i+((g_{i\beta}-m_{i\beta})+\omega_i(g_{0\beta}-m_{0\beta}))e_a^\beta)+\mathfrak{R}_{-2,0}=\mathfrak{R}_{-2,0},}
and
\fm{e_4((g_{i\beta}+\omega_ig_{0\beta})e_4^\beta)&=e_4(g_{i\beta}+\omega_ig_{0\beta})e_4^\beta+(g_{i\beta}+\omega_ig_{0\beta})e_4(e_4^\beta)\\
&=(e_4(g_{i\beta})+\omega_ie_4(g_{0\beta})+e_4(\omega_i)g_{0\beta})e_4^\beta+\eps\mathfrak{R}_{-2,0}\\
&=\mathfrak{R}_{-2,0}+\eps\mathfrak{R}_{-2,0}=\mathfrak{R}_{-2,0}.}
Since $(g_{i\beta}+\omega_ig_{0\beta})e_4^\beta=\mathfrak{R}_{-1,0}$ and $e_k(e_3(q))=\xi_{k3}^3e_3(q)=\mathfrak{R}_{-1,-1}$, we conclude from the Leibniz's rule that for $k\neq 3$,
\fm{e_k(q_i+\omega_iq_t)&=\frac{1}{2}e_k((g_{i\beta}+\omega_ig_{0\beta})e_4^\beta)e_3(q)+\frac{1}{2}(g_{i\beta}+\omega_ig_{0\beta})e_4^\beta e_k(e_3(q))\\
&=\mathfrak{R}_{-2,0}\cdot \mathfrak{R}_{0,0}+\mathfrak{R}_{-1,0}\cdot\mathfrak{R}_{-1,-1}=\mathfrak{R}_{-2,0}.} Besides,
\fm{u_i+\omega_iu_t&=r^{-1}\sum_j\omega_j\Omega_{ji}u+r^{-1}\omega_iSu+r^{-1}\omega_i(t+r)^{-1}(tSu-\sum_jx_j\Omega_{0j}u)=\mathfrak{R}_{-1,0}\cdot Zu.}
Note that $Zu=\eps\mathfrak{R}_{-1,0}$ and $e_k=\mathfrak{R}_{-1,0}\cdot V$ for $k\neq 3$. We conclude that \fm{e_k(u_i+\omega_iu_t)&=e_k(\mathfrak{R}_{-1,0})\cdot Zu+\mathfrak{R}_{-1,0}\cdot e_k(Zu)\\&=\mathfrak{R}_{-1,0}\cdot V_k(\mathfrak{R}_{-1,0})\cdot \eps\mathfrak{R}_{-1,0}+\mathfrak{R}_{-1,0}\cdot \mathfrak{R}_{-1,0}\cdot  V_k(\eps\mathfrak{R}_{-1,0})=\eps\mathfrak{R}_{-3,0}.}

(d) This part follows directly from \fm{\partial_t+\partial_r&=\sum\omega_i(\partial_i+\omega_i\partial_t),\hspace{2em}\partial_i-\omega_i\partial_r=\partial_i+\omega_i\partial_t-\sum\omega_i\omega_j(\partial_j+\omega_j\partial_t).}
\end{proof}
\rm

\prop{\label{l413}We have $e_4(e_3(u))+r^{-1}e_3(u)=\eps\mathfrak{R}_{-3,0}$ and $e_4(e_3(ru))=\eps\mathfrak{R}_{-2,0}$. }\rm
\begin{proof}
Note that
\fm{g^{\alpha\beta}(u)\partial_\alpha\partial_\beta u&=\sum_a e_a^\alpha e_a^\beta \partial_\alpha\partial_\beta u+\frac{1}{2}e_4^\alpha e_3^\beta \partial_\alpha\partial_\beta u+\frac{1}{2}e_3^\alpha e_4^\beta \partial_\alpha\partial_\beta u\\
&=\sum_a (e_a(e_a(u))-e_a(e_a^\alpha)\partial_\alpha u)+e_4(e_3(u))-e_4(e_3^\alpha)\partial_\alpha u.}
Here we have\fm{&\hspace{1.5em}e_a(e_a^\alpha)\partial_\alpha u\\&=-\xi_{aa'}^ae_{a'}(u)-\frac{1}{2}\chi_{aa}(e_3(u)+e_4(u))-\lra{e_a,e_a(g^{0\beta})\partial_\beta+g^{0\beta}e_a^\alpha\Gamma_{\alpha\beta}^\nu\partial_\nu}e_4(u)-e_a^\mu e_a^\nu\Gamma_{\mu\nu}^\alpha u_\alpha\\
&=-\xi_{aa'}^ae_{a'}(u)-\frac{1}{2}\chi_{aa}(e_3(u)+e_4(u))-(e_a^\alpha g_{\alpha\beta}e_a(g^{0\beta})+e_a^\mu g_{\mu\nu}g^{0\beta}e_a^\alpha\Gamma_{\alpha\beta}^\nu)e_4(u)-e_a^\mu e_a^\nu\Gamma_{\mu\nu}^\alpha u_\alpha\\
&=-\frac{1}{2}\chi_{aa}e_3(u)-e_a^\mu e_a^\nu\Gamma_{\mu\nu}^\alpha u_\alpha+\eps\mathfrak{R}_{-3,0}}
and 
\fm{e_4(e_3^\alpha)\partial_\alpha u&=\eps\mathfrak{R}_{-2,0}\cdot\eps\mathfrak{R}_{-1,-1}=\eps^2\mathfrak{R}_{-3,-1}.}
In addition, for $k,l\neq3$, we have
\fm{e_k^\mu e_l^\nu\Gamma_{\mu\nu}^\alpha u_\alpha&=\frac{1}{2}g^{\alpha\beta}(\partial_\mu g_{\nu\beta}+\partial_\nu g_{\mu\alpha}-\partial_\beta g_{\mu\nu})e_k^\mu e_l^\nu u_\alpha\\
&=\frac{1}{2}g^{\alpha\beta}e_k(g_{\nu\beta}) e_l^\nu u_\alpha+\frac{1}{2}g^{\alpha\beta}e_l( g_{\mu\alpha}) e_k^\mu  u_\alpha-\frac{1}{2}g^{\alpha\beta}\partial_\beta g_{\mu\nu}e_k^\mu e_l^\nu u_\alpha\\
&=\eps^2\mathfrak{R}_{-3,-1}-\frac{1}{2}\sum_c e_c (g_{\mu\nu})e_c(u)e_k^\mu e_l^\nu-\frac{1}{4}e_3(g_{\mu\nu})e_4(u)e_k^\mu e_l^\nu-\frac{1}{4}e_4(g_{\mu\nu})e_3(u)e_k^\mu e_l^\nu\\
&=\eps^2\mathfrak{R}_{-3,-1}.}
Since $\chi_{ab}-\delta_{ab}r^{-1}=\mathfrak{R}_{-2,0}$ and $e_3(u)=(3R-r+t)^{-1}V_3(u)=\eps\mathfrak{R}_{-1,-1}$, their product is $\eps \mathfrak{R}_{-3,-1}$. Thus we have
\fm{0&=\sum_a e_a(e_a(u))+e_4(e_3(u))+\frac{1}{2}\tr\chi e_3(u)+\eps\mathfrak{R}_{-3,0}\\&=\sum_a e_a(e_a(u))+e_4(e_3(u))+r^{-1} e_3(u)+\eps\mathfrak{R}_{-3,0}.}

Next, as in Lemma \ref{l3ca1}, we set 
\fm{h_i&:=r(\partial_i(ru)-q_iq_r^{-1}\partial_r(ru))=-r(u+ru_r)q_r^{-1}(q_i-\omega_iq_r)+r^2(u_i-\omega_iu_r).}
Recall from Lemma \ref{l3ca1} that \fm{e_a(ru)=\sum_i e_a(\omega_i)h_i.}
We claim that $h_i=\eps\mathfrak{R}_{0,0}$ and $e_a(h_i)=\eps\mathfrak{R}_{-1,0}$. In fact, note that $u+ru_r=\eps\mathfrak{R}_{-1,0}+\mathfrak{R}_{1,0}\cdot\eps\mathfrak{R}_{-1,-1}=\eps\mathfrak{R}_{0,-1}$. We also recall that $e_a(r)=\mathfrak{R}_{-1,0}$, so $e_a(r^{-1})=-r^{-2}e_a(r)=\mathfrak{R}_{-3,0}$. Thus by Lemma \ref{l414}, we have $h_i=\eps\mathfrak{R}_{0,0}$ and $e_a(h_i)=\eps\mathfrak{R}_{-1,0}$. We thus have
\fm{&\hspace{1.5em}e_a(e_a(u))\\&=e_a(r^{-1}e_a(ru)+e_a(r^{-1})ru)\\&=r^{-1}e_a(e_a(ru))+2e_a(r^{-1})e_a(ru)+e_a(e_a(r^{-1}))ru\\
&=r^{-1}e_a(e_a(ru))+\mathfrak{R}_{-3,0}\cdot  r^{-1}V_a(\eps\mathfrak{R}_{0,0})+V_a(\mathfrak{R}_{-3,0})\cdot\eps\mathfrak{R}_{-1,0}\\
&=r^{-1}\sum_ie_a(e_a(\omega_i))h_i+r^{-1}\sum_ie_a(\omega_i)e_a(h_i)+\eps\mathfrak{R}_{-4,0}\\
&=r^{-1}\sum_ie_a(r^{-1}(e_a^i-\omega_i\omega_je_a^j))h_i+\mathfrak{R}_{-1,0}\cdot r^{-1}V_a(\mathfrak{R}_{0,0})\cdot \eps \mathfrak{R}_{-1,0}+\eps\mathfrak{R}_{-4,0}\\
&=r^{-2}\sum_ie_a(e_a^i-\omega_i\omega_je_a^j)h_i+r^{-1}\sum_ie_a(r^{-1})(e_a^i-\omega_i\omega_je_a^j)h_i+\eps\mathfrak{R}_{-3,0}\\
&=r^{-2} e_a(\mathfrak{R}_{0,0})\cdot \eps\mathfrak{R}_{0,0}+r^{-1} \mathfrak{R}_{-3,0}\cdot\eps \mathfrak{R}_{0,0}+\eps\mathfrak{R}_{-3,0}=\eps\mathfrak{R}_{-3,0}.}
Thus,
\fm{0&=e_4(e_3(u))+r^{-1} e_3(u)+\eps\mathfrak{R}_{-3,0}.}
Finally, we have
\fm{e_4(e_3(ru))&=e_4(re_3(u))+e_4(e_3(r)u)=re_4(e_3(u))+e_4(r)e_3(u)+e_3(r)e_4(u)+e_4(e_3(r))u\\
&=-e_3(u)+e_4(r)e_3(u)+e_4(e_3^i\omega_i)u+\eps r\mathfrak{R}_{-3,0}+\eps\mathfrak{R}_{-2,0}\\
&=(e_4(r)-1)e_3(u)+t^{-1}V_4(1+(e_3^i-\omega_i)\omega_i) u+\eps\mathfrak{R}_{-2,0}\\
&=\mathfrak{R}_{-1,0}\cdot\eps\mathfrak{R}_{-1,-1}+\mathfrak{R}_{-1,0}\cdot V_4(\mathfrak{R}_{-1,0})\cdot\eps\mathfrak{R}_{-1,0}+\eps\mathfrak{R}_{-2,0}=\eps\mathfrak{R}_{-2,0}.}

\end{proof}\rm
\bigskip

Next we prove an estimate for $e_3(q)$. We start with the following lemma.
\lem{\label{l424im}Fix a function $f\in C^\infty(\R)$. Then, for $\eps\ll 1$, $f(u)-f(0)-f'(0)u=\eps^2\mathfrak{R}_{-2,0}$ where $u$ is a solution to \eqref{qwe}.}
\begin{proof}For $\eps\ll 1$, we have $f(u)-f(0)-f'(0)u=O(|u|^2)=O(\eps^2t^{-2+C\eps})$. Now, for each $I$ with $|I|>0$, we can write $V^I(f(u))-f'(u)(V^Iu)$ as a linear combination of terms of the form
\fm{f^{(s)}(u)V^{I_1}u\cdots V^{I_s}u,\hspace{2em}\sum|I_*|=|I|,\ s\geq 2,\ |I_*|>0.}
Since $u=\eps\mathfrak{R}_{-1,0}$, we can prove that each of these terms are $O((\eps t^{-1+C\eps})^s)=O(\eps^2t^{-2+C\eps})$. Finally, note that
$f'(u)V^Iu-f'(0)V^Iu=O(|u|\cdot|V^Iu|)=O(\eps^2t^{-1+C\eps})$. This finishes the proof.
\end{proof}\rm

Our main result is as follows.

\prop{\label{l422}In $\Omega\cap\{r-t<2R\}$, we have \fm{e_4(e_3(q))=-\frac{1}{4}e_3(u)G(\omega)e_3(q)+\eps\mathfrak{R}_{-2,0}.}}
\begin{proof}
We recall that
\fm{e_4(e_3(q))=-\Gamma_{\alpha\beta}^0 e_4^\alpha e_4^\beta e_3(q)=-\frac{1}{2}g^{0\nu}( e_4^\beta e_4( g_{\nu\beta})+e_4^\alpha e_4(g_{\nu\alpha})) e_3(q)+\frac{1}{2}g^{0\nu}\partial_\nu g_{\alpha\beta} e_4^\alpha e_4^\beta e_3(q).}
Here $e_3(q)=(3R-r+t)^{-1}V_3(q)=\mathfrak{R}_{0,0}$ and $e_4(g)=t^{-1}V_4(g)=\eps \mathfrak{R}_{-2,0}$. Thus, 
\fm{e_4(e_3(q))&=\frac{1}{2}g^{0\nu}\partial_\nu g_{\alpha\beta} e_4^\alpha e_4^\beta e_3(q)+\eps\mathfrak{R}_{-2,0}=\frac{1}{4}(e_3-e_4)(g_{\alpha\beta}) e_4^\alpha e_4^\beta e_3(q)+\eps\mathfrak{R}_{-2,0}\\
&=\frac{1}{4}e_3(g_{\alpha\beta}) e_4^\alpha e_4^\beta e_3(q)+\eps\mathfrak{R}_{-2,0}.}
Recall that the coefficients $(g^{\alpha\beta}(v))$ in \eqref{qwe} are known smooth functions, and that for all $|v|\ll 1$ the matrix $(g^{\alpha\beta}(v))$ has a smooth inverse $(g_{\alpha\beta}(v))$. We differentiate $g^{\alpha\sigma}(v)g_{\sigma\beta}(v)=\delta_{\alpha\beta}$ with respect to $v$ and then set $v=0$. Thus,
\fm{\frac{d}{dv}g^{\alpha\sigma}|_{v=0}\cdot m_{\sigma\beta}+m^{\alpha\sigma}\cdot\frac{d}{dv}g_{\sigma\beta}|_{v=0}=0.}
By setting $g_{\alpha\beta}^0=\frac{d}{dv}g_{\alpha\beta}|_{v=0}$ and $g^{\alpha\beta}_0=\frac{d}{dv}g^{\alpha\beta}|_{v=0}$, we conclude that \fm{g_{\alpha\beta}^0=-m_{\alpha\alpha}m_{\beta\beta}g^{\alpha\beta}_0.}
Here we do not take sum over $\alpha,\beta$. Thus we have \fm{g_{\alpha\beta}^0e_4^\alpha e_4^\beta&=-g^{00}e_4^0e_4^0+2g_0^{0i}e_4^0e_4^i-g^{ij}_0e_4^ie_4^j\\&=-G(\omega)+2g^{0i}_0(e_4^i-\omega_i)-g^{ij}_0e_4^i(e_4^j-\omega_j)-g^{ij}_0(e_4^i-\omega_i)\omega_j=-G(\omega)+\mathfrak{R}_{-1,0}.} By the previous lemma we have  \fm{e_4(e_3(q))&=\frac{1}{4}e_3(g_{\alpha\beta}^0u)e_4^\alpha e_4^\beta e_3(q)+\eps\mathfrak{R}_{-2,0}=-\frac{1}{4}e_3(u)G(\omega)e_3(q)+\eps\mathfrak{R}_{-2,0}.}
\end{proof}\rm

\section{The asymptotic equations and the scattering data}\label{s5}\rm
In Section \ref{s3}, we have constructed a global optical function $q(t,x)$ in $\Omega$ such that  $-q_t,q_r\geq C^{-1}t^{-C\eps}>0$. By setting 
\fm{\Omega':=\{(s,q,\omega):\ s>0,\ q>(\exp(\delta/\eps)-\exp((s+\delta)/\eps))/2+2R,\ \omega\in\mathbb{S}^2\},}
we have an invertible map from $\Omega$ to $\Omega'$, defined by \fm{\Phi(t,r,\omega)=(s,q,\omega):=(\eps\ln(t)-\delta,q(t,r\omega),\omega).} In fact,  we have $t=\exp((s+\delta)/\eps)$ and the map $r\mapsto q(t,r\omega)$ is strictly increasing for each fixed $(t,\omega)$. Thus, $\Phi$ is injective. Since $q=r-t$ when $r\geq t+2R$, we have $\lim_{r\to\infty}q(t,r\omega)=\infty$. Thus, $\Phi$  is surjective.  This gives us a new coordinate system $(s,q,\omega)$ on~$\Omega$.

In addition, $\Phi$ is  smooth since $q$ is a smooth function. Its inverse $\Phi^{-1}$ is also smooth, since we have $q_r>0$. So, any smooth function $F(t,x)$ induces a smooth function $F\circ \Phi^{-1}$. With an abuse of notation, we still write $F\circ \Phi^{-1}(s,q,\omega)$ as $F(s,q,\omega)$.

We define
\fm{(\mu,U)(t,x)=(q_t-q_r,\eps^{-1}ru)(t,x),\hspace{1em} (t,x)\in\Omega.}
Since $q$ and $u$ are both smooth, $\mu(t,x)$ and $U(t,x)$ are smooth. As discussed above, we also obtain two smooth functions $\mu(s,q,\omega)$ and  $U(s,q,\omega)$  in $\Omega'$. Our goal in this section is to derive a system of asymptotic equations for $(\mu,U)$ in the coordinate set $(s,q,\omega)$. Our main result is the following proposition.

\prop{\label{props5} Let $(\mu,U)(s,q,\omega)$ be defined as above. Then, by writing $t=\exp(\eps^{-1}(s+\delta))$ we have
\fm{\left\{\begin{array}{l}\displaystyle
\partial_s\mu=\frac{1}{4}G(\omega)\mu^2 U_q+\eps^{-1}\mathfrak{R}_{-1,0},\\[1em]
\displaystyle\partial_sU_q=-\frac{1}{4}G(\omega)\mu U_q^2+\eps^{-1}\mathfrak{R}_{-1,0}.
\end{array}\right.}

In addition, the following three limits exist for all $(q,\omega)\in\R\times\mathbb{S}^2$:
\fm{\left\{\begin{array}{l}
\displaystyle A(q,\omega):=-\frac{1}{2}\lim_{s\to\infty}(\mu U_q)(s,q,\omega),\\
\displaystyle A_1(q,\omega):=\lim_{s\to\infty}\exp(\frac{1}{2}G(\omega)A(q,\omega)s)\mu(s,q,\omega),\\
\displaystyle A_2(q,\omega):=\lim_{s\to\infty}\exp(-\frac{1}{2}G(\omega)A(q,\omega)s)U_q(s,q,\omega).
\end{array}\right.}
All of them are smooth functions of $(q,\omega)$ for $\eps\ll 1$. By setting
\fm{\left\{\begin{array}{l}
\displaystyle\wt{\mu}(s,q,\omega):=A_1\exp(-\frac{1}{2}GAs),\\[1em]
\displaystyle \wt{U}_q(s,q,\omega):=A_2\exp(\frac{1}{2}GAs).
\end{array}\right.}
we obtain an exact solution to our reduced system
\fm{\left\{\begin{array}{l}
\displaystyle \wt{\mu}_s=\frac{1}{4}G(\omega)\wt{\mu}^2\wt{U}_q,\\[1em]
\displaystyle \wt{U}_{sq}=-\frac{1}{4}G(\omega)\wt{\mu}\wt{U}_q^2,
\end{array}\right.}
which satisfies the following estimates:
\fm{\begin{array}{ll}(\lra{q}\partial_q)^m\partial_\omega^n (\mu U_q+2A)=O(t^{-1+C\eps}),&\hspace{2em} (\lra{q}\partial_q)^m\partial_\omega^n A=O(\lra{q}^{-1+C\eps});\\[.5em](\lra{q}\partial_q)^m\partial_\omega^n (\exp(\frac{1}{2}GAs)\mu-A_1)=O(t^{-1+C\eps}),&\hspace{2em} (\lra{q}\partial_q)^m\partial_\omega^n A_1=O(\lra{q}^{C\eps}),\\[.5em](\lra{q}\partial_q)^m\partial_\omega^n (\exp(-\frac{1}{2}GAs)U_q-A_2)=O(t^{-1+C\eps}),&\hspace{2em} (\lra{q}\partial_q)^m\partial_\omega^n A_2=O(\lra{q}^{-1+C\eps});\\[.5em]\partial^p_s(\lra{q}\partial_q)^{m}\partial_\omega^{n}(\wt{\mu}-\mu,\wt{U}_q-U_q)=O(\eps^{-p}t^{-1+C\eps}),&\hspace{2em}\partial^p_s \partial_\omega^{n}( \wt{U}-U)=O(\eps^{-p}\lra{q}t^{-1+C\eps}).\end{array}}
}\rm
\rmk{\rm \label{props5.1}Here $A$ is called the \emph{scattering data}.}\rm
\bigskip

After some preliminary computations in the new coordinate set $(s,q,\omega)$ in Section \ref{ss5.1}, we derive the asymptotic equations for $\mu$ and $U$ in Section \ref{ss5.2} and Section \ref{ss5.3}, respectively. Next, in Section \ref{ssdata}, we make use of the asymptotic equations to construct our scattering data. The main propositions in this subsection  are Proposition \ref{proplim} and Proposition \ref{proplim3}. Finally, in Section \ref{ssdata2}, we define an exact solution $(\wt{\mu},\wt{U})(s,q,\omega)$ to our reduced system and we show that it provides a good approximation of $(\mu,U)(s,q,\omega)$.

\subsection{Derivatives under the new coordinate}\label{ss5.1}
For convenience, from now on we make the following convention. For a function $F=F(s,q,\omega)$ where $\omega\in\mathbb{S}^2$, we extend it to all $\omega\neq 0$ by setting $F(s,q,\lambda\omega)=F(s,q,\omega)$ for each $\lambda>0$. Under such a setting, it is easy to compute the angular derivatives of $F$ since we can now define $\partial_{\omega_i}$. To avoid ambiguity, we will only use $\partial_{\omega_i}$ in the coordinate $(s,q,\omega)$ and will never use it in the coordinate $(t,r,\omega)$.

First we explain how to compute the derivatives of $U$ in $(s,q,\omega)$. Note by the chain rule, for any function $F=F(s,q,\omega)=F(t,r,\omega)$ we have
\fm{\left\{\begin{array}{l}F_t=\eps t^{-1}F_s+q_tF_q\\F_r=q_rF_q\end{array}\right.\Longrightarrow\left\{\begin{array}{l}F_s=\eps^{-1} t(F_t-q_tq_r^{-1}F_r)\\F_q=q_r^{-1}F_r\end{array}\right..}
In addition, by the homogeneity, we have $F(s,q,\omega)=F(s,q,\lambda\omega)$ and $\partial_{\omega_i}F(s,q,\omega)=\lambda\partial_{\omega_i}F(s,q,\lambda\omega)$ for each $\lambda>0$. At $(t,x)$, we set $\lambda=|x|$ which gives
\fm{F_i&=q_iF_q+r^{-1}F_{\omega_i}\Longrightarrow F_{\omega_i}=r(F_i-q_iq_r^{-1}F_r).}
Now we can explain the meaning of the function $h_i$ defined in Lemma \ref{l3ca1}; it is the derivative of $ru$ with respect to $\omega_i$ under the coordinate $(s,q,\omega)$.

To simplify our future computations, we note that $\partial_q$, $\partial_s$ and $\partial_{\omega_i}$ commute with each other. In fact,
\fm{\ [\partial_q,\partial_{\omega_i}]&=[q_r^{-1}\partial_r,r\partial_i-rq_iq_r^{-1}\partial_r]\\
&=q_r^{-1}\partial_i-q_r^{-1}\partial_r(rq_iq_r^{-1})\partial_r-r\partial_i(q_r^{-1}\omega_j)\partial_j+rq_iq_r^{-1}\partial_r(q_r^{-1})\partial_r\\
&=q_r^{-1}\partial_i-q_r^{-2}\partial_r(rq_i)\partial_r-r\partial_i(q_r^{-1})\partial_r-q_r^{-1}(\partial_i-\omega_i\partial_r)\\
&=-q_r^{-2}(q_i+r\partial_rq_i)\partial_r+rq_r^{-2}(\partial_r(q_i)+r^{-1}(q_i-\omega_iq_r))\partial_r+q_r^{-1}\omega_i\partial_r\\
&=0,}
\fm{\ [\partial_s,\partial_q]&=[\eps^{-1}t\partial_t-\eps^{-1}tq_tq_r^{-1}\partial_r,q_r^{-1}\partial_r]\\
&=\eps^{-1}t\partial_t(q_r^{-1})\partial_r-\eps^{-1}tq_tq_r^{-1}\partial_r(q_r^{-1})\partial_r+\eps^{-1}tq_r^{-1}\partial_r(q_tq_r^{-1})\partial_r\\
&=\eps^{-1}t\partial_t(q_r^{-1})\partial_r+\eps^{-1}tq_r^{-2}q_{tr}\partial_r=0,}
\fm{\ [\partial_s,\partial_{\omega_i}]&=[\eps^{-1}t\partial_t-\eps^{-1}tq_tq_r^{-1}\partial_r,r\partial_i-rq_iq_r^{-1}\partial_r]\\
&=-\eps^{-1}tr\partial_t(q_iq_r^{-1})\partial_r-\eps^{-1}tq_tq_r^{-1}(\partial_i-\partial_r(rq_iq_r^{-1})\partial_r)\\
&\hspace{1em}+\eps^{-1}tr\partial_i(q_tq_r^{-1}\omega_j)\partial_j-\eps^{-1}trq_iq_r^{-1}\partial_r(q_tq_r^{-1})\partial_r\\
&=-\eps^{-1}trq_{it}q_r^{-1}\partial_r-\eps^{-1}tq_tq_r^{-1}\partial_i+\eps^{-1}tq_tq_r^{-2}q_i\partial_r+\eps^{-1}trq_tq_r^{-2}\partial_r(q_i)\partial_r\\
&\hspace{1em}+\eps^{-1}tr q_{ti}q_r^{-1}\partial_r-\eps^{-1}trq_tq_r^{-2}\partial_i(q_r)\partial_r+\eps^{-1}tq_tq_r^{-1}(\partial_i-\omega_i\partial_r)\\
&=\eps^{-1}tq_tq_r^{-2}q_i \partial_r-\eps^{-1}tq_tq_r^{-2}(q_i-\omega_iq_r)\partial_r-\eps^{-1}tq_tq_r^{-1}\omega_i\partial_r=0.}

Moreover, we can express $(\partial_s,\partial_q,\partial_{\omega_i})$ in terms of the weighted null frame $\{V_k\}$.
\lem{\label{sqownf}We have
\fm{\partial_s&=\sum_{a}\eps^{-1}\mathfrak{R}_{-1,0}V_a+(\eps^{-1}+\mathfrak{R}_{-1,0})V_4,\\
\partial_{\omega_i}&=\sum_{k\neq 3}\mathfrak{R}_{-1,0}V_k+\sum_ae_a^iV_a=\sum_{k\neq 3}\mathfrak{R}_{0,0}V_k,\\
\partial_q&=\sum_{k}\mathfrak{R}_{0,-1}V_k.}}
\begin{proof}
We can express $\partial_s,\partial_{\omega_i}$ in terms of the null frame:
\fm{\partial_s&=\eps^{-1}t(g_{0\beta}e_a^\beta e_a+\frac{1}{2}g_{0\beta}e_4^\beta e_3+\frac{1}{2}g_{0\beta}e_3^\beta e_4)-\eps^{-1}tq_tq_r^{-1}(\omega_ig_{i\beta}e_a^\beta e_a+\frac{1}{2}\omega_ig_{i\beta}e_4^\beta e_3+\frac{1}{2}\omega_ig_{i\beta}e_3^\beta e_4)\\
&=\eps^{-1}t((g_{0\beta}-q_tq_r^{-1}\omega_ig_{i\beta})e_a^\beta e_a+\frac{1}{2}(g_{0\beta}-q_tq_r^{-1}\omega_ig_{i\beta})e_3^\beta e_4),\\
\partial_{\omega_i}&=r(g_{i\beta}e_a^\beta e_a+\frac{1}{2}g_{i\beta}e_4^\beta e_3+\frac{1}{2}g_{i\beta}e_3^\beta e_4)-rq_iq_r^{-1}(\omega_jg_{j\beta}e_a^\beta e_a+\frac{1}{2}\omega_jg_{j\beta}e_4^\beta e_3+\frac{1}{2}\omega_jg_{j\beta}e_3^\beta e_4)\\
&=r((g_{i\beta}-q_iq_r^{-1}\omega_jg_{j\beta})e_a^\beta e_a+\frac{1}{2}(g_{i\beta}-q_iq_r^{-1}\omega_jg_{j\beta})e_3^\beta e_4).}

We note that there is no term with $e_3$ in $\partial_s$ and $\partial_{\omega_i}$, since
\fm{(g_{0\beta}-q_tq_r^{-1}\omega_ig_{i\beta})e_4^\beta&=q_r^{-1}(q_rg_{0\beta}-q_t\omega_ig_{i\beta})e_4^\beta=\frac{1}{2}q_r^{-1}e_3(q)(\omega_ig_{i\nu}e_4^\nu g_{0\beta}e_4^\beta-g_{0\nu} e_4^\nu\omega_ig_{i\beta}e_4^\beta)=0,\\
(g_{i\beta}-q_iq_r^{-1}\omega_jg_{j\beta})e_4^\beta&=q_r^{-1}(q_rg_{i\beta}-q_i\omega_jg_{j\beta})e_4^\beta=\frac{1}{2}q_r^{-1}e_3(q)(\omega_j g_{j\nu}e_4^\nu g_{i\beta}e_4^\beta-g_{i\nu}e_4^\nu\omega_jg_{j\beta}e_4^\beta)=0.}
In these computations we use the equality  $q_\alpha=\frac{1}{2}g_{\alpha\beta}e_4^\beta e_3(q)$. In addition, we have
\fm{\eps^{-1}t(g_{0\beta}-q_tq_r^{-1}\omega_ig_{i\beta})e_a^\beta&=\eps^{-1}t((g_{0j}-m_{0j})-q_tq_r^{-1}\omega_i(g_{ij}-m_{ij}))e_a^j-\eps^{-1} tq_tq_r^{-1}e_a(r)\\
&= \mathfrak{R}_{0,0}+\eps^{-1}\mathfrak{R}_{0,0}=\eps^{-1}\mathfrak{R}_{0,0},\\
r(g_{i\beta}-q_iq_r^{-1}\omega_jg_{j\beta})e_a^\beta&=r((g_{ij'}-m_{ij'})-q_iq_r^{-1}\omega_j(g_{jj'}-m_{jj'}))e_a^{j'}+r(e_a^i-q_iq_r^{-1}e_a(r))\\
&=\mathfrak{R}_{0,0}+re_a^i.}
Besides, since $e_3^i\omega_i=2g^{0i}\omega_i+e_4^i\omega_i=1+\eps\mathfrak{R}_{-1,0}$, we have
\fm{\eps^{-1}t(g_{0\beta}-q_tq_r^{-1}\omega_ig_{i\beta})e_3^\beta&=\eps^{-1}t((g_{0\beta}-m_{0\beta})-q_tq_r^{-1}\omega_i(g_{i\beta}-m_{i\beta}))e_3^\beta+\eps^{-1}t(1-q_tq_r^{-1}e_3^i\omega_i)\\
&=\mathfrak{R}_{0,0}+\eps^{-1}tq_r^{-1}(2q_r-(q_t+q_r)-q_t(e_3^i\omega_i-1))=\mathfrak{R}_{0,0}+2\eps^{-1}t,\\
r(g_{i\beta}-q_iq_r^{-1}\omega_jg_{j\beta})e_3^\beta&=r((g_{i\beta}-m_{i\beta})-q_iq_r^{-1}\omega_j(g_{j\beta}-m_{j\beta}))e_3^\beta+r(e_3^i-q_iq_r^{-1}\omega_je_3^j)\\
&=\eps\mathfrak{R}_{0,0}+rq_r^{-1}((e_3^i-\omega_i)q_r-(q_i-\omega_iq_r)-q_i(e_3^j\omega_j-1))=\mathfrak{R}_{0,0}.}
Thus,
\fm{\partial_s&=\sum_{a}\eps^{-1}\mathfrak{R}_{0,0}e_a+(\eps^{-1}t+\mathfrak{R}_{0,0})e_4=\sum_{a}\eps^{-1}\mathfrak{R}_{-1,0}V_a+(\eps^{-1}+\mathfrak{R}_{-1,0})V_4,\\
\partial_{\omega_i}&=\sum_{k\neq 3}\mathfrak{R}_{0,0}e_k+\sum_are_a^ie_a=\sum_{k\neq 3}\mathfrak{R}_{-1,0}V_k+\sum_ae_a^iV_a=\sum_{k\neq 3}\mathfrak{R}_{0,0}V_k.}
It is also clear that
\fm{\partial_q&=\sum_k\mathfrak{R}_{0,0}e_k=\sum_{k}\mathfrak{R}_{0,-1}V_k.}
\end{proof}\rm

We end this subsection with the following estimates for $U$.
\lem{\label{lemmuq}We have \fm{(U,U_q,U_s,U_{\omega_i})=(\mathfrak{R}_{0,0},\mathfrak{R}_{0,-1},\eps^{-1}\mathfrak{R}_{0,0},\mathfrak{R}_{0,0}).} In conclusion, we have $\mu U_q=\mathfrak{R}_{0,-1}$.}
\begin{proof}
We have
\fm{U&=\eps^{-1}ru,\\
U_q&=q_r^{-1}\partial_r(\eps^{-1}ru)=\eps^{-1}q_r^{-1}(u+ru_r),\\
U_s&=\eps^{-2}tr(u_t+u_r-q_r^{-1}(q_t+q_r) u_r)-\eps^{-2}tq_tq_r^{-1}u,\\
U_{\omega_i}&=-\eps^{-1}r(q_i-\omega_iq_r)q_r^{-1}(u+ru_r)+\eps^{-1}r^2(u_i-\omega_iu_r).}
It follows directly from Lemma \ref{l5.6}, Lemma \ref{l414} and the proof of Proposition \ref{l413} that $(U,U_q,U_s,U_{\omega_i})=(\mathfrak{R}_{0,0},\mathfrak{R}_{0,-1},\eps^{-1}\mathfrak{R}_{0,0},\mathfrak{R}_{0,0})$. Finally, since $\mu=\mathfrak{R}_{0,0}$, we have $\mu U_q=\mathfrak{R}_{0,-1}$.
\end{proof}
\rmk{\rm  Note that we have $U_q=\mathfrak{R}_{0,-1}$ which is stronger than $(U,\eps U_s,U_{\omega_i})=\mathfrak{R}_{0,0}$. This is because we gain an additional factor $\lra{r-t}^{-1}$ when we estimate $u_r$ and its derivatives compared to $u$. We refer our readers to Lemma \ref{l2.1}.

We also remark that it is important for us to obtain $U_q=\mathfrak{R}_{0,-1}$ instead of $\mathfrak{R}_{0,0}$  here, because $U_q=\mathfrak{R}_{0,-1}$ is necessary for the scattering data to be defined later.}
\rm

\subsection{The asymptotic equation for $\mu$}\label{ss5.2}
We start with several estimates for $\mu=q_t-q_r$. By  Proposition \ref{l422}, we have
\fm{e_4(e_3(q))&=-\frac{1}{4}e_3(u)G(\omega)e_3(q)+\eps\mathfrak{R}_{-2,0}\\
&=-\frac{1}{4}(\eps r^{-1}e_3(U)-\eps r^{-2}e_3(r)U)G(\omega)e_3(q)+\eps\mathfrak{R}_{-2,0}\\
&=-\frac{\eps}{4 r}e_3(U)G(\omega)e_3(q)+\eps\mathfrak{R}_{-2,0}.}
Since $e_3^i-\omega_i=\mathfrak{R}_{-1,0}$, we have 
\fm{e_3(q)&=-\mu+\mathfrak{R}_{-1,0}\cdot\partial q=-\mu+\mathfrak{R}_{-1,0}.}
Moreover, 
\fm{e_4(e_3(q)+\mu)&=e_4((e_3^i-\omega_i)q_i)=e_4(e_3^i-\omega_i)q_i+(e_3^i-\omega_i)e_4(q_i)\\
&=-(e_4^j-\omega_j)r^{-1}(\delta_{ij}-\omega_i\omega_j)q_i+(e_3^i-\omega_i)e_4(\frac{1}{2}g_{i\beta}e_4^\beta e_3(q))+\eps\mathfrak{R}_{-2,0}\\
&=-r^{-1}(-q_t-q_r-q_r(e_4(r)-1))+\frac{1}{2}g_{i\beta}(e_3^i-\omega_i)e_4^\beta e_4( e_3(q))+\eps\mathfrak{R}_{-2,0}=\eps\mathfrak{R}_{-2,0}.}
To get the last equality, we use the following estimates: $e_4(r)-1=\eps\mathfrak{R}_{-1,0}$ by Lemma \ref{derqle4r}, $e_4(e_3(q))=\xi_{43}^3e_3(q)=\eps\mathfrak{R}_{-1,-1}$, and \fm{q_t+q_r&=\frac{1}{2}(g_{0\beta}+\omega_i g_{i\beta})e_4^\beta e_3(q)=\frac{1}{2}(-1+e_4^i\omega_i)e_3(q)+(g_{**}-m_{**})\cdot\mathfrak{R}_{0,0}=\eps\mathfrak{R}_{-1,0}.}

Besides, by the chain rule, we have \fm{ e_3(U)&=e_3(q)U_q-\eps t^{-1}U_s+\sum_ie_3(\omega_i)U_{\omega_i}=-\mu U_q+\mathfrak{R}_{-1,0}.}
 Here we apply Lemma \ref{lemmuq} and we note that $e_3(\omega_i)=(e_3^j-\omega_j)r^{-1}(\delta_{ij}-\omega_i\omega_j)=\mathfrak{R}_{-2,0}$. Thus, we have 
\fm{e_4(-\mu)+\eps\mathfrak{R}_{-2,0}&=-\frac{\eps}{4r} G(\omega)(-\mu U_q+\mathfrak{R}_{-1,0})(-\mu+\mathfrak{R}_{-1,0})+\eps \mathfrak{R}_{-2,0}\\
&=-\frac{\eps}{4r}G(\omega)\mu^2U_q+\eps\mathfrak{R}_{-2,0}.}
 Then,
\eq{\label{e4mucomp}e_4(\mu)&=\frac{\eps}{4r}G(\omega)\mu^2U_q+\eps\mathfrak{R}_{-2,0}.}
By Lemma \ref{sqownf} we have
\fm{\mu_s&=\eps^{-1}te_4(\mu)+\sum_{k\neq 3}\eps^{-1}\mathfrak{R}_{-1,0}V_k(\mu)=\eps^{-1}t(\frac{\eps}{4r}G(\omega)\mu^2U_q+\eps\mathfrak{R}_{-2,0})+\sum_{k\neq 3}\eps^{-1}\mathfrak{R}_{-1,0}V_k(\mathfrak{R}_{0,0})\\
&=\frac{ t}{4r}G(\omega)\mu^2U_q+\eps^{-1}\mathfrak{R}_{-1,0}=\frac{1}{4}G(\omega)\mu^2U_q+\frac{\eps (t-r)}{4r}G(\omega)\mu^2U_q+\eps^{-1}\mathfrak{R}_{-1,0}\\
&=\frac{1}{4}G(\omega)\mu^2U_q+\eps\mathfrak{R}_{-1,1}\cdot\mathfrak{R}_{0,0}\cdot\mathfrak{R}_{0,-1}+\eps^{-1}\mathfrak{R}_{-1,0}=\frac{1 }{4}G(\omega)\mu^2U_q+\eps^{-1}\mathfrak{R}_{-1,0}.}
We thus obtain the first asymptotic equation
\eq{\label{museqn} \mu_s=\frac{1}{4}G(\omega)\mu^2 U_q+\eps^{-1}\mathfrak{R}_{-1,0}.}

\subsection{The asymptotic equation for $U$}\label{ss5.3}
By Proposition \ref{l413}, we have \fm{e_4(e_3(U))=\eps^{-1}e_4(e_3(ru))=\mathfrak{R}_{-2,0}.}
Meanwhile, by Lemma \ref{lemmuq} we have \fm{e_4(e_3(U))&=e_4(e_3(q)U_q+\eps t^{-1}U_s+e_3(\omega_i)U_{\omega_i})\\
&=-e_4(\mu U_q)+e_4((e_3^i-\omega_i)q_iU_q+\eps t^{-1}U_s+(e_3^j-\omega_j)r^{-1}(\delta_{ij}-\omega_i\omega_j)U_{\omega_i})\\
&=-e_4(\mu U_q)+\mathfrak{R}_{-1,0}\cdot V_4(\mathfrak{R}_{-1,-1}+\eps t^{-1}\cdot\eps^{-1}\mathfrak{R}_{0,0}+\mathfrak{R}_{-1,0} \cdot r^{-1}\cdot \mathfrak{R}_{0,0})\\
&=-e_4(\mu U_q)+\mathfrak{R}_{-2,0}.}
Thus, $e_4(\mu U_q)=\mathfrak{R}_{-2,0}$.

Now, we compute $\partial_s(\mu U_q)$. By Lemma \ref{sqownf} we have
\fm{\partial_s(\mu U_q)&=\sum_a \eps^{-1}\mathfrak{R}_{-1,0}V_a(\mu U_q)+(\eps^{-1}+\mathfrak{R}_{-1,0})V_4(\mu U_q)\\
&=\sum_a \eps^{-1}\mathfrak{R}_{-1,0}V_a(\mathfrak{R}_{0,-1})+(\eps^{-1}+\mathfrak{R}_{-1,0})\mathfrak{R}_{-1,0}=\eps^{-1}\mathfrak{R}_{-1,0}.}
Thus, we have
\fm{\mu U_{sq}&=\partial_s(\mu U_q)-\mu_sU_q=\eps^{-1}\mathfrak{R}_{-1,0}-(\frac{1}{4}G(\omega)\mu^2U_q+\eps^{-1}\mathfrak{R}_{-1,-1}+\mathfrak{R}_{-1,0})U_q\\
&=-\frac{1}{4}G(\omega)\mu^2U_q^2+\eps^{-1}\mathfrak{R}_{-1,0}.}
Since $|\mu|>C^{-1}t^{-C\eps}$, we have $\mu^{-1}=\mathfrak{R}_{0,0}$. Thus we obtain the second asymptotic equation
\eq{\label{uqseqn} U_{sq}&=-\frac{1}{4}G(\omega)\mu U_q^2+\eps^{-1}\mathfrak{R}_{-1,0}.}

In summary, by \eqref{museqn} and \eqref{uqseqn}, we have proved the following proposition.
\prop{\label{propasyeqn}We have\eq{\label{asyeqn}\left\{\begin{array}{l}\displaystyle
\partial_s\mu=\frac{1}{4}G(\omega)\mu^2 U_q+\eps^{-1}\mathfrak{R}_{-1,0},\\[1em]
\displaystyle\partial_sU_q=-\frac{1}{4}G(\omega)\mu U_q^2+\eps^{-1}\mathfrak{R}_{-1,0}.
\end{array}\right.}
In other words, $(\mu,U_q)(s,q,\omega)$ is an apporximate solution to the reduced system of ODE's
\eq{\label{asyode}\left\{\begin{array}{l}\displaystyle
\partial_s\wt{\mu}=\frac{1}{4}G(\omega)\wt{\mu}^2 \wt{U}_q,\\[1em]
\displaystyle\partial_s\wt{U}_q=-\frac{1}{4}G(\omega)\wt{\mu} \wt{U}_q^2.
\end{array}\right.}}\rm

We remark that this proposition verifies the nonrigorous derivation in Section 3 of the author's previous paper \cite{yu2020}.

\subsection{The scattering data}\label{ssdata}
From the previous subsections, we have proved that $(\mu,U_q)(s,q,\omega)$ is an approximate solution to the reduced system \eqref{asyode}. In this subsection, we seek to construct an exact solution $(\wt{\mu},\wt{U}_q)$ to \eqref{asyode} which is a good approximation of $(\mu,U_q)$.

We start with the following key proposition. In this proposition, we define the \emph{scattering data} $A=A(q,\omega)$ for each $(q,\omega)\in\R\times\mathbb{S}^2$ and we show that it is a smooth function (in the sense defined in Section \ref{sec2.5}).

\prop{\label{proplim} In $\Omega'$, we have \fm{(\lra{q}\partial_q)^m\partial_\omega^n (\mu U_q)=O(\lra{q}^{-1}t^{C\eps}),\hspace{2em}\partial_s^p(\lra{q}\partial_q)^m\partial_\omega^n (\mu U_q)=O(\eps^{-p}t^{-1+C\eps}),\ p\geq 1.}

Moreover, for each $m,n$, the limit \fm{A_{m,n}(q,\omega):=-\frac{1}{2}\lim_{s\to\infty}(\lra{q}\partial_q)^m\partial_\omega^n (\mu U_q)(s,q,\omega)} exists for all $(q,\omega)\in\R\times\mathbb{S}^2$, and the convergence is uniform in $(q,\omega)$. So $A(q,\omega):=A_{0,0}(q,\omega)$ is a smooth function of $(q,\omega)$ in $\R\times\mathbb{S}^2$ such that $(\lra{q}\partial_q)^m\partial_\omega^nA=A_{m,n}$. We call this function $A$ the \emph{scattering data}. It is clear that $A\equiv 0$ for $q>R$.

Finally, we have \fm{(\lra{q}\partial_q)^m\partial_\omega^n (\mu U_q+2A)=O(t^{-1+C\eps}),\hspace{2em} (\lra{q}\partial_q)^m\partial_\omega^n A=O(\lra{q}^{-1+C\eps}).}}
\begin{proof}First we note that in the region $r-t>R$, we have $q=r-t$ and $u=0$. In this case, every estimate in the statement of this proposition is equal to  $0$, so there is nothing to prove. Thus, we can assume that $q<2R$ and $r-t<2R$ in the rest of this proof.

We need to derive an estimate for $\partial_s\partial_q^m\partial_{\omega}^n(\mu U_q)$. Here we apply Lemma \ref{sqownf}. Recall that  $\mu U_q=\mathfrak{R}_{0,-1}$ and $V_4(\mu U_q)=\mathfrak{R}_{-1,0}$. By the Leibniz's rule, we have
\fm{(\lra{q}\partial_q)^{m}\partial_\omega^{n}(\mu U_q)&=(\sum_k\mathfrak{R}_{0,0}V_k)^{m+n}(\mathfrak{R}_{0,-1})=O(\lra{q}^{-1+C\eps}t^{C\eps})=O(\lra{q}^{-1}t^{C\eps}).} In addition, for $p\geq 1$ we have \eq{\label{proplimf11}\partial_s^p(\lra{q}\partial_q)^{m}\partial_\omega^{n}(\mu U_q)&=\partial^{p-1}_s(\lra{q}\partial_q)^{m}\partial_\omega^{n}\partial_s(\mu U_q)\\&=\partial^{p-1}_s(\lra{q}\partial_q)^{m}\partial_\omega^{n}(\sum_{k\neq 3}\eps^{-1}\mathfrak{R}_{-1,0}\cdot V_k(\mu U_q)+\eps^{-1}V_4(\mu U_q))\\
&=\partial^{p-1}_s(\lra{q}\partial_q)^{m}\partial_\omega^{n}(\sum_a\eps^{-1}\mathfrak{R}_{-1,0}\cdot \mathfrak{R}_{0,-1}+\eps^{-1}\mathfrak{R}_{-1,0})\\
&=\eps^{1-p}(\sum\mathfrak{R}_{0,0}V_k)^{p+m+n-1}(\eps^{-1}\mathfrak{R}_{-1,0})=O(\eps^{-p}t^{-1+C\eps}).}
In both these estimates, we view $t$ as a function of $s$.

For fixed $q<2R$ and $\omega\in\mathbb{S}^2$, by the definition of $\Omega'$, we have $(s,q,\omega)\in\Omega'$ if and only if $s>0$ and 
\eq{\label{Opcond}\exp((s+\delta)/\eps)>\exp(\delta/\eps)-2q+4R.}We can write this condition as $s>s_{q,\delta,\eps}$ where $s_{q,\delta,\eps}\geq 0$ is a constant depending on its subscripts, such that  $(s_{q,\delta,\eps},q,\omega)\in\partial \Omega'$ corresponds with a point on $H$. Thus, for each fixed $(q,\omega)$ and $s_2>s_1\geq s_{q,\delta,\eps}=\exp(\delta/\eps)-2q+4R$, by \eqref{proplimf11} with $p=1$, we have
\fm{&\hspace{1.5em}|(\lra{q}\partial_q)^{m}\partial_\omega^{n}(\mu U_q)(s_2,q,\omega)-(\lra{q}\partial_q)^{m}\partial_\omega^{n}(\mu U_q)(s_1,q,\omega)|\\&\lesssim\int_{s_1}^{s_2}\eps^{-1}\exp((-1+C\eps)\eps^{-1}(s+\delta))\ ds\lesssim \exp((-1+C\eps)\eps^{-1}(s_1+\delta)).}In conclusion, $\{(\lra{q}\partial_q)^{m}\partial_\omega^{n}(\mu U_q)(s,q,\omega)\}_{s\geq s_{q,\delta,\eps}}$ is uniformly Cauchy for each $(q,\omega)$. Thus, the limit\fm{A_{m,n}(q,\omega):=-\frac{1}{2}\lim_{s\to\infty}(\lra{q}\partial_q)^{m}\partial_\omega^{n}(\mu U_q)(s,q,\omega)}exists, and the convergence is uniform in $(q,\omega)$. Besides, for each $s\geq s_{q,\delta,\eps}$, we have \eq{\label{proplimf1}|(\lra{q}\partial_q)^{m}\partial_\omega^{n}(\mu U_q)+2A_{m,n}|\lesssim t^{-1+C\eps}=\exp((-1+C\eps)\eps^{-1}(s+\delta)).}
By evaluating \eqref{proplimf1} at $(s_{q,\delta,\eps},q,\omega)$, we have \fm{|A_{m,n}(q,\omega)|&\lesssim |(\lra{q}\partial_q)^{m}\partial_\omega^{n}(\mu U_q)+2A_{m,n}|+|(\lra{q}\partial_q)^{m}\partial_\omega^{n}(\mu U_q)|\\&
\lesssim (\exp(\delta/\eps)-2q+4R)^{-1+C\eps}+\lra{q}^{-1}(\exp(\delta/\eps)-2q+4R)^{C\eps}\lesssim \lra{q}^{-1+C\eps}.}
In the last inequality, we note  that $(a+b)^{C\eps}\leq 2^{C\eps}\max\{a,b\}^{C\eps}\leq 2(a^{C\eps}+b^{C\eps})$ for each pairs $a,b\geq 0$. Since the convergence is uniform in $(q,\omega)$, if we define $A:=A_{0,0}$, then we have \fm{(\lra{q}\partial_q)^m\partial_\omega^n A=A_{m,n}=O(\lra{q}^{-1+C\eps}).}
\end{proof}\rm

Note that each function of $(s,q,\omega)$ can be viewed as a function of $(t,x)$. We then have the following lemma.
\lem{\label{proplim2}By viewing each function of $(s,q,\omega)$ as a function of $(t,x)\in\Omega\cap\{r-t<2R\}$, we have $(A,\partial_\omega A)=\mathfrak{R}_{0,-1}$, $\mu U_q+2A=\mathfrak{R}_{-1,0}$ and $\exp(\pm\frac{1}{2}G(\omega)As)-1=\mathfrak{R}_{0,-1}$.}
\begin{proof}
Note that $V^IA$ is a linear combination of terms of the form \fm{\partial_q^m\partial_\omega^{n}A\cdot V^{I_1}q\cdots V^{I_m}q\cdot V^{J_1}\omega\cdots V^{J_n}\omega,\hspace{2em} \sum|I_*|+\sum|J_*|=|I|.}
Each of these terms is $O(\lra{q}^{-1-m+C\eps}\cdot\lra{q}^{m}t^{C\eps})=O(\lra{q}^{-1}t^{C\eps})$, so $A=\mathfrak{R}_{0,-1}$. The proof of $\partial_\omega A=\mathfrak{R}_{0,-1}$ is essentially the same. 

Moreover, $V^I(\mu U_q+2A)$ is a linear combination of terms of the form \fm{\partial_q^m\partial_\omega^{n}(\mu U_q+A)\cdot V^{I_1}q\cdots V^{I_m}q\cdot V^{J_1}\omega\cdots V^{J_n}\omega,\hspace{2em} \sum|I_*|+\sum|J_*|=|I|;}
\fm{\partial_s^p\partial_q^m\partial_\omega^{n}(\mu U_q)\cdot V^{K_1}s\cdots V^{K_p}s\cdot V^{I_1}q\cdots V^{I_m}q\cdot \cdot V^{J_1}\omega\cdots V^{J_n}\omega,\\ \sum|I_*|+\sum|J_*|+\sum |K_*|=|I|,\ p>0.}
By applying \eqref{proplimf1} to the first row and \eqref{proplimf11} to the second row, we conclude that $V^I(\mu U_q+2A)=O(t^{-1+C\eps})$ and thus $\mu U_q+2A=\mathfrak{R}_{-1,0}$.

Finally, by the chain rule, for each $|I|>0$ we can write $V^I(\exp(\pm\frac{1}{2}G(\omega)As)-1)$ as a linear combination of terms of the form
\fm{\exp(\pm\frac{1}{2}G(\omega)As)\cdot V^{I_1}(\pm\frac{1}{2}G(\omega)As)\cdots V^{I_m}(\pm\frac{1}{2}G(\omega)As),\hspace{2em}\sum|I_*|=|I|,\ |I_*|>0.}
The first term in this product is $O(t^{C\eps})$, and each of the rest terms are $O(V^{I_*}(\mathfrak{R}_{0,-1}))=O(\lra{q}^{-1}t^{C\eps})$, so we conclude that $V^I(\exp(\pm\frac{1}{2}G(\omega)As)-1)=O(\lra{q}^{-1}t^{C\eps})$ for $|I|>0$. When $|I|=0$, since $|e^\rho-1|\lesssim |\rho|e^{|\rho|}$, we have \fm{|\exp(\pm\frac{1}{2}G(\omega)As)-1|\lesssim \lra{q}^{-1+C\eps}s\exp(C\lra{q}^{-1+C\eps}s)\lesssim \lra{q}^{-1}t^{C\eps}.}
Here we note that $s=\eps\ln(t)-\delta=O(t^{C\eps})$. In conclusion, $\exp(\pm\frac{1}{2}GAs)-1=\mathfrak{R}_{0,-1}$.
\end{proof}\rm

By \eqref{asyeqn} and Lemma \ref{proplim2}, we have
\fm{\left\{\begin{array}{l}\displaystyle
\partial_s\mu=-\frac{1}{2}G(\omega)A(q,\omega)\mu+\eps^{-1}\mathfrak{R}_{-1,0},\\[1em]
\displaystyle\partial_sU_q=\frac{1}{2}G(\omega)A(q,\omega)U_q+\eps^{-1}\mathfrak{R}_{-1,0}.
\end{array}\right.}
With the remainder terms omitted, we obtain two linear ODE's for $\mu$ and $U_q$. They motivate us to define
\eq{\label{defv12}\left\{\begin{array}{l}\displaystyle\wt{V}_1:=\exp(\frac{1}{2}G(\omega)A(q,\omega)s)\mu,\\[1em]
\displaystyle\wt{V}_2:=\exp(-\frac{1}{2}G(\omega)A(q,\omega)s)U_q.\end{array}\right.}
 Now we can prove the following proposition.

\prop{\label{proplim3}
We have 
\fm{(\lra{q}\partial_q)^m\partial_\omega^n \wt{V}_1=O(t^{C\eps}),\hspace{2em}\partial_s^p(\lra{q}\partial_q)^m\partial_\omega^n \wt{V}_1=O(\eps^{-p}t^{-1+C\eps}),\ p\geq 1;}
\fm{(\lra{q}\partial_q)^m\partial_\omega^n \wt{V}_2=O(\lra{q}^{-1}t^{C\eps}),\hspace{2em}\partial_s^p(\lra{q}\partial_q)^m\partial_\omega^n \wt{V}_2=O(\eps^{-p}t^{-1+C\eps}),\ p\geq 1.}

Moreover, for each $m,n$, the limit
\fm{A_{j,m,n}(q,\omega):=\lim_{s\to\infty}\wt{V}_j(s,q,\omega),\hspace{2em}j=1,2}
exists for all $(q,\omega)\in\R\times\mathbb{S}^2$, and the convergence is uniform in $(q,\omega)$. So, for $j=1,2$, $A_j:=A_{j,0,0}$  is smooth functions of $(q,\omega)$ in $\R\times\mathbb{S}^2$ such that $(\lra{q}\partial_q)^m\partial_\omega^nA_j=A_{j,m,n}$. It is clear that $A_1\equiv -2$ and $A_2\equiv 0$ for $q>R$. Besides, we have $A_1A_2=-2A$ everywhere.

Finally, we have  \fm{(\lra{q}\partial_q)^m\partial_\omega^n (\wt{V}_1-A_1)=O(t^{-1+C\eps}),\hspace{2em} (\lra{q}\partial_q)^m\partial_\omega^n A_1=O(\lra{q}^{C\eps}),}
\fm{(\lra{q}\partial_q)^m\partial_\omega^n (\wt{V}_2-A_2)=O(t^{-1+C\eps}),\hspace{2em} (\lra{q}\partial_q)^m\partial_\omega^n A_2=O(\lra{q}^{-1+C\eps}).}
\begin{proof} By \eqref{e4mucomp} and since $t/r=1+\mathfrak{R}_{-1,1}$, we have
\fm{V_4(\mu)&=\frac{\eps t}{4r}G(\omega)\mu^2U_q+\eps\mathfrak{R}_{-1,0}=\frac{\eps }{4}G(\omega)\mu^2U_q+\eps\mathfrak{R}_{-1,0}.}
Moreover, by viewing $(s,q,\omega)$ as functions of $(t,x)$, we have
\fm{e_4(G(\omega)A(q,\omega)s)&=\eps G(\omega)At^{-1}+e_4(\omega_j)\partial_{\omega_j} (G A)s=\eps G(\omega)At^{-1}+\mathfrak{R}_{-2,-1}.}
Here we note that $\partial_{\omega_j} (G A)=\mathfrak{R}_{0,-1}$ by Lemma \ref{proplim2} and $e_4(\omega_i)=(e_4^j-\omega_j)\partial_j\omega_i=\mathfrak{R}_{-2,0}$. Then, by Lemma \ref{proplim2}, we have $\wt{V}_1=\mathfrak{R}_{0,0}\cdot \mathfrak{R}_{0,0}=\mathfrak{R}_{0,0}$ and
\fm{V_4(\wt{V}_1)&=\frac{1}{2}V_4(GAs)\wt{V}_1+\exp(\frac{1}{2}GAs)V_4(\mu)\\
&=\frac{1}{4}(2\eps GA+\eps G\mu U_q+\mathfrak{R}_{-1,-1})\wt{V}_1+\eps\mathfrak{R}_{-1,0}\cdot \exp(\frac{1}{2}GAs)\\
&=\frac{1}{4}(\eps \mathfrak{R}_{-1,0}+\mathfrak{R}_{-1,-1})\cdot \mathfrak{R}_{0,0}+\eps\mathfrak{R}_{-1,0} \cdot \mathfrak{R}_{0,0}=\eps\mathfrak{R}_{-1,0}+\mathfrak{R}_{-1,-1}=\mathfrak{R}_{-1,0}.}
Next, we have $\wt{V}_1\wt{V}_2=\mu U_q$ and  $\mu U_q=\mathfrak{R}_{0,-1}$, $V_4(\mu U_q)=\mathfrak{R}_{-1,0}$ from Proposition \ref{proplim}. Since $\mu =q_t-q_r\leq -2C^{-1}t^{-C\eps}$ and $\exp(\frac{1}{2}GAs)\geq \exp(-Cs)=\exp(C\delta)t^{-C\eps}$, we have $|\wt{V}_1|=-\wt{V}_1\geq C^{-1}t^{-C\eps}$. We can express $V^I(\wt{V}_2)=V^I((\mu U_q)/\wt{V}_1)$ as a linear combination of terms of the form
\fm{\wt{V}_1^{-m-1}\cdot V^{I_1}(\wt{V}_1)\cdots V^{I_m}(\wt{V}_1)\cdot V^{I_0}(\mu U_q),\hspace{2em}\sum|I_*|=|I|.}
It is  easy to conclude that $\wt{V}_2=\mathfrak{R}_{0,-1}$ and $V_4(\wt{V}_2)=\mathfrak{R}_{-1,0}$. 

Now we can follow the proof in Proposition \ref{proplim}  to prove every estimate involving $A_2$ in the statement. As for $A_1$, we note that \fm{(\lra{q}\partial_q)^{m}\partial_\omega^{n}(\wt{V}_1)&=(\sum_k\mathfrak{R}_{0,0}V_k)^{m+n}(\mathfrak{R}_{0,0})=O(t^{C\eps}).} In addition,  for $p\geq 1$ we have \fm{\partial_s^p(\lra{q}\partial_q)^{m}\partial_\omega^{n}(\wt{V}_1)&=\partial^{p-1}_s(\lra{q}\partial_q)^{m}\partial_\omega^{n}\partial_s(\wt{V}_1)\\&=\partial^{p-1}_s(\lra{q}\partial_q)^{m}\partial_\omega^{n}(\sum_{k\neq 3}\eps^{-1}\mathfrak{R}_{-1,0}\cdot V_k(\wt{V}_1)+\eps^{-1}V_4(\wt{V}_1))\\
&=\partial^{p-1}_s(\lra{q}\partial_q)^{m}\partial_\omega^{n}(\sum_a\eps^{-1}\mathfrak{R}_{-1,0}\cdot \mathfrak{R}_{0,0}+\eps^{-1}\mathfrak{R}_{-1,0})\\
&=\eps^{1-p}(\sum\mathfrak{R}_{0,0}V_k)^{p+m+n-1}(\eps^{-1}\mathfrak{R}_{-1,0})=O(\eps^{-p}t^{-1+C\eps}).}
It is then clear that the estimates for $\wt{V}_1-A_1$ are the same as those for $\mu U_q+2A$. Finally, at $(s,q,\omega)=(s_{q,\delta,\eps},q,\omega)$ we have
\fm{|(\lra{q}\partial_q)\partial_\omega^nA_{1}(q,\omega)|&\lesssim |(\lra{q}\partial_q)\partial_\omega^n(\wt{V}_1-A_{1})(s,q,\omega)|+|(\lra{q}\partial_q)\partial_\omega^n(\wt{V}_1)(s,q,\omega)| \\
&\lesssim (\exp(\delta/\eps)-2q+4R)^{-1+C\eps}+(\exp(\delta/\eps)-2q+4R)^{C\eps}\lesssim \lra{q}^{C\eps}.}
In the last inequality, we note  that $(a+b)^{C\eps}\leq 2^{C\eps}\max\{a,b\}^{C\eps}\leq 2(a^{C\eps}+b^{C\eps})$ for each pairs $a,b\geq 0$.
\end{proof}
\rmk{\rm\label{proplim3.1} Following the proof of Lemma \ref{proplim2}, we can show that $(A_1,\partial_\omega A_1)=\mathfrak{R}_{0,0}$, $\wt{V}_1-A_1=\mathfrak{R}_{-1,0}$, $(A_2,\partial_\omega A_2)=\mathfrak{R}_{0,-1}$ and $\wt{V}_2-A_2=\mathfrak{R}_{-1,0}$. }\rm

\bigskip

Moreover, we note that $A_1\approx -2$ in the following sense.
\lem{\label{lema1-2} Fix $0<\kappa<1$. For $\eps\ll 1$ and for all $(q,\omega)\in\R\times\mathbb{S}^2$, we have $|A_1(q,\omega)+2|\leq \kappa\lra{q}^{-1+C\eps}$. The constant in the power may depend on $\kappa$. As a result, we have $A_1(q,\omega)<-1<0$.}
\begin{proof}
Since $A_1\equiv -2$ for $q>R$, we can assume $q<2R$ in the proof. Recall from the proof of Proposition \ref{proplim3} that
\fm{e_4(\wt{V}_1)&=\eps\mathfrak{R}_{-2,0}+\mathfrak{R}_{-2,-1}=O(\eps t^{-2+C\eps}+t^{-2+C\eps}\lra{q}^{-1}).}

Next we consider $\wt{V}_1|_H$. On $H$ we have $\mu=-2+O(|u|)=-2+O(\eps t^{-1+C\eps})$. As computed in Lemma \ref{proplim2}, on $H$ we have \fm{|(\exp(\frac{1}{2}GAs)-1)\mu|&\lesssim \lra{q}^{-1+C\eps}s\exp(C\lra{q}^{-1+C\eps}s)\cdot (2+O(\eps t^{-1+C\eps}))\\
&\lesssim \lra{q}^{-1+C\eps}s\exp(C\lra{q}^{-1+C\eps}s).} 
Thus, $\wt{V}_1|_H=-2+O(\eps t^{-1+C\eps}+\lra{q}^{-1+C\eps}s\exp(C\lra{q}^{-1+C\eps}s))$.

We integrate $e_4(\wt{V}_1)$ along the geodesic in $\mcl{A}$ passing through $(t,x)\in\Omega\cap\{r-t<2R\}$. Then, \fm{|\wt{V}_1(t,x)+2|&\lesssim \eps(x^0(0))^{-1+C\eps}+\lra{q}^{-1+C\eps}(\eps\ln x^0(0)-\delta)\exp(C\lra{q}^{-1+C\eps}(\eps\ln x^0(0)-\delta))\\&\hspace{2em}+(\eps+\lra{q}^{-1})\int_{x^0(0)}^t\tau^{-2+C\eps}\ d\tau\\
&\lesssim \eps(x^0(0))^{-1+C\eps}+\lra{q}^{-1+C\eps}(\eps\ln x^0(0)-\delta)\exp(C\lra{q}^{-1+C\eps}(\eps\ln x^0(0)-\delta))\\&\hspace{2em}+(\eps+\lra{q}^{-1})(x^0(0))^{-1+C\eps}.}
If $\eps\ln x^0(0)-\delta\leq c$ for some small constant $c>0$, we have
\fm{|\wt{V}_1(t,x)+2|&\leq C\eps\lra{q}^{-1+C\eps}+Cc\lra{q}^{-1+C\eps}\exp(Cc\lra{q}^{-1+C\eps})+C(\eps+\lra{q}^{-1})(\lra{q}+\exp(\delta/\eps))^{-1+C\eps}\\
&\leq C\eps \lra{q}^{-1+C\eps}+Cc\lra{q}^{-1+C\eps}.}
By choosing $c,\eps\ll_\kappa 1$, we can make $Cc+C\eps<\kappa$. Thus,
$|\wt{V}_1(t,x)+2|\leq \kappa\lra{q}^{-1+C\eps}$. If $\eps\ln(x^0(0))-\delta>c$, we have $x^0(0)>\exp((c+\delta)/\eps)$ and thus $q=(\exp(\delta/\eps)-x^0(0))/2+2R<-C^{-1}\exp((c+\delta)/\eps)$ for $\eps\ll 1$. Then we have $\lra{q}^{C'\eps}\geq C^{-C'\eps}\exp(C'(c+\delta))$ and thus \fm{|\wt{V}_1(t,x)+2|&\lesssim (\eps+\lra{q}^{-1})(x^0(0))^{-1+C\eps}+\lra{q}^{-1+C\eps}(x^0(0))^{C\eps}\\
&\lesssim (\eps+\lra{q}^{-1})\lra{q}^{-1}(\exp(\delta/\eps)+\lra{q})^{C\eps}+\lra{q}^{-1+C\eps}(\exp(\delta/\eps)+\lra{q})^{C\eps}\\&\lesssim \lra{q}^{-1+C\eps}\lesssim \lra{q}^{-1+(C+C')\eps}C^{C'\eps}\exp(-C'c).}
The second last inequality holds since $a^{C\eps}+b^{C\eps}\leq (2\max\{a,b\})^{C\eps}\leq 2^{C\eps}(a^{C\eps}+b^{C\eps})$ for $a,b>0$. By choosing $C'\gg_\kappa 1$ and $\eps\ll_\kappa 1$, again we have $|\wt{V}_1(t,x)+2|\leq \kappa\lra{q}^{-1+C\eps}$.

We finish the proof by sending $s\to\infty$.
\end{proof}\rm

\subsection{An exact solution to the reduced system}\label{ssdata2}
For each $(s,q,\omega)\in\R\times\R\times\mathbb{S}^2$, we define
\eq{\label{s5exactsol}\left\{\begin{array}{l}\displaystyle
\wt{\mu}(s,q,\omega)=A_1(q,\omega)\exp(-\frac{1}{2}G(\omega)A(q,\omega)s),\\[1em]
\displaystyle\wt{U}_q(s,q,\omega)=A_2(q,\omega)\exp(\frac{1}{2}G(\omega)A(q,\omega)s).
\end{array}\right.}
Since $\wt{\mu}\wt{U}_q=A_1A_2= -2A$, it is easy to show that $(\wt{\mu},\wt{U}_q)$ is indeed a solution to the reduced system \eqref{asyode}. To solve for $\wt{U}$ uniquely, we assume that $\lim_{q\to\infty}\wt{U}(s,q,\omega)=0$ (since $\lim_{q\to\infty}U(s,q,\omega)=0$). This also implies that $\wt{U}\equiv 0$ for $q\geq 2R$. At $(s,q,\omega)\in\Omega'\cap\{q<2R\}$ we have
\fm{
\wt{\mu}&=\mathfrak{R}_{0,0}\cdot(1+\mathfrak{R}_{0,-1})=\mathfrak{R}_{0,0},\hspace{2em}
\wt{U}_q&=\mathfrak{R}_{0,-1}(1+\mathfrak{R}_{0,0})=\mathfrak{R}_{0,-1},}
\fm{\wt{\mu}-\mu&=\exp(-\frac{1}{2}G(\omega)A(q,\omega)s)(A_1-\wt{V}_1)=\mathfrak{R}_{-1,0},\\
\wt{U}_q-U_q &=\exp(\frac{1}{2}G(\omega)A(q,\omega)s)(A_2-\wt{V}_2)=\mathfrak{R}_{-1,0}.}
Thus, for each $p,m,n$, we have
\eq{\label{sec5ef1}&\partial^p_s(\lra{q}\partial_q)^{m}\partial_\omega^{n}\wt{\mu}=\eps^{-p}(\sum_k\mathfrak{R}_{0,0}V_k)^{p+m+n}(\mathfrak{R}_{0,0})=O(\eps^{-p}t^{C\eps}),\\
&\partial^p_s(\lra{q}\partial_q)^{m}\partial_\omega^{n}\wt{U}_q=\eps^{-p}(\sum_k\mathfrak{R}_{0,0}V_k)^{p+m+n}(\mathfrak{R}_{0,-1})=O(\eps^{-p}\lra{q}^{-1}t^{C\eps}),\\&\partial^p_s(\lra{q}\partial_q)^{m}\partial_\omega^{n}(\wt{\mu}-\mu,\wt{U}_q-U_q) =\eps^{-p}(\sum_k\mathfrak{R}_{0,0}V_k)^{p+m+n}(\mathfrak{R}_{-1,0})=O(\eps^{-p}t^{-1+C\eps}).}
Moreover, since $U=\eps^{-1}ru=\mathfrak{R}_{0,0}$, we can also show that $\partial^p_s(\lra{q}\partial_q)^{m}\partial_\omega^{n}U=O(\eps^{-p}t^{C\eps})$. 
Now, by integrating  $\partial_s^p\partial_\omega^n(\wt{U}_q-U_q)$ with respect to $q$, we have
\eq{\label{sec5ef2}
\partial^p_s \partial_\omega^{n}( \wt{U}-U)=O(\eps^{-p}\lra{q}t^{-1+C\eps}),\hspace{2em} \partial^p_s \partial_\omega^{n} \wt{U}=O(\eps^{-p}\lra{q}t^{-1+C\eps}+\eps^{-p}t^{C\eps})=O(\eps^{-p}t^{C\eps}).}
Here we note that $\lra{q}\lesssim t$ in $\Omega'\cap\{q<2R\}$.
The estimates \eqref{sec5ef1} and \eqref{sec5ef2} will be used in Section \ref{sec7}.

\section{Gauge independence}\label{ss5.7}  At the beginning of Section \ref{s3}, we define a region $\Omega$ by \eqref{defomega} and then construct an optical function in $\Omega$. If we replace \eqref{defomega} with 
\fm{\Omega_{\kappa,\delta}:=\{(t,x):\ t>\exp(\delta/\eps),\ |x|-\exp(\delta/\eps)-2R>\kappa(t-\exp(\delta/\eps))\}}
for some fixed constants $\delta>0$ and $0<\kappa<1$, we are  still able to  construct an optical function in $\Omega_{\kappa,\delta}$ by following the proofs in Section \ref{s3} and Section \ref{s4}. We are also able to construct a scattering data by following the proofs in Section \ref{s5}. We do not expect that the scattering data to be independent of $(\kappa,\delta)$, but we have the next proposition.

\prop{\label{prop5.8gi}Suppose $q(t,x)$ and $\bar{q}(t,x)$ are two solutions to the same eikonal equation \fm{g^{\alpha\beta}(u)q_\alpha q_\beta=0} in different regions $\Omega_{\kappa,\delta}$ and $\Omega_{\bar{\kappa},\bar{\delta}}$, respectively, as constructed in Section \ref{s3} and Section \ref{s4}. Let $A(q,\omega)$ and $\bar{A}(\bar{q},\omega)$ be the corresponding scattering data constructed in Section \ref{ssdata}. Under the change of coordinates $(s,q,\omega)=(\eps\ln(t)-\delta,q(t,x),\omega)$, we can view $\bar{q}(t,x)$ as a function of $(s,q,\omega)$ in $\Omega_{\kappa,\delta}\cap\Omega_{\bar{\kappa}\cap\bar{\delta}}$. Then, the limit $\bar{q}_\infty(q,\omega):=\lim_{s\to\infty}\bar{q}(s,q,\omega)$ exists, and we have \fm{A(q,\omega)=\bar{A}(\bar{q}_\infty(q,\omega),\omega).}}
\begin{proof}We first recall several notations and estimates in Section \ref{s3}. For example, we have $\mu=q_t-q_r=O(t^{C\eps})$, $\nu=q_t+q_r=O(t^{-1+C\eps})$, and we have similar definitions and estimates for $\bar{\mu}$ and $\bar{\nu}$. By viewing $\bar{q}(t,x)$ as a function of $(s,q,\omega)=(\eps\ln(t)-\delta,q(t,x),\omega)$, we have
\fm{\partial_s\bar{q}&=\eps^{-1}t(\bar{q}_t-q_tq_r^{-1}\bar{q}_r)=t\eps^{-1}\bar{q}_r(\bar{\nu}\bar{q}_r^{-1}-\nu q_r^{-1}).}
By the eikonal equation, we have
\fm{0&=-(q_r-q_r)(q_r+q_r)+O(t^{-2+C\eps})+(g^{\alpha\beta}(u)-m^{\alpha\beta})q_\alpha q_\beta=-\nu\mu+\frac{1}{4}uG(\omega)\mu^2+O(t^{-2+C\eps}).}
Since $\mu\leq -C^{-1}t^{-C\eps}$, we have
\fm{\nu=\frac{1}{4}uG(\omega)\mu+O(t^{-2+C\eps})}
and thus
\fm{\frac{\nu}{q_r}&=\frac{1}{4}uG(\omega)\frac{\mu}{q_r}+O(t^{-2+C\eps})=\frac{1}{4}uG(\omega)(\frac{\nu}{q_r}-2)+O(t^{-2+C\eps})=-\frac{1}{2}uG(\omega)+O(t^{-2+C\eps}).}
We conclude that
\fm{\partial_s\bar{q}&=t\eps^{-1}\bar{q}_r^{-1}(-\frac{1}{2}uG(\omega)+O(t^{-2+C\eps})-(-\frac{1}{2}uG(\omega)+O(t^{-2+C\eps})))\\
&=O(\eps^{-1}t^{-1+C\eps})=O(\eps^{-1}\exp((-\eps^{-1}+C)(s+\delta))).}
As computed in Section \ref{ssdata}, we can show that $\bar{q}_\infty(q,\omega):=\lim_{s\to\infty}\bar{q}(s,q,\omega)$ exists for all $(q,\omega)$. Moreover, we can show that 
\fm{|\bar{q}(s,q,\omega)-\bar{q}_\infty(q,\omega)|\lesssim t^{-1+C\eps}.}

Since $\lim_{s\to\infty}(\mu U_q)(s,q,\omega)=-2A(q,\omega)$ and $\lim_{\bar{s}\to\infty}(\bar{\mu}\emph{\textbf{}} \bar{U}_q)(\bar{s},\bar{q},\omega)=-2\bar{A}(\bar{q},\omega)$ (recall that $\bar{s}+\bar{\delta}=s+\delta$), we have
\fm{\partial_r(\eps^{-1}ru)&=q_rU_q=-\frac{1}{2}\mu U_q+O(t^{-1+C\eps});\hspace{2em}
\partial_r(\eps^{-1}ru)=\bar{q}_r\bar{U}_{\bar{q}}=-\frac{1}{2}\bar{\mu} \bar{U}_{\bar{q}}+O(t^{-1+C\eps}).}
Then,
\fm{(\mu U_q)(s,q,\omega)&=(\bar{\mu} \bar{U}_q)(s+\delta-\bar{\delta},\bar{q}(s,q,\omega),\omega)+O(t^{-1+C\eps}).}
By sending $s$ (and thus $t$) to infinity, we conclude that $A(q,\omega)=\bar{A}(\bar{q}_\infty(q,\omega),\omega)$.
\end{proof}
\rm

\section{Approximation}\label{sec7}
Recall that we have constructed an exact solution to our reduced system in \eqref{s5exactsol}. In this section, we seek to prove that this exact solution gives a good approximation of the exact solution to \eqref{qwe}.

To state the result, we first recall the solution $(\wt{\mu},\wt{U})(s,q,\omega)$ to the reduced system defined in Proposition \ref{props5}, or in \eqref{s5exactsol}. We now solve
\fm{\wt{q}_t-\wt{q}_r=\wt{\mu}(\eps\ln(t)-\delta,\wt{q}(t,x),\omega)\quad\text{in }\Omega\cap\{r-t<2R\};\hspace{2em}\wt{q}=r-t\quad\text{when } r-t\geq 2R}
and set \fm{\wt{u}(t,x)=\eps r^{-1}\wt{U}(\eps\ln(t)-\delta,\wt{q}(t,x),\omega)\quad\text{in }\Omega\cap\{r-t<2R\}.}
We remark that the construction here is very similar to that in Section 4 of the author's previous paper \cite{yu2020}. We then have the following approximation result.

\prop{\label{props6} The function $\wt{u}=\wt{u}(t,x)$ is an approximate solution to \eqref{qwe} in the following sense:
\fm{|Z^I(g^{\alpha\beta}(\wt{u})\partial_{\alpha}\partial_\beta\wt{u})(t,x)|\lesssim \eps t^{-3+C\eps},\hspace{2em}\forall (t,x)\in\Omega,\ \forall I.}
Moreover, if we fix a constant $0<\gamma<1$ and a large integer $N$, then for $\eps\ll_{\gamma,N} 1$, at each $(t,x)\in\Omega$ such that $|r-t|\lesssim  t^{\gamma}$, we have $|Z^I(u-\wt{u})|\lesssim_\gamma \eps t^{-2+C\eps}\lra{r-t}$ for each $|I|\leq N$. }\rm\bigskip

The estimates for $u-\wt{u}$ in this proposition is better than the estimates for $u$ itself.

After making several definitions in Section \ref{sec6.1}, we introduce a simplification in Section \ref{sec6.2}. Instead of $(\wt{\mu},\wt{U}_q)$, the simplification in Section \ref{sec6.2} allows us work with $(\hat{\mu},\hat{U}_q)$ which is an exact solution to the reduced system \eqref{s5exactsol} with initial data $(-2,\hat{A})$. We thus get a new function $\hat{q}$  which is a solution to $\hat{q}_t-\hat{q}_r=\hat{\mu}$. In Section \ref{sec6.3}, we follow Section 4 of \cite{yu2020} to prove several estimates for $\hat{q}$ and $\hat{U}$. The most important result here is Proposition \ref{prop6.3main} which states that $\wt{u}=\hat{u}$ is indeed an approximate solution to \eqref{qwe}. In Section \ref{sec6.4}, we show that $\hat{q}$ approximates the optical function $q$ in a certain sense. Finally, in Section \ref{sec6.5}, we make use of the estimates in Section \ref{sec6.4} to prove Proposition \ref{props6}.

\subsection{Definitions}\label{sec6.1}
We first define a function $\wt{q}(t,x)$ in $\Omega$ by solving the following equation
\eq{\label{wtqeqn} \wt{q}_t-\wt{q}_r=\wt{\mu}(\eps\ln(t)-\delta,\wt{q}(t,x),\omega)\quad\text{in }\Omega\cap\{r-t<2R\};\hspace{2em}\wt{q}=r-t\quad\text{when } r-t\geq 2R.}
Recall that $\wt{\mu}$ is defined by\fm{\wt{\mu}(s,q,\omega):=A_1(q,\omega)\exp(-\frac{1}{2}G(\omega)A(q,\omega)s),\hspace{2em}\forall (s,q,\omega)\in\R\times\R\times\mathbb{S}^2.} In this section, when we write $q$, we usually mean a variable instead of the optical function $q(t,x)$.

As in \cite{yu2020}, we can use the method of characteristics to solve \eqref{wtqeqn}. We fix $(t,x)\in \Omega\cap\{r-t<2R\}$ and  set $z(\tau):=\wt{q}(\tau,r+t-\tau,\omega)$. Then, the function $z(\tau)$ is a solution to the autonomous system of ODE's
\fm{\dot{z}(\tau)=\wt{\mu}(\eps s(\tau)-\delta,z(\tau),\omega),\hspace{2em}\dot{s}(\tau)=\eps\tau^{-1}.}
The initial data is given by $(z,s)((r+t)/2-R)=(2R,\eps\ln((r+t)/2-R)-\delta)$. By Proposition \ref{proplim}, Proposition \ref{proplim3} and Lemma \ref{lema1-2}, we have $|A_1+2|=O(\lra{q}^{-1+C\eps})$, $(A_2,A)(q,\omega)=O(\lra{q}^{-1+C\eps})$ and $A_1<-1$ for all $(q,\omega)$. Thus,
\fm{0&\geq \mu(\eps s(\tau)-\delta,z(\tau),\omega)=A_1(z(\tau),\omega)\exp(-\frac{1}{2}G(\omega)A(z(\tau),\omega)(\eps s(\tau)-\delta))\\
&\geq -C\tau^{C\eps\lra{z(\tau)}^{-1+C\eps}}\geq -C\tau^{C\eps}.}
Then, $-C\tau^{C\eps}\leq\dot{z}(\tau)\leq 0$, so $|z(\tau)|$ cannot blow up in finite time. By the Picard's theorem, the system of ODE's above has a solution for all $ (r+t)/2-R\leq \tau<\frac{1}{3}(2(r+t)-4R-\exp(\delta/\eps))$. The upper bound here guarantees that $(\tau,r+t-\tau,\omega)\in\Omega$. Thus, \eqref{wtqeqn} has a solution $\wt{q}(t,x)$ in $\Omega$.

Next, we  define $\wt{U}(s,q,\omega)$ by\eq{\label{wtueqn}\wt{U}(s,q,\omega)&=-\int_q^\infty A_2(p,\omega)\exp(\frac{1}{2}G(\omega)A(p,\omega)s)\ dp.}Note that $A_2(q,\omega)=0$ whenever $q>R$, so when $q<R$, we can replace $\infty$ with $R$ in \eqref{wtueqn}. In $\Omega$ we set  \fm{\wt{u}(t,x)=\eps r^{-1}\wt{U}(\eps\ln(t)-\delta,\wt{q}(t,x),\omega).}
We seek to prove  that $\wt{u}(t,x)$ provides a good approximation of $u(t,x)$. 

\subsection{Simplification}\label{sec6.2}
We aim to introduce some simplification in this subsection. Define a new function $F(q,\omega)$ on $\R\times\mathbb{S}^2$ by  
\fm{F(q,\omega):=2R-\int^{q}_{2R} \frac{2}{A_1(p,\omega)}\ dp.}
Then, we have
\begin{enumerate}[a)]
\item $F$ is defined everywhere, and $2(q-R)\leq F(q,\omega)\leq 2(q+R)/3$ for all $q<2R$.  This is because  $A_1\in[-3,-1]$ by Lemma \ref{lema1-2}.
\item $F$ is a smooth function of $(q,\omega)$, in the sense that for each large integer $N$   and $\eps\ll_N1$, $F$ is in $C^N$. This is because $A_1\in[-3,-1]$ and by Proposition \ref{proplim3}.
\item $F(q,\omega)=q$ for $q>R$, and $\lra{F(q,\omega)}\sim\lra{q}$. This is because $A_1\equiv -2$ for $q>R$.
\item For each fixed $\omega$, the map $q\mapsto F(q,\omega)$ has an inverse denoted by $\hat{F}(q,\omega)$ which is also smooth (in the same sense as in a) above) in $\R\times\mathbb{S}^2$. This is because $F_q=-2/A_1\in[2/3,2]$.
\item  $\partial_q^a\partial_\omega^cF=O(\lra{q}^{1-a+C\eps})$. Recall that $A_1<-1$ and $\partial_q^a\partial_\omega^cA_1=O(\lra{q}^{-a+C\eps})$. If $a=0$, then  $|\partial_\omega^cF|\lesssim\int_{[q,2R]}\lra{p}^{C\eps}\ dp\lesssim \lra{q}^{1+C\eps}$. If $a\geq 1$, then $|\partial_q^a\partial_\omega^cF|=|\partial_q^{a-1}\partial_{\omega}^c(2/A_1)|\lesssim\lra{q}^{1-a+C\eps}$.
\end{enumerate}

For each $(s,q,\omega)$, we set
\fm{\hat{A}(q,\omega):=A(\hat{F}(q,\omega),\omega)}
and
\eq{\label{hatmuudef}\left\{\begin{array}{l}\displaystyle\hat{\mu}(s,q,\omega):=-2\exp(-\frac{1}{2}G(\omega)\hat{A}(q,\omega)s),\\[1em]
\displaystyle\hat{U}(s,q,\omega):=-\int_q^\infty\hat{A}(p,\omega)\exp(\frac{1}{2}G(\omega)\hat{A}(p,\omega)s)\ dp.
\end{array}\right.}
It is clear that $(\hat{\mu},\hat{U})$ is a solution to the reduced system \eqref{asyode}.

For each $(t,x)\in\Omega$, we set
\fm{\hat{q}(t,x):=F(\wt{q}(t,x),\omega),\hspace{2em}\hat{u}(t,x):=\eps r^{-1}\hat{U}(\eps\ln t-\delta,\hat{q}(t,x),\omega).}
We then have the next key lemma.

\lem{\label{lemsimplif} In $\Omega$, we have
\fm{\hat{q}_t-\hat{q}_r&=\hat{\mu}(\eps\ln t-\delta,\hat{q}(t,x),\omega)}
and $\hat{q}=r-t$ whenever $r-t>R$.  Moreover, we have $\hat{u}(t,x)=\wt{u}(t,x)$ everywhere.}
\begin{proof}
At $(t,x)\in\Omega$, we first have  \fm{\wt{q}(t,x)=\hat{F}(F(\wt{q}(t,x),\omega),\omega)=\hat{F}(\hat{q}(t,x),\omega).} Thus,
\fm{\hat{q}_t-\hat{q}_r&=(\partial_t-\partial_r)F(\wt{q}(t,x),\omega)=F_q(\wt{q}(t,x),\omega)\cdot\wt{\mu}(\eps\ln t-\delta,\wt{q}(t,x),\omega)\\
&=(-2/A_1\cdot A_1\exp(-\frac{1}{2}GAs))(\eps\ln t-\delta,\wt{q}(t,x),\omega)\\
&=-2\exp(-\frac{1}{2}G(\omega)A(\wt{q}(t,x),\omega)(\eps\ln t-\delta))\\
&=-2\exp(-\frac{1}{2}G(\omega)A(\hat{F}(\hat{q}(t,x),\omega),\omega)(\eps\ln t-\delta))\\
&=-2\exp(-\frac{1}{2}G(\omega)\hat{A}(\hat{q}(t,x),\omega)(\eps\ln t-\delta))=\hat{\mu}(\eps\ln t-\delta,\hat{q}(t,x),\omega).}
Since $F(q,\omega)=q$ for all $q>R$, we have $\hat{q}(t,x)=\wt{q}(t,x)=r-t$ whenever $r-t>R$.

Moreover, if $\rho=\hat{F}(p,\omega)$, then we have $p=F(\rho,\omega)$ and thus
\fm{A(\rho,\omega)&=A(\hat{F}(p,\omega),\omega)=\hat{A}(p,\omega).}
Then by the change of variables ($\rho=\hat{F}(p,\omega)$ and thus $p=F(\rho,\omega)$), we have
\fm{\hat{U}(s,\hat{q},\omega)&=-\int_{\hat{q}}^\infty\hat{A}(p,\omega)\exp(\frac{1}{2}G(\omega)\hat{A}(p,\omega)s)\ dp\\
&=-\int_{\wt{q}}^\infty A(\rho,\omega)\exp(\frac{1}{2}G(\omega)A(\rho,\omega)s)F_\rho(\rho,\omega)\ d\rho\\
&=-\int_{\wt{q}}^\infty A_2(\rho,\omega)\exp(\frac{1}{2}G(\omega)A(\rho,\omega)s)\ d\rho=\wt{U}(s,\wt{q},\omega).}
Here we note that $AF_q=-2A/A_1=A_2$.  That is, for each $(s,q,\omega)$ (not viewed as functions of $(t,x)$),
\eq{\label{lem7.2f1} \hat{U}(s,q,\omega)=\wt{U}(s,\hat{F}(q,\omega),\omega).}
We thus have $\wt{u}(t,x)=\hat{u}(t,x)$.
\end{proof}\rm

Because of Lemma \ref{lemsimplif}, we can work with $(\hat{u},\hat{q})$ instead of $(\wt{u},\wt{q})$.

We end this subsection with several useful estimates for $(\hat{A},\hat{\mu},\hat{U})$.
\prop{\label{proplim4} For each $(q,\omega)$, we have \fm{(\lra{q}\partial_q)^a\partial_\omega^c\hat{F}(q,\omega)=O(\lra{q}^{1+C\eps}),\hspace{2em}(\lra{q}\partial_q)^a\partial_\omega^c\hat{A}(q,\omega)=O(\lra{q}^{-1+C\eps}).} Besides, for each $(s,q,\omega)\in \Omega'\cap\{q<2R\}$, we have 
\fm{\partial_s^b(\lra{q}\partial_q)^a\partial_\omega^c\hat{U}=O(\eps^{-b}t^{C\eps}),\hspace{2em}\partial_s^b(\lra{q}\partial_q)^{a+1}\partial_\omega^c\hat{U}=O(t^{C\eps});}
\fm{\hat{\mu}=O(t^{C\eps}),\hspace{2em}\partial^b_s(\lra{q}\partial_q)^a\partial_\omega^c\hat{\mu}=O(\lra{q}^{-1+C\eps}t^{C\eps}|\hat{\mu}|),\hspace{1em}a+b+|c|>0.}}
\begin{proof}First, it is clear that $\lra{\hat{F}(q,\omega)}\sim \lra{q}$ and that $\hat{F}_q(q,\omega)=1/(F_q(\hat{F}(q,\omega),\omega))=-A_1(\hat{F}(q,\omega),\omega)/2\sim \lra{q}^{C\eps}$. In general we induct on $m+|n|$. By differentiating $q=F(\hat{F}(q,\omega),\omega)$, for $(a,c)\notin \{(0,0),(1,0)\}$, we have 
\fm{0&=F_q(\hat{F}(q,\omega),\omega)\cdot \partial_q^a\partial_\omega^c\hat{F}(q,\omega)+\sum C [(\partial_q^{m}\partial_\omega^{c'}F)(\hat{F}(q,\omega),\omega)\cdot\prod_{j=1}^m(\partial_q^{a_j}\partial_{\omega}^{c_j}\hat{F})(q,\omega)].}
Here  the sum on the right hand side is taken over all $(m,c',a_*,c_*)$ such that $\sum a_j=a$, $c'+\sum c_j=c$, $a_j+|c_j|<a+|c|$. We can now apply the induction hypotheses to conclude that 
\fm{0&=F_q(\hat{F}(q,\omega),\omega)\cdot \partial_q^a\partial_\omega^c\hat{F}(q,\omega)+\sum O(\lra{\hat{F}(q,\omega)}^{1-m+C\eps}\cdot \lra{\hat{q}}^{m-\sum a_j+C\eps})\\
&=F_q(\hat{F}(q,\omega),\omega)\cdot \partial_q^a\partial_\omega^c\hat{F}(q,\omega)+ O(\lra{q}^{1-a+C\eps}).}
And since $F_q\sim 1$, we conclude that $\partial_q^a\partial_\omega^c\hat{F}(q,\omega)=O(\lra{q}^{1-a+C\eps})$.   

Next, recall that 
\fm{\hat{A}(q,\omega)=A(\hat{F}(q,\omega),\omega),\hspace{2em}\hat{U}(s,q,\omega)=\wt{U}(s,\hat{F}(q,\omega),\omega).}
Then, $\partial_s^b\partial_q^a\partial_\omega^c\hat{U}(s,q,\omega)$ is a linear combination of terms of the form
\fm{\partial_s^b \partial_q^{m}\partial_\omega^{c'}\wt{U}(s,\hat{F}(q,\omega),\omega)\cdot\prod_{j=1}^m\partial_{q}^{a_j}\partial_\omega^{c_j}\hat{F}(q,\omega),\hspace{2em}\sum a_j=a,\ c'+\sum c_j=c.}
By \eqref{sec5ef1} and \eqref{sec5ef2}, we conclude that each of these terms are controlled by
\fm{\eps^{-b}\lra{\hat{F}(q,\omega)}^{-m}t^{C\eps}\cdot \lra{q}^{m-\sum a_j+C\eps}\lesssim \eps^{-b}\lra{q}^{-a}t^{C\eps}.} Thus, $\partial_s^b(\lra{q}\partial_q)^a\partial_\omega^c\hat{U}(s,q,\omega)=O(\eps^{-b}t^{C\eps})$.
Following the same proof, we can show that $(\lra{q}\partial_q)^a\partial_\omega^c\hat{A}(q,\omega)=O(\lra{q}^{-1+C\eps})$.

Finally, by \eqref{hatmuudef}, we can write $\partial_s^b\partial_q^a\partial_\omega^c\hat{U}_q(s,q,\omega)$ as a linear combination of terms of the form
\fm{\partial_q^{a'}\partial_{\omega}^{c'}\hat{A}(q,\omega)\cdot\exp(\frac{1}{2}G\hat{A}s)\prod_{j=1}^m \partial_s^{b_j}\partial_q^{a_j}\partial^{c_j}_\omega(\frac{1}{2}G\hat{A}s)} where $a'+\sum a_j=a,\ \sum b_j=b,\ c'+\sum c_j=c$. Each of these terms are controlled by
\fm{\lra{q}^{-1-a'+C\eps}\cdot t^{C\eps}\cdot \lra{q}^{-m-\sum a_j}t^{C\eps}\lesssim \lra{q}^{-1-a}t^{C\eps}.}
In conclusion, we have $\partial_s^b(\lra{q}\partial_q)^{a+1}\partial_\omega^c\hat{U}(s,q,\omega)=O(t^{C\eps})$. Here we do not have the factor $\eps^{-b}$ which is better. Moreover, we have $\hat{\mu}=O(t^{C\eps})$ and 
\fm{(\hat{\mu}_s,\lra{q}\hat{\mu}_q,\hat{\mu}_\omega)=-\frac{1}{2}(GA,\lra{q}GA_qs,\partial_\omega(GA)s)\hat{\mu}.} Following the same proof, we can show that $\partial^b_s(\lra{q}\partial_q)^a\partial_\omega^c\hat{\mu}(s,q,\omega)=O(\lra{q}^{-1+C\eps}t^{C\eps}|\hat{\mu}|)$ if $a+b+|c|>0$.
\end{proof}\rm

\subsection{Estimates for $\hat{q}$ and $\hat{U}$}\label{sec6.3}
We now follow Section 4 in \cite{yu2020} to prove several useful estimates. In this subsection, all functions of $(s,q,\omega)\in[0,\infty)\times\R\times\mathbb{S}^2$ are viewed as functions of $(t,x)\in \Omega$ by setting $(s,q,\omega)=(\eps\ln t-\delta,\hat{q}(t,x),\omega)$. This setting is different from that in the previous sections of this paper, where we take $q=q(t,x)$.
\lem{\label{wtqlem1} In $\Omega\cap\{r-t<2R\}$, we have $\lra{\hat{q}}/\lra{r-t}=t^{O(\eps)}$ and $\hat{q}(t,x)-r+t=O(\min\{\eps^{-1},\lra{\hat{q}}\}t^{C\eps})$.}
\begin{proof}
Fix $(t,x)\in\Omega\cap\{r-t<2R\}$. Then, we have 
\fm{|\hat{q}(t,x)-2R|&=\int_{(r+t)/2-R}^t(-\hat{\mu}(\eps\ln \tau-\delta,\hat{q}(\tau,r+t-\tau,\omega),\omega))\ d\tau\\
&\lesssim \int_{(r+t)/2-R}^t\exp(C\lra{\hat{q}}^{-1+C\eps}s)(\tau,r+t-\tau,\omega)\ d\tau\\
&\lesssim ((r-t)/2+R)t^{C\eps}\lesssim \lra{r-t}t^{C\eps};}
\fm{|\hat{q}(t,x)-2R|&=\int_{(r+t)/2-R}^t(-\hat{\mu}(\eps\ln \tau-\delta,\hat{q}(\tau,r+t-\tau,\omega),\omega))\ d\tau\\
&\gtrsim \int_{(r+t)/2-R}^t \exp(-C\lra{\hat{q}}^{-1+C\eps}s)(\tau,r+t-\tau,\omega)\ d\tau\\
&\gtrsim ((r-t)/2+R)t^{-C\eps}\gtrsim \lra{r-t}t^{-C\eps}.}
Thus, we have $t^{-C\eps}\lra{\hat{q}}\lesssim\lra{r-t}\lesssim t^{C\eps}\lra{\hat{q}}$. It follows that
\fm{|\hat{q}(t,x)-(r-t)|\leq |\hat{q}-2R|+|r-t-2R|\lesssim t^{C\eps}\lra{\hat{q}}+\lra{r-t}\lesssim \lra{\hat{q}}t^{C\eps}.}

To improve the estimate above, we note that
\fm{\hat{q}(t,x)&=2R+\int_{(r+t)/2-R}^t\hat{\mu}(\eps\ln \tau-\delta,\hat{q}(\tau,r+t-\tau,\omega),\omega)\ d\tau\\
&=r-t+\int_{(r+t)/2-R}^t(\hat{\mu}(\eps\ln \tau-\delta,\hat{q}(\tau,r+t-\tau,\omega),\omega)+2)\ d\tau.}
For each $(s,q,\omega)\in[0,\infty)\times\R\times\mathbb{S}^2$, by Proposition \ref{proplim2} and Lemma \ref{lema1-2} we have
\fm{&\hspace{1.5em}|\hat{\mu}(s,q,\omega)+2|\lesssim|1-\exp(-\frac{1}{2}GAs)|\lesssim \lra{q}^{-1+C\eps}|s|\exp(C\lra{q}^{-1+C\eps}s).}
By setting $(s,q,\omega)=(\eps\ln \tau-\delta,\hat{q}(\tau,r+t-\tau,\omega),\omega)$, we have \fm{|\hat{\mu}+2|(\tau)\lesssim \lra{r+t-2\tau}^{-1+C\eps}\tau^{C\eps}\lesssim (3R-r-t+2\tau)^{-1+C\eps}t^{C\eps}}
and then
\fm{|\hat{q}-r+t|\lesssim t^{C\eps}\int_{(r+t)/2-R}^t(3R-r-t+2\tau)^{-1+C\eps}\ d\tau\lesssim \eps^{-1}t^{C\eps}(3R-r+t)^{C\eps}.}
And since $0\leq 3R-r+t\lesssim 1+t\lesssim t$, we have $|\hat{q}-r+t|\lesssim \eps^{-1}t^{C\eps}$.
\end{proof}\rm

\lem{\label{lem7.5}In $\Omega$ we have
\fm{\hat{\nu}:=\hat{q}_t+\hat{q}_r=O( t^{-1+C\eps}),\hspace{2em}\hat{\lambda}_i:=\hat{q}_i-\omega_i\hat{q}_r=O((1+\ln\lra{r-t})t^{-1+C\eps}).}
It follows that $\hat{q}_r=(\hat{\nu}-\hat{\mu})/2>C^{-1}t^{-C\eps}$ and $\hat{q}_t=(\hat{\nu}+\hat{\mu})/2<-C^{-1}t^{-C\eps}$. Thus, for each fixed $(t,\omega)$ the function $r\mapsto \hat{q}(t,r\omega)$ is continuous and strictly increasing.}
\begin{proof}There is nothing to prove when $r-t>R$. Fix $(t,x)\in\Omega\cap\{r-t<2R\}$. Then,
\fm{(\partial_t-\partial_r)\hat{\nu}&=(\partial_t+\partial_r)\hat{\mu}=\hat{\mu}_q\hat{\nu}+\eps t^{-1}\hat{\mu}_s=\hat{\mu}_q\hat{\nu}-\frac{\eps}{2t}G(\omega)A(\hat{q},\omega)\hat{\mu}\\
&=-\frac{1}{2}G\hat{A}_qs\hat{\mu}\hat{\nu}-\frac{\eps}{2t}G\hat{A}\hat{\mu}.}
By setting $z(\tau):=\hat{q}(\tau,r+t-\tau,\omega)$, we have $\dot{z}=\hat{\mu}<0$ and thus
\fm{&\hspace{1.5em}\int_{(r+t)/2-R}^t|G\hat{A}_qs\hat{\mu}|(\tau,r+t-\tau,\omega)\ d\tau\lesssim \int_{(r+t)/2-R}^t(\eps\ln\tau+1)\lra{\hat{q}}^{-2+C\eps}(-\hat{\mu})\ d\tau\\
&\lesssim(\eps\ln t+1)\int_{(r+t)/2-R}^t\lra{z}^{-2+C\eps}(-\dot{z})\ d\tau\lesssim\eps\ln t+1,\\
&\hspace{1.5em}\int_{(r+t)/2-R}^t|\eps\tau^{-1}G\hat{A}\hat{\mu}|(\tau,r+t-\tau,\omega)\ d\tau\lesssim \eps ((r+t)/2-R)^{-1}\int_{(r+t)/2-R}^t\lra{\hat{q}}^{-1+C\eps}(-\hat{\mu})\ d\tau\\
&\lesssim \eps t^{-1}\int_{(r+t)/2-R}^t\lra{z}^{-1+C\eps}(-\dot{z})\ d\tau\lesssim t^{-1}\lra{\hat{q}}^{C\eps}\lesssim t^{-1+C\eps}.}
Here we note that $\lra{\hat{q}}\lesssim \lra{r-t}t^{C\eps}\lesssim t^{1+C\eps}$. Since $\hat{\nu}=0$ at $\tau=(r+t)/2-R$, by the Gronwall's inequality we conclude that $\hat{\nu}=O(t^{-1+C\eps})$.

Next, we have
\fm{(\partial_t-\partial_r)\hat{\lambda}_i&=(\partial_i-\omega_i\partial_r)\hat{\mu}+r^{-1}\hat{\lambda}_i=(\hat{\mu}_q+r^{-1})\hat{\lambda}_i+\sum_l(\partial_{\omega_l}\hat{\mu})(\partial_i\omega_l)\\
&=(\hat{\mu}_q+r^{-1})\hat{\lambda}_i-\frac{1}{2}\sum_l(\partial_{\omega_l}(G\hat{A}))(\eps\ln t-\delta)\hat{\mu} r^{-1}(\delta_{il}-\omega_i\omega_l)\\
&=(\hat{\mu}_q+r^{-1})\hat{\lambda}_i+O(\lra{\hat{q}}^{-1+C\eps}t^{-1+C\eps}|\hat{\mu}|).}
We have proved that $\int_{(r+t)/2-R}^t|\mu_q| \ d\tau\lesssim \eps \ln t+1$. Integrate along  the characteristic $(\tau,r+t-\tau,\omega)$ and we have
\fm{\int_{(r+t)/2-R}^t (r+t-\tau)^{-1}\ d\tau&=\ln\frac{(r+t)/2+R}{r}=O(1),}
\fm{\int_{(r+t)/2-R}^t \lra{\hat{q}}^{-1+C\eps}(-\hat{\mu})\tau^{-1+C\eps}\ d\tau&\lesssim \int_{(r+t)/2-R}^t \lra{\hat{q}}^{-1}(-\hat{\mu})\tau^{-1+C\eps}\ d\tau\\
&\lesssim t^{-1+C\eps}\int_{(r+t)/2-R}^t \lra{z}^{-1}(-\dot{z})\ d\tau\\
&\lesssim (1+\ln\lra{\hat{q}})t^{-1+C\eps}\lesssim (1+\ln\lra{r-t})t^{-1+C\eps}.}
Here note that $\lra{\hat{q}}\lesssim t^{1+C\eps}$ and $\ln\lra{\hat{q}}\lesssim \ln\lra{r-t}+C\eps\ln t$ in $\Omega\cap\{r-t<2R\}$.
Since $\hat{\lambda}_i=0$ at $\tau=(r+t)/2-R$, by Gronwall's inequality we conclude that  $\hat{\lambda}_i=O((1+\ln\lra{r-t})t^{-1+C\eps})$.
\end{proof}\rm

\lem{\label{lem7.6}In $\Omega$, we have
\fm{\hat{\nu}=\frac{\eps G(\omega)}{4t}\hat{\mu} \hat{U}+O(\eps t^{-2+C\eps}\lra{r-t}),\hspace{1.5em} \hat{\nu}_q=\frac{\eps G(\omega)}{4t} (\hat{\mu} \hat{U}_q+\hat{\mu}_q \hat{U})+O(\eps(1+\ln\lra{r-t})t^{-2+C\eps}).}}
\begin{proof}We have
\fm{&\hspace{1.5em}(\partial_t-\partial_r)(\hat{\nu}-\frac{\eps G(\omega)}{4t}\hat{\mu} \hat{U})\\
&=\hat{\mu}_q\hat{\nu}-\frac{\eps}{2t}G\hat{A}\hat{\mu}+\frac{\eps G\hat{\mu}\hat{U}}{4t^2}-\frac{\eps G}{4t}(\hat{\mu}_q\hat{U}+\hat{\mu}\hat{U}_q)\hat{\mu}-\frac{\eps G}{4t}(\hat{\mu}_s\hat{U}+\hat{\mu}\hat{U}_s)\eps t^{-1}\\
&=\hat{\mu}_q\hat{\nu}-\frac{\eps}{2t}G\hat{A}\hat{\mu}+\frac{\eps G\hat{\mu}\hat{U}}{4t^2}-\frac{\eps G}{4t}(\hat{\mu}_q\hat{U}-2\hat{A})\hat{\mu}-\frac{\eps G}{4t}(-\frac{1}{2}G\hat{A}\hat{\mu}\hat{U}+\hat{\mu}\hat{U}_s)\eps t^{-1}\\
&=\hat{\mu}_q(\hat{\nu}-\frac{\eps G}{4t}\hat{\mu} \hat{U})+\frac{\eps G\hat{\mu}\hat{U}}{4t^2}-\frac{\eps^2 G}{4t^2}(-\frac{1}{2}G\hat{A}\hat{U}+\hat{U}_s)\hat{\mu}.}
Since $\hat{U}=O(t^{C\eps})$ and $\hat{U}_s=O(\eps^{-1}t^{C\eps})$ by Proposition \ref{proplim4}, we have
\fm{|\frac{\eps G\hat{\mu}\hat{U}}{4t^2}-\frac{\eps^2 G}{4t^2}(-\frac{1}{2}G\hat{A}\hat{U}+\hat{U}_s)\hat{\mu}|&\lesssim \eps t^{-2+C\eps}.}
Besides, we have
\fm{\int_{(r+t)/2-R}^t\eps\tau^{-2+C\eps}&\lesssim ((r+t)/2-R)^{-2+C\eps}\cdot \eps((t-r)/2-R)\lesssim \eps t^{-2+C\eps}\lra{r-t}.}
And since $\hat{\nu}-\frac{\eps G}{4t}\hat{\mu} \hat{U}=0$ at $\tau=(r+t)/2-R$, by Gronwall's inequality we conclude that \fm{\hat{\nu}-\frac{\eps G}{4t}\hat{\mu} \hat{U}=O(\eps t^{-2+C\eps}\lra{r-t}).}

Next, we have
\fm{&\hspace{1.5em}(\partial_t-\partial_r)\partial_r(\hat{\nu}-\frac{\eps G(\omega)}{4t}\hat{\mu} \hat{U})=\partial_r(\partial_t-\partial_r)(\hat{\nu}-\frac{\eps G(\omega)}{4t}\hat{\mu} \hat{U})\\
&=\partial_r(\hat{\mu}_q(\hat{\nu}-\frac{\eps G}{4t}\hat{\mu} \hat{U})+\frac{\eps G\hat{\mu}\hat{U}}{4t^2}-\frac{\eps^2 G}{4t^2}(-\frac{1}{2}G\hat{A}\hat{U}+\hat{U}_s)\hat{\mu})\\
&=\hat{\mu}_q\partial_r(\hat{\nu}-\frac{\eps G}{4t}\hat{\mu} \hat{U})+\hat{q}_r\hat{\mu}_{qq}(\hat{\nu}-\frac{\eps G}{4t}\hat{\mu} \hat{U})+\frac{\eps G\hat{q}_r\partial_q(\hat{\mu}\hat{U})}{4t^2}\\
&\hspace{1em}-\frac{\eps^2 G}{4t^2}(-\frac{1}{2}G\hat{A}\hat{U}+\hat{U}_s)\hat{\mu}_q\hat{q}_r-\frac{\eps^2 G}{4t^2}(-\frac{1}{2}G\partial_q(\hat{A}\hat{U})+\hat{U}_{sq})\hat{\mu}\hat{q}_r.}
By  Proposition \ref{proplim4}, we have
\fm{|\hat{\mu}_{qq}(\hat{\nu}-\frac{\eps G}{4t}\hat{\mu} \hat{U})|&\lesssim |\partial_q(G\hat{A}_qs\hat{\mu})|\cdot  \eps t^{-2+C\eps}\lra{r-t}\lesssim \eps t^{-2+C\eps}\lra{\hat{q}}^{-2+C\eps},}
\fm{|\partial_q(\hat{\mu}\hat{U})|&\lesssim |\hat{\mu}_q\hat{U}|+|2\hat{A}|\lesssim t^{C\eps}\lra{\hat{q}}^{-2+C\eps}+\lra{\hat{q}}^{-1+C\eps}\lesssim t^{C\eps}\lra{\hat{q}}^{-1+C\eps},}
\fm{|(-\frac{1}{2}G\hat{A}\hat{U}+\hat{U}_s)\hat{\mu}_q|\lesssim (\lra{q}^{-1+C\eps} t^{C\eps}+ \eps^{-1}t^{C\eps})\cdot \lra{q}^{-2+C\eps}t^{C\eps}\lesssim \eps^{-1}\lra{\hat{q}}^{-2+C\eps}t^{C\eps},}
\fm{|(-\frac{1}{2}G\partial_q(\hat{A}\hat{U})+\hat{U}_{sq})\hat{\mu}|&\lesssim |(-\frac{1}{2}G\partial_q(\hat{A}\hat{U})+\frac{1}{2}G\hat{A}\hat{U}_q)\hat{\mu}|\lesssim |(-\frac{1}{2}G\hat{A}_q\hat{U})\hat{\mu}|\lesssim \lra{\hat{q}}^{-2+C\eps}t^{C\eps}.}
In conclusion,
\fm{&\hspace{1.5em}(\partial_t-\partial_r)\partial_r(\hat{\nu}-\frac{\eps G(\omega)}{4t}\hat{\mu} \hat{U})\\&=\hat{\mu}_q\partial_r(\hat{\nu}-\frac{\eps G}{4t}\hat{\mu} \hat{U})+O(|\hat{q}_r|\eps\lra{\hat{q}}^{-1+C\eps}t^{-2+C\eps})\\
&=\hat{\mu}_q\partial_r(\hat{\nu}-\frac{\eps G}{4t}\hat{\mu} \hat{U})+O((-\hat{\mu})\eps\lra{\hat{q}}^{-1+C\eps}t^{-2+C\eps}+|\hat{\nu}|\eps\lra{\hat{q}}^{-1+C\eps}t^{-2+C\eps})\\
&=\hat{\mu}_q\partial_r(\hat{\nu}-\frac{\eps G}{4t}\hat{\mu} \hat{U})+O((-\hat{\mu})\eps\lra{\hat{q}}^{-1+C\eps}t^{-2+C\eps}+\eps\lra{\hat{q}}^{-1+C\eps}t^{-3+C\eps}).}
Take integral of the remainder terms along a charactersitic $(\tau,r+t-\tau,\omega)$ for $(r+t)/2-R\leq\tau\leq t$. We have
\fm{\int_{(r+t)/2-R}^t&\tau^{-2+C\eps}\eps\lra{z}^{-1+C\eps} (-\dot{z})+\eps\tau^{-3+C\eps}\ d\tau\lesssim \eps(1+\ln\lra{r-t})t^{-2+C\eps}.}
The proof of this estimate can be found in the proof of Lemma \ref{lem7.5}.
Since $\hat{\nu}-\frac{\eps G(\omega)}{4t}\hat{\mu} \hat{U}=0$ whenever $r-t>R$, we have $\partial_r(\hat{\nu}-\frac{\eps G(\omega)}{4t}\hat{\mu} \hat{U})=0$ at $\tau=(r+t)/2-R$. By Gronwall's inequality, we conclude that $\partial_r(\hat{\nu}-\frac{\eps G(\omega)}{4t}\hat{\mu} \hat{U})=O(\eps(1+\ln\lra{r-t})t^{-2+C\eps})$. To end the proof, we recall that $\partial_r=\hat{q}_r\partial_q$ where $\hat{q}_r>C^{-1}t^{-C\eps}$ in $\Omega\cap\{r-t<2R\}$.
\end{proof}\rm

Before we state the next lemma, we introduce the following notation.
\defn{\rm Fix $s,p\in\R$. We say that a function $F=F(t,x)$ with domain $\Omega\cap\{r-t<2R\}$ belongs to $S^{s,p}$, if for $\eps\ll_{s,p}1$, we have $Z^I(F)=O_I(t^{s+C_I\eps}\lra{r-t}^{p})$ for all multiindices $I$ in $\Omega\cap\{r-t<2R\}$.

It follows directly that $S^{s,p}+ S^{s',p'}\in S^{\max\{s,s'\},\max\{p,p'\}}$, that $S^{s,p}\cdot S^{s',p'}\in S^{s+s',p+p'}$, and that $Z^IS^{s,p}\in S^{s,p}$.

Following the proof of Corollary \ref{cor4.21}, we can show that $\mathfrak{R}_{s,p}\in S^{s,p}$. Here we prefer this new notation $S^{*,*}$ since it does not rely on the optical function $q(t,x)$ and the  corresponding null frames.}
\lem{\label{lemqhats}We have $\hat{q}\in S^{0,1}$. We also have $\Omega_{kk'}\hat{q}\in S^{0,\gamma}$ for each $1\leq k<k'\leq 3$ and $0<\gamma<1$. In other words, in $\Omega\cap\{r-t<2R\}$, for each $I$ we have
\eq{\label{zqhatest} |Z^I\hat{q}|\lesssim_I \lra{r-t}t^{C_I\eps};}
\eq{\label{zqhatest2} |Z^I\Omega_{kk'}\hat{q}|\lesssim_I t^{C_I\eps}\lra{r-t}^\gamma.}
As a result, we have $\partial_q^m\partial_\omega^n\hat{A}\in S^{0,-1-m}$, $\hat{\mu}\in S^{0,0}$, $\partial_s^p\partial_q^m\partial_\omega^n\hat{\mu} \in S^{0,-1-m}$ for $m+n+p>0$, $\partial_s^p\partial_\omega^n\hat{U}\in \eps^{-p}S^{0,0}$ and $\partial_s^p\partial_q^m\partial_\omega^n\hat{U}_q\in S^{0,-1-m}$. All functions here are of $(s,q,\omega)=(\eps\ln t-\delta,\hat{q}(t,x),\omega)$.}
\begin{proof}
We prove \eqref{zqhatest} by induction on $|I|$. The case $|I|=0$ has been proved in Lemma \ref{wtqlem1}. In general, suppose \eqref{zqhatest} holds for all $|I|\leq k$, and fix a multiindex $I$ with $|I|=k+1$. By the chain rule and Leibniz's rule, we express $Z^I\hat{\mu}$ as a linear combination of terms of the form
\eq{\label{zqclf1}  (\partial_s^b\partial_q^a\partial_\omega^c\hat{\mu})\cdot Z^{I_1}\hat{q}\cdots Z^{I_a}\hat{q}\cdot Z^{J_1}(\eps\ln t-\delta)\cdots Z^{J_b}(\eps\ln t-\delta)\cdot\prod_l Z^{K_{l,1}}\omega_l\cdots Z^{K_l,c_l}\omega_l}
where $a+b+|c|>0$, $|I_*|,|J_*|,|K_{*,*}|$ are nonzero, and the sum of all these multiindices is $k+1$. The only term with some $|I_*|>k$ is $\hat{\mu}_qZ^I\hat{q}$. All the other terms have an upper bound
\fm{\lra{\hat{q}}^{-1-a+C\eps}t^{C\eps}|\hat{\mu}|\cdot (\lra{r-t}t^{C\eps})^{a}\cdot \eps^b\cdot 1\lesssim \lra{\hat{q}}^{-1}t^{C\eps}|\hat{\mu}|.}
Here we apply Proposition \ref{proplim4} and  the induction hypotheses to control $Z^{I_*}\hat{q}$. In summary, we have $Z^I\hat{\mu}=\hat{\mu}_qZ^I\hat{q}+O(\lra{\hat{q}}^{-1+C\eps}t^{C\eps}|\hat{\mu}|)$. Following the same proof, we also have \fm{\sum_{0<|J|\leq k}|Z^J\hat{\mu}|=O(\lra{\hat{q}}^{-1+C\eps}t^{C\eps}|\hat{\mu}|).}

In addition, by the induction hypotheses and Lemma \ref{l2.1}, we have
\fm{\sum_{|J|<|I|} |(\partial_i+\omega_i\partial_t)Z^J\hat{q}|&\lesssim \sum_{|J|<k+1}(1+t+r)^{-1}|ZZ^J\hat{q}|\\&\lesssim (1+t+r)^{-1}\sum_{|J|=k+1}|Z^J\hat{q}|+t^{-1+C\eps}\lra{r-t}.}
In summary, by \eqref{lemtrcomf1} in Lemma \ref{lemtrcom} we have
\fm{|(\partial_t-\partial_r)Z^I\hat{q}|&\lesssim |\hat{\mu}_qZ^I\hat{q}|+(1+t+r)^{-1}\sum_{|J|=k+1}|Z^J\hat{q}|+ t^{C\eps}(-\hat{\mu})+t^{-1+C\eps}\lra{r-t}.}
Here we note that \fm{\sum_{|J|\leq k}|Z^J\hat{\mu}|\lesssim |\hat{\mu}|+\lra{\hat{q}}^{-1+C\eps}t^{C\eps}|\hat{\mu}|\lesssim t^{C\eps}(-\hat{\mu}).}
Now, we fix $(t,x)\in\Omega\cap\{r-t<2R\}$, integrate $(\partial_t-\partial_r)Z^I\hat{q}$ along the characteristic $(\tau,r+t-\tau,\omega)$ for $(t+r)/2-R\leq \tau\leq t$, and sum over all $|I|=k+1$. We then have
\fm{&\hspace{1.5em}\sum_{|I|=k+1}|Z^I\hat{q}(t,x)-Z^I\hat{q}|_{\tau=(r+t)/2-R}|\\&\lesssim \int_{(r+t)/2-R}^t(|\hat{\mu}_q|+(1+t+r)^{-1})\sum_{|I|=k+1}|Z^I\hat{q}|(\tau)+ \tau^{C\eps}(-\hat{\mu})+\tau^{C\eps}\ d\tau\\
&\lesssim \int_{(r+t)/2-R}^t(|\hat{\mu}_q|+(1+t+r)^{-1})\sum_{|I|=k+1}|Z^I\hat{q}|(\tau)\ d\tau+t^{C\eps}\lra{\hat{q}}+ \lra{r-t}t^{C\eps}.}
Moreover, we have $\hat{q}=r-t$ for $r-t>R$ and $\hat{q}=2R$ at $\tau=(r+t)/2-R$, so \fm{|Z^I\hat{q}|_{\tau=(r+t)/2-R}|=|Z^I(r-t)|_{\tau=(r+t)/2-R}|\lesssim t^{C\eps}.}
By Gronwall's inequality, we conclude that $\sum_{|I|=k+1}|Z^I\hat{q}(t,x)|\lesssim \lra{r-t}t^{C\eps}$.

Fix $\gamma>0$. Now we prove \eqref{zqhatest2} by induction on $|I|$.  First,  in Lemma \ref{lem7.5} we have proved $\hat{\lambda}_i=O((1+\ln\lra{r-t})t^{-1+C\eps})=O_\gamma(\lra{r-t}^\gamma t^{-1+C\eps})$. So we have $\Omega_{kk'}\hat{q}=x_k\lambda_{k'}-x_{k'}\lambda_k=O(\lra{r-t}^\gamma rt^{-1+C\eps})=O(\lra{r-t}^\gamma t^{C\eps})$, so the case $|I|=0$ is proved. In general, we fix $I$ with $|I|>0$. As computed above, we have
\fm{Z^I\Omega_{kk'}\hat{\mu}=\hat{\mu}_qZ^I\Omega_{kk'}\hat{q}+O(\lra{\hat{q}}^{-1+C\eps}t^{C\eps}|\hat{\mu}|),\hspace{2em}\sum_{|J|\leq |I|}|Z^J \hat{\mu}|=O(\lra{\hat{q}}^{-1+C\eps}t^{C\eps}|\hat{\mu}|);}
\fm{\sum_{|J|<|I|}|(\partial_i+\omega_i\partial_t)Z^J\Omega_{kk'}\hat{q}|&\lesssim (1+t+r)^{-1}\sum_{|J|\leq |I|}|Z^J\Omega_{kk'}\hat{q}|\\
&\lesssim (1+t+r)^{-1}\sum_{|J|=|I|}|Z^J\Omega_{kk'}\hat{q}|+t^{-1+C\eps}\lra{r-t}^\gamma.}
Thus, by \eqref{lemtrcomf2}, we have
\fm{|(\partial_t-\partial_r)Z^I\Omega_{kk'}\hat{q}|&\lesssim|\hat{\mu}_qZ^I\Omega_{kk'}\hat{q}|+(1+t+r)^{-1}\sum_{|J|=|I|}|Z^J\Omega_{kk'}\hat{q}|\\&\hspace{1em}+\lra{\hat{q}}^{-1+C\eps}t^{C\eps}(-\hat{\mu})+t^{-1+C\eps}\lra{r-t}^\gamma.}
Fix $(t,x)\in\Omega\cap\{r-t<2R\}$ and take integrals along a geodesic $(\tau,r+t-\tau,\omega)$. We note that
\fm{&\hspace{1.5em}\int_{(r+t)/2-R}^t\lra{\hat{q}(\tau)}^{-1+C\eps}\tau^{C\eps}(-\hat{\mu}(\tau))+\tau^{-1+C\eps}\lra{r+t-2\tau}^\gamma\ d\tau\\&\lesssim t^{C\eps}\int_{(r+t)/2-R}^t\lra{z(\tau)}^{-1}(-\dot{z}(\tau))\ d\tau+t^{-1+C\eps}\lra{r-t}^{1+\gamma}\\
&\lesssim (1+\ln\lra{r-t})t^{C\eps}+t^{C\eps}\lra{r-t}^\gamma\lesssim t^{C\eps}\lra{r-t}^\gamma.}
In addition, recall that $Z^Iq|_{\tau=(r+t)/2-R}=O(t^{C\eps})$. We finish the proof by applying Gronwall.

Finally, if $Q=Q(s,q,\omega)$ is a given function of $(s,q,\omega)$ and if we take $(s,q,\omega)=(\eps\ln t-\delta,\hat{q}(t,x),\omega)$, then $Z^IQ$ is a linear combination of terms of the form \eqref{zqclf1} with $\hat{\mu}$ replaced by $Q$. Thus,
\fm{|Z^IQ|\lesssim \sum_{a+b+|c|\leq |I|}\eps^b\lra{r-t}^{a} t^{C\eps}|\partial_s^b\partial_q^a\partial_\omega^cQ|.}
We combine this inequality with Proposition \ref{proplim4}.
As a result, we have $\partial_q^m\partial_\omega^n\hat{A}\in S^{0,-1-m}$, $\hat{\mu}\in S^{0,0}$, $\partial_s^p\partial_q^m\partial_\omega^n\hat{\mu} \in S^{0,-1-m}$ for $m+n+p>0$, $\partial_s^p\partial_\omega^n\hat{U}\in \eps^{-p}S^{0,0}$ and $\partial_s^p\partial_q^m\partial_\omega^n\hat{U}_q\in S^{0,-1-m}$.
\end{proof}\rm

\lem{Fix $\gamma\in(0,1)$. We have $\hat{\nu}\in \eps S^{-1,0}$, $\hat{\nu}_q\in \eps S^{-1,-1}$, $\hat{\lambda}_i\in S^{-1,\gamma}$ and \fm{\hat{\nu}-\frac{\eps}{4t}G(\omega)\hat{\mu} \hat{U}\in \eps S^{-2,1},\hspace{2em}\hat{\nu}_q-\frac{\eps}{4t}G(\omega)(\hat{\mu}_q \hat{U}-2\hat{A})\in \eps S^{-2,0}.}
All functions here are of $(s,q,\omega)=(\eps\ln t-\delta,\hat{q}(t,x),\omega)$.}
\begin{proof}First, we have \fm{\hat{\lambda}_i=\sum_jr^{-1}\omega_j\Omega_{ji}\hat{q}\in S^{-1,0}\cdot S^{0,\gamma}\subset S^{-1,\gamma}.}

Next, we set $Q:=\hat{\nu}-\eps G(\omega)\hat{\mu}\hat{U}/(4t)$. We have proved $Q=O(\eps t^{-2+C\eps}\lra{r-t})$ in Lemma \ref{lem7.6}. In general, we fix $I$ with $|I|>0$ and suppose $Z^JQ=O(\eps t^{-2+C\eps}\lra{r-t})$ whenever $|J|<|I|$. As computed in Lemma \ref{lem7.6}, we have
\fm{Q_t-Q_r&=\hat{\mu}_qQ+\frac{\eps G\hat{\mu}\hat{U}}{4t^2}-\frac{\eps^2 G}{4t^2}(-\frac{1}{2}G\hat{A}\hat{U}+\hat{U}_s)\hat{\mu}=\hat{\mu}_qQ+\eps S^{-2,0}.}
By \eqref{lemtrcomf1} in Lemma \ref{lemtrcom}, we have
\fm{&\hspace{1.5em}|(\partial_t-\partial_r)Z^IQ| \lesssim |Z^I(\hat{\mu}_qQ+\eps S^{-2,0})|+\sum_{|J|<|I|}[|Z^J(\hat{\mu}_qQ+\eps S^{-2,0})|+(1+t+r)^{-1}|ZZ^JQ|]\\
&\lesssim |\hat{\mu}_qZ^IQ|+(1+t+r)^{-1}\sum_{|J|=|I|}|Z^JQ|+\sum_{|K_1|+|K_2|\leq |I|\atop |K_2|<|I|}(|Z^{K_1}\hat{\mu}_q|+t^{-1})|Z^{K_2}Q|+\eps t^{-2+C\eps}\\
&\lesssim|\hat{\mu}_qZ^IQ|+(1+t+r)^{-1}\sum_{|J|=|I|}|Z^JQ|+\eps t^{-2+C\eps}\lra{r-t}^{-1}+\eps t^{-2+C\eps}.}
The last estimate follows from $\hat{\mu}_q\in S^{0,-2}$ and the induction hypotheses. Since $Q\equiv 0$ near $\tau=(r+t)/2-R$, and since \fm{\int_{(r+t)/2-R}^t \eps\tau^{-2+C\eps}\ d\tau\lesssim  \eps t^{-2+C\eps}\lra{r-t},}
we conclude by Gronwall that $Z^IQ=O(\eps t^{-2+C\eps}\lra{r-t})$. So $Q\in \eps S^{-2,1}$.

Since $\hat{\mu},\hat{U}\in S^{0,0}$ and since $\lra{r-t}\lesssim t$ in $\Omega\cap\{r-t<2R\}$, we have $\hat{\nu}=Q+\eps G(\omega)\hat{\mu}\hat{U}/(4t)\in \eps S^{-2,1}+\eps S^{-1,0}\subset \eps S^{-1,0}$. Moreover, for each $I$ we have
\fm{|Z^IQ_q|&\lesssim |Z^I(\hat{q}_r^{-1}\omega\cdot\partial Q)|\lesssim \sum_{|J|\leq |I|}t^{C\eps}|Z^J\partial Q|\\&\lesssim \sum_{|J|\leq |I|}t^{C\eps}|\partial Z^J Q|\lesssim \lra{r-t}^{-1}t^{C\eps}\sum_{|J|\leq |I|+1}|Z^J Q|\lesssim \eps t^{-2+C\eps}.}
Here we use the estimate $\hat{q}_r^{-1}\in S^{0,0}$ which follows from $\hat{q}_r\in S^{0,0}$ and $\hat{q}_r>C^{-1}t^{-C\eps}$. Thus,
\fm{Q_q&=\hat{\nu}_q-\frac{\eps}{4t}G(\omega)(\hat{\mu}_q \hat{U}-2\hat{A})\in \eps S^{-2,0}.}
Since $\hat{\mu}_q\hat{U}\in S^{0,-2}$ and $\hat{A}\in S^{0,-1}$, we conclude that $\hat{\nu}_q\in \eps S^{-1,-1}+\eps S^{-2,0}=\eps S^{-1,-1}$.
\end{proof}\rm

Now we prove that $\hat{q}$ is an approximate optical function.
\prop{We have
\fm{g^{\alpha\beta}(\hat{u})\hat{q}_\alpha \hat{q}_\beta\in S^{-2,1}.}}
\begin{proof}Fix $\gamma\in(0,1/2)$ and suppose we have obtained $\hat{\lambda}_i\in S^{-1,\gamma}$ from the pervious lemma. We note that $\hat{q}_t=\frac{1}{2}(\hat{\mu}+\hat{\nu})\in S^{0,0}$ and $\hat{q}_i=\frac{1}{2}(-\hat{\mu}+\hat{\nu})\omega_i+\hat{\lambda}_i\in S^{0,0}$. Thus,
\fm{g^{\alpha\beta}_0\hat{q}_\alpha \hat{q}_\beta
&=\frac{1}{4}g^{00}_0(\hat{\mu}+\hat{\nu})^2+\frac{1}{2}g^{0i}(\hat{\mu}+\hat{\nu})((-\hat{\mu}+\hat{\nu})\omega_i+2\hat{\lambda}_i)\\&\hspace{1em}+\frac{1}{4}g^{ij}_0((-\hat{\mu}+\hat{\nu})\omega_i+2\hat{\lambda}_i)((-\hat{\mu}+\hat{\nu})\omega_j+2\hat{\lambda}_j)\\
&=\frac{1}{4}G(\omega)\hat{\mu}^2+\frac{1}{2}g^{00}_0\hat{\mu}\hat{\nu}+\frac{1}{4}g^{00}_0\hat{\nu}^2+\frac{1}{2}g^{0i}_0(2\hat{\mu} \hat{\lambda}_i+\hat{\nu}^2\omega_i+2\hat{\nu}\hat{\lambda}_i)\\&\hspace{1em}+\frac{1}{4}g^{ij}_0(-\hat{\mu}(2\hat{\nu}\omega_j\omega_i+2\hat{\lambda}_j\omega_i+2\hat{\lambda}_i\omega_j)+(\hat{\nu}\omega_i+2\hat{\lambda}_i)(\hat{\nu}\omega_j+2\hat{\lambda}_j)).}
Since $\hat{\nu}\in \eps S^{-1,0}$ and $\hat{\lambda}_i\in S^{-1,\gamma}$, we have $\hat{\nu}^2,\hat{\nu}\hat{\lambda}_i,\hat{\lambda}_i\hat{\lambda}_j\in S^{-2,2\gamma}$ and thus
\fm{g^{\alpha\beta}_0\hat{q}_\alpha \hat{q}_\beta
&=\frac{1}{4}G(\omega)\hat{\mu}^2+\frac{1}{2}(g^{00}_0-g_{0}^{ij}\omega_i\omega_j)\hat{\mu}\hat{\nu}+g^{0i}_0\hat{\mu} \hat{\lambda}_i-\frac{1}{2}g^{ij}_0\hat{\mu}(\hat{\lambda}_j\omega_i+\hat{\lambda}_i\omega_j)\mod S^{-2,2\gamma}\\
&=\frac{1}{4}G(\omega)\hat{\mu}^2\mod S^{-1,\gamma}.}
If we replace $(g^{\alpha\beta}_0)$ with $(m^{\alpha\beta})$ in the computations, we have
\fm{-\hat{q}_t^2+\sum_j\hat{q}_j^2&=-\hat{\mu}\hat{\nu}-\frac{1}{2}m^{ij}\hat{\mu}(\hat{\lambda}_j\omega_i+\hat{\lambda}_i\omega_j)\mod S^{-2,2\gamma}=-\hat{\mu}\hat{\nu}\mod S^{-2,2\gamma}.}
Here we note that $m^{ij}\hat{\lambda}_j\omega_i=m^{ij}\hat{\lambda}_i\omega_j=\sum_j\omega_j(\hat{q}_j-\omega_j\hat{q}_r)=0$.

Moreover, note that $\hat{u}=\eps r^{-1}\hat{U}\in \eps S^{-1,0}$. Following the proof of Lemma \ref{l424im} with $V$ replaced by $Z$, we can prove that $f(\hat{u})-f(0)-f'(0)\hat{u}\in \eps^2 S^{-2,0}$ for each smooth function $f$.  Thus,
\fm{g^{\alpha\beta}(\hat{u})\hat{q}_\alpha\hat{q}_\beta&=-\hat{q}_t^2+\sum_j\hat{q}_j^2+g^{\alpha\beta}_0\hat{u}\hat{q}_\alpha \hat{q}_\beta+(g^{\alpha\beta}(\hat{u})-g^{\alpha\beta}_0\hat{u}-m^{\alpha\beta})\hat{q}_\alpha \hat{q}_\beta\\
&=-\hat{\mu}(\hat{\nu}-\frac{\eps}{4r}G(\omega)\hat{\mu}\hat{U})\mod  S^{-2,2\gamma}\\
&=-\hat{\mu}(\hat{\nu}-\frac{\eps}{4t}G(\omega)\hat{\mu}\hat{U})+\frac{\eps(t-r)}{4rt}G(\omega)\hat{\mu}^2\hat{U}\mod  S^{-2,2\gamma}\\&=\eps S^{-2,1}\mod  S^{-2,2\gamma}.}
Since $\gamma\in(0,1/2)$, we have $\eps S^{-2,1}\subset S^{-2,1}$ and $S^{-2,2\gamma}\subset S^{-2,1}$.
\end{proof}\rm

In order to prove that $\hat{u}$ is an approximate solution to \eqref{qwe}, we need the following lemma.
\lem{For each $\gamma\in(0,1/2)$, we have \fm{g^{\alpha\beta}(\hat{u})\partial_\alpha\partial_\beta\hat{q}=-r^{-1}\hat{\mu}+\frac{\eps}{2t}G\hat{A}\hat{\mu}\mod S^{-2,\gamma}.}}
\begin{proof}
Fix $\gamma\in(0,1/2)$ and suppose we have obtained $\hat{\lambda}_i\in S^{-1,\gamma}$. First we note that 
\fm{\eps t^{-1}\hat{\nu}_s&=\hat{\nu}_t-\hat{\nu}_q\hat{q}_t=\hat{\nu}_t+\hat{\nu}_r-\hat{\nu}\hat{\nu}_q,\\
\sum_j(\partial_i\omega_j)\hat{\nu}_{\omega_j}&=\hat{\nu}_i-\hat{\nu}_q\hat{q}_i=\hat{\nu}_i-\omega_i\hat{\nu}_r-\hat{\lambda}_i\hat{\nu}_q.}
Note that
\fm{\partial_t+\partial_r&=\frac{\sum_j\omega_j\Omega_{0j}+S}{r+t},\hspace{2em}\partial_i-\omega_i\partial_r=r^{-1}\sum_{j}\omega_j\Omega_{ji},}
and that $\hat{\nu}\in \eps S^{-1,0}$. Thus, we conclude that $\hat{\nu}_t+\hat{\nu}_r,\hat{\nu}_i-\omega_i\hat{\nu}_r\in \eps S^{-2,0}$. Besides, we have $\hat{\nu}\hat{\nu}_q\in \eps^2 S^{-2,-1}$ and $\hat{\lambda}_i\hat{\nu}_q\in \eps S^{-2,-1+\gamma}$. We conclude that $\eps t^{-1}\hat{\nu}_s,\sum_j(\partial_i\omega_j)\hat{\nu}_{\omega_j}\in \eps S^{-2,0}$.

Now, we have
\fm{\hat{q}_{tt}&=\partial_t(\frac{1}{2}(\hat{\mu}+\hat{\nu}))=\frac{1}{2}((\hat{\mu}_q+\hat{\nu}_q)\cdot\frac{1}{2}(\hat{\mu}+\hat{\nu})+\eps t^{-1}\hat{\mu}_s+\eps t^{-1}\hat{\nu}_s)\\
&=\frac{1}{4}\hat{\mu}_q \hat{\mu}+\frac{1}{4}\hat{\mu}_q\hat{\nu}+\frac{1}{4}\hat{\nu}_q \hat{\mu}  +\frac{\eps}{2t} \hat{\mu}_s\mod \eps S^{-2,0}=\frac{1}{4}\hat{\mu}_q \hat{\mu} \mod  \eps S^{-1,-1},}
\fm{\hat{q}_{ti}&=\partial_i(\frac{1}{2}(\hat{\mu}+\hat{\nu}))=\frac{1}{2}((\hat{\mu}_q+\hat{\nu}_q)\cdot(\frac{1}{2}(\hat{\nu}-\hat{\mu})\omega_i+\hat{\lambda}_i)+\sum_j(\partial_i\omega_j)\hat{\mu}_{\omega_j}+\sum_j(\partial_i\omega_j)\hat{\nu}_{\omega_j})\\
&=-\frac{1}{4}\hat{\mu}\hat{\mu}_q\omega_i\mod S^{-1,-1},}
\fm{\hat{q}_{ij}&=\partial_i(\frac{1}{2}(\hat{\nu}-\hat{\mu})\omega_j+\hat{\lambda}_j)\\
&=\frac{1}{2}(\hat{\nu}_q-\hat{\mu}_q)(\frac{1}{2}(\hat{\nu}-\hat{\mu})\omega_i+\hat{\lambda}_i)\omega_j+\frac{1}{2}\sum_k(\hat{\nu}_{\omega_k}-\hat{\mu}_{\omega_k})(\partial_i\omega_k)\omega_j+\frac{1}{2}(\hat{\nu}-\hat{\mu})\partial_i\omega_j+\partial_i\hat{\lambda}_j\\
&=\frac{1}{4}(\hat{\mu}\hat{\mu}_q-\hat{\mu}_q\hat{\nu}-\hat{\nu}_q\hat{\mu})\omega_i\omega_j-\frac{1}{2}\hat{\mu}_q\hat{\lambda}_i\omega_j-\frac{1}{2}\sum_k\hat{\mu}_{\omega_k}(\partial_i\omega_k)\omega_j-\frac{1}{2}\hat{\mu}\partial_i\omega_j+\partial_i\hat{\lambda}_j\mod \eps S^{-2,0}\\
&=\frac{1}{4}\hat{\mu}\hat{\mu}_q\omega_i\omega_j\mod S^{-1,0}.}
In the last estimate, we note that  $\partial_i\hat{\lambda}_j\in S^{-1,0}$ since for each $I$, \fm{|Z^I\partial_i\hat{\lambda}_j|&\lesssim\sum_{|J|\leq |I|}|\partial Z^{J}\hat{\lambda}_j|\lesssim \lra{r-t}^{-1}\sum_{|J|\leq |I|+1}|Z^{J}\hat{\lambda}_j|\\&\lesssim\lra{r-t}^{-1}\cdot t^{-1+C\eps}\lra{r-t}^\gamma\lesssim t^{-1+C\eps}\lra{r-t}^{1-\gamma}.}
Thus, we have $\partial^2\hat{q}\in S^{0,-2}+S^{-1,-1}= S^{0,-2}$ and
\fm{g^{\alpha\beta}_0\hat{q}_{\alpha\beta}&=\frac{1}{4}G(\omega)\hat{\mu}_q\hat{\mu}\mod S^{-1,0}.}
In addition,
\fm{\Box \hat{q}&=-(\frac{1}{4}\hat{\mu}_q \hat{\mu}+\frac{1}{4}\hat{\mu}_q\hat{\nu}+\frac{1}{4}\hat{\nu}_q \hat{\mu}  +\frac{\eps}{2t} \hat{\mu}_s)+[\frac{1}{4}(\hat{\mu}\hat{\mu}_q-\hat{\mu}_q\hat{\nu}-\hat{\nu}_q\hat{\mu})-r^{-1}\hat{\mu}+\sum_i\partial_i\hat{\lambda}_i]\mod \eps S^{-2,0}\\
&=-(\frac{1}{2}\hat{\mu}_q\hat{\nu}+\frac{1}{2}\hat{\nu}_q \hat{\mu}  +\frac{\eps}{2t} \hat{\mu}_s)-r^{-1}\hat{\mu}+\sum_i\partial_i\hat{\lambda}_i\mod \eps S^{-2,0}.}
Since $\sum_i\omega_i\hat{\lambda}_i=0$, we have $0=\partial_r(\sum_i\omega_i\hat{\lambda}_i)=\sum_i\omega_i\partial_r\hat{\lambda}_i$. And since $\hat{\lambda}_i\in  S^{-1,\gamma}$, we have
\fm{\sum_i\partial_i\hat{\lambda}_i=\sum_i(\partial_i-\omega_i\partial_r)\hat{\lambda}_i=\sum_{i,j}r^{-1}\omega_i\Omega_{ji}\hat{\lambda}_i\in S^{-2,\gamma}}
Finally, we have
\fm{g^{\alpha\beta}(\hat{u})\partial_\alpha\partial_\beta\hat{q}&=\Box\hat{q}+g^{\alpha\beta}_0\hat{u}\partial_\alpha\partial_\beta\hat{q}+(g^{\alpha\beta}(\hat{u})-g^{\alpha\beta}_0\hat{u}-m^{\alpha\beta})\partial_\alpha\partial_\beta\hat{q}\\
&=-(\frac{1}{2}\hat{\mu}_q\hat{\nu}+\frac{1}{2}\hat{\nu}_q \hat{\mu}  +\frac{\eps}{2t} \hat{\mu}_s)-r^{-1}\hat{\mu}+\frac{\eps}{4r}G(\omega)\hat{\mu}\hat{\mu}_q\hat{U}\mod  S^{-2,\gamma}\\
&=-\frac{1}{2}\hat{\mu}_q\cdot\frac{\eps}{4t}G\hat{\mu}\hat{U}-\frac{1}{2}\hat{\mu} \cdot\frac{\eps}{4t}G(\hat{\mu}_q\hat{U}-2\hat{A})+\frac{\eps}{4t}G\hat{A}\hat{\mu}-r^{-1}\hat{\mu}\\&\hspace{1em}+\frac{\eps}{4t}G\hat{\mu}\hat{\mu}_q\hat{U}+\frac{\eps(t-r)}{4tr}G\hat{\mu}\hat{\mu}_q\hat{U}\mod  S^{-2,\gamma}\\
&=-r^{-1}\hat{\mu}+\frac{\eps}{2t}G\hat{A}\hat{\mu}\mod  S^{-2,\gamma}.}
\end{proof}\rm

Now we claim that $\hat{u}=\eps r^{-1}\hat{U}(\eps\ln t-\delta,\hat{q}(t,x),\omega)$ is an approximate solution to \eqref{qwe}.
\prop{\label{prop6.3main}We have
\fm{g^{\alpha\beta}(\hat{u})\partial_\alpha\partial_\beta \hat{u}\in \eps S^{-3,0} .}}
\begin{proof}We have
\fm{\hat{u}_t&=\eps r^{-1}(\eps t^{-1}\hat{U}_s+\hat{q}_t\hat{U}_q),\hspace{2em}\hat{u}_i=-\eps r^{-2}\omega_i\hat{U}+\eps r^{-1}(\hat{U}_q\hat{q}_i+\sum_k\hat{U}_{\omega_k}\partial_i\omega_k).}
By Lemma \ref{lemqhats}, we have $\partial_s^b\partial_\omega^c\hat{U}\in \eps^{-b}S^{0,0}$. Thus we have
\fm{\hat{u}_{tt}&=\eps r^{-1}(-\eps t^{-2}\hat{U}_s+\eps^2t^{-2}\hat{U}_{ss}+2\eps t^{-1}\hat{q}_t\hat{U}_{sq}+\hat{q}_{tt}\hat{U}_q+\hat{q}_t^2\hat{U}_{qq})\\
&=\eps r^{-1}(2\eps t^{-1}\hat{q}_t\hat{U}_{sq}+\hat{q}_{tt}\hat{U}_q+\hat{q}_t^2\hat{U}_{qq})\mod \eps S^{-3,0}\\
&=\eps r^{-1}(\hat{q}_{tt}\hat{U}_q+\hat{q}_t^2\hat{U}_{qq})\mod \eps S^{-2,-1},}
\fm{\hat{u}_{ti}&=-\eps r^{-2}\omega_i(\eps t^{-1}\hat{U}_s+\hat{q}_t\hat{U}_q)\\&\hspace{1em}+\eps r^{-1}(\eps t^{-1}\hat{U}_{sq}\hat{q}_i+\eps t^{-1}\sum_k\hat{U}_{s\omega_k}\partial_i\omega_k+\hat{q}_{ti}\hat{U}_q+\hat{q}_t\hat{U}_{qq}\hat{q}_i+\hat{q}_t\sum_k\hat{U}_{q\omega_k}\partial_i\omega_k)\\
&=\eps r^{-1}(\hat{q}_{ti}\hat{U}_q+\hat{q}_t\hat{U}_{qq}\hat{q}_i)\mod \eps S^{-2,-1},}
\fm{\hat{u}_{ij}&=-\eps \partial_i(r^{-2}\omega_j)\hat{U}-\eps r^{-2}\omega_j(\hat{U}_q\hat{q}_i+\sum_k\hat{U}_{\omega_k}\partial_i\omega_k)-\eps r^{-2}\omega_i(\hat{U}_q\hat{q}_j+\sum_k\hat{U}_{\omega_k}\partial_j\omega_k)\\&\hspace{1em}+\eps r^{-1}[\hat{U}_{qq}\hat{q}_i\hat{q}_j+\sum_k\hat{U}_{q\omega_k}(\partial_i\omega_k)\hat{q}_j+\hat{U}_q\hat{q}_{ij}\\&\hspace{5em}+\sum_k(\hat{U}_{\omega_kq}\hat{q}_i\partial_j\omega_k+\hat{U}_{\omega_k}\partial_i\partial_j\omega_k)+\sum_{k,k'}\hat{U}_{\omega_k\omega_{k'}}(\partial_i\omega_k)(\partial_j\omega_{k'})]\\
&=-\eps r^{-2}\omega_j\hat{U}_q\hat{q}_i-\eps r^{-2}\omega_i\hat{U}_q\hat{q}_j\\&\hspace{1em}+\eps r^{-1}[\hat{U}_{qq}\hat{q}_i\hat{q}_j+\sum_k\hat{U}_{q\omega_k}((\partial_i\omega_k)\hat{q}_j+(\partial_j\omega_k)\hat{q}_i)+\hat{U}_q\hat{q}_{ij}]\mod \eps S^{-3,0}\\
&=\eps r^{-1}(\hat{U}_{qq}\hat{q}_i\hat{q}_j+\hat{U}_q\hat{q}_{ij})\mod \eps S^{-2,-1}.}
Since $g^{\alpha\beta}(\hat{u})-m^{\alpha\beta}=g^{\alpha\beta}_0\hat{u}\bmod \eps^2S^{-2,0}\in \eps S^{-1,0}$, we have
\fm{&\hspace{1.5em}g^{\alpha\beta}(\hat{u})\partial_\alpha\partial_\beta\hat{u}=\Box\hat{u}+(g^{\alpha\beta}(\hat{u})-m^{\alpha\beta})\partial_\alpha\partial_\beta\hat{u}\\
&=-\eps r^{-1}(2\eps t^{-1}\hat{q}_t\hat{U}_{sq}+\hat{q}_{tt}\hat{U}_q+\hat{q}_t^2\hat{U}_{qq})-2\eps r^{-2}\hat{U}_q\hat{q}_r\\&\hspace{1em}+\eps r^{-1}\sum_i[\hat{U}_{qq}\hat{q}_i^2+\sum_k2\hat{U}_{q\omega_k}(\partial_i\omega_k)\hat{q}_i+\hat{U}_q\hat{q}_{ii}]\\&\hspace{1em}+(g^{\alpha\beta}(\hat{u})-m^{\alpha\beta})\cdot\eps r^{-1}(\hat{q}_{\alpha\beta}\hat{U}_q+\hat{q}_\alpha\hat{q}_\beta\hat{U}_{qq})\mod \eps S^{-3,0}\\
&=-\eps^2 (tr)^{-1}\hat{q}_tG\hat{A}\hat{U}_{q}-2\eps r^{-2}\hat{U}_q\hat{q}_r+\eps r^{-1}\sum_i\sum_k2\hat{U}_{q\omega_k}(\partial_i\omega_k)(\hat{\lambda}_i+\omega_i\hat{q}_r)\\&\hspace{1em}+\eps r^{-1}(g^{\alpha\beta}(\hat{u})\hat{q}_{\alpha\beta}\hat{U}_q+g^{\alpha\beta}(\hat{u})\hat{q}_\alpha\hat{q}_\beta\hat{U}_{qq})\mod \eps S^{-3,0}\\
&=-\eps^2 (rt)^{-1}\hat{q}_tG\hat{A}\hat{U}_{q}-2\eps r^{-2}\hat{U}_q\hat{q}_r-\eps r^{-2}\hat{\mu}\hat{U}_{q}+\eps^2(2tr)^{-1}G\hat{A}\hat{\mu}\hat{U}_q\mod \eps S^{-3,0}\\
&=-\frac{1}{2}\eps^2 r^{-2}\hat{\nu}G\hat{A}\hat{U}_{q}-\eps r^{-2}\hat{\nu}\hat{U}_{q}\mod \eps S^{-3,0}\in \eps S^{-3,0}.}
In the third equality, we note that \fm{\eps r^{-1}[g^{\alpha\beta}(\hat{u})\hat{q}_{\alpha\beta}+r^{-1}\hat{\mu}-\frac{\eps}{2t}G\hat{A}\hat{\mu}]\hat{U}_{q}\in \eps S^{-1,0}\cdot S^{-2,\gamma}\cdot S^{0,-1}\subset \eps S^{-3,0},\\ \eps r^{-1}g^{\alpha\beta}(\hat{u})\hat{q}_{\alpha}\hat{q}_{\beta}\hat{U}_{qq}\in \eps S^{-1,0}\cdot S^{-2,1}\cdot S^{0,-2}\subset \eps S^{-3,0}}
and that
\fm{\eps r^{-1}\sum_i\sum_k2\hat{U}_{q\omega_k}(\partial_i\omega_k)(\hat{\lambda}_i+\omega_i\hat{q}_r)&=\eps r^{-1}\sum_i\sum_k2\hat{U}_{q\omega_k}(\partial_i\omega_k)\hat{\lambda}_i+\eps r^{-1}\sum_k2\hat{U}_{q\omega_k}(\partial_r\omega_k)\hat{q}_r\\
&\in\eps S^{-1,0}\cdot S^{0,-1}\cdot S^{-1,0}\cdot S^{-1,\gamma}+0\subset \eps S^{-3,0}.}
\end{proof}\rm

\subsection{Approximation of the optical function}\label{sec6.4}

We set $p(t,x):=F(q(t,x),\omega)-\hat{q}(t,x)$ in $\Omega$, where $q(t,x)$ is the optical function constructed in Section \ref{s3}.  
\prop{\label{prop6.6}Fix a constant $\gamma\in(0,1)$. Then, for $\eps\ll_\gamma 1$, at each $(t,x)\in\Omega$ such that $|r-t|\lesssim t^{\gamma}$, we have $|p(t,x)|\lesssim_\gamma t^{-1+C\eps}\lra{r-t}$.}
\begin{proof} It is clear that $p\equiv 0$ in the region $\{r-t>R\}$. In $\Omega\cap\{r-t<2R\}$, by setting $s=\eps\ln t-\delta$ we have
\eq{\label{prop6.6f1}p_t-p_r&=F_q\mu(s,q(t,x),\omega)-\hat{\mu}(s,\hat{q}(t,x),\omega)\\&=[F_q\mu(s,q(t,x),\omega)-\hat{\mu}(s,F(q(t,x),\omega),\omega)]+[\hat{\mu}(s,F(q(t,x),\omega),\omega)-\hat{\mu}(s,\hat{q}(t,x),\omega)]\\
&=:\mathcal{R}_1+\mathcal{R}_2.}
Since $\hat{A}(F(q,\omega),\omega)=A(q,\omega)$, we have
\eq{\label{prop6.6f2}\mathcal{R}_1&=-\frac{2}{A_1(q(t,x),\omega)}\wt{V}_1(s,q(t,x),\omega)\exp(-\frac{1}{2}G(\omega)A(q(t,x),\omega)s)\\
&\hspace{1em}+2\exp(-\frac{1}{2}G(\omega)\hat{A}(F(q(t,x),\omega),\omega)s)\\
&=(-\frac{2}{A_1(q(t,x),\omega)}\wt{V}_1(s,q(t,x),\omega)+2)\exp(-\frac{1}{2}G(\omega)A(q(t,x),\omega)s)\\
&=-\frac{2}{A_1(q(t,x),\omega)}(\wt{V}_1(s,q(t,x),\omega)-A_1(q(t,x),\omega))\exp(-\frac{1}{2}G(\omega)A(q(t,x),\omega)s).}
By Proposition \ref{proplim3}, we have
\fm{|\mathcal{R}_1|\lesssim |\wt{V}_1(s,q(t,x),\omega)-A_1(q(t,x),\omega)|\exp(C\lra{q}^{-1+C\eps}s)\lesssim t^{-1+C\eps}.}
Moreover, \fm{|\mathcal{R}_2|&=|\int_{\hat{q}}^{F(q,\omega)}\hat{\mu}_\rho(s,\rho,\omega)\ d\rho|\lesssim \left|\int_{\hat{q}}^{F(q,\omega)}\lra{\rho}^{-2+C\eps}s|\hat{\mu}(s,\rho,\omega)|\ d\rho\right|\\
&\lesssim (\eps\ln t-\delta)|p|\cdot\max_{\kappa\in[0,1]}\big[\lra{\hat{q}+\kappa p}^{-2+C\eps}\exp(-\frac{1}{2}G(\omega)\hat{A}(\hat{q}+\kappa p,\omega)s)\big].}

We now use a continuity argument to end the proof. Fix $(t,x)\in\Omega\cap\{r-t<2R,\ |r-t|\lesssim t^{\gamma}\}$. Suppose that for some $t_0\in[(r+t)/2-R,t)$, we have 
\eq{\label{prop6.6ca}|p(\tau,r+t-\tau,\omega)|\leq \frac{\delta}{10\eps\ln \tau},\hspace{2em}\forall\tau\in[(r+t)/2-R,t_0].} 
Note that $\eqref{prop6.6ca}$ holds for $t_0=(r+t)/2-R$, since $p((r+t)/2-R,(r+t)/2+R,\omega)=0$. At $(\tau,r+t-\tau,\omega)$ for $(r+t)/2-R\leq\tau\leq t_0$ and for each $\kappa\in[0,1]$, we have \fm{\lra{\hat{q}+\kappa p}\sim 1+|\hat{q}+\kappa p|\geq 1+|\hat{q}|-|\kappa p|\geq 1+|\hat{q}|-\frac{1}{10}\gtrsim \lra{\hat{q}}.}
In the second last inequality we note that $\tau>\exp(\delta/\eps)$, so $\eps\ln\tau>\delta$ and thus $|p|\leq 1/10$. Moreover, 
\fm{\exp(-\frac{1}{2}G(\omega)(\hat{A}(\hat{q}+\kappa p,\omega)-\hat{A}(\hat{q},\omega))s)\lesssim\exp(C\kappa |p|s)\lesssim\exp(\delta/10)\lesssim 1.}
In conclusion, at $(\tau,r+t-\tau,\omega)$ for $(r+t)/2-R\leq\tau\leq t_0$, we have
\fm{|\mathcal{R}_2|&\lesssim(\eps\ln \tau-\delta)[|p|\lra{\hat{q}}^{-2+C\eps}\exp(-\frac{1}{2}G(\omega)\hat{A}(\hat{q},\omega)s)](\tau,r+t-\tau,\omega)\\
&\lesssim (\eps\ln \tau-\delta)[|p|\lra{\hat{q}}^{-2+C\eps}(-\hat{\mu})](\tau,r+t-\tau,\omega).}
If we fix any $t_1\in[(r+t)/2-R,t_0]$, then \fm{\int_{(r+t)/2-R}^{t_1}(\eps\ln \tau-\delta)\lra{\hat{q}}^{-2+C\eps}(-\hat{\mu})(\tau,r+t-\tau,\omega)\ d\tau&\lesssim\eps\ln t_1\int_{(r+t)/2-R}^{t_1} \lra{z}^{-2+C\eps}(-\dot{z})\ d\tau\\
&\lesssim \eps\ln t_1}
and\fm{\int_{(r+t)/2-R}^{t_1}|\mathcal{R}_1|(\tau,r+t-\tau,\omega)\ d\tau&\lesssim\int_{(r+t)/2-R}^{t_1} \tau^{-1+C\eps}\ d\tau\\&\lesssim ((r+t)/2-R)^{-1+C\eps}(t_1-(r+t)/2+R)\\
&\lesssim t_1^{-1+C\eps}\lra{r-t}.}
Here we recall that $[(r+t)/2-R]\sim t\sim t_1$. And since $p=0$ at $\tau=(r+t)/2-R$, by applying the Gronwall's inequality to $p_t-p_r=\mathcal{R}_1+\mathcal{R}_2$, we conclude that 
\eq{\label{prop6.6ca2}|p(t_1,r+t-t_1,\omega)|\lesssim t_1^{-1+C\eps}\lra{r-t}\cdot\exp(C\eps\ln (Ct_1))&\lesssim t_1^{-1+C\eps}\lra{r-t},\\ &\forall t_1\in[(r+t)/2-R,t_0].}
For $\eps\ll_\gamma 1$ (where $\eps$ does not depend on $(t,x)$) and $t_1\in[(r+t)/2-R,t_0]$, we have $|r-t|\lesssim t^\gamma\sim t_1^\gamma$ and thus \fm{t_1^{-1+C\eps}\lra{r-t}\lesssim t_1^{-1+\gamma+C\eps}\leq t_1^{(\gamma-1)/2}\leq\delta/(20\eps\ln t_1).}  And since $\tau\mapsto\eps(\ln\tau) p(\tau,r+t-\tau,\omega)$ is a continuous function, \eqref{prop6.6ca} holds with $t_0$ replaced by some $t_0'>t_0$. By the continuity argument we conclude that $|p(t,x)|\lesssim t^{-1+C\eps}\lra{r-t}$. The constants here do not depend on $(t,x)$.
\end{proof}\rm

Next we consider $Z^Ip$. We need the following  lemma.
\lem{\label{lem6.8r12} Let $\mathcal{R}_1$ and $\mathcal{R}_2$ be defined as in \eqref{prop6.6f1}. Then, we have $\mathcal{R}_1\in S^{-1,0}$ and for $|I|>0$ we have
\fm{|Z^I\mathcal{R}_2|&\lesssim \lra{r-t}^{-2}t^{C\eps}\sum_{|J|<|I|}|Z^Jp|+|\hat{\mu}_qZ^Ip|.}}
\begin{proof} By \eqref{prop6.6f2}, Remark \ref{proplim3.1} and Lemma \ref{proplim2}, and since $A_1<-1$ everywhere, we have $\mathcal{R}_{1}=\mathfrak{R}_{0,0}\cdot\mathfrak{R}_{-1,0}\cdot \mathfrak{R}_{0,0}=\mathfrak{R}_{-1,0}\in S^{-1,0}$. 

To estimate $\mathcal{R}_2$,  we fix an arbitrary multiindex $I$ with $|I|>0$. By the chain rule and Leibniz's rule, we can express $Z^I\hat{\mu}(s,F(q(t,x),\omega),\omega)-Z^I\hat{\mu}(s,\hat{q}(t,x),\omega)$ as a linear combination of terms of the form
\eq{\label{zqclf2} [(\partial_s^b\partial_q^a\partial_\omega^c\hat{\mu})(s,F(q,\omega),\omega)\cdot \prod_{i=1}^aZ^{I_i}(F(q,\omega))-(\partial_s^b\partial_q^a\partial_\omega^c\hat{\mu})(s,\hat{q},\omega)\cdot \prod_{i=1}^aZ^{I_i}\hat{q}]\\\cdot \prod_{j=1}^bZ^{J_j}(\eps\ln t-\delta)\cdot\prod_{l=1}^c Z^{K_{l,1}}\omega_l\cdots Z^{K_l,c_l}\omega_l}
where $|I_*|,|J_*|,|K_{*,*}|$ are nonzero, and the sum of all these multiindices is $|I|$. The only term with $|I_j|=|I|$ for some $j$ is $\hat{\mu}_qZ^Ip$, so from now on we assume $|I_j|<|I|$ for each $j$ in \eqref{zqclf2}. Here the second row in \eqref{zqclf2} is $O(\eps^b)$. The first row is equal to the sum of 
\eq{\label{zqclf3} [(\partial_s^b\partial_q^a\partial_\omega^c\hat{\mu})(s,F(q,\omega),\omega)-(\partial_s^b\partial_q^a\partial_\omega^c\hat{\mu})(s,\hat{q},\omega)]\cdot \prod_{i=1}^aZ^{I_i}(F(q,\omega))}
and for each $j=1,2,\dots,a$
\eq{\label{zqclf4}(\partial_s^b\partial_q^a\partial_\omega^c\hat{\mu})(s,\hat{q},\omega)\cdot \prod_{i=1}^{j-1}Z^{I_i}(F(q,\omega))\cdot Z^{I_j}p\cdot\prod_{i=j+1}^aZ^{I_i}\hat{q}.}
Since $|I|>0$, we must have $a>0$ if \eqref{zqclf4} does appear. 

To control \eqref{zqclf3} and \eqref{zqclf4}, we first recall from Lemma \ref{lemqhats} and Proposition \ref{proplim4} that \fm{Z^{I_*}(\hat{q}(t,x),F(q(t,x),\omega))&=O(\lra{r-t}t^{C\eps});}
\fm{(\partial_s^b\partial_q^a\partial_\omega^c\hat{\mu})(s,\hat{q},\omega)&=O(\lra{\hat{q}}^{-a-1+C\eps}t^{C\eps})=O(\lra{r-t}^{-a-1}t^{C\eps}),\hspace{2em} \text{when }a+b+|c|>0.} It follows immediately that \eqref{zqclf4} is $O(\sum_{|J|<|I|}t^{C\eps}\lra{r-t}^{-2}|Z^Jp|)$. In addition, we have $\lra{F(q,\omega)}/\lra{r-t}\sim\lra{q}/\lra{r-t}=t^{O(\eps)}$ and $\lra{\hat{q}}/\lra{r-t}=t^{O(\eps)}$. Thus, for each $\tau\in[0,1]$, 
\eq{\label{lem6.8r12f1}\lra{\tau\hat{q}+(1-\tau)F(q,\omega)}\sim \tau\lra{\hat{q}}+(1-\tau)\lra{F(q,\omega)}\gtrsim\lra{r-t}t^{-C\eps}.} Then, we have
\fm{&\hspace{1.5em}|(\partial_s^b\partial_q^a\partial_\omega^c\hat{\mu})(s,F(q,\omega),\omega)-(\partial_s^b\partial_q^a\partial_\omega^c\hat{\mu})(s,\hat{q},\omega)|=|\int_{\hat{q}}^{F(q,\omega)}(\partial_s^b\partial_q^{a+1}\partial_\omega^c\hat{\mu})(s,\rho,\omega)\ d\rho|\\
&\lesssim|\int_{\hat{q}}^{F(q,\omega)}\lra{\rho}^{-2-a+C\eps}\exp(Cs)\ d\rho|\lesssim |p(t,x)|t^{C\eps}\lra{r-t}^{-a-2}.}
Thus, \eqref{zqclf3} is $O(|p|t^{C\eps}\lra{r-t}^{-2})$. 

In conclusion, for $|I|>0$ we have
\fm{|Z^I\mathcal{R}_2|&\lesssim \lra{r-t}^{-2}t^{C\eps}\sum_{|J|<|I|}|Z^Jp|+|\hat{\mu}_qZ^Ip|.}
\end{proof}\rm

\prop{Fix a constant $\gamma\in(0,1/2)$ and a large integer $N$. Then, for $\eps\ll_{\gamma,N} 1$, at each $(t,x)\in\Omega$ such that $|r-t|\lesssim t^{\gamma}$, we have $|Z^Ip(t,x)|\lesssim_\gamma t^{-1+C\eps}\lra{r-t}$ for each $|I|\leq N$.}
\begin{proof} We prove by induction on $|I|$. The case $|I|=0$ has been proved in Proposition \ref{prop6.6}.
Fix a multiindex $I$ with $|I|>0$, and suppose that we have proved the proposition for all $|J|<|I|$. By Lemma \ref{lemtrcom}, we have
\fm{(\partial_t-\partial_r)Z^Ip&=Z^I(p_t-p_r)+\sum_{|J|<|I|}[f_0Z^J(p_t-p_r)+\sum_if_0(\partial_i+\omega_i\partial_t)Z^Jp].}
By Lemma \ref{lem6.8r12} and our induction hypotheses, in $\Omega\cap\{r-t<2R,\ |r-t|\lesssim t^{\gamma}\}$ we have
\fm{&\hspace{1.5em}|(\partial_t-\partial_r)Z^Ip|\lesssim |Z^I(\mcl{R}_1+\mcl{R}_2)|+\sum_{|J|<|I|}|Z^J(\mcl{R}_1+\mcl{R}_2)|+t^{-1}|ZZ^Jp|]\\&\lesssim t^{-1+C\eps}+\lra{r-t}^{-2}t^{C\eps}\sum_{|J|<|I|}|Z^Jp|+|\hat{\mu}_qZ^Ip|+\sum_{|J|\leq|I|}t^{-1}|Z^Jp|\\
&\lesssim t^{-1+C\eps}+\lra{r-t}^{-2}\cdot t^{-1+C\eps}\lra{r-t}+|\hat{\mu}_qZ^Ip|+\sum_{|J|=|I|}t^{-1}|Z^Jp|+t^{-2+C\eps}\lra{r-t}\\
&\lesssim t^{-1+C\eps}+|\hat{\mu}_qZ^Ip|+\sum_{|J|=|I|}t^{-1}|Z^Jp|.}
The integral of $|\hat{\mu}_q|$ and $t^{-1}$ along a characteristic $(\tau,r+t-\tau,\omega)$, $\tau\in[(r+t)/2-R,t]$, is $O(\eps\ln t+1)$. Moreover,
\fm{&\hspace{1.5em}\int_{(r+t)/2-R}^t \tau^{-1+C\eps}\ d\tau\lesssim ((r+t)/2-R)^{-1+C\eps}((t-r)/2+R)\lesssim t^{-1+C\eps}\lra{r-t}.}
Since  $Z^Ip\equiv 0$ in the region $\Omega\cap\{r-t>R\}$, by  Gronwall's inequality we conclude that $|Z^Ip|\lesssim t^{-1+C\eps}\lra{r-t}$. 
\end{proof}
\rm

\subsection{Approximation of the solution to \eqref{qwe}}\label{sec6.5}
We can now discuss the difference $u-\hat{u}$ where $u$ is a solution to \eqref{qwe} and $\hat{u}$ is defined in Section \ref{sec6.2}. Again, we fix a point in region $\Omega\cap\{|r-t|\lesssim t^{\gamma}\}$ for some $0<\gamma<1$. Note that
\fm{u-\hat{u}&=\eps r^{-1}U(s,q(t,x),\omega)-\eps r^{-1}\hat{U}(s,\hat{q}(t,x),\omega)\\&=\eps r^{-1}U(s,q(t,x),\omega)-\eps r^{-1}\hat{U}(s,F(q(t,x),\omega),\omega)\\&\hspace{1em}+\eps r^{-1}\hat{U}(s,F(q(t,x),\omega),\omega)-\eps r^{-1}\hat{U}(s,\hat{q}(t,x),\omega)\\
&=:\mathcal{R}_3+\mathcal{R}_4.}
Now we estimate $\mathcal{R}_3$ and $\mathcal{R}_4$ separately.

\lem{Fix a constant $0<\gamma<1$ and a large integer $N$. Then, for $\eps\ll_{\gamma,N} 1$, at each $(t,x)\in\Omega$ such that $|r-t|\lesssim  t^{\gamma}$, we have $|Z^I\mathcal{R}_3|\lesssim_\gamma \eps t^{-2+C\eps}\lra{r-t}$ for each $|I|\leq N$.}
\begin{proof}
As computed in Lemma \ref{lemsimplif}, by change of variables  we can prove that
\fm{\hat{U}(s,F(q(t,x),\omega),\omega)&=\wt{U}(s,q(t,x),\omega).}
Thus,
\fm{\mathcal{R}_3&=\eps r^{-1}(U(s,q(t,x),\omega)-\wt{U}(s,q(t,x),\omega)).}
By \eqref{sec5ef2}, we have $|U-\wt{U}|\lesssim \lra{q}t^{-1+C\eps}$ at $(s,q,\omega)=(\eps\ln t-\delta,q(t,x),\omega)$, so 
\fm{|\mathcal{R}_3|\lesssim \eps t^{-2+C\eps}\lra{q}\lesssim \eps t^{-2+C\eps}\lra{r-t}.}

Next we fix a multiindex $I$ with $|I|>0$. Then, $Z^I\mathcal{R}_3$ can be expressed as a linear combination of terms of the form
\eq{\label{6.9f1}Z^{I'}(\eps r^{-1})\cdot(\partial_s^b\partial_q^a\partial_\omega^c (U-\wt{U}))(s,q,\omega)\cdot\prod_{i=1}^a Z^{I_i}q\cdot \prod_{i=1}^b Z^{J_i}s\cdot \prod_{i=1}^c Z^{K_i}\omega.}
The sum of all the $|I'|,|I_*|,|J_*|,|K_*|$ is $|I|$. If $a\geq 1$, by \eqref{sec5ef1}, we have\fm{|\partial_s^b\partial_q^{a-1}\partial_\omega^c(U_q-\wt{U}_q)|&\lesssim\eps^{-b}\lra{q}^{1-a}t^{-1+C\eps}.}
Thus, the terms \eqref{6.9f1} with $a>0$ have an upper bound \fm{\eps t^{-1}\cdot\eps^{-b}\lra{q}^{1-a}t^{-1+C\eps}\cdot (\lra{q}t^{C\eps})^a\cdot \eps^b\lesssim\eps\lra{q} t^{-2+C\eps}\lesssim \eps\lra{r-t}t^{-2+C\eps}.}
Moreover, by \eqref{sec5ef2}, we have
\fm{|\partial_s^b\partial_\omega^c(U-\wt{U})|\lesssim \eps^{-b}\lra{q}t^{-1+C\eps}.}
Thus, the terms \eqref{6.9f1} with $a=0$ have an upper bound \fm{\eps t^{-1}\cdot\eps^{-b}\lra{q}t^{-1+C\eps}\cdot \eps^b\lesssim\eps\lra{q} t^{-2+C\eps}\lesssim \eps\lra{r-t} t^{-2+C\eps}.}
In conclusion, $|Z^I\mathcal{R}_3|\lesssim\eps t^{-2+C\eps}\lra{r-t}$ for $|I|>0$.
\end{proof}\rm

\lem{Fix a constant $0<\gamma<1$ and a large integer $N$. Then, for $\eps\ll_{\gamma,N} 1$, at each $(t,x)\in\Omega$ such that $|r-t|\lesssim  t^{\gamma}$, we have $|Z^I\mathcal{R}_4|\lesssim_\gamma \eps t^{-2+C\eps}\lra{r-t}$ for each $|I|\leq N$. }
\begin{proof}First we consider the case  $|I|=0$. We have 
\fm{|\mathcal{R}_4|&\lesssim\eps r^{-1}|\hat{U}(s,F(q(t,x),\omega),\omega)-\eps r^{-1}\hat{U}(s,\hat{q}(t,x),\omega)|\\
&\lesssim\eps t^{-1}|\int_{\hat{q}}^{F(q,\omega)}|\hat{U}_\rho(s,\rho,\omega)|\ d\rho|\lesssim \eps t^{-1}|\int_{\hat{q}}^{F(q,\omega)}(|\partial_\rho A_2|+|A_2||\partial_\rho A|)t^{C\eps}\ d\rho|\\
&\lesssim \eps \lra{r-t}^{-2}t^{-1+C\eps}|p(t,x)|\lesssim \eps t^{-2+C\eps}\lra{r-t}^{-1}.}
In the second last inequality, we apply \eqref{lem6.8r12f1} to see that the integrand is $O(\lra{r-t}^{-2}t^{C\eps})$. In the last inequality we apply Proposition \ref{prop6.6}.

In general, fix a multiindex $I$ with $|I|>0$. Then, we can express $Z^I\mathcal{R}_4$ as a linear combination of terms of the form \eq{\label{6.10f1} [(\partial_s^b\partial_q^a\partial_\omega^c\hat{U})(s,F(q,\omega),\omega)\cdot \prod_{i=1}^aZ^{I_i}(F(q,\omega))-(\partial_s^b\partial_q^a\partial_\omega^c\hat{U})(s,\hat{q},\omega)\cdot \prod_{i=1}^aZ^{I_i}\hat{q}]\\\cdot Z^{I'}(\eps r^{-1})\cdot \prod_{j=1}^bZ^{J_j}(\eps\ln t-\delta)\cdot\prod_{l=1}^c Z^{K_l}\omega}
where the sum of all these multiindices is $|I|$. The estimates for such terms are similar to those for \eqref{zqclf2}. The second row is $O(\eps^{b+1} t^{-1+C\eps})$ while the first row is equal to the sum of
 \eq{\label{6.10f2} [(\partial_s^b\partial_q^a\partial_\omega^c\hat{U})(s,F(q,\omega),\omega)-(\partial_s^b\partial_q^a\partial_\omega^c\hat{U})(s,\hat{q},\omega)]\cdot \prod_{i=1}^aZ^{I_i}(F(q,\omega))}
and for each $j=1,2,\dots,a$
\eq{\label{6.10f3}(\partial_s^b\partial_q^a\partial_\omega^c\hat{U})(s,\hat{q},\omega)\cdot \prod_{i=1}^{j-1}Z^{I_i}(F(q,\omega))\cdot Z^{I_j}p\cdot\prod_{i=j+1}^aZ^{I_i}\hat{q}.}
Since $|I|>0$, we must have $a+b+|c|>0$ if \eqref{6.10f3} appears.

Note that  \fm{Z^{I_*}(\hat{q},F(q,\omega))&=O(\lra{r-t}t^{C\eps}),\hspace{2em}Z^{I_*}p=O(t^{-1+\gamma+C\eps});}
\fm{(\partial_s^b\partial_q^a\partial_\omega^c\hat{U})(s,\hat{q},\omega)&=O(\eps^{-b}\lra{\hat{q}}^{1-a+C\eps}t^{C\eps})=O(\eps^{-b}\lra{r-t}^{1-a}t^{C\eps}),\hspace{2em} \text{when }a+b+|c|>0.} 
So \eqref{6.10f3} has an upper bound
\fm{\eps^{-b}\lra{r-t}^{1-a}t^{C\eps}\cdot (\lra{r-t}t^{C\eps})^{a-1}\cdot t^{-1+C\eps}\lra{r-t}\lesssim \eps^{-b}t^{-1+C\eps}\lra{r-t}.}
Besides, by applying Proposition \ref{proplim4} and \eqref{lem6.8r12f1}, we have
\fm{&\hspace{1.5em}|(\partial_s^b\partial_q^a\partial_\omega^c\hat{U})(s,F(q,\omega),\omega)-(\partial_s^b\partial_q^a\partial_\omega^c\hat{U})(s,\hat{q},\omega)|\lesssim|\int_{\hat{q}}^{F(q,\omega)} |\partial_s^b\partial_q^{a+1}\partial_\omega^c\hat{U}|(s,\rho,\omega)\ d\rho|\\
&\lesssim|\int_{\hat{q}}^{F(q,\omega)} \lra{\rho}^{-a-1+C\eps}t^{C\eps}\ d\rho|\lesssim |p(t,x)|\cdot\lra{r-t}^{-a-1+C\eps}t^{C\eps}\lesssim   t^{-1+C\eps}\cdot \lra{r-t}^{-a}.}
In conclusion, \eqref{6.10f2} has an upper bound
\fm{t^{-1+C\eps}\lra{r-t}^{-a}\cdot (\lra{r-t}t^{C\eps})^{a}\lesssim t^{-1+C\eps}.}
Combine all the estimates above and we conclude that
$|Z^I\mathcal{R}_4|\lesssim \eps t^{-2+C\eps}\lra{r-t}$.
\end{proof}\rm

\bigskip
We thus conclude the following approximation result. 
\prop{Fix a constant $0<\gamma<1$ and a large integer $N$. Then, for $\eps\ll_{\gamma,N} 1$, at each $(t,x)\in\Omega$ such that $|r-t|\lesssim  t^{\gamma}$, we have $|Z^I(u-\wt{u})|\lesssim_\gamma \eps t^{-2+C\eps}\lra{r-t}$ for each $|I|\leq N$. }\rm

\bibliography{paperas}{}
\bibliographystyle{plain}

\end{document}